\documentclass[11pt]{article}

\usepackage[T1]{fontenc}
\usepackage[utf8]{inputenc}
\usepackage[margin=1in]{geometry}
\usepackage[numbers,sort&compress]{natbib}
\usepackage{amsthm,amsmath,amsfonts,amssymb}
\usepackage{mathtools,dsfont,mathrsfs,bm,enumitem,upgreek}
\usepackage{graphicx,booktabs,tabularx,longtable,array,multirow}
\usepackage{algorithm,algorithmic}
\usepackage{bbm,mathdots,yhmath,cancel,siunitx,extarrows}
\usepackage{tikz}
\usepackage{setspace}
\usepackage{xcolor}
\usepackage[hidelinks,hypertexnames=false]{hyperref}
\hypersetup{
  pdftitle={Wasserstein Projection Tests: Higher-Order Asymptotics and Connections to Empirical Likelihood},
  pdfauthor={Sirui Lin, Jose Blanchet, Peter Glynn, and Viet Anh Nguyen}
}
\usetikzlibrary{arrows.meta,fadings,patterns,positioning,shadows.blur,shapes}
\allowdisplaybreaks

\DeclareMathOperator{\Tr}{tr}

\DeclareMathOperator*{\argmax}{arg\,max}

\newcommand{\norm}[1]{\left\|#1\right\|}
\newcommand{\abs}[1]{\left|#1\right|}
\newcommand{\brak}[1]{\left\langle#1\right\rangle}
\newcommand{\brac}[1]{\left[#1\right]}
\newcommand{\pare}[1]{\left(#1\right)}
\newcommand{\cure}[1]{\left\{#1\right\}}

\newcommand{\be}{\begin{equation}}
\newcommand{\ee}{\end{equation}}

\newcommand{\mf}{\mathbf}
\newcommand{\Wass}{\mathds W}

\newcommand{\RR}{\mathds R}

\newcommand{\opt}{^\star}

\newcommand{\Cov}{\mathrm{Cov}}

\newcommand{\PD}{\mathbb{S}_{++}}

\newcommand{\mc}{\mathcal}

\newcommand{\EE}{\mathds{E}}
\newcommand{\QQ}{\mathbb{Q}}
\newcommand{\PP}{\mathbb{P}}

\newcommand{\D}{\mathcal{D}}

\newcommand{\Let}{\triangleq}

\newcommand{\diff}{\mathrm{d}}
\newcommand{\eps}{\varepsilon}
\newcommand{\half}{\frac{1}{2}}

\theoremstyle{plain}
\newtheorem{theorem}{Theorem}[section]
\newtheorem{proposition}[theorem]{Proposition}
\newtheorem{corollary}[theorem]{Corollary}
\newtheorem{lemma}[theorem]{Lemma}
\theoremstyle{definition}
\newtheorem{definition}[theorem]{Definition}
\newtheorem{assumption}[theorem]{Assumption}
\newtheorem{remark}[theorem]{Remark}
\newtheorem{example}[theorem]{Example}

\graphicspath{{figures/}}

\newcommand{\EL}{\mathrm{EL}}
\renewcommand{\arraystretch}{1.2}
\DeclareMathOperator{\inte}{int}
\DeclareMathOperator{\conv}{conv}

\newcommand{\suppAlterEdgeworthNumber}{S2.23}

\newcommand{\MainConvexHullNumber}{2.4}
\newcommand{\MainBasicNumber}{2.8}
\newcommand{\MainPowerExpansionNumber}{3.5}
\newcommand{\MainLocationShiftNumber}{3.7}
\newcommand{\MainPowerComparisonNumber}{3.8}
\newcommand{\MainSmoothnessNumber}{3.9}
\newcommand{\MainHigherExpansionNumber}{3.10}
\newcommand{\MainBartlettCoverageNumber}{3.11}
\newcommand{\MainBartlettOneNumber}{3.12}
\newcommand{\MainBartlettTwoNumber}{3.13}
\newcommand{\MainLocalizationNumber}{4.2}
\newcommand{\MainComputableTestNumber}{4.3}
\newcommand{\MainApproxRWPINumber}{6}
\newcommand{\MainHigherExpansionEquationNumber}{12}

\title{Wasserstein Projection Tests: Higher-Order Asymptotics and Connections to Empirical Likelihood}
\author{Sirui Lin$^{1}$ \quad Jos\'{e} Blanchet$^{1}$ \quad Peter Glynn$^{1}$ \quad Viet Anh Nguyen$^{2}$\\[0.5em]
\small $^1$Stanford University \quad $^2$The Chinese University of Hong Kong\\
\small \texttt{siruilin@stanford.edu} \quad \texttt{jose.blanchet@stanford.edu}\\
\small \texttt{glynn@stanford.edu} \quad
\texttt{nguyen@se.cuhk.edu.hk}}
\date{}

\begin{document}
\maketitle

\begin{abstract}
Tests of moment restrictions can be constructed by projecting the empirical distribution onto the set of laws satisfying the null. Wasserstein projection (WP) performs this projection by transporting mass in the sample space, rather than only reweighting observed atoms as in empirical likelihood (EL), and therefore incorporates both the ground geometry and the local variation of the moment function. Existing WP theory provides a first-order null limit but does not quantify finite-sample calibration or distinguish tests with the same limiting local power. We derive stochastic expansions of the rescaled WP statistic under the null and under Pitman-type local alternatives, identifying its order-$n^{-1/2}$ correction with a uniform remainder of order $\widetilde{O}_p(n^{-1})$. The resulting Edgeworth theory gives an order-$\widetilde{O}(n^{-1})$ size error for the plug-in WP test and an explicit order-$n^{-1/2}$ correction to local power. For scalar moments under local location shifts, WP, EL, and Hotelling's $T^2$ have the same first-order power, whereas their second-order ranking is determined by derivative-based curvature for WP and value-based skewness for EL. Extending the expansion by one further order, we construct two Bartlett-type corrections that improve coverage accuracy to $\widetilde{O}(n^{-3/2})$. Finally, we develop a certified localized dual algorithm that controls numerical error in the test decision. A theory-guided smooth fairness experiment using a common asymptotic chi-square reference quantile finds WP power gains over the corresponding EL and $T^2$ tests in the derivative-curvature regime predicted by the comparison theory.
\end{abstract}



\tableofcontents
\clearpage

\section{Introduction}\label{sec:introduction}

Moment restrictions underlie estimating-equation methods~\cite{hansen1982large} and empirical likelihood (EL) inference~\cite{owen2001empirical,qin1994empirical,hjort2009extending}.  Given an empirical distribution \(\QQ_n\) and a null model
\[
    \mc P_0=\{\PP:\EE_{\PP}[\mf h(X)]=\mf 0\},
\]
EL measures the discrepancy from \(\QQ_n\) to \(\mc P_0\) by reweighting the observed sample.  Wasserstein projection (WP) provides a geometric alternative:
\[
    R_n(\mf h)=\inf_{\PP\in\mc P_0}\Wass_c(\QQ_n,\PP),
\]
the least transportation cost required to bring the empirical distribution into agreement with the null moment restriction.  The associated test rejects for large values of \(nR_n(\mf h)\).  This paper uses higher-order asymptotics to quantify finite-sample size distortion and to identify when transport geometry yields power gains.

WP is built from optimal transport (OT), which compares distributions through a ground cost on the sample space.  Although estimation of an unrestricted Wasserstein distance suffers from the curse of dimensionality~\cite{dudley1969speed,fournier2015rate,weed2019sharp}, projection onto a model defined by finitely many moment restrictions retains a parametric rate~\cite{blanchet2019robust}; related results extend to certain infinite-dimensional restrictions~\cite{si2020quantifying}.  This distinction has made WP useful in hypothesis testing, fairness, and distributionally robust optimization.  Other approaches, including Sinkhorn, smoothed, and sliced Wasserstein distances, address dimensionality by modifying or reducing the transport geometry~\cite{cuturi2013sinkhorn,genevay2019sample,nietert2021smooth,goldfeld2024limit,bonneel2015sliced,xi2022distributional}; our focus is instead the inferential behavior of the projection statistic.

The close connection between WP and EL makes their comparison especially informative.  Both project the empirical law onto the same moment model, but they use different discrepancies.  EL changes only the weights assigned to observed atoms, whereas WP can split and transport mass to new locations.  WP therefore incorporates the ground geometry and the local variation of \(\mf h\); it also remains finite whenever the convex hull of the range of \(\mf h\) contains the origin, including settings in which the EL constraint can be infeasible without adjustment~\cite{chen2008adjusted,liu2010adjusted}.  At higher order, this structural difference becomes explicit: the WP power correction is governed by derivative-based curvature of the moment function, while the EL correction is governed by value-based moments of \(h(X)\).

Existing WP theory provides its first-order null distribution and hence an asymptotically valid test~\cite{blanchet2019robust}.  First-order theory, however, does not quantify the error from approximating the finite-sample law by this limit.  Nor can first-order theory distinguish WP from its main competitors under Pitman alternatives when the moment dimension is one (\(d_h=1\)), because WP, EL, and Hotelling's \(T^2\) have the same limiting local-power function.  Determining whether transport geometry helps therefore requires the next term in the size and power expansions.  This is the central motivation for the higher-order analysis developed here.

Fairness testing gives a concrete example.  Moment restrictions can encode equality of opportunity or equality of odds, and \(R_n(\mf h)\) measures the smallest modification of the empirical distribution that makes a fixed classifier satisfy the chosen criterion~\cite{taskesen2021statistical,ref:si2021testing}.  The projected law also has a counterfactual interpretation as the least costly fair modification of the observed population.  For nonlinear fairness restrictions the WP statistic generally has no closed form, so an implementable test additionally requires a reliable numerical procedure.  Section~\ref{sec:computation_algorithm} develops such a procedure and illustrates it in this setting.

More broadly, WP quantiles calibrate Wasserstein ambiguity radii and confidence regions in distributionally robust optimization~\cite{blanchet2019robust,blanchet2022confidence}, while related projections quantify model sensitivity and stability~\cite{blanchet2024stability,il2024quantile}.  In these applications, replacing an unavailable finite-sample quantile by its limiting approximation can produce systematic coverage or conservatism errors.  Thus the higher-order expansion and Bartlett-type corrections below also improve uncertainty calibration beyond hypothesis testing; related higher-order questions for DRO are studied in~\cite{he2025higher}.

Throughout, $\widetilde{O}(\cdot)$ and $\widetilde{O}_p(\cdot)$ denote deterministic and
stochastic orders, respectively, with logarithmic factors in $n$ suppressed.

Our contributions are as follows.
\begin{enumerate}[label=(\arabic*)]
    \item We derive stochastic expansions of \(nR_n(\mf h)\) under the null and under Pitman-type local alternatives.  The expansion identifies the leading quadratic term, its order-\(n^{-1/2}\) cubic correction, and a uniform remainder of order \(\widetilde{O}(n^{-1})\); see Theorems~\ref{thm:main_expan} and~\ref{thm:alter_expand}.

    \item We translate the statistic expansion into inferential guarantees.  The resulting Edgeworth theory yields an order-\(\widetilde{O}(n^{-1})\) size error for the plug-in WP test and an explicit order-\(n^{-1/2}\) correction to its local power; see Theorems~\ref{thm:edgeworth1}, \ref{thm:wp_coverage_accuracy}, and~\ref{thm:wp_power_expansion}.

    \item For scalar moments under local location shifts, we compare WP with EL and Hotelling's \(T^2\).  The tests have the same first-order power, while their second-order ranking is determined by a derivative-based WP curvature coefficient and a value-based EL skewness coefficient; see Proposition~\ref{prop:powercomparison}, Table~\ref{tab:power_comparison_scores}, and Supplementary Sections~S2.9.4 and~S4.2.

    \item We extend the statistic expansion by one further order and construct two Bartlett-type corrections for scalar moment restrictions.  These corrections remove the leading order-\(n^{-1}\) calibration error and improve coverage accuracy to \(\widetilde{O}(n^{-3/2})\); see Theorems~\ref{thm:higherexp_rwp}--\ref{thm:bartlettOT} and Propositions~\ref{prop:BartlettI}--\ref{prop:BartlettII}.

    \item We develop a certified localized dual algorithm for WP testing when the projection lacks a closed form.  The procedure controls the numerical error in the test decision and is illustrated through a nonlinear fairness experiment; see Theorem~\ref{thm:wp_computable_test}, Algorithm~1 in Supplementary Section~S2.10.4, and Section~\ref{sec:comp_fairness}.
\end{enumerate}

Section~\ref{sec:problem_formulation} formulates the test, records the first-order benchmark, and summarizes the regularity conditions.  Section~\ref{sec:MainResults} develops the higher-order expansions and their testing implications.  Section~\ref{sec:computation_algorithm} treats computation, and Section~\ref{sec:discuss} discusses the WP--EL connection and broader extensions.  The notation appendix and Supplementary Sections~S2--S4 provide a comprehensive notation reference, all technical proofs, and computational checks.

\section{Testing framework, first-order benchmark, and assumptions}\label{sec:problem_formulation}

We first fix the transport notation used throughout the paper.
\begin{definition}[Optimal transport cost]\label{def:ot}
For a lower semicontinuous ground cost
$c:\RR^{d_X}\times\RR^{d_X}\to[0,\infty]$ and probability laws
$\QQ,\PP\in\mc P(\RR^{d_X})$, define
\begin{align}\label{def:wass}
    \Wass_c(\QQ,\PP)
    \Let \inf_{\pi}
    \left\{\EE_{\pi}[c(\bar X,X)]:
    \pi|_{\bar X}=\QQ,\ \pi|_X=\PP\right\},
\end{align}
where the infimum ranges over probability laws on
$\RR^{d_X}\times\RR^{d_X}$ and $\pi|_{\bar X},\pi|_X$ denote their marginals.
\end{definition}
Throughout, we use the squared Mahalanobis ground cost
\begin{align}\label{eq:mahala_norm}
    c(\bar x,x)=\norm{\bar x-x}_{\Sigma}^2
    \Let (\bar x-x)^\top\Sigma^{-1}(\bar x-x),
    \qquad \Sigma\in\PD.
\end{align}
Mathematically, $\Sigma$ specifies the local transport geometry.
In practice, it may encode scientifically meaningful relative perturbation costs, with covariance-based choices measuring transport in standardized, correlation-adjusted units and $\Sigma=I_{d_X}$ recovering the ordinary squared Euclidean cost.
We use ``distance'' for
$\Wass_c$ even though the squared cost does not itself satisfy the triangle
inequality.

Let $X_1,\ldots,X_n$ be i.i.d.~observations in $\RR^{d_X}$ with common law
$\PP^\star$, and let $\mf h:\RR^{d_X}\to\RR^{d_h}$ be a user-specified moment
function.  We test whether the data-generating distribution satisfies the moment
restriction
\begin{align}\label{eq:moment_testing_problem}
    \mc H_0: \EE_{\PP^\star}[\mf h(X)] = \mf 0
    \qquad\text{vs.}\qquad
    \mc H_1: \EE_{\PP^\star}[\mf h(X)] \neq \mf 0.
\end{align}
Equivalently, if
\[
    \mc P_0 \Let \{\PP\in\mc P(\RR^{d_X}): \EE_{\PP}[\mf h(X)] = \mf 0\},
\]
then $\mc H_0$ is the membership statement $\PP^\star\in\mc P_0$.  This setup
includes ordinary moment tests, inference based on estimating equations
$\mf h(x)=\mf f(\theta,x)$, and fairness tests in which the moment restriction
encodes the target fairness criterion~\cite{taskesen2021statistical}.

Let $\QQ_n=n^{-1}\sum_{i=1}^n\delta_{X_i}$ denote the empirical law.
\begin{definition}[Wasserstein projection]\label{def:rwp}
The WP distance from $\QQ_n$ to the null model is
\begin{align}\label{eq:rwp}
    R_n(\mf h)
    \Let \inf_{\PP\in\mc P_0}\Wass_c(\QQ_n,\PP)
    = \inf_{\PP\in\mc P(\RR^{d_X})}
    \left\{\Wass_c(\QQ_n,\PP): \EE_{\PP}[\mf h(X)]=\mf 0\right\}.
\end{align}
\end{definition}
This is a projection test: the null determines a set of distributions, and the
statistic measures the distance from the empirical law to that set.  Hence
$R_n(\mf h)$ is the least cost required to transport the data into agreement
with the null, and a large value is evidence against $\mc H_0$.  Under the
null, the empirical moment fluctuates at order $n^{-1/2}$ and the transport
cost is locally quadratic, so $R_n(\mf h)=O_p(n^{-1})$ and $nR_n(\mf h)$ has a
nondegenerate limit.

For nominal level $\varrho$, the WP test rejects when
\begin{align}\label{eq:wp_test_overview}
    \text{reject } \mc H_0
    \quad\text{if}\quad
    nR_n(\mf h)>\hat z_{1-\varrho},
\end{align}
where $\hat z_{1-\varrho}$ is the plug-in $(1-\varrho)$-quantile of the
limiting null law.  Hotelling's $T^2$, empirical likelihood (EL), and WP all
measure violation of the same moment restriction but impose different
geometries.  Hotelling's statistic measures the empirical moment directly; EL
projects by reweighting the observed atoms; WP transports mass in the sample
space and therefore responds to both the ground metric $\Sigma$ and local
changes in $\mf h$.  Their comparison isolates whether this transport geometry
improves finite-sample calibration or power after their common first-order
behavior is taken into account.

\subsection{First-order benchmark}

Near the null, the transport perturbations required to satisfy the moment restriction are
small, so linearizing the moment restriction reduces the WP problem to a quadratic program
that determines the test's first-order approximation.  The exact rescaled dual representation is given in Supplementary Section~S2.1.1.  We use
$[d]\Let\{1,\ldots,d\}$ and let $\D\mf h(x)\in\RR^{d_h\times d_X}$ denote the
Jacobian of $\mf h$, with $\D h^\beta(x)$ its $\beta$th row.  For a vector $u$,
$u^{\otimes k}$ denotes its $k$th outer power, and
$\brak{A,u^{\otimes k}}$ denotes full contraction with a $k$th-order tensor
$A$.  Define
\begin{align*}
    H_n&\Let\frac{1}{\sqrt n}\sum_{i=1}^n\mf h(X_i),
    &
    \mc V_n&\Let\frac1n\sum_{i=1}^n
    \D\mf h(X_i)\Sigma\D\mf h(X_i)^\top.
\end{align*}
Expanding the inner dual objective then gives the compact approximation
\[
    nR_n(\mf h)=H_n^\top\mc V_n^{-1}H_n+o_p(1).
\]
This yields the standard first-order null law below; see
\cite[Propositions~3 and~5]{blanchet2019robust}.

\begin{proposition}[First-order null limit]\label{prop:first_order_null_limit}
Suppose that $\mf h\in C^1(\RR^{d_X})$ is bounded, the moment constraint is
feasible, and the following two matrices are positive definite:
\[
    W\Let \EE[\mf h(X)^{\otimes 2}],
    \qquad
    V\Let \EE[\D\mf h(X)\Sigma\D\mf h(X)^\top].
\]
Then, under $\mc H_0$,
\[
    nR_n(\mf h)\overset{d}{\longrightarrow} H^\top V^{-1}H,
    \qquad
    H \sim \mc N(\mf 0,W).
\]
\end{proposition}
Boundedness is used here only as a convenient sufficient condition for this
first-order benchmark.  The higher-order results below instead accommodate
unbounded moment functions through the stated derivative-envelope and moment
conditions.

The limit justifies the plug-in critical value in~\eqref{eq:wp_test_overview},
but it leaves its finite-sample accuracy unresolved.  Both the statistic and
its critical value contain estimated quantities, so asymptotic validity alone
does not identify the leading size distortion or indicate how it should be
corrected.  This requires a distributional expansion beyond the quadratic
limit.

First-order theory also cannot determine when transport geometry improves
power.  For scalar moments under Pitman location shifts, WP, EL, and Hotelling's
$T^2$ have the same limiting local power; their first meaningful ranking appears
in the order-$n^{-1/2}$ term.  That term distinguishes the derivative-based
geometry used by WP from the value-based information used by EL.  Extending the
expansion by one further order then identifies the leading calibration error and
makes a Bartlett-type correction possible.  These three gaps motivate the null,
local-power, and higher-order developments in Section~\ref{sec:MainResults}.

\subsection{Assumptions}\label{sec:prelim}
We state the regularity conditions here to make the scope of the main results
explicit.  Assumptions~\ref{a:convexhull}--\ref{a:finite_moments} support the
stochastic expansion of the WP statistic, while Assumption~\ref{a:moments_II}
provides the additional moment control needed for its Edgeworth expansion.
Further motivation and representative bounds are collected in
Supplementary Section~S2.1.2.

\subsubsection{\texorpdfstring{Conditions on the moment function $\mf h$}{}}
\begin{assumption}[Feasibility of the moment constraint]\label{a:convexhull}
    The origin lies in the interior of the convex hull of the range of $\mf h$:
    \[
        \mf 0 \in \inte\left(\conv\{\mf h(x):x\in\RR^{d_X}\}\right).
    \]
\end{assumption}

This interiority requirement is a Slater-type condition ensuring that the WP projection
admits an attained finite-dimensional dual representation; see Supplementary Section~S2.7.2.
Expanding that dual formulation is the central device driving our higher-order theory.

\begin{assumption}[Lipschitz control of the Jacobian]\label{a:D1}
    We assume $\mf h\in C^3(\RR^{d_X})$. There exists a function $\kappa_1:\RR^{d_X}\to\RR_+$ such that, for all $(x,\Delta,\zeta)\in\RR^{d_X}\times\RR^{d_X}\times\RR^{d_h}$,
    \[
        \Big\| \sum_{\beta = 1}^{d_h} \zeta^{(\beta)}
        \left(\D h^{\beta}(x+\Delta)-\D h^{\beta}(x)\right) \Big\|_2
        \leq \kappa_1(x)\norm{\Delta}_2\norm{\zeta}_2,
    \]
    where $\mf h=(h^\beta)_{\beta\in[d_h]}$ and $\zeta=(\zeta^{(\beta)})_{\beta\in[d_h]}$.
\end{assumption}

In particular, $\kappa_1$ bounds directional Hessians and enforces growth compatible with the
squared ground cost, preventing projection degeneracy; a counterexample and alternative cost
choices are given in Supplementary Section~S2.1.2.

\begin{assumption}[Local Lipschitz control of the Hessian]\label{a:D2}
    There exist $\kappa_2:\RR^{d_X}\to\RR_+$ and $\hat\delta>0$ such that, for all $(x,\Delta,\zeta)\in\RR^{d_X}\times\RR^{d_X}\times\RR^{d_h}$ with $\norm{\Delta}_2\leq\hat\delta$,
    \[
        \Big\| \sum_{\beta=1}^{d_h} \zeta^{(\beta)}
        \left(\D^2 h^{\beta}(x+\Delta)-\D^2 h^{\beta}(x)\right) \Big\|_2
        \leq \kappa_2(x)\norm{\Delta}_2\norm{\zeta}_2.
    \]
\end{assumption}

This condition controls the third-order Taylor remainder; without such smoothness, the
projection may enter a different boundary-layer regime, as discussed in Supplementary Section~S2.1.2.

\subsubsection{\texorpdfstring{Conditions on the data-generating distribution $\PP^\star$}{}}
\begin{assumption}[Absolute continuity and Cram\'er regularity]\label{a:cont_P}
    The distribution $\PP^\star$ has a density with respect to Lebesgue
    measure on $\RR^{d_X}$.  Moreover, after deterministic affine redundancies
    have been removed, the functions of a single observation $X$ whose sample
    averages enter the Edgeworth and signed-root expansions jointly satisfy
    the standard multivariate nonlattice (Cram\'er) regularity condition.  The
    precise collection of functions and the formal condition are given in
    Supplementary Section~S2.6.
\end{assumption}

The nonlattice requirement is used only to justify the Edgeworth expansions.

\begin{assumption}[Null hypothesis and nondegeneracy]\label{a:basic}
    Under $\PP^\star$,
    \[
        \EE[\mf h(X)]=\mf 0,
        \qquad
        \mathrm{Cov}(\mf h(X),\mf h(X))\in\PD,
        \qquad
        \EE\brac{\D\mf h(X)\Sigma\D\mf h(X)^\top}\in\PD.
    \]
\end{assumption}

The first condition is the null hypothesis.  The two positive-definiteness
conditions rule out degenerate moment variation and a degenerate local
transport metric, respectively.

\begin{assumption}[Moment envelope for the asymptotic expansion]\label{a:finite_moments}
    Under $\PP^\star$,
    \[
        \EE\left[\left(\norm{\mf h(X)}_2+\norm{\D \mf h(X)}_2+\kappa_1(X)+\kappa_2(X)\right)^8\right]<\infty.
    \]
\end{assumption}

These moments control the stochastic Taylor remainder uniformly with probability
at least $1-O(n^{-1})$.

\begin{assumption}[Moment condition for the Edgeworth expansion]\label{a:moments_II}
    Under $\PP^\star$,
    \[
        \EE\left[\norm{\D \mf h(X)}_2^8\kappa_1(X)^4\right]<\infty.
    \]
\end{assumption}

This product moment supplies the fourth-moment bound for the summands entering
$\mc K_n$, as required for an order-$n^{-1}$ Edgeworth expansion.


\section{Main results}\label{sec:MainResults}
\subsection{Second-order asymptotic expansion}\label{sec:main_result}
We first derive an expansion of the WP statistic in powers of $n^{-\half}$ and use it to obtain a second-order Edgeworth approximation based on classical normal-approximation theory~\cite{bhattacharya2010normal}.
 In this section, sample-level probabilities are evaluated under the product
 law $(\PP^\star)^{\otimes n}$, and $\EE$ denotes expectation under that law;
 for a single observation, $\PP^\star$ retains its usual meaning.

\begin{theorem}[Second-order asymptotic expansion]
\label{thm:main_expan}
    If Assumptions~\ref{a:convexhull}--\ref{a:finite_moments} hold, then there is a deterministic sequence $\delta_n = \widetilde{O}\left(n^{-1}\right)$ and an integer $N$ such that when $n \geq N$, we have
    \begin{align}\label{eq:approx_rwpi}
         nR_n(\mf h) = 
         \brak{ \mc V_n, \xi_n^{\otimes 2} }  + \frac{1}{
         \sqrt{n}} \brak{ \mc K_n, \xi_n^{\otimes 3}}  + \eps_n.
    \end{align}
    The remainder satisfies $\PP\opt\pare{\abs{\eps_n}\leq\delta_n}=1-O(n^{-1})$.
    Here $\mc V_n$, $\xi_n$, and $\mc K_n$ are empirical quantities dependent on $\mf h$ via
    \begin{align*}
        \mc V_n ~&\Let~ \left(\frac{1}{n}\sum_{i=1}^{n} \D h^{\beta}(X_i)\Sigma\D h^{\gamma}(X_i)^\top\right)_{\beta, \gamma \in [d_h]} \in \RR^{d_h\times d_h},\\
        \xi_n ~&\Let~ \frac{1}{\sqrt{n}}\sum_{i=1}^n \mc V_n^{-1}\mf h(X_i) \Let \left(\xi_n^{(\beta)}\right)_{\beta \in [d_h]} \in \RR^{d_h},\\
         \mc K_n ~&\Let \left(\frac{1}{n}\sum_{i=1}^n \D h^{\beta}(X_i)\Sigma\D^2 h^{\gamma}(X_i)\Sigma\D h^{\omega}(X_i)^\top\right)_{\beta, \gamma, \omega \in [d_h]} \in (\RR^{d_h})^{\otimes 3}.
    \end{align*} 

\end{theorem}
The proof, the explicit choice of $N$, and exact checks for linear and quadratic moments are
given, respectively, in Supplementary Sections~S2.7, S2.12, and~S2.5.  The latter checks
verify the expansion against exact formulas for linear and quadratic moment functions.

To pass from Theorem~\ref{thm:main_expan} to a distributional approximation, its leading
terms can be represented as the squared norm of a signed-root statistic.  Under the null,
the order-$n^{-1/2}$ Edgeworth correction is odd and therefore integrates to zero over the
symmetric acceptance region.  The signed-root construction, its moment expansion, and this
cancellation are given in Supplementary Sections~S2.8.1--S2.8.2 and~S2.13.

We next establish a valid Edgeworth expansion under the regularity assumptions imposed in Section~\ref{sec:prelim}.

\begin{theorem}[Second-order Edgeworth expansion]\label{thm:edgeworth1}
    Under Assumptions~\ref{a:convexhull}--\ref{a:moments_II}, for every fixed $z_0>0$,
    \begin{align}\label{eq:edgeworth_ninf}
        \sup_{z\geq z_0}\left|\PP^\star\left(nR_n(\mf h) \leq z\right)-\Psi(z)\right|
        = \widetilde{O}\left(\frac{1}{n}\right),
    \end{align}
    where $\Psi$ is the CDF of the random variable $H^\top V^{-1} H$ with $H\sim \mc N(\mf 0,W)$. In particular, the expansion holds at every fixed $z>0$.
\end{theorem}

The cancellation of the $n^{-\half}$ term under the null is common for chi-squared-like
statistics, including likelihood-ratio~\cite[Equation~(30)]{Box1994} and empirical-likelihood
statistics~\cite[Section~3.2]{hall1990methodology}.  Under local alternatives, this term is
generally nonzero and determines the second-order power correction derived in
Theorem~\ref{thm:wp_power_expansion}.

\subsection{Asymptotic expansion under local alternatives}\label{sec:expan_alter}

We now let the one-observation law $\PP_n\opt$ vary with $n$, draw
$X_1,\ldots,X_n$ independently from this law, and write $\EE_n$ for the
corresponding expectation.  For sample-level events, $\PP_n\opt$ denotes the
product law $(\PP_n\opt)^{\otimes n}$.  We consider Pitman-type local alternatives satisfying
\[
    \EE_n[\mf h(X)]=O(n^{-1/2}),
\]
so the violation of the null moment condition is of the same order as sampling
noise.  The uniform nondegeneracy, moment, and Cram\'{e}r-type conditions needed
for the expansion are stated and explained in
Supplementary Section~S2.8.3.  Their role is to keep the
stochastic expansion and its remainder uniform along the sequence
$(\PP_n\opt)_{n\geq1}$.

\begin{theorem}[Asymptotic expansion under local alternatives]\label{thm:alter_expand}
    Suppose Assumptions~\ref{a:convexhull}--\ref{a:D2} and the uniform
    local-alternative conditions in Assumption~\suppAlterEdgeworthNumber{} of
    Supplementary Section~S2.8.3 hold.  Then there is a
    deterministic sequence $\delta_n=\widetilde{O}(n^{-1})$ and an integer $N$ such that, for
    $n\geq N$,
    \begin{align*}
         nR_n(\mf h) = 
         \brak{ \mc V_n, \xi_n^{\otimes 2} }  + \frac{1}{
         \sqrt{n}} \brak{ \mc K_n, \xi_n^{\otimes 3}}  + \eps_n.
    \end{align*}
    Here $\xi_n$, $\mc V_n$, and $\mc K_n$ are defined immediately after
    \eqref{eq:approx_rwpi} in Theorem~\ref{thm:main_expan},
    and $\PP_{n}\opt\pare{\abs{\eps_n}\leq\delta_n}=1-O(n^{-1})$.
\end{theorem}

The proof follows the null expansion with uniform bounds in place of fixed-law
bounds and is given in Supplementary Section~S2.8.3.  A stronger set
of conditions, needed only for the Edgeworth power expansion, is stated beside
that proof in Supplementary Section~S2.9.3.

\subsection{Implications for WP-based hypothesis testing}\label{sec:implication}
For confidence level $1-\varrho\in(0,1)$, the ideal asymptotic rule rejects
$\mc H_0$ when $nR_n(\mf h)>z_{1-\varrho}$, where $z_{1-\varrho}$ is the
$(1-\varrho)$-quantile of the limiting law of $nR_n(\mf h)$. Under suitable
assumptions, this law is that of the quadratic form $H^\top V^{-1}H$, with
$H\sim\mc N(\mf 0,W)$. Because the population matrices $W$ and $V$ are
unknown, however, this quantile is not directly available.
To address this issue, we use the sample analogues
\begin{align}\label{eq:WnVn}
    \mc W_n \Let  \frac{1}{n}\sum_{i=1}^n \mf h(X_i)^{\otimes 2},\quad \mc V_n \Let \frac{1}{n}\sum_{i=1}^n \D \mf h(X_i)\Sigma\D\mf h(X_i)^\top,    
\end{align}
and consider the plug-in CDF $\hat F_n$ defined by
\begin{align*}
    \hat F_n(z) \Let \int_{\left\{v \in \RR^{d_h}:~v^\top \mc W_n^{\half} \mc V_n^{-1} \mc W_n^{\half} v \leq z\right\}} \phi(v) \diff v,
\end{align*}
where $\phi$ is the density function of the standard normal distribution on $\RR^{d_h}$. By the strong law of large numbers and Assumption~\ref{a:basic}, $\mc W_n$ and $\mc V_n$ are eventually invertible almost surely. Thus $\mc W_n^{\half}\mc V_n^{-1}\mc W_n^{\half}$ and the following plug-in quantile are well defined on an event whose probability tends to one.

We therefore estimate $z_{1-\varrho}$ by the $(1-\varrho)$-quantile of $\hat F_n$:
\begin{equation} \label{eq:hat-z-def}
    \hat z_{1-\varrho} \Let \hat F_n^{-1}(1-\varrho).
\end{equation}
The WP decision rule is therefore
\begin{equation} \label{def:rwptest}
    \text{reject $\mc H_0$ ~if ~$n R_n(\mf h) > \hat z_{1-\varrho}$.} \tag{WP test}
\end{equation}
By Slutsky's theorem, under the null the acceptance probability of the decision
rule in~\eqref{def:rwptest} converges to $1-\varrho$, or equivalently its Type-I
error converges to $\varrho$.  This is the desired first-order property.

\subsubsection{Confidence level expansion}

The confidence level of the WP-based test is quantified by
\[
\underbrace{\PP^\star\left(n R_n(\mf h) \leq \hat z_{1-\varrho}\right)}_{\text{Actual confidence level}}~~=~~ \underbrace{1-\varrho}_{\text{Desired confidence level}} ~~+~~ \text{error}.
\]
The error term characterizes the order at which the actual confidence level
converges to the target. This error has two sources: (i) the finite-sample
$(1-\varrho)$-quantile of $nR_n(\mf h)$ is unavailable, so we use the quantile
of its \textit{limiting} distribution; and (ii) this limiting distribution
depends on unknown population quantities that must be estimated. To quantify
the error, we need the second-order expansion of $nR_n(\mf h)$ from
Theorem~\ref{thm:main_expan} and the corresponding expansion of the plug-in
quantile $\hat z_{1-\varrho}$ derived below.
The statistic expansion and the corresponding expansion of the plug-in critical value can be
combined through a signed-root reduction.  The quantile expansion and this reduction are
given with the proof in Supplementary Sections~S2.9.1--S2.9.2.
\begin{theorem}[Coverage accuracy of the WP test I]
\label{thm:wp_coverage_accuracy}
Under Assumptions~\ref{a:convexhull}--\ref{a:moments_II}, we have
\begin{align}\label{eq:expan_error}
    \PP^\star\left(n R_n(\mf h) \leq \hat z_{1-\varrho}\right) = 1-\varrho + \widetilde{O}\left(\frac{1}{n}\right).
\end{align}
\end{theorem}

According to~\cite[Section 23.3]{van2000asymptotic}, the above result shows that the confidence level of~\eqref{def:rwptest} is accurate up to order $\widetilde{O}(n^{-1})$. The rate $n^{-1}$ is typical for the coverage error of both parametric likelihood and nonparametric methods (see, e.g., the discussions of empirical likelihood in~\cite[Section~2.6]{owen2001empirical} and the bootstrap in~\cite[Section~3.5]{hall2013bootstrap}).

When $\mf h(x)=\mf f(\theta,x)$, inversion of the WP test gives the confidence region
\begin{equation}\label{eq:confidence_region_Theta}
    \hat\Theta_n=\{\theta:nR_n(\mf f(\theta,\cdot))\leq\hat z_{1-\varrho}(\theta)\},
\end{equation}
where $\hat z_{1-\varrho}(\theta)$ is the plug-in quantile computed using the
moment function $\mf f(\theta,\cdot)$.  For each fixed true value $\theta_0$
for which the assumptions of Theorem~\ref{thm:wp_coverage_accuracy}
hold with $\mf h(x)=\mf f(\theta_0,x)$, its pointwise coverage satisfies
$\PP^\star(\theta_0\in\hat\Theta_n)=1-\varrho+\widetilde{O}(n^{-1})$.

\subsubsection{Local power expansion}\label{sec:localter}

At confidence level $1-\varrho$, the power is
$\PP_n\opt(nR_n(\mf h)>\hat z_{1-\varrho})$.  We expand this probability under
the Pitman-type alternatives introduced above.  The first-order term gives the
limiting power, while the $n^{-1/2}$ term distinguishes tests that share that
limit.  The latter requires the stronger uniform moment and convergence
conditions in Supplementary Section~S2.9.3; unlike classical
contiguity assumptions, these conditions are formulated for a sequence of
individual-observation distributions converging in total variation.

The following result gives the power of the WP test~\eqref{def:rwptest}, namely, the probability that $n R_n(\mf h) > \hat z_{1-\varrho}$. Recall that $\EE_n$ denotes expectation under $\PP_n\opt$. Define $W_{n} = \EE_{n}[\mc W_n]$ and $V_{n} = \EE_{n}[\mc V_n]$, so that
\[
\pare{W_{n}}_{\beta\gamma} = \EE_{n}\brac{h^{\beta}(X) h^{\gamma}(X)},~~ \pare{V_{n}}_{\beta\gamma} = \EE_{n}\brac{\D h^{\beta}(X) \Sigma \D h^{\gamma}(X)^\top}~~\forall \beta, \gamma \in [d_h].
\]

\begin{theorem}[Power expansion of WP-based test]\label{thm:wp_power_expansion}
    Under Assumptions~\ref{a:convexhull}--\ref{a:D2} and the additional Edgeworth conditions stated in Supplementary Section~S2.9.3, we have
    \begin{align*}
        \PP_{n}\opt\left(n R_n(\mf h) > \hat z_{1-\varrho}\right)
        ={}& \underbrace{\int_{\norm{v + \tau_n}_2^2 > z_{n,1-\varrho}} \bar \phi(v) \diff v}_{\text{first-order}} \\
        &{}+ \underbrace{\frac{1}{\sqrt{n}} \int_{\norm{v + \tau_n}_2^2 > z_{n,1-\varrho}} \bar p(v) \bar \phi(v) \diff v}_{\text{second-order}}
        + \widetilde{O}\left(\frac{1}{n}\right),
    \end{align*}
    where $\tau_n = \sqrt{n} V_{n}^{-\half} \EE_{n}\brac{\mf h(X)} \in \RR^{d_h}$, $\bar \phi$ is the density of $H \sim \mc N\big(0,V_{n}^{-\half} W_{n} V_{n}^{-\half} \big)$, and $z_{n,1-\varrho}$ is the $(1-\varrho)$-quantile of $\norm{H}_2^2$.  The explicit correction polynomial $\bar p$, including the cumulant maps and Hermite polynomials, is given in Supplementary Section~S4.1.
\end{theorem}
The additional conditions strengthen those of
Theorem~\ref{thm:alter_expand} only enough to validate this second-order Edgeworth
expansion.  Its first-order consequence is the following limiting-power result.

\begin{corollary}[Asymptotic power of the WP-based test]\label{cor:limiting_power}
    Under Assumptions~\ref{a:convexhull}--\ref{a:D2} and the additional Edgeworth conditions stated in Supplementary Section~S2.9.3, if $\lim_{n\rightarrow\infty} \tau_n = \tau \in \RR^{d_h}$, $\lim_{n\rightarrow\infty} V_n = V\in \RR^{d_h\times d_h},$ and $\lim_{n\rightarrow\infty} W_n = W\in \RR^{d_h\times d_h},$ then we have
    \begin{align*}
        \lim_{n\rightarrow \infty} \PP_{n}\opt\left(n R_n(\mf h) > \hat z_{1-\varrho}\right) = \int_{\norm{v + \tau}_2^2 > z_{1-\varrho}} \bar \phi(v) \diff v,
    \end{align*}
    where $\bar \phi$ is the density function of
    $H \sim \mc N\big(0,V^{-\half} W V^{-\half} \big)$ and
    $z_{1-\varrho}$ is the $(1-\varrho)$-quantile of $\norm{H}_2^2$.
\end{corollary}

\subsection{Power comparison under location shifts}\label{sec:localshift_power}

We now specialize the local alternatives in Section~\ref{sec:localter} to
location shifts of the data-generating distribution.

\begin{assumption}[Location-shift alternatives]\label{a:location_shift}
    Let $Y\sim\PP_0\opt\in\mc P_0$ and, for a fixed
    $\tau_0\in\RR^{d_X}$, define
    \[
        X=Y+\frac{\tau_0}{\sqrt n}
        \qquad\text{under }\PP_n\opt.
    \]
    We also impose the mild density regularity stated in
    Supplementary Section~S2.9.4.1.
\end{assumption}

The target of this section is to compare the power of the WP-based test with the empirical likelihood test and Hotelling's $T^2$ test under the local alternatives in Assumption~\ref{a:location_shift}.
The WP-based test is the test in \eqref{def:rwptest}.
For a general $d_h$-dimensional moment function $\mf h$, the empirical likelihood test is defined by
\begin{equation} \label{def:EL}
    \text{reject $\mc H_0$ ~if ~$nR^{EL}_n(\mf h) > \chi^2_{d_h;1-\varrho}$,} \tag{EL test}
\end{equation}
where $\mc R$ is Owen's support-restricted empirical likelihood ratio, obtained
by reweighting the observed atoms \cite[Chapters~2--3]{owen2001empirical}.
Equivalently,
\[
    R^{EL}_n(\mf h)
    =2\inf_{\substack{\PP\in\mc P_0\\ \PP\ll\QQ_n}}
       D_{\mathrm{KL}}(\QQ_n\Vert\PP)
    =-\frac{2}{n}\log\mc R,
\]
where $\PP\ll\QQ_n$ restricts $\PP$ to the empirical support.\footnote{See
Section~\ref{sec:ELvsWP} for further discussion.}
The Hotelling-type $T^2$ test used in the comparison is defined by
\begin{equation}\label{def:hotelling}
    \text{reject $\mc H_0$ ~if ~$\brak{\mc W_n^{-1}, \pare{\frac{1}{\sqrt{n}}\sum_{i=1}^n \mf h(X_i)}^{\otimes 2}} > \chi^2_{d_h;1-\varrho}$,} \tag{$T^2$ test}
\end{equation}
where $\mc W_n$ is the uncentered empirical second-moment matrix defined in
\eqref{eq:WnVn}.  Under the null and the local alternatives considered below,
this statistic differs from the conventional centered Hotelling statistic by
$O_p(n^{-1})$ and is therefore equivalent at the order-$n^{-1/2}$ comparison
studied here \cite{hotelling1947}.
Fix a confidence level $1-\varrho\in(0,1)$.
The rejection thresholds for the three tests are chosen so that, under the null hypothesis, their Type-I error probabilities are asymptotically equal to $\varrho$.
Then, for a test $s$, we want to compare the local power, defined by the rejection probability under this sequence of alternatives,
\[
    \PP_n^\star\pare{\text{Test $s$ rejects } \mc H_0},
    \qquad
    s\in\{\mathrm{WP},\mathrm{EL},T^2\}.
\]

This section uses $\EE_0$ for the expectation under $\PP_0^{\star}$.

\subsubsection{Vector moment restrictions}
We first discuss the comparison when the moment function $\mf h$ is multi-dimensional ($d_h > 1$).
Let
\[
    W_0=\EE_0[\mf h(X)\mf h(X)^\top],
    \qquad
    V_0=\EE_0[\D\mf h(X)\Sigma\D\mf h(X)^\top],
    \qquad
    \Delta_0=\EE_0[\D\mf h(X)]\tau_0.
\]
Here $W_0$ and $V_0$ are deterministic population quantities under the null distribution $\PP_0^\star$; they are the null counterparts of the population matrices $W_n$ and $V_n$ in Theorem~\ref{thm:wp_power_expansion}, and should not be confused with the random plug-in matrices $\mc W_n$ and $\mc V_n$ in \eqref{eq:WnVn}.
The vector $\Delta_0$ denotes the first-order drift of the moment mean induced by the local shift direction $\tau_0$.
Under Assumption~\ref{a:location_shift},
\[
    \EE_n[\mf h(X)]
    =
    \frac{\Delta_0}{\sqrt n}+O(n^{-1}),
\]
so the local shift produces a Pitman drift in the moment condition.

We choose the rejection thresholds so that each test has asymptotic confidence level $1-\varrho$ under $\mc H_0$.
For the WP-based test, this means using the plug-in threshold $\hat z_{1-\varrho}$ in \eqref{eq:hat-z-def}; by Theorem~\ref{thm:wp_coverage_accuracy},
\[
    \PP_0^\star\left(nR_n(\mf h)\leq \hat z_{1-\varrho}\right)
    =
    1-\varrho+\widetilde{O}(n^{-1})
    =
    1-\varrho+o(n^{-1/2}).
\]
For EL and Hotelling's $T^2$, the usual chi-square threshold $\chi^2_{d_h;1-\varrho}$ is used in \eqref{def:EL} and \eqref{def:hotelling}; under the corresponding Wilks and studentized-mean approximations,
\[
    \PP_0^\star\left(nR_n^{EL}(\mf h)\leq \chi^2_{d_h;1-\varrho}\right)
    =
    1-\varrho+o(n^{-1/2}),
\]
and
\[
    \PP_0^\star\left(
    \brak{\mc W_n^{-1},\pare{\frac1{\sqrt n}\sum_{i=1}^n\mf h(X_i)}^{\otimes2}}
    \leq \chi^2_{d_h;1-\varrho}
    \right)
    =
    1-\varrho+o(n^{-1/2}).
\]

We now compare the first-order Type-II errors under the local alternatives.
Let
\[
    A_0=V_0^{-1/2}W_0V_0^{-1/2},
    \qquad
    \tau_{\mathrm{WP}}=V_0^{-1/2}\Delta_0,
    \qquad
    \tau_{\mathrm{EL}}=W_0^{-1/2}\Delta_0,
\]
and let $\phi_{d_h}$ denote the density of $\mc N(0,I_{d_h})$.
The vector $\tau_{\mathrm{WP}}$ is the limit of the standardized drift $\tau_n=\sqrt n\,V_n^{-1/2}\EE_n[\mf h(X)]$ in Theorem~\ref{thm:wp_power_expansion}, while $\tau_{\mathrm{EL}}$ is the corresponding covariance-standardized drift for EL and Hotelling's $T^2$.
Let $z^{\mathrm{WP}}_{1-\varrho}$ be the $(1-\varrho)$-quantile of $Z^\top A_0Z$ for $Z\sim \mc N(0,I_{d_h})$.
Then Theorem~\ref{thm:wp_power_expansion} implies that the Type-II error of the WP-based test satisfies
\[
    \PP_n^\star\left(nR_n(\mf h)\leq \hat z_{1-\varrho}\right)
    =
    \int_{\norm{A_0^{1/2}z+\tau_{\mathrm{WP}}}_2^2\leq z^{\mathrm{WP}}_{1-\varrho}}
    \phi_{d_h}(z)\diff z
    +o(1).
\]
For EL and Hotelling's $T^2$, the first-order Type-II errors are both
\[
    \int_{\norm{z+\tau_{\mathrm{EL}}}_2^2\leq \chi^2_{d_h;1-\varrho}}
    \phi_{d_h}(z)\diff z
    +o(1).
\]
Thus, in the multi-dimensional case, WP and EL/Hotelling's $T^2$ generally need not have the same first-order local power.
The WP limit involves the non-standard quadratic form $Z^\top A_0Z$, whereas EL and Hotelling's $T^2$ reduce to the usual noncentral chi-square form.
Consequently, a universal ordering is not available in general; for a fixed local shift direction $\tau_0$, the test with the smaller first-order Type-II error above has the larger local power.
This first-order comparison issue disappears when $d_h=1$: after scalar normalization, the three tests share the same first-order local power, so the comparison can be made through the second-order terms.

\subsubsection{Scalar moment restriction}
We specialize to a single scalar moment restriction $d_h=1$ (with $d_X$
arbitrary) and write the moment function as $h$.  The ordering below applies
only to the local location shifts and regularity conditions of
Proposition~\ref{prop:powercomparison}, not to vector moments or other
alternatives.
In this case, the non-standard first-order comparison from the previous subsection collapses to a scalar normalization.
The first-order power term in Theorem~\ref{thm:wp_power_expansion} becomes
\[
\int_{\norm{v + \tau_n}_2^2 > z_{n,1-\varrho}} \bar \phi(v) \diff v = \int_{\abs{u + \sqrt{n}\EE_{n}\brac{h(X)}/\sqrt{\alpha_{2,n}}}^2 > \chi^2_{1;1-\varrho}} \phi(u) \diff u,
\]
where $\alpha_{2,n} = \EE_{n}\brac{h(X)^2}$, $\phi$ is the $\mc N(0,1)$ density, and
$\chi^2_{1;1-\varrho}$ is the $(1-\varrho)$-quantile of $\chi_1^2$.  EL and
Hotelling's $T^2$ share this first-order power, so only the $n^{-1/2}$ terms
distinguish the three tests.

The full second-order WP expansion and its cumulant coefficients are given in
Supplementary Sections~S2.9.4 and~S4.2; here we retain only the comparison with
EL and Hotelling's $T^2$.
For EL, set $Z_{n,i}=h(X_i)$.  The EL constraint depends on $X_i$ only
through $Z_{n,i}$, so the statistic is exactly the scalar mean EL statistic for
the transformed observations.  A uniform signed-root derivation for the
triangular array $(Z_{n,i})_{i=1}^n$ is given in Supplementary
Section~S2.9.4; its coefficients agree with the fixed-law scalar mean formula
in \cite[Theorem~3.1]{chen1994comparing}.
For Hotelling's $T^2$, the test in \eqref{def:hotelling} is equivalent to
\begin{equation*}
    \text{reject $\mc H_0$ ~if ~$\brak{\mc V_n, \xi_n^{\otimes 2}} > \hat z_{1-\varrho}$,}
\end{equation*}
where $\xi_n$ and $\mc V_n$ are defined in \eqref{eq:approx_rwpi}, and
$\hat z_{1-\varrho}$ is defined in \eqref{eq:hat-z-def}.
Thus its power calculation is the WP calculation with the cubic tensor
$\mc K_n$ set to zero.

Let $s \in \{\mathrm{WP}, \mathrm{EL}, T^2\}$ index the tests in
\eqref{def:rwptest}, \eqref{def:EL}, and \eqref{def:hotelling}.  The following
proposition gives their pairwise second-order comparison.
\begin{proposition}[Pairwise second-order power comparison]\label{prop:powercomparison}
    Suppose Assumption~\ref{a:location_shift} and the conditions of
    Theorem~\ref{thm:wp_power_expansion} hold with $d_h=1$.
    Write $\alpha_j=\EE_0[h(X)^j]$ for $j=2,3$ and
    $\tilde\alpha_2=\EE_0[\D h(X)\Sigma\D h(X)^\top]$.
    Also set $\tilde\alpha_3=
    \EE_0[\D h(X)\Sigma\D^2h(X)\Sigma\D h(X)^\top]$; the standing
    nondegeneracy conditions give $\alpha_2,\tilde\alpha_2>0$.
    For a confidence level $1-\varrho \in (0, 1)$, define
    \[
        \pi_{n,s}\Let
        \PP_{n}\opt\pare{\text{Test $s$ rejects the null hypothesis $\mc H_0$}}.
    \]
    Then, for every $s,t\in\{\mathrm{WP},\mathrm{EL},T^2\}$,
    \[
        \pi_{n,s}-\pi_{n,t}
        =\frac{B(s)-B(t)}{\sqrt n}+o(n^{-1/2}),
    \]
    where
    \[
    B(s) = \begin{cases}
        \frac{\tilde \alpha_{3}}{\tilde \alpha_{2}^2} \times \sqrt{\alpha_2}\,I(\varrho, \tau)  & \text{if } s = \mathrm{WP}, \\[1.5ex]
        \frac{ 2\alpha_{3} }{3\alpha_{2}^2} \times \sqrt{\alpha_2}\,I(\varrho, \tau)  & \text{if } s= \mathrm{EL}, \\
        0 &\text{if } s = T^2.
    \end{cases}
    \]
    In the comparison score $B(s)$ above, $\phi$ denotes the standard-normal
    density.  Set $\tau=\EE_0[\D h(X)]\tau_0/\sqrt{\alpha_2}$ and
    $w_1=\sqrt{\chi^2_{1;1-\varrho}}-\tau$,
    $w_2=\sqrt{\chi^2_{1;1-\varrho}}+\tau$, and define
    $I:(0,1)\times\RR\to\RR$ by
    \[
    \begin{aligned}
        I(\varrho,\tau)
        &\Let \int_{|v+\tau|^2>\chi^2_{1;1-\varrho}}
        \left\{\frac{1+\tau^2}{2}v+\tau(v^2-1)
        +\frac12(v^3-3v)\right\}\phi(v)\diff v\\
        &=\frac{\chi^2_{1;1-\varrho}}{2}
        \{\phi(w_1)-\phi(w_2)\}.
    \end{aligned}
    \]
    Moreover, all three tests have the common limiting local power
    \[
        \lim_{n\to\infty}\pi_{n,s}
        =\int_{|v+\tau|^2>\chi^2_{1;1-\varrho}}\phi(v)\diff v.
    \]
\end{proposition}

The leading local power is common to all three tests and therefore does not
determine a preference among WP, EL, and Hotelling's $T^2$; their pairwise
power differences are governed by the second-order scores $B(s)-B(t)$.
Because $\PP_n\opt$ differs from $\PP_0\opt$ by an $O(n^{-1/2})$ location shift,
the relevant empirical averages are $O_p(n^{-1/2})$-close to their
$\PP_0\opt$ counterparts: both their sampling error and their shift bias have
that order.  Away from the nondegenerate denominators in
Proposition~\ref{prop:powercomparison}, plug-in therefore changes $B(s)$ by
$O_p(n^{-1/2})$ and the estimated $n^{-1/2}$ power gap by $O_p(n^{-1})$.
This statement concerns estimation of the comparison scores, not replacement
of a population critical value by $\hat z_{1-\varrho}$.  Thus, for a specified
$\tau_0$, the scores have a direct plug-in implementation; an unspecified
$\tau_0$ requires structural information or an external estimator.  A
data-driven ranking is stable away from $\tau_h=0$ and pairwise score ties.

Let $\tau_h=\EE_0[\D h(X)]\tau_0$ and define the sign-adjusted scores
\[
    C_{\mathrm{WP}}
    =\mathrm{sign}(\tau_h)\frac{\tilde\alpha_3}{\tilde\alpha_2^2},
    \qquad
    C_{\mathrm{EL}}
    =\mathrm{sign}(\tau_h)\frac{2\alpha_3}{3\alpha_2^2}.
\]
If $c=\sqrt{\chi^2_{1;1-\varrho}}$, then
\[
    \phi(c-\tau)-\phi(c+\tau)
    =2\phi(c)e^{-\tau^2/2}\sinh(c\tau),
    \]
so $I(\varrho,\tau)$ has the same sign as $\tau_h$.  The second-order
ranking is determined by Table~\ref{tab:power_comparison_scores}.  When
$\tau_h=0$, both scores vanish and all three tests agree to this order.

\begin{table}[H]
\centering
\small
\caption{Second-order local-power scores (larger is better; ties indicate equal power).}
\label{tab:power_comparison_scores}
\begin{tabular}{@{}lll@{}}
\toprule
Test & Second-order score & Preferred when \\
\midrule
WP
& $C_{\mathrm{WP}}$
& $C_{\mathrm{WP}}\geq\max\{C_{\mathrm{EL}},0\}$ \\
EL
& $C_{\mathrm{EL}}$
& $C_{\mathrm{EL}}\geq\max\{C_{\mathrm{WP}},0\}$ \\
Hotelling's $T^2$
& $0$
& $0\geq\max\{C_{\mathrm{WP}},C_{\mathrm{EL}}\}$ \\
\bottomrule
\end{tabular}
\end{table}

\begin{figure}[!t]
    \centering
    \includegraphics[width=0.90\textwidth]{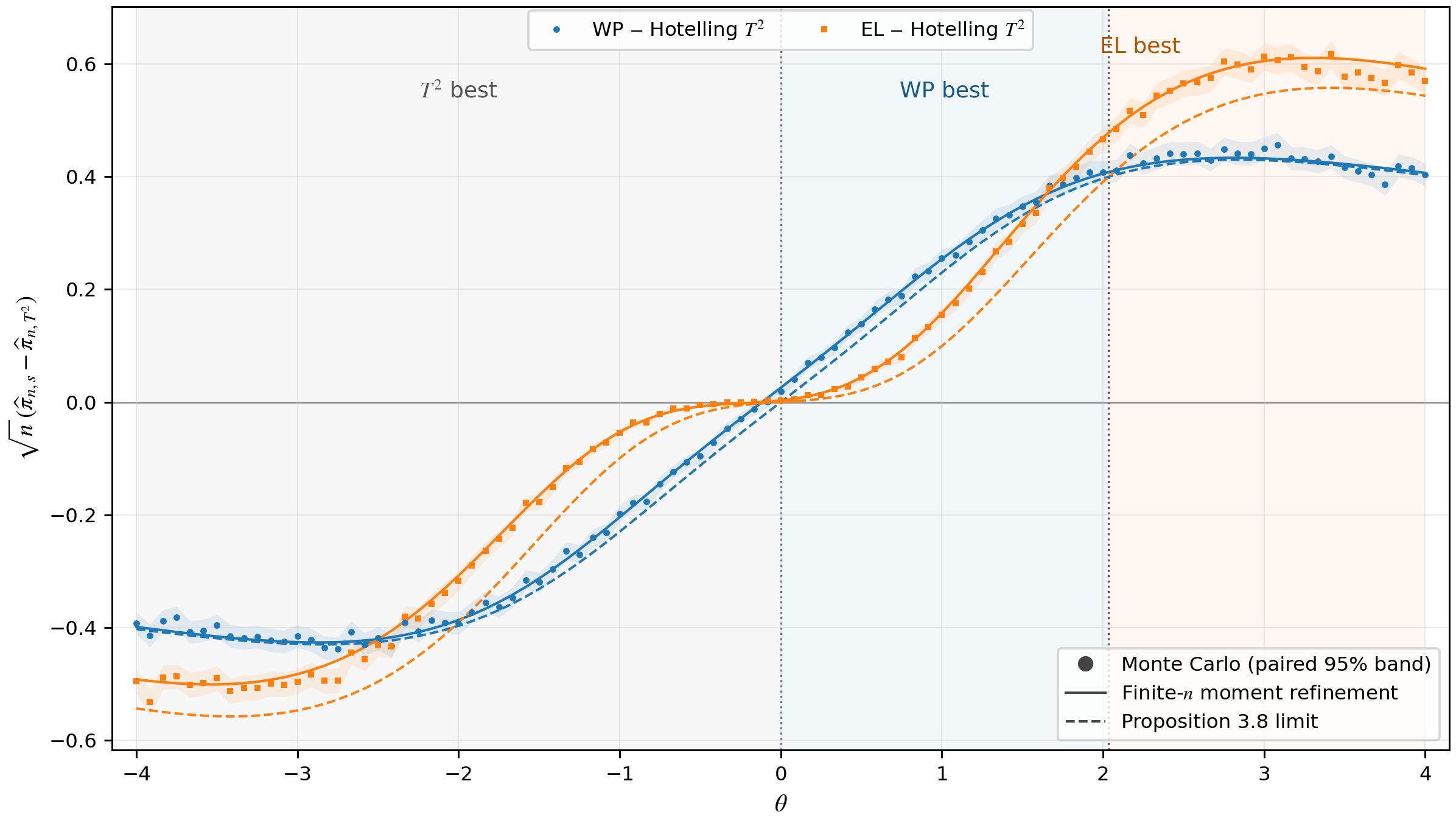}
    \par\vspace{-2mm}
    {\small (a) Second-order local-power comparison.}

    \vspace{2mm}
    \includegraphics[width=0.96\textwidth]{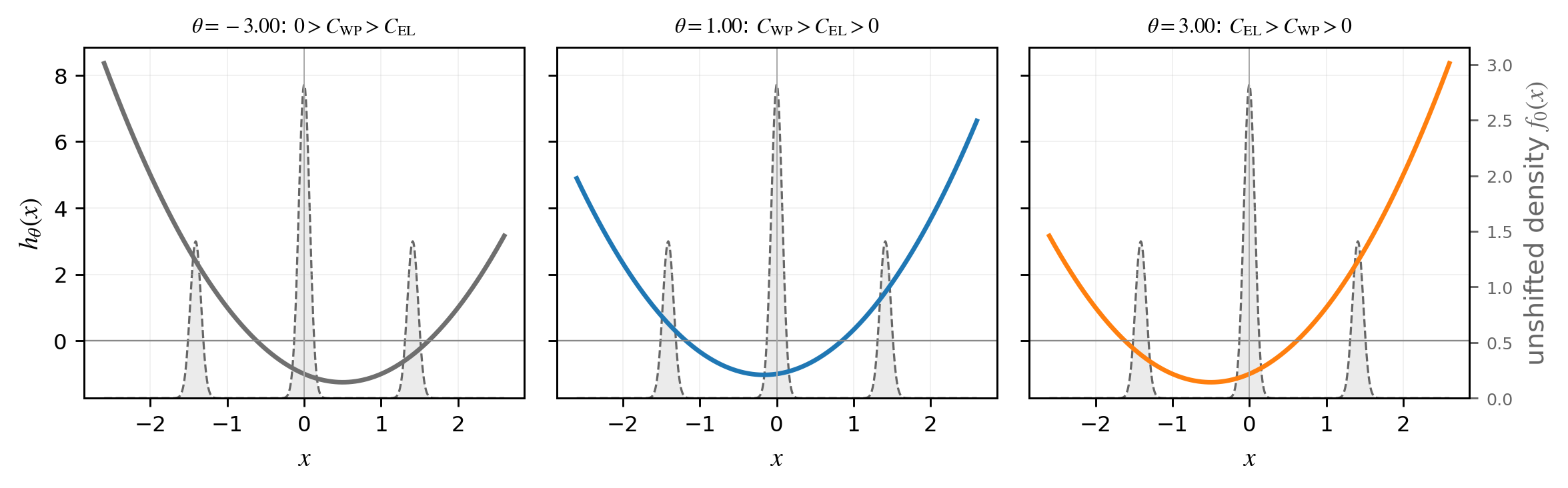}
    \par\vspace{-1mm}
    {\small (b) Moment shapes and unshifted null density.}
    \caption{Pairwise local-power comparison.  Panel~(a) plots
    $\sqrt n(\widehat\pi_{n,s}-\widehat\pi_{n,T^2})$ for $s\in\{\mathrm{WP},
    \mathrm{EL}\}$ using common random samples, $n=6000$, $\tau_0=3$, and
    $300{,}000$ replications.  Shading gives paired pointwise 95\% intervals.
    Dashed curves are the limiting scores $B(s)-B(T^2)$ from
    Proposition~\ref{prop:powercomparison}, solid curves are finite-$n$ moment
    refinements that account for the exact shifted moments, and the backgrounds
    mark the asymptotically predicted winner.  Panel~(b) overlays three representative moment
    functions with the unshifted null density (gray fill and dashed curve);
    each title gives the ordering of the corresponding comparison scores.
    Full details are in Supplementary Section~S2.2.}
    \label{fig:local_power_numeric}
\end{figure}

The WP score measures the gradient-weighted curvature $\tilde\alpha_3$
relative to the squared gradient energy $\tilde\alpha_2^2$, whereas the EL
score measures value skewness $\alpha_3$ relative to the squared moment
variance $\alpha_2^2$.  The selection rule compares their sign-adjusted ratios
$\tilde\alpha_3/\tilde\alpha_2^2$ and $2\alpha_3/(3\alpha_2^2)$, using zero as
the $T^2$ benchmark.
Figure~\ref{fig:local_power_numeric} makes the mechanism concrete under a
smooth symmetric three-component Gaussian-mixture null with
$h_\theta(x)=x^2+(\theta/3)x-1$ and $\tau_0=3$.  Here
$\tau_h=\theta$ and $h_\theta''\equiv2$.
For $\theta<0$, the induced moment drift is negative and both sign-adjusted
ratios are negative, so $T^2$ is preferred; at $\theta=0$, all scores vanish.  For
$0<\theta<\theta_\star$, the fixed curvature and the gradient energy
$\tilde\alpha_2=4+\theta^2/9$, which remains close to its minimum, make the
WP gradient-curvature ratio substantially larger; at $\theta=1$ the WP ratio
is more than twice the EL ratio.  As $\theta$ increases, the transformed
moment becomes more value-skewed; after the crossing at
$\theta_\star=2.03423$, the EL value-skewness ratio is larger (its standardized
skewness is $1.084$ at $\theta=3$), so EL is preferred.  All calculations and
simulation details are in Supplementary Section~S2.2.

\subsection{Higher-order asymptotic expansion and Bartlett-type correction}\label{sec:higher_expan_bartlett}

\subsubsection{Higher-order asymptotic expansion} \label{sec:higher_expan}

As a natural extension of the second-order expansion in Theorem~\ref{thm:main_expan},
this section expands the WP statistic through order $n^{-1}$, with an
$\widetilde{O}(n^{-3/2})$ remainder, using the proof strategy detailed in
Supplementary Sections~S2.9.6--S2.9.9 under suitable regularity conditions.

The expansion to this order and the resulting Edgeworth expansion refine~\eqref{eq:edgeworth_ninf} by identifying the $n^{-1}$ term. This result is closely related to the Bartlett-type correction~\cite{cribari1996bartlett}, which removes that term and therefore improves on the approximation in Theorem~\ref{thm:edgeworth1}.

Higher-order expansions require additional smoothness and moment conditions beyond those in Section~\ref{sec:prelim}. To focus on demonstrating the proof strategy and avoid cumbersome assumptions, we impose fourth-order smoothness together with global at-most-quadratic growth and exponential-moment control of the relevant envelopes.

\begin{assumption}[Global growth, smoothness, and exponential moments]\label{a:bddfuncs}
    We assume that $\mf h\in C^4(\RR^{d_X})$ and that the global
    Jacobian-Lipschitz condition in Assumption~\ref{a:D1} holds with envelope
    $\kappa_1$.  Define the fourth-order envelope, for $x\in\RR^{d_X}$, by
    \[
        \mc E_4(x) \Let
        \norm{\mf h(x)}_2+\kappa_1(x)
        + \sum_{j=1}^4 \sup_{\norm{u-x}_2\leq 1}\norm{\D^j\mf h(u)}_{\mathrm{op}} .
    \]
    We assume that there exists $\bar \eta>0$ such that
    \[
        \EE\left[\exp\left\{\bar \eta\mc E_4(X)\right\}\right] < \infty .
    \]
\end{assumption}

The global part of this assumption is essential for a squared transport cost:
it prevents a vanishing amount of mass from being sent to infinity to repair a
superquadratic moment at arbitrarily small cost.  The local derivative terms
provide the higher-order Taylor envelopes, while inclusion of $\kappa_1(X)$
supplies the concentration bounds used in the localization argument.

We next present the higher-order asymptotic expansion of $nR_n(\mf h)$.
\begin{theorem}[Higher-order expansion of {$nR_n(\mf h)$}]\label{thm:higherexp_rwp}
    Under Assumptions~\ref{a:convexhull}, \ref{a:cont_P}, \ref{a:basic}, and \ref{a:bddfuncs}, there is a deterministic sequence $\hat \delta_n = \widetilde{O}\left(n^{-\frac{3}{2}}\right)$ and an integer $N$ such that, when $n \geq N$,
    \begin{align}\label{eq:higher_expan}
        n R_n(\mf h) 
        =  \brak{ \mc V_n, \xi_n^{\otimes 2} }  + \frac{1}{
         \sqrt{n}} \brak{ \mc K_n, \xi_n^{\otimes 3}} + \frac{1}{n} \brak{ \mc L_n, \xi_n^{\otimes 4}} + \hat \eps_n.
    \end{align}
    The remainder satisfies
    $\PP\opt\pare{\abs{\hat \eps_n}\leq\hat \delta_n}=1-O(n^{-3/2})$.
    The explicit fourth-order tensor $\mc L_n$ is defined in
    Supplementary Section~S2.9.6.
\end{theorem}

With the expansion~\eqref{eq:higher_expan} in hand, we can apply the same analysis as in the discussion before Theorem~\ref{thm:edgeworth1}. Thus, the higher-order Edgeworth expansion and the related Bartlett-type correction are applicable, as we will demonstrate for the one-dimensional case $d_h = 1$ in the next section.

\subsubsection{Bartlett-type correction}\label{sec:bartlett}
In the WP-based test, we use an estimated quantile of the limiting distribution of
$nR_n(\mf h)$ as the rejection threshold for the null hypothesis $\mc H_0$.
As noted earlier, the resulting confidence level is correct only asymptotically,
which can produce size distortions \cite[Chapter 2]{cordeiro2014introduction}.
Expanding the actual confidence level in powers of $n^{-1}$ reveals its rate of
convergence to the target level.  Once the second-order term is identified, a
Bartlett correction \cite{bartlett1937properties} can adjust the statistic and
improve this rate.  In particular, under the null, the distribution of the
corrected statistic is generally approximated more accurately by the original
limiting distribution, thereby reducing small-sample size distortion.

To obtain an explicit second-order expansion term, in contrast to Theorem~\ref{thm:wp_coverage_accuracy}, we need to expand $nR_n(\mf h)$ (available in Theorem~\ref{thm:higherexp_rwp}) and the estimated quantile $\hat z_{1-\varrho}$ to higher orders. For simplicity, we present the result when $d_h = 1$, and thus, we write $\mf h$ as $h$. Note that in this case, $\hat z_{1-\varrho} = (\mc W_n/\mc V_n) \chi^2_{1;1-\varrho}$, for $\mc W_n, \mc V_n \in \RR$.

\begin{theorem}[Coverage accuracy of the WP test II]\label{thm:bartlettOT}
    Under Assumptions~\ref{a:convexhull}, \ref{a:cont_P}, \ref{a:basic}, and \ref{a:bddfuncs}, we have
    \begin{align*}
        \PP^\star\left(n R_n(h) \leq \hat z_{1-\varrho}\right) = 1-\varrho + \frac{1}{n}  g_{1}(\chi^2_{1;1-\varrho}) q(\chi^2_{1;1-\varrho}) + \widetilde{O}\left(n^{-\frac{3}{2}}\right),
    \end{align*}
    where $g_1(\cdot)$ and $\chi^2_{1;1-\varrho}$ are the density and
    $(1-\varrho)$-quantile of the chi-square distribution with one degree of
    freedom, and $q(t)=\sum_{k=1}^3 C_kt^k$.  The coefficients $C_k$ and their
    moment and cumulant definitions are collected in
    Supplementary Sections~S2.9.10 and~S4.3.
\end{theorem}

Theorem~\ref{thm:bartlettOT} implies that a simple correction can be applied to the quantile $\hat z_{1-\varrho}$ to remove the $n^{-1}$ term. This leads to the first correction presented in the following proposition.

\begin{proposition}[Bartlett-type correction I]\label{prop:BartlettI}
    Under Assumptions~\ref{a:convexhull}, \ref{a:cont_P}, \ref{a:basic}, and \ref{a:bddfuncs}, we have
    \begin{align*}
        \PP^\star\left(n R_n(h) \leq \pare{1 - C}\hat z_{1-\varrho}\right) = 1-\varrho + \widetilde{O}\left(n^{-\frac{3}{2}}\right),
    \end{align*}    
    where $C = \sum_{k = 1}^{3} C_{k} \pare{\chi^2_{1;1-\varrho}}^{k-1}/n$.
\end{proposition}

The correction in Proposition~\ref{prop:BartlettI} is generally level-dependent, since $C$ depends on the target quantile $\chi^2_{1;1-\varrho}$. The classical Bartlett correction \cite{bartlett1937properties}, by contrast, rescales the test statistic by a single factor that is the same for all confidence levels. This classical form is recovered from Proposition~\ref{prop:BartlettI} when $B_2 = B_3 = 0$, where $B_2$ and $B_3$ are defined in Supplementary Section~S2.9.10. In this case $C_2 = C_3 = 0$, so $C = C_1/n$ is independent of $\varrho$, and the correction can be written directly as
\begin{align}\label{eq:bartlett_correctionI}
    \PP^\star\left(\pare{1 - C}^{-1} n R_n(h) \leq \hat z_{1-\varrho}\right) = 1-\varrho + \widetilde{O}\left(n^{-\frac{3}{2}}\right).
\end{align}
This Bartlett correction to the test statistic uses the same constant $C$ for all confidence levels $1-\varrho \in (0,1)$, and is valid only when $B_2 = B_3 = 0$ (see the discussion in \cite[Section~2.3]{diciccio1991empirical}).

Nevertheless, even when the original Bartlett correction is invalid, another general Bartlett-type correction proposed by \cite{cordeiro1991modified} can be applied to the test statistic. Specifically, we adopt their construction to get
\begin{proposition}[Bartlett-type correction II]\label{prop:BartlettII}
    For $\mc V_n, \mc W_n \in \RR$ defined in~\eqref{eq:WnVn}, let \[S = \frac{\mc V_n}{\mc W_n} n R_n(h).\]
    Under Assumptions~\ref{a:convexhull}, \ref{a:cont_P}, \ref{a:basic}, and \ref{a:bddfuncs}, we have
    \begin{align}\label{eq:bartlett_correctionII}
        \PP\opt\pare{\pare{1 + \frac{1}{n} \sum_{k=1}^3 C_k S^{k-1}} nR_n(h) \leq \hat z_{1-\varrho}} = 1-\varrho + \widetilde{O}\left(n^{-\frac{3}{2}}\right).
    \end{align}
\end{proposition}

Note that the correction in Proposition~\ref{prop:BartlettII} involves the random variable $S$, in contrast to the original Bartlett correction in~\eqref{eq:bartlett_correctionI}; nevertheless, the same corrected statistic is used for all confidence levels $1-\varrho\in(0,1)$. In practice, the moments entering $C_k$, $k=1,2,3$, may be replaced by their empirical counterparts under $\QQ_n$. Under additional concentration conditions ensuring that the resulting estimators satisfy $\widehat C_k-C_k=\widetilde{O}_p(n^{-1/2})$, together with suitable local anti-concentration bounds, the induced perturbation of the corrected statistic is of order $\widetilde{O}_p(n^{-3/2})$. Accordingly, Propositions~\ref{prop:BartlettI}--\ref{prop:BartlettII} are stated for the population coefficients $C_k$, while the simulation below uses their empirical counterparts.

To conclude this section, we present a simulation to illustrate the improvement of coverage using the corrections in Propositions~\ref{prop:BartlettI} and~\ref{prop:BartlettII}.
\begin{example}[Gaussian quadratic moment function]
Consider the one-dimensional quadratic moment function
$
    h(x)=x^2-1, \PP^{\star} = \mc N(0,1),
$
with squared Euclidean ground cost and nominal coverage $1-\varrho=0.95$.
The normal design makes the centered moment $h(X)=X^2-1$ substantially skewed, so the first-order chi-square approximation exhibits a visible finite-sample coverage error. In each replication, every moment entering $C_1,C_2,C_3$ is estimated from the same sample used to compute the test statistic. Table~\ref{tab:bartlett_sim_quadratic_normal} reports Monte Carlo coverage probabilities based on $2\times 10^6$ replications; the numbers in parentheses are Monte Carlo standard errors.

\begin{table}[t]
\centering
\small
\caption{Coverage comparison with empirically estimated Bartlett coefficients for $h(x)=x^2-1$, $X\sim \mc N(0,1)$, nominal coverage $0.95$, and $2\times 10^6$ replications. Entries are percentages; parentheses report Monte Carlo standard errors.}
\label{tab:bartlett_sim_quadratic_normal}
\begin{tabular}{rrrr}
\toprule
$n$ & Uncorrected & Correction I & Correction II \\
\midrule
50 & $93.19\%\ (0.018\%)$ & $93.61\%\ (0.017\%)$ & $94.02\%\ (0.017\%)$ \\
100 & $93.96\%\ (0.017\%)$ & $94.45\%\ (0.016\%)$ & $94.59\%\ (0.016\%)$ \\
150 & $94.25\%\ (0.016\%)$ & $94.69\%\ (0.016\%)$ & $94.78\%\ (0.016\%)$ \\
200 & $94.41\%\ (0.016\%)$ & $94.80\%\ (0.016\%)$ & $94.86\%\ (0.016\%)$ \\
250 & $94.54\%\ (0.016\%)$ & $94.89\%\ (0.016\%)$ & $94.93\%\ (0.016\%)$ \\
\bottomrule
\end{tabular}
\end{table}

Both plug-in corrections reduce undercoverage at every reported sample size, with Correction II closer to nominal coverage throughout. The remaining undercoverage is most pronounced at $n=50$ and decreases as $n$ grows.
\end{example}

\section{Computation for WP tests}\label{sec:computation_algorithm}

We develop a certified computational scheme for the WP test.  At nominal significance level
$\varrho\in(0,1)$, the test rejects when $nR_n(\mf h)>\hat z_{1-\varrho}$.

A direct approach would compute the exact value of $R_n(\mf h)$ from its unscaled dual
representation derived in Supplementary Section~S2.7.2. It can be written compactly as
\begin{align}
    R_n(\mf h)
    =-\inf_{\zeta\in\RR^{d_h}}\frac{1}{n}\sum_{i=1}^n\sup_{x\in\RR^{d_X}}
    \left[\zeta^\top\mf h(x)-\norm{x-X_i}_{\Sigma}^2\right].
    \label{eq:unscaled_wp_dual_computation}
\end{align}
This is a convex--nonconcave minimax problem: the outer objective is convex in $\zeta$, whereas
the inner maximization need not be concave in $x$ for a nonlinear moment function $\mf h$.  For
the WP test, however, computing the exact projection distance is unnecessary.  Given the
rejection threshold $\delta=\hat z_{1-\varrho}/n$, it suffices to decide whether
$R_n(\mf h)>\delta$.  We first recast this comparison as a distributionally robust optimization
(DRO) problem over the Wasserstein ball centered at $\QQ_n$
(Proposition~\ref{prop:wp_dro_reformulation}).  We then show that, when $\delta$ falls below an
explicit observable threshold, the dual of the DRO problem localizes to a compact box on which
the inner minimizations are strongly convex and the maximization over the direction variable is
strongly concave (Proposition~\ref{prop:wp_localization}).  These two facts certify the nested
first-order scheme detailed as Algorithm~1 in Supplementary Section~S2.10.4, whose total cost is $\widetilde{O}(n)$ with probability $1-O(n^{-1})$ under the regularity and accuracy conditions stated below and whose statistical
behavior provably matches that of the oracle test up to $\widetilde{O}(n^{-1})$
(Theorem~\ref{thm:wp_computable_test}).  When higher-order accuracy is required, as with a
Bartlett correction, the same scheme can be adapted by tightening the optimization tolerance
while retaining its $\widetilde{O}(n)$ total cost with high probability; details are given in Supplementary Section~S2.10.7.  In short,
the proposed simplification reduces the computational burden without compromising the
statistical validity of the test.

We begin by reformulating the threshold comparison.  For $\delta>0$, let
\begin{align}
    \mc B_\delta(\QQ_n)
    \Let\left\{\PP\in\mc P(\RR^{d_X}):\Wass_c(\QQ_n,\PP)\leq\delta\right\},
    \qquad
    \mathbb B_2\Let\left\{\xi\in\RR^{d_h}:\norm{\xi}_2\leq1\right\},
    \label{eq:comp_ball_defs}
\end{align}
and define, for $\alpha\geq0$ and $\xi\in\RR^{d_h}$,
\begin{align}
    F_\delta(\alpha,\xi)
    \Let
    -\alpha\delta
    +\frac{1}{n}\sum_{i=1}^n\inf_{x\in\RR^{d_X}}
    \left[\xi^\top\mf h(x)+\alpha\norm{x-X_i}_{\Sigma}^{2}\right],
    \qquad
    \Phi_n(\delta)
    \Let
    \sup_{\substack{\alpha\geq0\\ \xi\in\mathbb B_2}}F_\delta(\alpha,\xi).
    \label{eq:def_F_delta}
\end{align}
The function $F_\delta$ is possibly extended-valued and jointly concave in $(\alpha,\xi)$,
because it is an average of pointwise infima of affine functions.  It is also positively
homogeneous along rays:
\begin{align}
    F_\delta(t\alpha,t\xi)=t\,F_\delta(\alpha,\xi)
    \qquad\text{for all }t>0.
    \label{eq:comp_F_homogeneity}
\end{align}
Since $F_\delta(0,\mf 0)=0$, we always have $\Phi_n(\delta)\geq0$.  The following proposition
shows that the sign of $\Phi_n(\delta)$ exactly determines the WP testing decision.

\begin{proposition}[DRO reformulation of the WP threshold comparison]
\label{prop:wp_dro_reformulation}
Suppose Assumptions~\ref{a:convexhull} and~\ref{a:D1} hold, and let $\delta>0$.  Then:
\begin{enumerate}[label=\textup{(\roman*)}]
    \item The following strong-duality identity holds:
    \begin{align}
        \Phi_n(\delta)
        =\inf_{\PP\in\mc B_\delta(\QQ_n)}\,\max_{\xi\in\mathbb B_2}\,
        \EE_{\PP}\!\left[\xi^\top\mf h(X)\right]
        =\inf_{\PP\in\mc B_\delta(\QQ_n)}\norm{\EE_{\PP}[\mf h(X)]}_2
        \geq 0,
        \label{eq:wp_dro_strong_duality}
    \end{align}
    \item The threshold comparison is equivalent to the sign of the dual value:
    \begin{align}
        R_n(\mf h)>\delta
        \ \Longleftrightarrow\
        \Phi_n(\delta)>0.
        \label{eq:wp_threshold_sign}
    \end{align}
\end{enumerate}
\end{proposition}

The proof is given in Supplementary Section~S2.10.1.  To turn this reformulation
into an algorithm, we next introduce the curvature and nondegeneracy quantities used to localize
the dual problem.  Assumption~\ref{a:D1} provides the Hessian envelope $\kappa_1$; we impose the
uniform strengthening
\begin{align}
    \overline\kappa_1
    \Let\sup_{x\in\RR^{d_X}}\kappa_1(x)<\infty,
    \qquad
    K_{\mathrm{curv}}\Let\norm{\Sigma}_2\,\overline\kappa_1,
    \label{eq:comp_uniform_hessian}
\end{align}
which implies the uniform directional Hessian bound
\begin{align}
    \sup_{\substack{x\in\RR^{d_X}\\ \norm{u}_2=1}}
    \norm{\Sigma^{1/2}\D^2\!\left(u^\top\mf h\right)(x)\Sigma^{1/2}}_2
    \leq K_{\mathrm{curv}}.
    \label{eq:comp_hessian_bound}
\end{align}
For nondegeneracy, recall the empirical Jacobian Gram matrix $\mc V_n$ from
Theorem~\ref{thm:main_expan} and define its extreme eigenvalues by
\begin{align}
    \mc V_n
    \Let
    \frac{1}{n}\sum_{i=1}^n\D\mf h(X_i)\Sigma\D\mf h(X_i)^\top,
    \qquad
    \ell_n\Let\lambda_{\min}(\mc V_n),
    \qquad
    U_n\Let\lambda_{\max}(\mc V_n),
    \label{eq:comp_empirical_eigenvalues}
\end{align}
and write $A(x)\Let\D\mf h(x)\Sigma^{\half}\in\RR^{d_h\times d_X}$, so that
$\mc V_n=\frac1n\sum_{i=1}^nA(X_i)A(X_i)^\top$.  The deterministic results below are stated on
the observable event $\ell_n>0$, whose probability approaches one under
Assumptions~\ref{a:basic} and~\ref{a:finite_moments}.

Although $\Phi_n(\delta)$ is a concave maximization, its multiplier range $\alpha\geq0$ is
unbounded.  The next proposition localizes the maximizers to a compact box.  On the event
$\ell_n>0$, define
\begin{align}
\begin{aligned}
    \varkappa_n&\Let\frac{U_n}{\ell_n}\geq1,
    &\underline\alpha_\delta&\Let
    \frac12\pare{\sqrt{\frac{\ell_n}{\delta}}-K_{\mathrm{curv}}},\\
    \overline\alpha_\delta&\Let
    \frac12\pare{\sqrt{\frac{U_n}{\delta}}+K_{\mathrm{curv}}},
    &\lambda(\alpha)&\Let2\alpha-K_{\mathrm{curv}},\\
    \delta_0&\Let\frac{\ell_n}{K_{\mathrm{curv}}^2\pare{3+2\sqrt{\varkappa_n}}^2},
    &\mc K_\delta&\Let
    \brac{\underline\alpha_\delta,\overline\alpha_\delta}\times\mathbb B_2,
\end{aligned}
\label{eq:comp_box}
\end{align}
with $\delta_0\Let+\infty$ when $K_{\mathrm{curv}}=0$.
The box endpoints scale as $\delta^{-1/2}$: in the testing regime
$\delta=\hat z_{1-\varrho}/n$, the relevant transport multipliers $\alpha$ are
bracketed between constant multiples of
$\sqrt{n\ell_n/\hat z_{1-\varrho}}$ and
$\sqrt{nU_n/\hat z_{1-\varrho}}$.  In particular, when the empirical condition
number $\varkappa_n=U_n/\ell_n$ remains bounded, they are of exact order
$\sqrt{n\ell_n/\hat z_{1-\varrho}}$, a fact that drives both the conditioning of
the dual problem and the accuracy accounting below.

\begin{proposition}[Localization and curvature of the dual problem]
\label{prop:wp_localization}
Suppose Assumptions~\ref{a:convexhull} and \ref{a:D1} hold,
\eqref{eq:comp_uniform_hessian} is satisfied, and $\ell_n>0$.  Let $0<\delta\leq\delta_0$, and
define the localized dual value
\begin{align}
    \Phi^{\mathrm{loc}}_n(\delta)
    \Let
    \max_{(\alpha,\xi)\in\mc K_\delta}F_\delta(\alpha,\xi).
    \label{eq:def_Phi_loc}
\end{align}
Then:
\begin{enumerate}[label=\textup{(\roman*)}]
    \item \textup{(Localization.)}  The supremum defining $\Phi_n(\delta)$ is attained, and
    \begin{align}
        \Phi_n(\delta)=\max\left\{\Phi^{\mathrm{loc}}_n(\delta),\,0\right\}.
        \label{eq:wp_localization_identity}
    \end{align}
    If $\Phi_n(\delta)>0$, then every maximizer $(\alpha^\star,\xi^\star)$ of $\Phi_n(\delta)$
    satisfies $\norm{\xi^\star}_2=1$ and $(\alpha^\star,\xi^\star)\in\mc K_\delta$, so that
    $\Phi^{\mathrm{loc}}_n(\delta)=\Phi_n(\delta)$ in this case.  Consequently,
    \begin{align}
        R_n(\mf h)>\delta
        \quad\Longleftrightarrow\quad
        \Phi^{\mathrm{loc}}_n(\delta)>0.
        \label{eq:comp_box_decision}
    \end{align}
    \item \textup{(Strong convexity in $x$.)}  For every $\alpha\geq\underline\alpha_\delta$,
    $\xi\in\mathbb B_2$, and $i\in[n]$, the inner objective
    $x\mapsto\xi^\top\mf h(x)+\alpha\norm{x-X_i}_{\Sigma}^{2}$ is $\lambda(\alpha)$-strongly
    convex and $(2\alpha+K_{\mathrm{curv}})$-smooth with respect to $\norm{\cdot}_{\Sigma}$, with
    \begin{align}
        \lambda(\alpha)
        \geq\sqrt{\frac{\ell_n}{\delta}}-2K_{\mathrm{curv}}
        \geq\frac35\sqrt{\frac{\ell_n}{\delta}}>0 ;
        \label{eq:wp_lambda_lower_bound}
    \end{align}
    \item \textup{(Strong concavity in $\xi$.)}  For every $\alpha\geq\underline\alpha_\delta$,
    the function $\xi\mapsto F_\delta(\alpha,\xi)$ is $\mu(\alpha)$-strongly concave on
    $\mathbb B_2$, where
    \begin{align}
        \mu(\alpha)
        \Let
        \frac{\left(\sqrt{\ell_n}
        -\dfrac{K_{\mathrm{curv}}\sqrt{U_n}}{\lambda(\alpha)}\right)^{2}}{2\alpha+K_{\mathrm{curv}}}
        \geq
        \frac{\ell_n}{\left(1+\sqrt{\varkappa_n}\right)^{2}(2\alpha+K_{\mathrm{curv}})}
        >0.
        \label{eq:wp_xi_strong_concavity}
    \end{align}
\end{enumerate}
\end{proposition}

The proof is given in Supplementary Section~S2.10.2;
Assumption~\ref{a:convexhull} enters only through part~(i), where it rules out maximizers at
$\alpha=0$.  The localization identity~\eqref{eq:wp_localization_identity} is also important
computationally.  Because $\Phi_n(\delta)$ is nonnegative, the entire event
$\{R_n(\mf h)\leq\delta\}$ is mapped to $\Phi_n(\delta)=0$; thus, under the null, the
distribution of $\Phi_n(\delta)$ has an atom at zero corresponding to the non-rejection event.
This atom makes a direct numerical comparison of $\Phi_n(\delta)$ with zero delicate.

The localized value $\Phi_n^{\mathrm{loc}}(\delta)$ retains the sign information lost through
this truncation.  If $R_n(\mf h)<\delta$, choose a null distribution $\PP_0$ satisfying
$\EE_{\PP_0}[\mf h(X)]=\mf0$ and a coupling with $\QQ_n$ whose transportation cost $r$ is
less than $\delta$.  For every
$(\alpha,\xi)\in\mc K_\delta$, weak duality gives
\[
    F_\delta(\alpha,\xi)
    \leq-\alpha\delta+\xi^\top\EE_{\PP_0}[\mf h(X)]+\alpha r
    =-\alpha(\delta-r)
    \leq-\underline\alpha_\delta(\delta-r)<0,
\]
and hence $\Phi_n^{\mathrm{loc}}(\delta)<0$.  Conversely,
\eqref{eq:comp_box_decision} gives $\Phi_n^{\mathrm{loc}}(\delta)>0$ whenever
$R_n(\mf h)>\delta$.  Consequently,
\[
    \{\Phi_n^{\mathrm{loc}}(\delta)=0\}
    \subseteq
    \{R_n(\mf h)=\delta\}.
\]
Under the standard no-ties condition $\PP^\star(R_n(\mf h)=\delta)=0$, the localized value
therefore has no atom at zero, making it the more useful quantity for computation.

\subsection{A certified nested algorithm}\label{sec:comp_nested}

The homogeneity in~\eqref{eq:comp_F_homogeneity} prevents joint strong concavity on
full-dimensional domains, so the
algorithm uses the partial curvature in Proposition~\ref{prop:wp_localization}.  It searches
over the scalar multiplier $\alpha$, applies projected gradient ascent in $\xi$, and solves the
$n$ strongly convex inner problems in parallel.  On $\mc K_\delta$, the resulting
strong-concavity-to-smoothness condition number in $\xi$ can be taken as
$30\varkappa_n^{3/2}$ and therefore has no explicit dependence on $n$ when the empirical condition
number $\varkappa_n=U_n/\ell_n$ remains bounded.  Explicit moduli, data-scale constants, stopping
rules, and the complete complexity ledger are given in Supplementary
Sections~S2.10.3--S2.10.6.

\begin{theorem}[Computable WP test]\label{thm:wp_computable_test}
Suppose Assumptions~\ref{a:convexhull} and~\ref{a:D1} hold and
\eqref{eq:comp_uniform_hessian} is satisfied.  Set
$\delta=\hat z_{1-\varrho}/n$ and let $\eps_{\mathrm{alg},n}>0$.  On the observable event
\begin{align}
    \mathfrak E_n:\qquad
    \ell_n>0,\qquad
    \delta\leq\delta_0/2,\qquad
    \eps_{\mathrm{alg},n}\leq\sqrt{\ell_n\delta}/8,
    \label{eq:comp_event}
\end{align}
run Algorithm~1 in Supplementary Section~S2.10.4 and let $\hat\Phi_n$ be its
returned value.  Define the computable test $\hat{\mathsf T}_n$ to reject
when $\mathfrak E_n$ holds and $\hat\Phi_n>\eps_{\mathrm{alg},n}$, and to fail to reject otherwise.  Define
$w_n=5\eps_{\mathrm{alg},n}\sqrt{2n\hat z_{1-\varrho}/\ell_n}$ on $\mathfrak E_n$ and set $w_n=0$ on
$\mathfrak E_n^c$.  Then:
\textup{(i)} On $\mathfrak E_n$,
$|\hat\Phi_n-\Phi_n^{\mathrm{loc}}(\delta)|\leq\eps_{\mathrm{alg},n}$.
\textup{(ii)} Deterministically,
\begin{align}
    \{\text{computable test rejects}\}
    &\subseteq \{nR_n(\mf h)>\hat z_{1-\varrho}\}, \notag\\
    \{nR_n(\mf h)>\hat z_{1-\varrho}\}\cap\mathfrak E_n
    &\subseteq \{\text{computable test rejects}\}\notag\\
    &\quad\cup\bigl(\{0<nR_n(\mf h)-\hat z_{1-\varrho}\leq w_n\}
       \cap\mathfrak E_n\bigr).
       \label{eq:comp_decision_band}
\end{align}
Equivalently, the second inclusion holds globally after adjoining $\mathfrak E_n^c$ to its
right-hand side.  Thus failure of the observable certification event is explicitly retained
in the decision comparison.
\textup{(iii)} Suppose additionally that $d_h$ is fixed,
Assumptions~\ref{a:convexhull}--\ref{a:moments_II} hold, and
$\eps_{\mathrm{alg},n}^{-1}=O(n^c)$ for some fixed $c>0$.  There exist events
$\mathfrak C_n$ with $\PP^\star(\mathfrak C_n^c)=O(n^{-1})$ such that, on
$\mathfrak E_n\cap\mathfrak C_n$, the algorithm uses
$O(n\log^3 n)=\widetilde{O}(n)$ first-order oracle calls, where one call jointly evaluates
$(\mf h,\D\mf h)$ at a single point.  For the accuracy specified in part~\textup{(iv)},
$\PP^\star(\mathfrak E_n\cap\mathfrak C_n)=1-O(n^{-1})$.
\textup{(iv)} For the statistical conclusions, take
$\eps_{\mathrm{alg},n}=1/(n^{3/2}\log^2n)$.  Under
Assumptions~\ref{a:convexhull}--\ref{a:moments_II} and
\eqref{eq:comp_uniform_hessian}, the computable test inherits the oracle test's size
expansion up to $\widetilde{O}(n^{-1})$.  For each fixed true parameter value at
which these assumptions hold, its inverted confidence region likewise inherits
the oracle pointwise coverage expansion up to $\widetilde{O}(n^{-1})$.  Under the assumptions of
Theorem~\ref{thm:wp_power_expansion}, together with
\eqref{eq:comp_uniform_hessian}, it likewise inherits the oracle local-power expansion up to
    $\widetilde{O}(n^{-1})$.  Finally, when $d_h=1$, under the assumptions of
    Proposition~\ref{prop:BartlettI}, together with
    \eqref{eq:comp_uniform_hessian}, taking
$\eps_{\mathrm{alg},n}=1/(n^2\log^2n)$ preserves the corrected test's
$\widetilde{O}(n^{-3/2})$ coverage accuracy.  Under the same assumptions, the
corrected algorithm also uses $O(n\log^3n)$ first-order oracle calls with
probability $1-O(n^{-3/2})$.
\end{theorem}

Thus computation is simplified without altering the statistical guarantees at the order
studied in the paper.  The factor $\varkappa_n^{3/2}$ also has a statistical interpretation:
when $\ell_n$ is nearly zero, both the optimization problem and the limiting WP law are nearly
degenerate.  The proof, complete complexity ledger, stopping rules, and pseudocode are given
in Supplementary Sections~S2.10.4 and~S2.10.6.

Supplementary Section~S2.3 provides a complementary finite-dimensional stress test in
which the data dimension, number of moment restrictions, and sample size all vary.  Across
the tested designs, the literal certified algorithm maintains null rejection rates close to
the nominal level, local power close to its fixed-dimensional limit, and certification
frequencies of at least \(0.996\); its decisions also closely track the leading quadratic
WP statistic.  These results support both the finite-sample calibration predicted by the
theory and the reliability of the certified implementation beyond the scalar experiment below.

\subsection{Numerical experiment: fairness testing}\label{sec:comp_fairness}

We now implement the certified nested procedure of Section~\ref{sec:comp_nested}
(Algorithm~1 in the Supplementary Material) in a nonlinear fairness test.  Let
$A\in\{-1,1\}$ be a protected attribute and $X\in\RR$ a score.  For $k>0$, the
null requires parity of the groups' mean transformed scores:
\[
    \EE[\tanh(kX)\mid A=1]=\EE[\tanh(kX)\mid A=-1].
\]
The bounded monotone transformation interpolates between smooth mean-score
parity for small $k$ and parity of positive decisions near the threshold zero
for large $k$.

We consider a prospective balanced audit of a fixed scoring system.  Before
observing scores, the auditor draws $n$ cases independently from each of two
fixed group-specific sampling frames and independently randomizes the order
within each group.  Combining observations with the same index gives
\[
    Z_j=(X_{j,+},X_{j,-}),\qquad j=1,\ldots,n.
\]
This pairing is only bookkeeping and does not assert covariate matching.  The
vectors $Z_j$ are i.i.d. from $P_+\otimes P_-$, and the group roles are fixed
coordinate labels.  Thus the sample consists of $n$ i.i.d. paired vectors.  Define
\[
\begin{aligned}
    \widetilde h_k(Z_j)
       &=\frac{\tanh(kX_{j,+})-\tanh(kX_{j,-})}{2k},\\
    \widetilde c_{\mathrm{aniso}}(z,z')
       &=\frac12(z_+-z'_+)^2+5(z_--z'_-)^2,\\
    \widetilde c_{\mathrm{id}}(z,z')
       &=(z_+-z'_+)^2+(z_--z'_-)^2.
\end{aligned}
\]
The factor $1/k$ does not change the null or the WP projection, but prevents
the derivative scale from collapsing as $k$ decreases.  To evaluate the
geometry-dependent second-order prediction of
Proposition~\ref{prop:powercomparison}, we compare
$\Sigma_{\mathrm{aniso}}=\operatorname{diag}(2,0.2)$ with
$\Sigma_{\mathrm{id}}=I_2$.  Both procedures use the same moment restriction,
samples, and asymptotic chi-square reference quantile and differ only in
transport geometry; the fixed group-coordinate labels cannot be exchanged.
The identity geometry is
therefore a matched control, while the anisotropic geometry is chosen because
the theory predicts a larger favorable second-order correction without a
change in first-order power.  The two matrices are thus controlled design
choices for evaluating the theoretical prediction.  Both costs are of the
Mahalanobis form in~\eqref{eq:mahala_norm}.  Under either geometry,
$R_n(\widetilde h_k)$ is the least average quadratic score movement needed to
attain transformed-score parity, so the theory applies with i.i.d. sample size
$n$.

The redesigned experiment uses $k=0.5$ and
$n\in\{75,100,125,150,250,500\}$ paired observations.  Under the null, the coordinates
are independent with
\begin{equation}
\begin{aligned}
  Y_+&\sim0.338\mc N(-7.60,0.1^2)
          +0.662\mc N(-1.319,0.1^2),\\
  Y_-&\sim \mc N(-1.817,0.1^2).
\end{aligned}
  \label{eq:fairness_certified_null}
\end{equation}
The displayed parameters are rounded.  The negative-group location was
obtained by Gaussian quadrature so that
$\EE[\tanh(0.5Y_+)]=\EE[\tanh(0.5Y_-)]\approx-0.720$.  The far-left mixture
component lies in the nearly flat lower tail of the transformation, where
$\tanh(0.5y)$ is close to $-1$, whereas the other component lies near the point
of maximal positive curvature.  Numerical quadrature gives
the second-order comparison scores
\[
    C_{\mathrm{WP,aniso}}=2.216,\qquad
    C_{\mathrm{WP,id}}=0.228,\qquad
    C_{\mathrm{EL}}=-2.018,\qquad C_{T^2}=0,
\]
so Proposition~\ref{prop:powercomparison} predicts a pronounced advantage for
anisotropic WP and only a small advantage for identity WP over $T^2$.
Under the local alternative,
\begin{equation}
  X_+=Y_++\frac{\tau}{2\sqrt n},\qquad
  X_-=Y_--\frac{\tau}{2\sqrt n},\qquad
  \tau\approx1.759.
  \label{eq:fairness_certified_alt}
\end{equation}
Writing
$\alpha_2=\EE[\widetilde h_{0.5}(Y_+,Y_-)^2]$, this value makes the
standardized first-order drift
\[
 \frac{\tau\{\EE[\operatorname{sech}^2(0.5Y_+)]
       +\EE[\operatorname{sech}^2(0.5Y_-)]\}}
      {4\sqrt{\alpha_2}}=2.
\]
Consequently the experiment targets common first-order limiting power
$\PP(|Z+2|>1.96)=0.516$ and uses curvature, rather than different first-order
signal strengths, to distinguish the methods.

To complement this local comparison with a fixed-alternative consistency
check, we also use
\begin{equation}
  X_+=Y_++\frac{\Delta}{2},\qquad
  X_-=Y_--\frac{\Delta}{2},\qquad
  \Delta=\frac{\tau}{\sqrt{150}}\approx0.144.
  \label{eq:fairness_certified_fixed_alt}
\end{equation}
Thus the fixed and local alternatives coincide at $n=150$; the fixed
alternative is weaker when $n<150$ and stronger when $n>150$.  At each $n$,
the two alternatives use the same underlying paired draws $(Y_+,Y_-)$.

For the certified test, let $\widetilde h_j=\widetilde h_{0.5}(Z_j)$ and define
\[
 \widehat W=\frac1n\sum_{j=1}^n\widetilde h_j^2,\qquad
 \widehat V_{\Sigma}=\frac1{4n}\sum_{j=1}^n
 \{\sigma_+\operatorname{sech}^4(0.5X_{j,+})
   +\sigma_-\operatorname{sech}^4(0.5X_{j,-})\},
\]
where $(\sigma_+,\sigma_-)=(2,0.2)$ for anisotropic WP and $(1,1)$
for identity WP.  The resulting plug-in radius is
$\delta=\widehat W\chi^2_{1;0.95}/(n\widehat V_{\Sigma})$.  The required
derivatives are available in closed form, and
the global curvature bound
$K_{\mathrm{curv}}=2k\max(\sigma_+,\sigma_-)/(3\sqrt3)$ makes the certification
event $\mathfrak E_n$ fully observable.  It held in $17{,}953$ of $18{,}000$
anisotropic-WP replications and all $18{,}000$ identity-WP replications;
uncertified runs were conservatively counted as non-rejections.

Table~\ref{tab:fairness_power_redesigned} reports the common-reference power
comparison.  Each design uses $1000$ common samples under the null and each
alternative.  All four methods use the common asymptotic reference quantile
$\chi^2_{1;0.95}=3.841$; WP converts it to the appropriate sample-specific
radius and uses the literal certified algorithm with
$\eps_{\mathrm{alg},n}=1/(n^{3/2}\log^2n)$, while EL-infeasible samples count as rejections.
This common reference directly evaluates the finite-sample ordering predicted
by Proposition~\ref{prop:powercomparison}, with the observed null rejection
rates reporting any finite-sample calibration differences.

Panel~A treats the local alternative and Panel~B the fixed alternative.  Null
rejection rates are similar across methods, while anisotropic WP has the
highest estimated power throughout.  For $n=75$ through $150$ in Panel~A,
its power exceeds EL by 3.9--5.0 percentage points, $T^2$ by 2.9--4.0 points,
and identity WP by 2.6--3.5 points; identity WP exceeds $T^2$ by only
0.1--0.6 points across all six sizes.  The local gaps diminish with $n$, as predicted by the
order-$n^{-1/2}$ correction, whereas fixed-alternative powers approach one.
The two panels coincide at $n=150$, providing an implementation check.  These
design-specific results do not imply general dominance; Figure~\ref{fig:local_power_numeric}
shows geometries that instead favor EL or $T^2$.

\begin{table}[t]
\centering
\scriptsize
\caption{Certified common-reference fairness comparisons. $\mathrm{WP}_{\mathrm{aniso}}$ uses $\Sigma=\operatorname{diag}(2,0.2)$ and $\mathrm{WP}_{\mathrm{id}}$ uses $\Sigma=I_2$; both use literal Algorithm~1. Panel~A uses local separation $\tau/\sqrt n$, $\tau=1.759$, and Panel~B uses fixed separation $\tau/\sqrt{150}=0.144$, so the panels coincide at $n=150$.  Each entry uses $B=1000$ common samples; parentheses give Monte Carlo standard errors.}
\label{tab:fairness_power_redesigned}
\begin{minipage}[t]{0.49\textwidth}
\centering
\textit{Panel A: Local alternative}\par\smallskip
\begin{tabular}{@{}rlcc@{}}
\toprule
$n$ & Method & Null rejection & Power\\
\midrule
75 & $\mathrm{WP}_{\mathrm{aniso}}$ & 0.038 (0.006) & \textbf{0.499 (0.016)}\\
 & $\mathrm{WP}_{\mathrm{id}}$ & 0.045 (0.007) & 0.466 (0.016)\\
 & EL & 0.052 (0.007) & 0.456 (0.016)\\
 & $T^2$ & 0.048 (0.007) & 0.461 (0.016)\\
\addlinespace[2pt]
100 & $\mathrm{WP}_{\mathrm{aniso}}$ & 0.037 (0.006) & \textbf{0.499 (0.016)}\\
 & $\mathrm{WP}_{\mathrm{id}}$ & 0.044 (0.006) & 0.464 (0.016)\\
 & EL & 0.042 (0.006) & 0.449 (0.016)\\
 & $T^2$ & 0.044 (0.006) & 0.459 (0.016)\\
\addlinespace[2pt]
125 & $\mathrm{WP}_{\mathrm{aniso}}$ & 0.048 (0.007) & \textbf{0.480 (0.016)}\\
 & $\mathrm{WP}_{\mathrm{id}}$ & 0.050 (0.007) & 0.454 (0.016)\\
 & EL & 0.052 (0.007) & 0.441 (0.016)\\
 & $T^2$ & 0.051 (0.007) & 0.449 (0.016)\\
\addlinespace[2pt]
150 & $\mathrm{WP}_{\mathrm{aniso}}$ & 0.057 (0.007) & \textbf{0.493 (0.016)}\\
 & $\mathrm{WP}_{\mathrm{id}}$ & 0.055 (0.007) & 0.466 (0.016)\\
 & EL & 0.048 (0.007) & 0.452 (0.016)\\
 & $T^2$ & 0.055 (0.007) & 0.464 (0.016)\\
\addlinespace[2pt]
250 & $\mathrm{WP}_{\mathrm{aniso}}$ & 0.054 (0.007) & \textbf{0.481 (0.016)}\\
 & $\mathrm{WP}_{\mathrm{id}}$ & 0.053 (0.007) & 0.469 (0.016)\\
 & EL & 0.057 (0.007) & 0.452 (0.016)\\
 & $T^2$ & 0.053 (0.007) & 0.463 (0.016)\\
\addlinespace[2pt]
500 & $\mathrm{WP}_{\mathrm{aniso}}$ & 0.051 (0.007) & \textbf{0.527 (0.016)}\\
 & $\mathrm{WP}_{\mathrm{id}}$ & 0.051 (0.007) & 0.522 (0.016)\\
 & EL & 0.051 (0.007) & 0.509 (0.016)\\
 & $T^2$ & 0.051 (0.007) & 0.521 (0.016)\\
\bottomrule
\end{tabular}
\end{minipage}\hfill
\begin{minipage}[t]{0.49\textwidth}
\centering
\textit{Panel B: Fixed alternative}\par\smallskip
\begin{tabular}{@{}rlcc@{}}
\toprule
$n$ & Method & Null rejection & Power\\
\midrule
75 & $\mathrm{WP}_{\mathrm{aniso}}$ & 0.038 (0.006) & \textbf{0.308 (0.015)}\\
 & $\mathrm{WP}_{\mathrm{id}}$ & 0.045 (0.007) & 0.279 (0.014)\\
 & EL & 0.052 (0.007) & 0.258 (0.014)\\
 & $T^2$ & 0.048 (0.007) & 0.275 (0.014)\\
\addlinespace[2pt]
100 & $\mathrm{WP}_{\mathrm{aniso}}$ & 0.037 (0.006) & \textbf{0.371 (0.015)}\\
 & $\mathrm{WP}_{\mathrm{id}}$ & 0.044 (0.006) & 0.350 (0.015)\\
 & EL & 0.042 (0.006) & 0.331 (0.015)\\
 & $T^2$ & 0.044 (0.006) & 0.347 (0.015)\\
\addlinespace[2pt]
125 & $\mathrm{WP}_{\mathrm{aniso}}$ & 0.048 (0.007) & \textbf{0.425 (0.016)}\\
 & $\mathrm{WP}_{\mathrm{id}}$ & 0.050 (0.007) & 0.392 (0.015)\\
 & EL & 0.052 (0.007) & 0.372 (0.015)\\
 & $T^2$ & 0.051 (0.007) & 0.386 (0.015)\\
\addlinespace[2pt]
150 & $\mathrm{WP}_{\mathrm{aniso}}$ & 0.057 (0.007) & \textbf{0.493 (0.016)}\\
 & $\mathrm{WP}_{\mathrm{id}}$ & 0.055 (0.007) & 0.466 (0.016)\\
 & EL & 0.048 (0.007) & 0.452 (0.016)\\
 & $T^2$ & 0.055 (0.007) & 0.464 (0.016)\\
\addlinespace[2pt]
250 & $\mathrm{WP}_{\mathrm{aniso}}$ & 0.054 (0.007) & \textbf{0.686 (0.015)}\\
 & $\mathrm{WP}_{\mathrm{id}}$ & 0.053 (0.007) & 0.673 (0.015)\\
 & EL & 0.057 (0.007) & 0.661 (0.015)\\
 & $T^2$ & 0.053 (0.007) & 0.671 (0.015)\\
\addlinespace[2pt]
500 & $\mathrm{WP}_{\mathrm{aniso}}$ & 0.051 (0.007) & \textbf{0.932 (0.008)}\\
 & $\mathrm{WP}_{\mathrm{id}}$ & 0.051 (0.007) & 0.928 (0.008)\\
 & EL & 0.051 (0.007) & 0.925 (0.008)\\
 & $T^2$ & 0.051 (0.007) & 0.928 (0.008)\\
\bottomrule
\end{tabular}
\end{minipage}
\end{table}


\section{Discussion}\label{sec:discuss}

The results reveal two roles for transport geometry: locally, it enters the
higher-order distribution and power expansions through derivatives of the
moment function; computationally, it yields the localized nested formulation
of Section~\ref{sec:computation_algorithm}.  We next relate WP to empirical
likelihood, discuss nonlocal power, and outline extensions beyond i.i.d. sampling.

\subsection{WP and EL: geometry versus reweighting}\label{sec:ELvsWP}
Empirical likelihood (EL) also projects onto the null: it minimizes
$D_{\mathrm{KL}}(\QQ_n\Vert\PP)$ over laws supported on the observed atoms,
equivalently maximizing empirical likelihood.  With the normalization in
\eqref{def:EL},
\[
    R_n^{\mathrm{EL}}(\mf h)
    =2\inf_{\substack{\PP\in\mc P_0\\ \PP\ll\QQ_n}}
       D_{\mathrm{KL}}(\QQ_n\Vert\PP)
    =-\frac{2}{n}\log\mc R;
\]
see \cite[Chapters~2--3]{owen2001empirical}.  The local dual expansions have
the same algebraic form but different tensors: EL uses the observed values of
$\mf h$ and generates $W$, whereas WP uses derivatives of $\mf h$ and the
ground-cost matrix $\Sigma$ and generates $V$.

At the next order, the curvature scores $\tilde\alpha_3/\tilde\alpha_2^2$ for
WP and $2\alpha_3/(3\alpha_2^2)$ for EL determine their relative power; see
Proposition~\ref{prop:powercomparison}.  The complete dual expansion, its
quadratic through quartic coefficients, and the consistency check against
\cite[Theorem~3.1]{chen1994comparing} appear in Supplementary Section~S2.4.1.

The interpolation proposed by \cite{blanchet2024stability} suggests a broader
family that tunes the balance between geometric and value-based information;
its higher-order analysis is a natural direction for future work.

\subsection{From local to nonlocal power}\label{sec:non-local_power}
The local comparison and fairness experiment reflect the same mechanism: EL
and Hotelling's $T^2$ use only $\{\mf h(X_i)\}$, whereas WP also uses the
sample-space geometry.  This distinction persists nonlocally.  Supplementary
Sections~S2.4.1--S2.4.2 construct smooth alternatives along which the population
EL and $T^2$ discrepancies converge to zero while the WP discrepancy stays bounded away
from zero.  All three tests are pointwise consistent, but only WP retains
uniform separation along the sequence.

\subsection{Extensions beyond i.i.d. sampling and computational scope}
The deterministic expansion and optimization do not intrinsically require
independence, but inferential extensions need a joint Edgeworth expansion and
high-probability control of the regularity event; relevant tools include
\cite{rinott2003edgeworth,skovgaard1981transformation,jirak2016berry,banna2016bernstein}.
The computation remains valid conditional on the observable localization
event, whereas dependent-data inference also needs a valid critical-value
approximation.  The certified algorithm decides only whether $R_n(\mf h)$
exceeds the testing threshold; it does not reconstruct a minimizing projected
law or provide statistical inference for that optimizer.  We therefore make no condition-free claim for arbitrary
non-i.i.d. sequences.

\section{Conclusions}
We developed higher-order theory and certified computation for WP tests of
moment restrictions.  The expansions yield
$\widetilde{O}(n^{-1})$ size accuracy, an order-$n^{-1/2}$ local-power
correction, and scalar Bartlett-type corrections with
$\widetilde{O}(n^{-3/2})$ coverage accuracy.  WP's second-order power reflects
derivative-based curvature and geometry, whereas EL's reflects moment
skewness.  The fairness experiment illustrates this contrast without implying
dominance.  Finally, we propose a certified localized dual algorithm that
implements the WP test with oracle accuracy and high-probability
$\widetilde{O}(n)$ complexity.

\section*{Funding}
The material in this paper is based upon work supported by the Air Force Office of Scientific Research under award number FA9550-20-1-0397. Additional support was provided by NSF awards 1915967, 2118199, and 2229012. J. Blanchet gratefully acknowledges support from ONR under award N00014-24-1-2655, and the National Science Foundation (NSF) under grants 2312204 and 2403007. Viet Anh Nguyen acknowledges support from UGC ECS Grant 24210924 and the CUHK Improvement on Competitiveness in Hiring New Faculties Funding Scheme.

\section*{Code availability}
\addcontentsline{toc}{section}{Code availability}
The code and compact numerical outputs supporting the computational results are
available in the public repository \emph{Reproducibility code for Wasserstein
Projection Tests}~\cite{lin2026wpcode}.

\clearpage
\part*{Supplementary Material}
\addcontentsline{toc}{part}{Supplementary Material}
\begin{quote}
\small This supplement provides the notation reference, supplementary examples and derivations,
complete proofs, and computational details for the main paper.
\end{quote}

\renewcommand{\thesection}{S\arabic{section}}
\renewcommand{\theequation}{S.\arabic{equation}}
\renewcommand{\thefigure}{S.\arabic{figure}}
\renewcommand{\thetable}{S.\arabic{table}}
\renewcommand{\theHsection}{supp.\arabic{section}}
\renewcommand{\theHsubsection}{supp.\arabic{section}.\arabic{subsection}}
\renewcommand{\theHsubsubsection}{supp.\arabic{section}.\arabic{subsection}.\arabic{subsubsection}}
\renewcommand{\theHequation}{supp.\arabic{equation}}
\renewcommand{\theHfigure}{supp.\arabic{figure}}
\renewcommand{\theHtable}{supp.\arabic{table}}
\renewcommand{\theHtheorem}{supp.\arabic{section}.\arabic{theorem}}
\setcounter{section}{1}
\setcounter{equation}{0}
\setcounter{figure}{0}
\setcounter{table}{0}

\phantomsection
\section*{Appendix S1: Notation reference}
\label{sec:notation_guide}
\addcontentsline{toc}{section}{Appendix S1: Notation reference}
Table~\ref{tab:notation_guide} collects the global and recurring notation used in the main paper
and this supplement.  Dummy indices and
one-line proof-local abbreviations are defined where they occur.

\vspace{0.75em}
\begingroup
\footnotesize
\setstretch{1.0}
\setlength{\LTpre}{0.5em}
\setlength{\LTpost}{0.5em}
\setlength{\tabcolsep}{4pt}
\renewcommand{\arraystretch}{1.05}
\begin{longtable}{>{\raggedright\arraybackslash}p{0.18\textwidth}
                  >{\raggedright\arraybackslash}p{0.56\textwidth}
                  >{\raggedright\arraybackslash}p{0.18\textwidth}}
\caption{Notation used throughout the paper.}\label{tab:notation_guide}\\
\toprule
\textbf{Symbol} & \textbf{Definition or meaning} & \textbf{Defined at}\\
\midrule
\endfirsthead
\multicolumn{3}{l}{\small\itshape Table~\ref{tab:notation_guide} continued from the previous page.}\\
\toprule
\textbf{Symbol} & \textbf{Definition or meaning} & \textbf{Defined at}\\
\midrule
\endhead
\midrule
\multicolumn{3}{r}{\small\itshape Continued on the next page.}\\
\endfoot
\bottomrule
\endlastfoot

\multicolumn{3}{l}{\textit{Basic objects and probability}}\\*
$n$, $[n]$ & Sample size and index set $[n]=\{1,\ldots,n\}$. & Section~1\\
$d_X$, $d_h$ & Dimensions of the sample space and moment-function codomain, respectively. & Definition~2.2\\
$\RR^d$ & $d$-dimensional Euclidean space. & Definition~2.1\\
$\mc P(\RR^{d_X})$ & Set of probability distributions on $\RR^{d_X}$. & Definition~2.1\\
$X$, $X_1,\ldots,X_n$ & Generic observation and i.i.d. sample in $\RR^{d_X}$. & Definition~2.2\\
$\PP^\star$, $\EE_{\PP^\star}$ & Data-generating distribution and expectation under it; $\EE$ is used when the governing distribution is clear. & Eq.~(3)\\
$\PP_n^\star$, $\EE_n$ & Sample-size-dependent one-observation law under local alternatives and expectation under that law; for events involving the full sample, $\PP_n^\star$ denotes the corresponding product law $(\PP_n^\star)^{\otimes n}$. & Section~3.2\\
$\QQ_n$ & Empirical distribution $n^{-1}\sum_{i=1}^n\delta_{X_i}$. & Definition~2.2\\
$\delta_x$ & Dirac probability measure concentrated at $x$. & Definition~2.2\\
$\pi$, $\pi|_{\bar X}$, $\pi|_X$ & Transport coupling and its two marginals. & Definition~2.1\\
$\EE_{\PP}[\cdot]$, $\Cov(\cdot)$ & Expectation under $\PP$ and covariance. & Definitions~2.1--2.2; Assumption~\MainBasicNumber\\
$\mc N(\mu,W)$, $\phi$ & Gaussian law with mean $\mu$ and covariance $W$; $\phi$ denotes the standard-normal density of the dimension in context. & Proposition~2.3; Eq.~(8)\\

\multicolumn{3}{l}{\textit{Wasserstein projection and testing}}\\*
$c(\bar x,x)$ & Ground cost $\norm{\bar x-x}_{\Sigma}^{2}=(\bar x-x)^\top\Sigma^{-1}(\bar x-x)$. & Eq.~(2)\\
$\Sigma$ & Positive-definite geometry matrix in the Mahalanobis ground cost. & Eq.~(2)\\
$\norm{x}_{\Sigma}$, $\norm{x}_2$ & Mahalanobis norm $(x^\top\Sigma^{-1}x)^{1/2}$ and Euclidean norm; for $x\in\RR^d$, $\norm{x}_2=\norm{x}_{I_d}$. & Eq.~(2)\\
$\Wass_c(\QQ,\PP)$ & Optimal transport cost between $\QQ$ and $\PP$ under ground cost $c$. & Definition~2.1, Eq.~(1)\\
$\mf h=(h^\beta)_{\beta\in[d_h]}$ & Vector-valued moment function $\RR^{d_X}\to\RR^{d_h}$ and its scalar components. & Definition~2.2\\
$\mc P_0$ & Null model $\{\PP:\EE_{\PP}[\mf h(X)]=\mf0\}$. & Eq.~(3)\\
$\mc H_0$, $\mc H_1$ & Null and alternative moment hypotheses $\EE_{\PP^\star}[\mf h(X)]=\mf0$ and $\ne\mf0$. & Eq.~(3)\\
$R_n(\mf h)$ & Wasserstein projection distance from $\QQ_n$ to $\mc P_0$. & Definition~2.2, Eq.~(4)\\
$nR_n(\mf h)$ & Rescaled WP test statistic. & Eq.~(5)\\
$\varrho$, $1-\varrho$ & Significance level and confidence level. & Eq.~(5)\\
$z_{1-\varrho}$, $z_{n,1-\varrho}$, $\hat z_{1-\varrho}$ & Limiting, local-alternative, and plug-in $(1-\varrho)$ critical quantiles for $nR_n(\mf h)$, respectively. & Eq.~(9), Theorem~\MainPowerExpansionNumber, and the WP decision rule\\
$\hat F_n$ & Plug-in CDF induced by $\mc W_n^{1/2}\mc V_n^{-1}\mc W_n^{1/2}$ and a standard Gaussian vector. & Section~3.3, before Eq.~(9)\\
$\hat\Theta_n$ & WP confidence region obtained by inverting the moment test over $\theta$. & Eq.~(11)\\
$\EL$, $T^2$ & Empirical-likelihood and Hotelling tests used as comparators. & EL and $T^2$ decision rules\\

\multicolumn{3}{l}{\textit{Linear algebra, derivatives, and tensors}}\\*
$\D h^\beta$, $\D\mf h$ & Row gradient of $h^\beta$ and Jacobian $(h^\beta_\gamma)_{\beta,\gamma}$. & Notation reference\\
$\D^2h^\beta$, $\D^j\mf h$ & Hessian of $h^\beta$ and $j$th derivative tensor of $\mf h$. & Assumptions~2.5--2.6\\
$h^\beta_{\alpha_1\cdots\alpha_k}$ & Partial derivative $\partial^kh^\beta/\partial x_{\alpha_1}\cdots\partial x_{\alpha_k}$. & Notation reference\\
$\kappa_1$, $\kappa_2$ & Envelopes controlling, respectively, directional Hessians and local Hessian variation. & Assumptions~2.5--2.6\\
$A^\top$, $A^{-1}$, $A^{1/2}$ & Transpose, inverse, and principal positive-semidefinite square root when $A\succcurlyeq0$ (positive definite when $A\succ0$). & Notation reference\\
$\norm{A}_2$, $\norm{A}_F$ & Matrix operator norm and matrix/tensor Frobenius norm. & Notation reference\\
$\norm{\D^j\mf h(x)}_{\mathrm{op}}$ & Induced multilinear operator norm $\sup_{\norm{u_1}_2\leq1,\ldots,\norm{u_j}_2\leq1}\norm{\D^j\mf h(x)[u_1,\ldots,u_j]}_2$. & Assumption~\MainSmoothnessNumber\\
$A\succ0$, $A\succcurlyeq0$, $\PD$ & Positive definite, positive semidefinite, and the cone of symmetric positive-definite matrices. & Notation reference\\
$\lambda_{\min}(A)$, $\lambda_{\max}(A)$, $\Tr(A)$ & Smallest eigenvalue, largest eigenvalue, and trace. & Eq.~(23)\\
$u^{\otimes\ell}$ & $\ell$-fold outer product, with entries $u_{j_1}\cdots u_{j_\ell}$. & Notation reference\\
$\langle A,B\rangle$ & Euclidean/Frobenius tensor inner product, with contraction over all displayed tensor indices. & Notation reference\\
$\inte(\cdot)$, $\conv(\cdot)$ & Interior and convex hull. & Assumption~\MainConvexHullNumber\\
$\mathbb B_2$ & Closed Euclidean unit ball in $\RR^{d_h}$. & Eq.~(16)\\

\multicolumn{3}{l}{\textit{Asymptotic notation and core expansion quantities}}\\*
$O(\cdot)$, $o(\cdot)$ & Deterministic big-$O$ and little-$o$ asymptotic orders as $n\to\infty$. & Notation reference\\
$O_p(\cdot)$, $o_p(\cdot)$ & Corresponding stochastic orders; $O_p(1)$ denotes tightness. & Notation reference\\
$\widetilde{O}(\cdot)$, $\widetilde{O}_p(\cdot)$ & Deterministic and stochastic orders with logarithmic factors in $n$ suppressed. & Notation reference\\
$H_n$ & Scaled empirical moment $n^{-1/2}\sum_i\mf h(X_i)$. & Eq.~\eqref{eq:prelim_rescale_dual}\\
$M_n(\zeta)$ & Rescaled nonlinear dual term obtained by optimizing over the displacement $\Delta$. & Eq.~\eqref{eq:prelim_rescale_dual}\\
$\zeta$, $\Delta$ & Dual moment multiplier and rescaled transport displacement in the WP dual representation. & Eq.~\eqref{eq:prelim_rescale_dual}\\
$W$ & Population moment matrix $\EE[\mf h(X)^{\otimes2}]$. & Proposition~2.3\\
$V$ & Population Jacobian Gram matrix $\EE[\D\mf h(X)\Sigma\D\mf h(X)^\top]$. & Proposition~2.3\\
$\mc W_n$ & Empirical moment matrix $n^{-1}\sum_i\mf h(X_i)^{\otimes2}$. & Eq.~(8)\\
$\mc V_n$ & Empirical Jacobian Gram matrix $n^{-1}\sum_i\D\mf h(X_i)\Sigma\D\mf h(X_i)^\top$. & Eq.~(6); Eq.~(8)\\
$\xi_n$ & Empirical signed direction $n^{-1/2}\sum_i\mc V_n^{-1}\mf h(X_i)$. & Eq.~(6)\\
$\widetilde\xi_n$ & Population-matrix analogue $n^{-1/2}\sum_iV^{-1}\mf h(X_i)$. & Eq.~\eqref{eq:nRn_expan_gamma}\\
$\mc K_n$ & Empirical cubic derivative tensor appearing in the $n^{-1/2}$ expansion term. & Eq.~(6)\\
$\mc L_n$ & Empirical quartic tensor appearing in the $n^{-1}$ expansion term; distinct from the scalar $\ell_n$ in Section~4. & Eq.~(\MainHigherExpansionEquationNumber); Supplementary Section~S2.9.6\\
$\gamma$ & Signed-root vector satisfying $nR_n(\mf h)=\norm{\gamma}_2^2+\widetilde{O}_p(n^{-1})$. & Eqs.~\eqref{eq:nRn_expan_gamma}--\eqref{eq:defgamma_n}\\
$m_k^{\alpha_1\cdots\alpha_k}$ & Raw $k$th moment tensor $\EE[\gamma_{\alpha_1}\cdots\gamma_{\alpha_k}]$. & Eq.~\eqref{eq:gamma_moment}\\
$\Psi$ & CDF of $H^\top V^{-1}H$ for $H\sim\mc N(\mf0,W)$. & Theorem~3.2\\
$\mathrm h_1$, $\mathrm h_2$, $\mathrm h_3$ & First three multivariate Hermite polynomial tensors used in the Edgeworth expansion. & Eq.~\eqref{eq:hermit_poly}\\
$\eps_n$, $\delta_n$; $\hat\eps_n$, $\hat\delta_n$ & Random remainders and corresponding deterministic high-probability bounds in the second- and higher-order expansions, respectively. & Eqs.~(\MainApproxRWPINumber) and~(\MainHigherExpansionEquationNumber)\\

\multicolumn{3}{l}{\textit{Local alternatives and power comparison}}\\*
$W_n$, $V_n$ & Population moment and Jacobian Gram matrices under $\PP_n^\star$; calligraphic versions are empirical. & Section~3.3.2\\
$\tau_n$, $\tau$ & Standardized local drift $\sqrt n\,V_n^{-1/2}\EE_n[\mf h(X)]$ and its limit. & Theorem~\MainPowerExpansionNumber\\
$\bar\phi$, $\bar p$ & Gaussian density under the local alternative and its Edgeworth correction polynomial. & Theorem~\MainPowerExpansionNumber\\
$K_1,K_2,K_3$ & Homogeneous cumulant-correction forms, with coefficient tensors $\mathsf K_1,\mathsf K_2,\mathsf K_3$, in the local-power Edgeworth polynomial. & Eq.~\eqref{eq:K_II}\\
$\PP_0^\star$, $\EE_0$, $\theta_n$ & Baseline distribution, its expectation, and local location shift in the location-shift model. & Assumption~\MainLocationShiftNumber; Supplementary Section~S2.9.4\\
$\tau_0$ & Fixed vector determining the location shift $\theta_n=\tau_0/\sqrt n$. & Assumption~\MainLocationShiftNumber; Supplementary Section~S2.9.4\\
$\alpha_j$ & Scalar raw moment $\EE_0[h(X)^j]$ (or the analogous null moment in the Bartlett section). & Eq.~\eqref{eq:k11k21}; Theorem~\MainBartlettCoverageNumber\\
$\tilde\alpha_2$, $\tilde\alpha_3$, $\tilde\alpha_4$ & Derivative-based population contractions governing WP's higher-order terms. & Eq.~\eqref{eq:k11k21}; Theorem~\MainBartlettCoverageNumber\\
$k_1,k_2,k_3$ & Scalar cumulant coefficients in the one-dimensional local-power expansion. & Eq.~\eqref{eq:k11k21}\\
$E_2(\varrho,\tau_0)$ & Second-order coefficient in WP power under a location shift. & Eq.~\eqref{eq:power_expan_one_dim}\\
$B(s)$ & Test-specific second-order power coefficient for $s\in\{\mathrm{WP},\mathrm{EL},T^2\}$. & Proposition~\MainPowerComparisonNumber\\

\multicolumn{3}{l}{\textit{Bartlett correction}}\\*
$\mc E_4(X)$ & Fourth-order envelope including the global growth control and local derivatives of $\mf h$. & Assumption~\MainSmoothnessNumber\\
$g_1$, $\chi^2_{1;1-\varrho}$ & Density and $(1-\varrho)$ quantile of a chi-square variable with one degree of freedom. & Theorem~\MainBartlettCoverageNumber\\
$k_{11},k_{22},k_{31},k_{42}$ & Cumulant-expansion coefficients for the signed-root statistic. & Theorem~\MainBartlettCoverageNumber\\
$B_0,\ldots,B_3$ & Polynomial coefficients formed from $k_{11},k_{22},k_{31},k_{42}$. & Theorem~\MainBartlettCoverageNumber\\
$u_k$ & Gamma-function normalization $2^k\Gamma(k+\tfrac12)/\Gamma(\tfrac12)$. & Theorem~\MainBartlettCoverageNumber\\
$C_1,C_2,C_3$ & Bartlett polynomial coefficients $C_k=-2u_k^{-1}\sum_{r=k}^3B_r$. & Theorem~\MainBartlettCoverageNumber\\
$q(t)$ & Cubic coverage-error polynomial $\sum_{k=1}^3C_kt^k$. & Theorem~\MainBartlettCoverageNumber\\
$C$ & Level-dependent correction $n^{-1}\sum_{k=1}^3C_k(\chi^2_{1;1-\varrho})^{k-1}$. & Proposition~\MainBartlettOneNumber\\
\end{longtable}
\clearpage
\begin{longtable}{>{\raggedright\arraybackslash}p{0.18\textwidth}
                  >{\raggedright\arraybackslash}p{0.56\textwidth}
                  >{\raggedright\arraybackslash}p{0.18\textwidth}}
\caption*{Table~\ref{tab:notation_guide} (continued): Notation used throughout the paper.}\\
\toprule
\textbf{Symbol} & \textbf{Definition or meaning} & \textbf{Defined at}\\
\midrule
\endfirsthead
\multicolumn{3}{l}{\small\itshape Table~\ref{tab:notation_guide} continued from the previous page.}\\
\toprule
\textbf{Symbol} & \textbf{Definition or meaning} & \textbf{Defined at}\\
\midrule
\endhead
\midrule
\multicolumn{3}{r}{\small\itshape Continued on the next page.}\\
\endfoot
\bottomrule
\endlastfoot

$S$ & Studentized scalar WP statistic $(\mc V_n/\mc W_n)nR_n(h)$. & Proposition~\MainBartlettTwoNumber\\

\multicolumn{3}{l}{\textit{Computation section}}\\*
$\varrho$, $1-\varrho$ & Significance and confidence levels, using the paper-wide convention. & Section~4\\
$\delta$ & Threshold radius $\hat z_{1-\varrho}/n$ in the WP decision problem. & Section~4\\
$\alpha$, $\xi$ & Nonnegative transport multiplier and unit-ball direction variable in the computational dual; this $\alpha$ is not the significance level. & Eq.~(17)\\
$\mc B_\delta(\QQ_n)$ & Wasserstein ball of radius $\delta$ centered at $\QQ_n$. & Eq.~(16)\\
$\mc M_\delta$ & Convex image of $\mc B_\delta(\QQ_n)$ under the moment map $\PP\mapsto\EE_{\PP}[\mf h(X)]$; it need not be closed under quadratic growth. & Eq.~\eqref{eq:comp_moment_image}\\
$F_\delta(\alpha,\xi)$ & Possibly extended-valued concave dual objective for the threshold problem. & Eq.~(17)\\
$\Phi_n(\delta)$ & Global dual value $\sup_{\alpha\geq0,\xi\in\mathbb B_2}F_\delta(\alpha,\xi)$. & Eq.~(17)\\
$\overline\kappa_1$, $K_{\mathrm{curv}}$ & Uniform Hessian envelope and scaled curvature bound $\norm{\Sigma}_2\overline\kappa_1$. & Eq.~(21)\\
$A(x)$ & Whitened Jacobian $\D\mf h(x)\Sigma^{1/2}$. & Eq.~(23)\\
$\ell_n$, $U_n$ & Smallest and largest eigenvalues of $\mc V_n$ in the computation section. & Eq.~(23)\\
$\varkappa_n$ & Empirical condition ratio $U_n/\ell_n$. & Eq.~(24)\\
$\underline\alpha_\delta$, $\overline\alpha_\delta$ & Lower and upper localized bounds for the transport multiplier. & Eq.~(24)\\
$\lambda(\alpha)$ & Strong-convexity modulus $2\alpha-K_{\mathrm{curv}}$ of each inner problem. & Eq.~(24)\\
$\delta_0$ & Explicit certified radius threshold for the localization and curvature results. & Eq.~(24)\\
$\mc K_\delta$ & Compact search box $[\underline\alpha_\delta,\overline\alpha_\delta]\times\mathbb B_2$. & Eq.~(24)\\
$\Phi_n^{\mathrm{loc}}(\delta)$ & Maximum of $F_\delta$ over the localized box $\mc K_\delta$. & Eq.~(25)\\
$g_i(\alpha,\xi)$ & Unique minimizer of the $i$th inner objective. & Proof of Proposition~\MainLocalizationNumber(ii), before Eq.~\eqref{eq:comp_per_i_displacement}\\
$\mu(\alpha)$ & Strong-concavity modulus of $F_\delta(\alpha,\cdot)$. & Eq.~(29)\\
$G_\delta(\alpha)$ & Profile objective $\max_{\xi\in\mathbb B_2}F_\delta(\alpha,\xi)$. & Supplementary Section~S2.10.3\\
$\mu_\delta$, $\mathsf L_\delta$ & Box-uniform strong-concavity and smoothness moduli in $\xi$. & Eq.~\eqref{eq:comp_constants}\\
$\Lambda_\delta$ & Lipschitz bound for the scalar profile $G_\delta$. & Eq.~\eqref{eq:comp_constants}\\
$\mathsf s_n$, $\bar{\mathsf m}_n$ & Maximum empirical Jacobian scale and root-mean-square empirical moment scale. & Eq.~\eqref{eq:comp_constants}\\
$\mathsf M_\delta$ & Uniform bound on $\norm{\D_\xi F_\delta}_2$. & Eq.~\eqref{eq:comp_constants}\\
$\mathfrak E_n$ & Observable event on which the certified numerical procedure is run. & Eq.~(30)\\
$\mathfrak C_n$ & High-probability event controlling the empirical scales in the runtime bound. & Theorem~\MainComputableTestNumber(iii)\\
$\eps_{\mathrm{alg},n}$ & Requested solver accuracy in Section~4. & Theorem~\MainComputableTestNumber\\
$\eps_G$, $\eps_x$ & Profile-evaluation and inner-solve accuracy budgets. & Theorem~\MainComputableTestNumber; Algorithm~\ref{alg:wp_nested}\\
$k^\star$, $T$, $k_{\mathrm{in}}$ & Numbers of trisection rounds, projected-gradient steps, and inner gradient steps. & Algorithm~\ref{alg:wp_nested}\\
$\hat\Phi_n$, $\hat{\mathsf T}_n$ & Returned approximate dual value and resulting computable WP test. & Theorem~\MainComputableTestNumber\\
$w_n$ & Width of the oracle/computable-test disagreement band. & Theorem~\MainComputableTestNumber\\

\multicolumn{3}{l}{\textit{Recurring proof machinery}}\\*
$\mc Z_n$ & Score-adaptive dual ball $\{\zeta:\norm{\zeta}_2\leq8(\norm{H_n}_2\vee n^{-1/2})/\lambda_{\min}(\mc V_n)\}$. & Eq.~\eqref{eq:Z}\\
$\Delta_{n,i}(\zeta)$ & Optimizing rescaled displacement in the $i$th term of $M_n(\zeta)$. & Lemma~\ref{lem:DelinMn}, Eq.~\eqref{eq:opt_Delta}\\
$F_n(\zeta)$, $L_n(\zeta)$ & Quadratic--cubic approximation to $M_n$ and its cubic correction. & Eqs.~\eqref{eq:expanL} and~\eqref{eq:def_Fn}\\
$F_n^\dagger(\zeta)$ & Higher-order approximation to $M_n(\zeta)$. & Proposition~\ref{prop:expanMn2}\\
$\zeta_n^\dagger$ & Optimizer of the quadratic--cubic dual approximation $F_n$. & Eq.~\eqref{eq:zetan_dagger}\\
$\mc L_{W,V}$ & Fr\'echet derivative/linear functional governing the plug-in critical-quantile expansion. & Eq.~\eqref{eq:mcLn}\\
$\gamma^\dagger$, $\gamma^\ddagger$ & Modified signed-root vectors used in the coverage and local-power proofs. & Eq.~\eqref{eq:defgammada_n}; Lemma~\ref{lem:alter_cumulantII}\\
\end{longtable}
\endgroup


\section{Technical development}\label{sec:supple}

\subsection*{Relation to statistical limit theory for optimal transport}
Classical empirical-transport limit theory typically fixes one or two target
laws and studies fluctuations of a full transport distance such as
$\Wass_c(\QQ_n,\PP)$ or $\Wass_c(\QQ_n,\PP_m)$.  The resulting behavior is
sensitive to the support, dimension, cost, and whether the two population laws
coincide.  On finite spaces, Sommerfeld and Munk~\cite{sommerfeld2018inference} obtained generally
non-Gaussian limits from directional Hadamard differentiability; related work
covers countable metric spaces and Gaussian models
\cite{tameling2019empirical,rippl2016limit}.  For quadratic transport in
general Euclidean dimension, del Barrio and Loubes~\cite{del2019central} established Gaussian
central limit theorems, centered by the expected empirical cost, using
uniqueness and stability of Kantorovich potentials.  This centering is
material: the bias between the expected empirical cost and the population
transport cost can decay at a dimension-dependent empirical-transport rate
and need not be negligible on the $n^{-1/2}$ scale
\cite{fournier2015rate,weed2019sharp}.  Hundrieser, Klatt, Munk, and
Staudt~\cite{hundrieser2022unifying} subsequently developed a unifying
empirical-process limit theory for optimal transport.  Under Donsker
conditions on the relevant class of Kantorovich potentials, their theory
centers directly at the population transport cost and expresses the limits as
suprema of Gaussian processes.  These results describe first-order
fluctuations of the unrestricted transport functional around a fixed pair of
laws.

WP concerns a different functional,
\[
    R_n(\mf h)=\inf_{\PP:\,\EE_{\PP}[\mf h(X)]=\mf0}
    \Wass_c(\QQ_n,\PP).
\]
Under the null, the population law itself belongs to the target set, so the
first variation vanishes and the leading term comes from a local quadratic
projection onto the finite-dimensional moment constraint.  Consequently,
$R_n(\mf h)$ is of order $n^{-1}$ and $nR_n(\mf h)$ has a quadratic-form limit
whose effective dimension is $d_h$, rather than the ambient sample-space
dimension $d_X$ \cite{blanchet2019robust,si2020quantifying}.  The expansions in
the main paper refine this projection regime through Edgeworth and
Bartlett-type corrections; they are therefore complementary to, rather than a
replacement for, limit theory for the unrestricted empirical Wasserstein
distance.

Other dimension-robust transport discrepancies change the geometry itself.
Entropic regularization yields Sinkhorn divergences with a parametric
$n^{-1/2}$ estimation rate for fixed regularization, with constants that
deteriorate as the regularization is removed \cite{genevay2019sample}.
Gaussian-smoothed Wasserstein distances convolve both laws before transport
and admit parametric rates and functional limit theory in any fixed dimension
\cite{nietert2021smooth,goldfeld2024limit}.  Sliced Wasserstein distances
instead aggregate one-dimensional projected transport problems, for which a
uniform process limit gives distributional results for several slicing rules
\cite{xi2022distributional}.  WP retains the original ground cost and obtains
its parametric behavior from projection onto finitely many restrictions, while
Sinkhorn, smoothing, and slicing trade some of the original transport geometry
for statistical or computational regularity.

This contrast can also be seen from the dual representations.  In the
unifying theory of Hundrieser, Klatt, Munk, and
Staudt~\cite{hundrieser2022unifying}, unrestricted transport is a supremum over
a class of Kantorovich potentials, and a population-centered limit theorem
requires the corresponding empirical process to be Donsker.  Whether this
holds depends on the dimension and the complexity of the potential class; in
high dimension its metric entropy can be too large, which is an
empirical-process manifestation of the curse of dimensionality.  The modified
transport discrepancies reduce this dual complexity in different ways.  For
fixed regularization, Sinkhorn confines the relevant dual optimizers to a
bounded Sobolev--RKHS ball~\cite{genevay2019sample}.  Gaussian smoothing moves
the convolution to the dual test functions, thereby regularizing the indexed
class sufficiently for parametric empirical-process behavior under suitable
conditions~\cite{nietert2021smooth,goldfeld2024limit}.  Slicing replaces the
ambient-dimensional transport problem by a uniform process of one-dimensional
transport problems indexed by projection directions~\cite{xi2022distributional}.
By contrast, the local WP dual is indexed by the finite-dimensional multiplier
$\zeta\in\mathbb R^{d_h}$ associated with the prescribed moment map $\mf h$.
Thus its null fluctuations are controlled by the number and regularity of the
moment restrictions, rather than by the full class of ambient-dimensional
Kantorovich potentials.  In this paper, this finite-dimensional dual
representation underlies both the stochastic expansions and the computational
algorithm developed below.

\subsection{Supplementary details for Section~2}
\label{app:section2_details}

\subsubsection{Exact rescaled dual representation}
The following exact formula underlies the local quadratic approximation in
Section~2.

\begin{proposition}[Rescaling of the WP function]\label{prop:prelim_rescale}
Suppose that $\mf h \in C^1(\RR^{d_X})$ is bounded and that the moment constraint is feasible:
\[
    \mf 0\in \inte\left(\conv\{\mf h(x):x\in\RR^{d_X}\}\right).
\]
For $\QQ_n=n^{-1}\sum_{i=1}^n\delta_{X_i}$ and
$c(\bar x,x)=\norm{\bar x-x}_{\Sigma}^{2}$,
\begin{align}\label{eq:prelim_rescale_dual}
    nR_n(\mf h)
    =
    \sup_{\zeta\in\RR^{d_h}}
    \left\{
        -\zeta^\top H_n-M_n(\zeta)
    \right\},
    \qquad
    H_n\Let \frac{1}{\sqrt n}\sum_{i=1}^n \mf h(X_i),
\end{align}
where
\begin{align*}
    M_n(\zeta)
    \Let \frac{1}{n}\sum_{i=1}^n
    \sup_{\Delta\in\RR^{d_X}}
    \left\{
        \zeta^\top
        \int_0^1
        \D\mf h(X_i+n^{-1/2}\Delta u)\Delta\,\diff u
        -\norm{\Delta}_{\Sigma}^{2}
    \right\}.
\end{align*}
\end{proposition}

\subsubsection{Comments on the regularity conditions}
Letting $\Delta\to0$ in Assumption~2.5 gives
\[
    \Big\|\sum_{\beta\in[d_h]}\zeta^{(\beta)}\D^2h^\beta(x)\Big\|_2
    \leq \kappa_1(x)\norm{\zeta}_2,
\]
so $\kappa_1$ controls the Hessian of $\mf h$ in every scalar direction.

Assumption~2.5 is compatible with moment functions whose growth at
infinity is at most quadratic, matching the squared Mahalanobis cost
in~(2).  This compatibility prevents the WP distance from
degenerating.  For example, suppose $d_X=d_h=1$, $\QQ_n=\delta_1$, and
$h(x)=x^p$ for odd $p>2$.  For
\[
    \PP_\varepsilon=(1-\varepsilon)\delta_1
    +\varepsilon\delta_{-((1-\varepsilon)/\varepsilon)^{1/p}},
\]
we have $\EE_{\PP_\varepsilon}[h(X)]=0$, whereas
\[
    \Wass_c(\QQ_n,\PP_\varepsilon)
    \leq \varepsilon\left(1+\left((1-\varepsilon)/\varepsilon\right)^{1/p}\right)^2
    =O\left(\varepsilon^{1-2/p}\right)\to0.
\]
Thus a superquadratic moment can be repaired by transporting an arbitrarily
small amount of mass over an arbitrarily large distance at vanishing quadratic
cost.  One may instead take
$c(\bar x,x)=\ell(\norm{\bar x-x}_2)$, with $\ell(r)$ behaving like $r^p$ at
infinity and like $r^2$ near zero.  Because the expansion is driven by local
perturbations of the empirical atoms, the same local analysis applies where
this modified cost remains quadratic.

The smoothness requirement in Assumption~2.6 is structural rather than merely technical.
For example, the fairness constraint studied by Si et al.~\cite{ref:si2021testing} has an
estimating function that is discontinuous at a decision boundary.  In one dimension, when
the sampling density is positive near the boundary, $O_p(n^{1/2})$ observations lie within
distance $O_p(n^{-1/2})$ of the discontinuity.  Their boundary-layer mechanism moves only
these observations by $O_p(n^{-1/2})$, rather than perturbing the entire sample.  With the
squared ground cost used in the main paper, the resulting heuristic scale is
$O_p(n^{1/2}\cdot n^{-1}\cdot n^{-1})=O_p(n^{-3/2})$, instead of the
$O_p(n^{-1})$ scale in the smooth case.  Si et al.~\cite{ref:si2021testing} use a first-order
transport cost and correspondingly obtain an $O_p(n^{-1})$ projection.  Thus dropping
smoothness changes both the scaling and the geometry of the projection and requires a
different analysis.

The eighth-moment envelope in Assumption~2.9 yields the
uniform bounds used in the stochastic expansion.  For example, Markov's
inequality and a union bound give
\[
    \PP^{\star}\left(\max_{i\in[n]}\kappa_1(X_i)>n^{1/4}\right)
    \leq n\frac{\EE[\kappa_1(X)^8]}{n^2}
    =O(n^{-1}),
\]
and hence $\max_{i\in[n]}\kappa_1(X_i)/\sqrt n=o_p(1)$.

Assumption~2.10 is tailored to the Edgeworth expansion of
$\sqrt n(\mc K_n-\EE[\mc K_n])$.  Each summand contributing to $\mc K_n$ is
bounded, up to constants, by
$\norm{\D\mf h(X_i)}_2^2\kappa_1(X_i)$, so its fourth moment is controlled by
the displayed product moment in that assumption.  This is the moment order
needed for an order-$n^{-1}$ Edgeworth expansion; see
\cite[Theorem~20.1]{bhattacharya2010normal}.  The Cram\'er regularity in
Assumption~2.7 supplies the accompanying nonlattice condition; useful
sufficient conditions appear in Section~\ref{sec:valid_Edgeworth}.

\newif\ifcompactdisplay
\newcommand{\proofdisplaysize}{%
  \ifcompactdisplay\scriptsize\else\normalsize\fi}
\AtBeginEnvironment{align}{\proofdisplaysize}
\AtBeginEnvironment{align*}{\proofdisplaysize}
\AtBeginEnvironment{equation}{\proofdisplaysize}
\AtBeginEnvironment{equation*}{\proofdisplaysize}
\AtBeginEnvironment{gather}{\proofdisplaysize}
\AtBeginEnvironment{gather*}{\proofdisplaysize}
\AtBeginEnvironment{multline}{\proofdisplaysize}
\AtBeginEnvironment{multline*}{\proofdisplaysize}

\subsection{Simulation design for Figure~1 of the main paper}
\label{app:local_power_simulation}

The numerical comparison uses the smooth symmetric three-component
Gaussian-mixture null
\begingroup\compactdisplaytrue
\[
    \PP_0^\star
    =\tfrac14\mc N(-\mu,\sigma^2)+\tfrac12\mc N(0,\sigma^2)
      +\tfrac14\mc N(\mu,\sigma^2),
    \qquad
    \sigma^2=0.005,
    \qquad
    \mu=\sqrt{2(1-\sigma^2)}=\sqrt{1.99}.
\]
\endgroup
It is symmetric with mean zero and variance one.  The moment function is
\[
    h_\theta(x)=x^2+\frac{\theta}{3}x-1,
\]
so $\EE_0[h_\theta(X)]=0$ and
\[
    h_\theta'(x)=2x+\frac{\theta}{3},
    \qquad h_\theta''(x)=2,
    \qquad h_\theta^{(3)}(x)=h_\theta^{(4)}(x)=0.
\]
Thus the smoothness and moment conditions used for the local-power expansion
hold for every fixed $\theta$.  Direct calculation gives
\[
\begin{aligned}
    \EE_0[X^4]&=2+2\sigma^2-\sigma^4,\\
    \EE_0[X^6]&=4+18\sigma^2-3\sigma^4-4\sigma^6,\\
    \kappa
    &\Let\EE_0[(X^2-1)^2]
      =1+2\sigma^2-\sigma^4,\\
    \lambda
    &\Let\EE_0[(X^2-1)^3]
      =12\sigma^2-4\sigma^6,
\end{aligned}
\]
and hence
\[
\begin{aligned}
    \alpha_2
    &=\EE_0[h_\theta(X)^2]=\kappa+\frac{\theta^2}{9},\\
    \alpha_3
    &=\EE_0[h_\theta(X)^3]
      =\lambda+\frac{\kappa\theta^2}{3},\\
    \tilde\alpha_2
    &=\EE_0[(h_\theta'(X))^2]=4+\frac{\theta^2}{9},\\
    \tilde\alpha_3
    &=\EE_0[(h_\theta'(X))^2h_\theta''(X)]
      =2\left(4+\frac{\theta^2}{9}\right).
\end{aligned}
\]
For the local shift below,
$\tau_h=\EE_0[h_\theta'(X)]\tau_0=\theta$.  The sign-adjusted scores in
Table~1 of the main paper therefore reduce to
\[
    C_{\mathrm{WP}}
    =\operatorname{sign}(\theta)\frac{2}{4+\theta^2/9},
    \qquad
    C_{\mathrm{EL}}
    =\operatorname{sign}(\theta)
      \frac{2(\lambda+\kappa\theta^2/3)}
      {3(\kappa+\theta^2/9)^2}.
\]
When $\theta<0$, both are negative and Hotelling's $T^2$ is best.  For
$\theta>0$, the sign of $C_{\mathrm{EL}}-C_{\mathrm{WP}}$ is the sign of
\[
    d(t)
    =3(\kappa-1)t^2+(\lambda+6\kappa)t
      +(4\lambda-3\kappa^2),
    \qquad t=\frac{\theta^2}{9}.
\]
The polynomial has a unique positive root
$t_\star\approx0.45979$, giving $\theta_\star\approx2.03423$.
Thus WP is best for $0<\theta<\theta_\star$: throughout this interval
$h_\theta''=2$ is fixed and
$\tilde\alpha_2=4+\theta^2/9$ stays close to its minimum $4$, so
$C_{\mathrm{WP}}$ remains large.  At $\theta=1$, for example,
\[
    C_{\mathrm{WP}}=0.48649>0.21040=C_{\mathrm{EL}},
\]
making the WP advantage pronounced rather than marginal.  Beyond the crossing,
the EL score is larger.  At $\theta=3$ the conventional skewness of the
transformed moment is
\[
    \frac{\alpha_3}{\alpha_2^{3/2}}=1.08433,
    \qquad
    C_{\mathrm{EL}}=0.50989>C_{\mathrm{WP}}=0.40000,
\]
which exhibits the value-skewness mechanism favoring EL.  At $\theta=0$,
$\tau_h=0$ and all three tests tie to order $n^{-1/2}$.

We take $n=6000$ and $\tau_0=3$, so the observations under the alternative
are shifted by $\tau_0/\sqrt n\approx0.03873$.  We evaluate $\theta$ on an equally
spaced grid from $-4$ to $4$ in increments of $1/12$, including the tie point
zero.  Every WP sample is calibrated with the feasible threshold
\[
    \hat z_{1-\varrho}
    =\frac{n^{-1}\sum_{i=1}^n h_\theta(X_i)^2}
    {n^{-1}\sum_{i=1}^n h_\theta'(X_i)^2}
      \chi^2_{1;1-\varrho}
\]
analyzed in Proposition~\MainPowerComparisonNumber. Panel~(a) reports
$\sqrt n(\widehat\pi_{n,s}-\widehat\pi_{n,T^2})$ for
$s\in\{\mathrm{WP},\mathrm{EL}\}$ and overlays the corresponding theoretical
scores $B(s)-B(T^2)$ from that proposition.  The rejection probabilities are
estimated from $300{,}000$ Monte Carlo replications using common random samples
across all values of $\theta$ and random seed 20260806.  The WP projection is
evaluated by its exact quadratic-moment formula
\[
    R_n(h_\theta)
    =\left\{
      \sqrt{\frac1n\sum_{i=1}^n(X_i+\theta/6)^2}
      -\sqrt{1+\theta^2/36}
      \right\}^2,
\]
and the EL statistic by its exact scalar dual equation.  Raw estimates are
shown as markers.  Shaded regions are pointwise
paired 95\% Monte Carlo intervals, computed from the sample variance of the
within-replication difference of rejection indicators; the largest standard
error on the displayed $\sqrt n$ scale is $0.0126$.

The dashed comparison scores are limits and omit the order-$n^{-1}$ term in
the power difference.  That term is unusually visible here because the
quadratic moment produces the exact drift
\[
    m_{1,n}\Let\EE_n[h_\theta(X)]
    =\frac{\theta}{\sqrt n}+\frac{9}{n},
    \qquad
    \sqrt n\,m_{1,n}=\theta+\frac{9}{\sqrt n}.
\]
Thus at $n=6000$ the finite-sample drift on the standardized-mean scale
contains the addition $9/\sqrt n\approx0.116$.  This is not Monte Carlo
or EL root-solving error: the largest residual in the EL dual equation is
$5.0\times10^{-12}$.  To display the resulting pre-asymptotic effect, let
\[
    m_{j,n}=\EE_n[h_\theta(X)^j],\qquad
    \widetilde m_{2,n}=\EE_n[(h_\theta'(X))^2],\qquad
    \tau_n^{\mathrm{ref}}=\frac{\sqrt n\,m_{1,n}}{\sqrt{m_{2,n}}},
\]
and put $q=\chi^2_{1;1-\varrho}$ and $c=\sqrt q$.  The solid curves use the
finite-$n$ moment refinements
\[
\begin{aligned}
    B_{\mathrm{WP},n}^{\mathrm{ref}}
    &=\frac{q\sqrt{m_{2,n}}}{\widetilde m_{2,n}}
      \left\{\phi(c-\tau_n^{\mathrm{ref}})
             -\phi(c+\tau_n^{\mathrm{ref}})\right\},\\
    B_{\mathrm{EL},n}^{\mathrm{ref}}
    &=\frac{q m_{3,n}}{3m_{2,n}^{3/2}}
      \left\{\phi(c-\tau_n^{\mathrm{ref}})
             -\phi(c+\tau_n^{\mathrm{ref}})\right\}.
\end{aligned}
\]
They evaluate the same score functionals at the exact shifted moments and
converge to the dashed Proposition~\MainPowerComparisonNumber{} scores. Their difference from those
limits is $O(n^{-1/2})$ on the displayed scale, hence $O(n^{-1})$ on the
unscaled power scale covered by the proposition's remainder.  For example, at
$\theta=1$ the limiting EL score is $0.09942$, the refinement is $0.15799$,
and the simulated scaled gap is $0.15518$.  At $\theta=3$ the corresponding
values are $0.54708$, $0.60513$, and $0.61271$.  The refinement therefore
explains the apparent upward displacement of EL without changing the
asymptotic ranking rule.  The background shading continues to mark the three
limiting regimes derived above.  Panel~(b) uses $\theta=-3$, $1$, and $3$ and
overlays each $h_\theta$ with the density of the unshifted law
$\PP_0^\star$.  The gray fill and dashed curve show this density, while each
panel title reports the ordering of $C_{\mathrm{WP}}$, $C_{\mathrm{EL}}$, and
the $T^2$ benchmark zero.

\subsection{Certified finite-dimensional stress test}
\label{app:comp_dimension_stress}

We examine the finite-sample behavior of the WP test as the data dimension
$d_X$, the number of moment restrictions $d_h$, and the sample size $n$ vary.
Every decision reported below uses the literal fixed-loop implementation of
Algorithm~1, including its observable certification event and its
fail-to-reject convention off that event.  No direct projection, KKT solver,
adaptive stopping rule, or warm start is used.

The design is
\begin{align}
    X&\sim \mc N(\mu_n,I_{d_X}),
    &
    \mf h_k(X)&=\frac{\tanh(kBX)}{k},
    &
    BB^\top&=I_{d_h},
    &
    k&=0.5,
    \label{eq:dimension_stress_design}
\end{align}
where $\tanh$ acts componentwise.  Under the null, $\mu_n=0$.  Under the
local alternative,
\begin{align}
    \mu_n=\frac{2.5}{\sqrt{nd_h}}B^\top\mf 1_{d_h},
    \label{eq:dimension_stress_local}
\end{align}
so $\sqrt n\norm{\mu_n}_2=2.5$ for every $d_h$.  We use
$d_X\in\{5,10,20\}$, $d_h\in\{1,2,5\}$, and
$n\in\{150,250,500\}$.  Hence all nine dimension pairs are nondegenerate.
The scaling $k=0.5$ retains a smooth nonlinear moment while reducing the
global Jacobian Lipschitz constant to
$K_{\mathrm{curv}}=4k/(3\sqrt3)=0.3849$, which makes the observable certification condition
informative at these sample sizes.

For every sample, the plug-in critical value is computed in the general form
of equations~(8)--(9).  In particular,
\begin{align*}
    \mc W_n&=\frac1n\sum_{i=1}^n\mf h_k(X_i)\mf h_k(X_i)^\top,
    &
    \mc V_n&=\frac1n\sum_{i=1}^n
      \D\mf h_k(X_i)\D\mf h_k(X_i)^\top,
\end{align*}
and $\hat z_{0.95}$ is the $0.95$ quantile of
$\sum_{j=1}^{d_h}\hat\lambda_jZ_j^2$, where the $Z_j$ are independent standard
normals and the $\hat\lambda_j$ are the eigenvalues of
$\mc V_n^{-1/2}\mc W_n\mc V_n^{-1/2}$.  We evaluate this quantile
deterministically by a positive-mixture expansion with coefficient-tail
tolerance $10^{-13}$.

We set $\eps_{\mathrm{alg},n}=1/\{n^{3/2}\log^2(n)\}$.  Algorithm~1 first checks
\begin{align*}
    \mathfrak E_n:\qquad
    \ell_n>0,\qquad
    \hat z_{0.95}/n\leq\delta_0/2,\qquad
    \eps_{\mathrm{alg},n}\leq\sqrt{\ell_n\hat z_{0.95}/n}/8,
\end{align*}
and fails to reject when this event does not hold.  On $\mathfrak E_n$, it
uses exactly the prescribed $k^\star$ trisection rounds, $T$
projected-gradient steps, and $k_{\mathrm{in}}$ preconditioned-gradient
steps, restarting every inner problem at $X_i$ and every profile query at
$\xi=0$.  It rejects precisely when the returned
    $\widehat\Phi_n>\eps_{\mathrm{alg},n}$.

The implementation stores the score $z=Bx$ during the inner iterations.  This
is an exact coordinate representation of Algorithm~1: starting from $x=X_i$,
every preconditioned-gradient displacement lies in the row space of $B$, and
multiplication by $B$ maps the ambient iterate bijectively to the score
iterate.  Consequently, decisions and loop counts are invariant to $d_X$ for
fixed $(d_h,n)$; only the arithmetic cost of an ambient first-order oracle
call depends on $d_X$.  Dense deterministic orthonormal-row matrices for all
nine $(d_X,d_h)$ pairs are retained with the simulation archive.

For orientation, let $Z\sim \mc N(0,1)$ and define
\begin{align*}
    a&=\EE[\operatorname{sech}^2(kZ)],
    &
    w&=\EE\left[\{\tanh(kZ)/k\}^2\right],
    &
    v&=\EE[\operatorname{sech}^4(kZ)].
\end{align*}
Here $a=0.8265$, $w=0.6941$, and $v=0.7174$.  Thus the null matrices are
$W=wI_{d_h}$ and $V=vI_{d_h}$, and the WP limiting local-power calculation
reduces to a noncentral chi-square distribution with common noncentrality
\begin{align*}
    \lambda=2.5^2a^2/w=6.1510.
\end{align*}
The limiting local power is therefore
$\PP\{\chi^2_{d_h}(\lambda)>\chi^2_{d_h;0.95}\}$, which equals $0.6985$,
$0.5954$, and $0.4433$ for $d_h=1,2,5$, respectively.

Each null and local-alternative cell uses $1{,}000$ Monte Carlo samples.  The
Gaussian score arrays are nested across $n$ within each $(d_h,\text{stage})$
block, so changes across sample sizes can also be assessed with paired
comparisons.  Table~\ref{tab:dimension_stress_inference} reports each
$d_X$-invariant result once.

\begin{table}[H]
    \centering
    \scriptsize
    \caption{Finite-sample inference from literal Algorithm~1 in the nonlinear
    stress test.  Every null and local-alternative entry is based on $1{,}000$
    samples.  The largest Monte Carlo standard errors are $0.0074$ for null
    rejection and $0.0157$ for local power.  The limiting-power column gives
    $\PP\{\chi^2_{d_h}(6.1510)>\chi^2_{d_h;0.95}\}$.  The event columns report
    the frequency of the observable certification event $\mathfrak E_n$; off
    that event Algorithm~1 fails to reject.  Results are invariant to
    $d_X\in\{5,10,20\}$.}
    \label{tab:dimension_stress_inference}
    {\setlength{\tabcolsep}{3.5pt}
    \begin{tabular}{rrccccc}
        \toprule
        $d_h$ & $n$ & Null rejection & Local power & Limiting power
        & \shortstack{Null\\$\mathfrak E_n$ frequency}
        & \shortstack{Alternative\\$\mathfrak E_n$ frequency}\\
        \midrule
        1 & 150 & 0.042 & 0.679 & 0.6985 & 1.000 & 1.000\\
        1 & 250 & 0.044 & 0.664 & 0.6985 & 1.000 & 1.000\\
        1 & 500 & 0.054 & 0.684 & 0.6985 & 1.000 & 1.000\\
        \addlinespace
        2 & 150 & 0.044 & 0.564 & 0.5954 & 1.000 & 1.000\\
        2 & 250 & 0.058 & 0.576 & 0.5954 & 1.000 & 1.000\\
        2 & 500 & 0.045 & 0.560 & 0.5954 & 1.000 & 1.000\\
        \addlinespace
        5 & 150 & 0.044 & 0.406 & 0.4433 & 0.999 & 0.996\\
        5 & 250 & 0.048 & 0.438 & 0.4433 & 1.000 & 1.000\\
        5 & 500 & 0.051 & 0.408 & 0.4433 & 1.000 & 1.000\\
        \bottomrule
    \end{tabular}
    }
\end{table}

The null rejection rates differ from $0.05$ by at most $0.008$ and by at most
$1.27$ Monte Carlo standard errors.  Local powers differ from their limits by
at most $0.038$ and by at most $2.41$ cellwise Monte Carlo standard errors.
Because the mean shift in~\eqref{eq:dimension_stress_local} shrinks at the
Pitman rate, power should stabilize rather than increase with $n$.  In the
paired data, the power changes from $n=150$ to $n=500$ are $0.005$, $-0.004$,
and $0.002$ for $d_h=1,2,5$, with paired standard errors $0.0161$, $0.0174$,
and $0.0181$.  Thus the size and local-power results are stable across the
tested sample sizes and are close to the fixed-dimensional limits.

As an additional link to the asymptotic expansion, we recomputed the leading
quadratic WP decision based on
$H_n^\top\mc V_n^{-1}H_n$ using the same plug-in critical value.  Its decisions
disagree with literal Algorithm~1 in at most $0.006$ of any cell, and in at
most $0.002$ of any null cell.  Hence the reported behavior is not driven by
an optimization shortcut or by certification failures: the literal certified
algorithm closely tracks the leading statistic in these samples.

Table~\ref{tab:dimension_stress_computation} reports the amortized wall-clock
time per decision for each design cell.

\begin{table}[H]
    \centering
    \scriptsize
    \caption{Amortized wall-clock time per decision for literal Algorithm~1.
    Each entry gives the null/local-alternative time in seconds, obtained by
    dividing the recorded runtime of each $1{,}000$-decision batch by $1{,}000$.
    The batches used eight worker threads, so the entries measure throughput
    rather than single-decision latency.  Timings include the certification
    check and all replications, including those that fail the certification
    event.  The score-coordinate implementation is invariant to $d_X$, so the
    time is reported once for each $(d_h,n)$ cell.}
    \label{tab:dimension_stress_computation}
    {\setlength{\tabcolsep}{6pt}
    \begin{tabular}{rccc}
        \toprule
        $d_h$ & $n=150$ & $n=250$ & $n=500$\\
        \midrule
        1 & $0.215/0.231$ & $0.456/0.467$ & $1.130/1.171$\\
        2 & $0.568/0.571$ & $1.091/1.082$ & $2.626/2.619$\\
        5 & $1.671/1.648$ & $3.038/3.057$ & $7.277/7.326$\\
        \bottomrule
    \end{tabular}
    }
\end{table}

All original score arrays, replication-level decisions, seeds, critical
values, certification outcomes, loop counts, and design matrices are
retained.  An independent audit
recomputed all critical values, events, decisions, and summaries; it found
zero discrepancies, a minimum certification rate of $0.996$, and maximum
design-matrix orthonormality error $2.44\times10^{-16}$.


\subsection{Supplementary details for Section~5}
\label{app:discussion_details}

\subsubsection{Dual coefficients in the WP--EL comparison}
With the Wilks normalization used by the EL test, empirical likelihood can be written as
$R_n^{\mathrm{EL}}(\mf h)=2\inf_{\substack{\PP\in\mc P_0\\ \PP\ll\QQ_n}}
D_{\mathrm{KL}}(\QQ_n\Vert\PP)
=-2n^{-1}\log\mc R(\theta)$ when $\mf h(x)=\mf f(\theta,x)$; see
\cite[Chapters~2--3]{owen2001empirical}.  Here $\PP\ll\QQ_n$ is Owen's
support restriction to laws that reweight the observed atoms.  The WP and EL dual representations
(see \cite[Section~12.3]{owen2001empirical} for EL) have local expansions of the form
\begin{align*}
    nR_n^{(s)}(\mf h)
    &=\{F_n^{(s)}\}^*(-H_n)+\widetilde{O}_p(n^{-1}),
    \qquad s\in\{\mathrm{WP},\mathrm{EL}\},
\end{align*}
where $R_n^{(\mathrm{WP})}=R_n$, $H_n=n^{-1/2}\sum_{i=1}^n\mf h(X_i)$, and
\begin{align*}
    F_n^{(s)}(\zeta)
    &=\frac1n\sum_{i=1}^n\bigg\{
      A^{s}_{\omega^1\omega^2}(X_i)\zeta^{(\omega^1)}\zeta^{(\omega^2)}
      +\frac{A^{s}_{\omega^1\omega^2\omega^3}(X_i)}{\sqrt n}
       \zeta^{(\omega^1)}\zeta^{(\omega^2)}\zeta^{(\omega^3)}\\
    &\hspace{42mm}
      +\frac{A^{s}_{\omega^1\omega^2\omega^3\omega^4}(X_i)}{n}
       \prod_{j=1}^4\zeta^{(\omega^j)}\bigg\}.
\end{align*}
The quadratic and cubic coefficients are compared in Table~\ref{tab:supp_compare_OT_EL}.
\begin{table}[!t]
    \centering
    \caption{Quadratic and cubic coefficients in the WP and EL dual expansions.  All
    derivatives and moment values are evaluated at $X_i$.}
    \label{tab:supp_compare_OT_EL}
    \begin{tabular}{lll}
        \toprule
        Coefficient & WP & EL\\
        \midrule
        $A^{s}_{\omega^1\omega^2}$
        & $\frac14h^{\omega^1}_{\beta}h^{\omega^2}_{\gamma}\Sigma_{\beta\gamma}$
        & $\frac14h^{\omega^1}h^{\omega^2}$\\[1mm]
        $A^{s}_{\omega^1\omega^2\omega^3}$
        & $\frac18h^{\omega^1}_{\beta}h^{\omega^2}_{\lambda\eta}
          h^{\omega^3}_{\gamma}\Sigma_{\beta\lambda}\Sigma_{\eta\gamma}$
        & $\frac1{12}h^{\omega^1}h^{\omega^2}h^{\omega^3}$\\
        \bottomrule
    \end{tabular}
\end{table}

The quartic coefficients, which enter the $n^{-1}$ coverage term and the Bartlett-type
corrections, are
\begin{align*}
    A^{\mathrm{WP}}_{\omega^1\omega^2\omega^3\omega^4}
    &=\frac{1}{16}h^{\omega^1}_{\alpha\alpha'}
      h^{\omega^2}_{\beta'\gamma'}h^{\omega^3}_{\beta}h^{\omega^4}_{\omega'}
      \Sigma_{\omega'\gamma'}\Sigma_{\beta\alpha}\Sigma_{\beta'\alpha'}\\
    &\quad+\frac{1}{48}h^{\omega^1}_{\alpha\alpha'\alpha''}
      h^{\omega^2}_{\beta}h^{\omega^3}_{\beta'}h^{\omega^4}_{\beta''}
      \Sigma_{\beta\alpha}\Sigma_{\beta'\alpha'}\Sigma_{\beta''\alpha''},\\
    A^{\mathrm{EL}}_{\omega^1\omega^2\omega^3\omega^4}
    &=\frac{1}{32}h^{\omega^1}h^{\omega^2}h^{\omega^3}h^{\omega^4}.
\end{align*}
All quantities are evaluated at the corresponding observation $X_i$.  Replacing
$\tilde\alpha_2$ by $\alpha_2$ and $\tilde\alpha_3$ by $2\alpha_3/3$ in the WP
expansion gives its EL analogue.  Substitution into the scalar power expansion recovers
\cite[Theorem~3.1]{chen1994comparing} through order $n^{-1/2}$.

\begin{example}[Uniform separation by transport geometry]\label{exp:example_nonlocalpower}
For squared Euclidean ground cost, there exist a smooth even function
$h:\RR\to[-1,2]$, points $x_0,x_1,\ldots$, and
distributions $\PP^{(k)}=\frac12\delta_{x_0}+\frac12\delta_{x_k}$, $k\geq1$, such that, for
$\mc P_0=\{\PP:\EE_{\PP}[h(X)]=0\}$,
\begin{subequations}
\begin{align}
    \abs{\EE_{\PP^{(k)}}[h(X)]}^2
    &=\frac1{4^{k+1}}\longrightarrow0,
    \label{eq:hotelling}\\
    \inf_{\PP\in\mc P_0}D_{\mathrm{KL}}(\PP^{(k)}\Vert\PP)
    &=\frac12\log\left(1+\frac1{2^{k+1}}\right)
      +\frac12\log\left(1-\frac1{2(2^k+1)}\right)
      \longrightarrow0,
    \label{eq:ELDPk}\\
    \inf_{\PP\in\mc P_0}\Wass_c(\PP^{(k)},\PP)
    &\geq c_\varepsilon>0,
    \label{eq:RPkh}
\end{align}
\end{subequations}
where $c_\varepsilon$ does not depend on $k$. Thus the value-based population separations
underlying Hotelling's $T^2$ and EL can become arbitrarily small along this class, while the
transport separation remains uniformly positive. For each fixed $k$, all three tests are
consistent; WP alone retains uniform separation as $k\to\infty$.
\end{example}

\begin{figure}[t]
    \centering
    \includegraphics[width=0.72\textwidth]{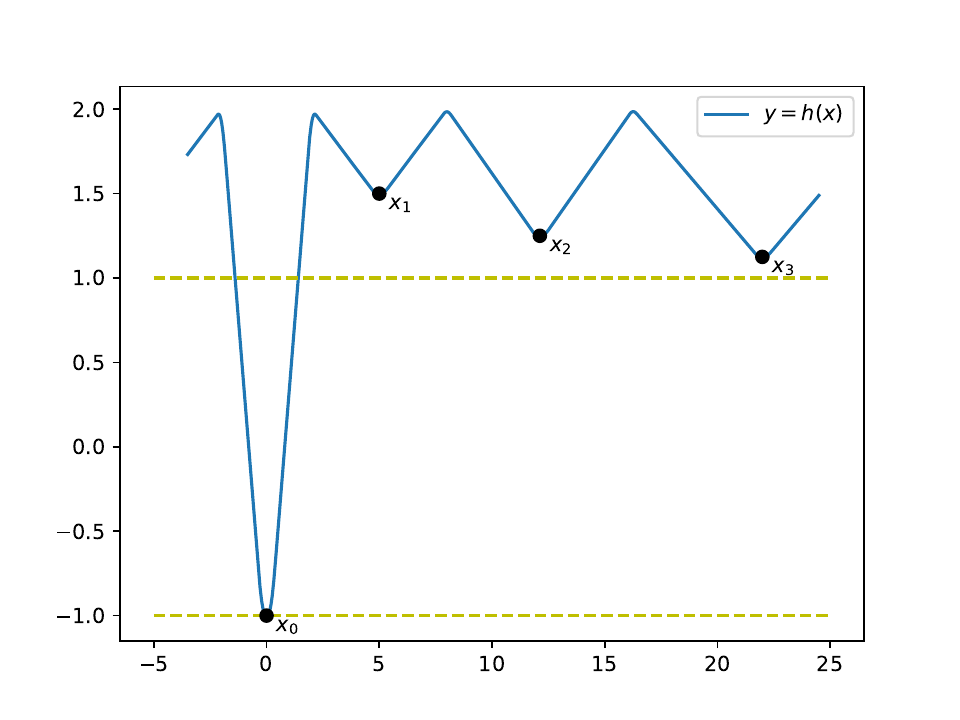}
    \caption{The smooth moment function in Example~\ref{exp:example_nonlocalpower}. Its
    increasingly distant shallow minima make the moment values approach the null while the
    required transport remains bounded away from zero.}
    \label{fig:counterexample}
\end{figure}

\subsubsection{Construction for Example~\ref{exp:example_nonlocalpower}}
\label{app:nonlocal_separation}
Set $\varepsilon=0.01$, $d_0=2$, $x_0=0$, and $y_0=-1$.  For $k\geq1$, let
\begin{align*}
    d_k&=\sqrt{2(2^{k+1}+1)},
    &x_k&=2+2\sum_{j=1}^{k-1}d_j+d_k,
    &y_k&=1+2^{-k}.
\end{align*}
The intervals $[x_k-d_k,x_k+d_k]$ meet at their endpoints.  Define a continuous even
piecewise linear function $f:\RR\to[-1,2]$ by, for every $k\geq0$,
\begin{align*}
f(x)=
\begin{cases}
-\dfrac{2-y_k}{d_k-\varepsilon}(x-x_k+\varepsilon)+y_k,
    &x_k-d_k\leq x<x_k-\varepsilon,\\[1mm]
y_k,
    &x_k-\varepsilon\leq x<x_k+\varepsilon,\\[1mm]
\dfrac{2-y_k}{d_k-\varepsilon}(x-x_k-\varepsilon)+y_k,
    &x_k+\varepsilon\leq x\leq x_k+d_k,
\end{cases}
\end{align*}
and by symmetry on the remaining negative intervals.  Let
\begin{align*}
    K_\varepsilon(t)
    &=C_\varepsilon
      \exp\left\{-\frac{1}{1-t^2/\varepsilon^2}\right\}
      \mathbbm1\{|t|<\varepsilon\},
    &\int_\RR K_\varepsilon(t)\diff t&=1,
\end{align*}
and set $h=f*K_\varepsilon$.  Then $h$ is smooth and even, $-1\leq h\leq2$, and the constant
pieces of $f$ give
\begin{align*}
    h(x_0)=-1,
    \qquad
    h(x_k)=1+2^{-k},
    \qquad k\geq1.
\end{align*}
Consequently, for $\PP^{(k)}=\frac12\delta_{x_0}+\frac12\delta_{x_k}$,
$\EE_{\PP^{(k)}}[h(X)]=2^{-(k+1)}$, which proves~\eqref{eq:hotelling}.

Write $a_k=2^{-k}$.  For any feasible $\PP$, let $p_0$ and $p_k$ denote its masses at
$x_0$ and $x_k$.  Since $h\geq-1$, feasibility implies $p_k\leq1/(2+a_k)$; hence
$p_0p_k\leq p_k(1-p_k)\leq(1+a_k)/(2+a_k)^2$.  Equality is attained by the two-point law
with weights
\begin{align*}
    q_0=\frac{1+a_k}{2+a_k},
    \qquad
    q_k=\frac{1}{2+a_k};
\end{align*}
these weights solve $-q_0+(1+a_k)q_k=0$.  Substitution into
$\frac12\log\{(1/2)/q_0\}+\frac12\log\{(1/2)/q_k\}$ gives
\eqref{eq:ELDPk}.

It remains to prove the uniform transport bound.  Let $\PP\in\mc P_0$ and disintegrate any
coupling from $\PP^{(k)}$ to $\PP$ into conditional destination laws $Q_0$ and $Q_k$ given
the source points $x_0$ and $x_k$.  Since $h\geq-1$ and
$\frac12\EE_{Q_0}[h]+\frac12\EE_{Q_k}[h]=0$, we have $\EE_{Q_k}[h]\leq1$.
Moreover,
\begin{align*}
    h(y)\geq1+a_k
    \qquad\text{whenever}\qquad
    \abs{y-x_k}<d_k-\varepsilon,
\end{align*}
because the support of the mollifier then remains inside the $k$th linear segment, where
$f\geq y_k$.  If
$p_k=Q_k\{\abs{Y-x_k}\geq d_k-\varepsilon\}$, the global bound $h\geq-1$ therefore yields
\begin{align*}
    1\geq\EE_{Q_k}[h]
    \geq(1-p_k)(1+a_k)-p_k,
    \qquad\text{so}\qquad
    p_k\geq\frac{a_k}{2+a_k}.
\end{align*}
The transport cost of the mass originating at $x_k$ is consequently at least
\begin{align*}
    \frac12p_k(d_k-\varepsilon)^2
    &\geq\frac{(d_k-\varepsilon)^2}{2(2^{k+1}+1)}
    =\left(1-\frac{\varepsilon}{d_k}\right)^2
    \geq\left(1-\frac{\varepsilon}{d_1}\right)^2
    \Let c_\varepsilon>0.
\end{align*}
Taking the infimum over $\PP\in\mc P_0$ and over couplings proves~\eqref{eq:RPkh}.

\subsection{Exact checks for Theorem~3.1}
\label{app:checks_main_expansion}

\begin{example}[Linear and quadratic moment functions]
The expansion in Theorem~3.1 agrees with the following two exactly
solvable projection problems.
\begin{enumerate}[label=(\roman*)]
    \item If $\mf h(x)=x$, then
    \[
        nR_n(\mf h)
        =\left\langle\Sigma^{-1},
          \left(\frac1{\sqrt n}\sum_{i=1}^n\mf h(X_i)\right)^{\otimes2}
          \right\rangle.
    \]
    Here $\mc V_n=\Sigma$ and $\mc K_n=0$, so
    (6) reproduces the exact value.

    \item Suppose $\mf h(x)=\norm{x}_2^2-1$, $\Sigma=I_{d_X}$,
    $\EE\norm{X}_2^2=1$, and $\EE\norm{X}_2^4<\infty$.  With
    \[
        \Delta_n=\frac1{\sqrt n}\sum_{i=1}^n(\norm{X_i}_2^2-1),
    \]
    the exact projection value satisfies
    \begin{equation}
        nR_n(\mf h)
        =n\left(\sqrt{1+\frac{\Delta_n}{\sqrt n}}-1\right)^2
        =n\left(\frac{\Delta_n}{2\sqrt n}-\frac{\Delta_n^2}{8n}
          +O_p(n^{-3/2})\right)^2.
        \label{eq:19a}
    \end{equation}
    On the other hand, Theorem~3.1 gives
    \begin{equation}
        \text{RHS of (6)}
        =\frac{\Delta_n^2}{4}-\frac{\Delta_n^3}{8\sqrt n}
          +O_p(n^{-1}).
        \label{eq:approx_rwpi_quad}
    \end{equation}
    Expanding the exact value in~\eqref{eq:19a} verifies agreement through the
    stated $O_p(n^{-1})$ remainder.
\end{enumerate}
\end{example}

\subsection{Validity of Edgeworth expansion}\label{sec:valid_Edgeworth}

The Cram\'er condition for $Y\in \RR^l$ is
\[
    \sup_{\norm{t}_2\geq b}
    \abs{\EE\brac{\exp\pare{\mathrm{i}t^\top Y}}}<1
    \qquad\text{for every }b>0.
\]

\paragraph*{Precise formulation of Assumption~2.7}
For one observation $X$, define the vectorized summands
\begingroup\compactdisplaytrue
\begin{align*}
    U_H(X)
    &\Let \mf h(X),\\
    U_V(X)
    &\Let \mathrm{vec}\!\left(\D\mf h(X)\Sigma
        \D\mf h(X)^\top\right),\\
    U_K(X)
    &\Let \mathrm{vec}\!\left(
        \left(
        \D h^{\beta}(X)\Sigma\D^2h^{\gamma}(X)\Sigma
        \D h^{\omega}(X)^\top
        \right)_{\beta,\gamma,\omega\in[d_h]}
        \right),\\
    U_W(X)
    &\Let \mathrm{vec}\!\left(\mf h(X)^{\otimes2}\right),
\end{align*}
\endgroup
and stack them into the finite-dimensional master vector
\[
    F_0(X)\Let
    \left(U_H(X)^\top,U_V(X)^\top,U_K(X)^\top,U_W(X)^\top\right)^\top.
\]
Write $F_0=(f_1,\ldots,f_k)^\top$.  Select a maximal subcollection
$\tilde F_0=(\tilde f_1,\ldots,\tilde f_r)^\top$ for which
$\{1,\tilde f_1,\ldots,\tilde f_r\}$ is linearly independent as a family of
functions on $\RR^{d_X}$.  By maximality, there are a deterministic vector $a$
and a deterministic matrix $L$ such that
\[
    F_0(x)=a+L\tilde F_0(x)\qquad\text{for every }x\in\RR^{d_X}.
\]
The full content of Assumption~2.7 is that $\PP^\star$ is absolutely
continuous with respect to Lebesgue measure and that the reduced joint vector
$\tilde F_0(X)$ satisfies
\begin{equation*}\tag{C}\label{eq:assumption_cramer_full}
    \sup_{\norm{t}_2\geq b}
    \abs{\EE\brac{\exp\pare{\mathrm{i}t^\top\tilde F_0(X)}}}<1
    \qquad\text{for every }b>0.
\end{equation*}
This single joint condition covers the projections generating $H_n$,
$\mc V_n$, $\mc K_n$, and $\mc W_n$, as well as the smooth signed-root
statistics constructed from their empirical averages.  Polynomial products
that appear only in moment and cumulant calculations are not additional
coordinates in~\eqref{eq:assumption_cramer_full}; their integrability is
supplied by Assumptions~2.9--2.10 and, for the higher-order results, by
Assumption~3.9.

More generally, when $Y$ satisfies the displayed Cram\'er condition,
\cite[Theorem 20.1, Corollary 20.4]{bhattacharya2010normal} implies that $\frac{1}{\sqrt{n}}\sum_{i=1}^n (Y_i - \EE[Y])$ admits an Edgeworth expansion up to order $n^{-\frac{s-2}{2}}$ provided that $\EE\brac{\norm{Y}_2^s} < \infty$ for $s \geq 3$, where $(Y_i)_{i \in [n]}$ are i.i.d.~copies of $Y$. Further, for a smooth vector-valued function $G\in C^s\pare{\RR^l}$, \cite[Theorem 2(b)]{bhattacharya1978validity} asserts that
\begin{align}\label{eq:Sn}
    S_n = \sqrt{n} \pare{G\pare{\frac{1}{n}\sum_{i=1}^n Y_i} - G\pare{\EE[Y]}}
\end{align}
also admits an Edgeworth expansion when its asymptotic covariance matrix is non-singular.

In the following development, $Y$ is frequently a smooth function of $X$.
The next lemma gives a useful sufficient condition for the required Cram\'er
property.
\begin{lemma}[Sufficient condition for Cram\'er regularity {\cite[Lemma 1.4]{bhattacharya1977refinements}}]\label{lem:cramer}
    Let $X$ be a random vector with values in $\RR^{d_X}$ whose distribution has
    a nonzero absolutely continuous component $\bar \PP$ relative to Lebesgue measure on $\RR^{d_X}$. Let $(f_i)_{i\in[k]}$ be Borel measurable real-valued functions on $\RR^{d_X}$. Assume that there is an open ball $\mc B\subset\RR^{d_X}$ on which the density of $\bar \PP$ is positive almost everywhere and the functions $f_i$ are continuously differentiable. If $\{1,f_1,\ldots,f_k\}$ is linearly independent on $\mc B$, then the joint distribution of $(f_i(X))_{i\in[k]}$ satisfies the Cram\'er condition.
\end{lemma}

Absolute continuity by itself does not imply
\eqref{eq:assumption_cramer_full} for a nonlinear transformation of $X$.
Lemma~\ref{lem:cramer} provides a concrete verification route: it is enough
to find a ball on which the density is positive almost everywhere and
$\{1,\tilde f_1,\ldots,\tilde f_r\}$ is linearly independent.  In that case,
for every smooth map $G$ used below,
\[
G\pare{\frac{1}{n}\sum_{i=1}^n F_0\pare{X_i}}
=G\pare{a+L\pare{\frac{1}{n}\sum_{i=1}^n\tilde F_0(X_i)}}.
\]
Consequently, the Edgeworth expansion in~\cite[Theorem 2(b)]{bhattacharya1978validity} can be applied to $S_n$ in~\eqref{eq:Sn} using the reduced vector $\tilde F_0(X)$.

\subsection{Proof of Theorem~3.1}\label{sec:pf_mainthm}

\subsubsection{Proof roadmap}\label{sec:proofroadmap}

The proof proceeds in five steps, summarized in Figure~\ref{fig:proofmap}.
\begin{figure}[H]
    \centering
\begin{tikzpicture}[
    >=Latex,
    roadmapbox/.style={
        rounded corners=6pt,
        fill={rgb,255:red,235;green,241;blue,222},
        draw=none,
        align=center,
        inner xsep=9pt,
        inner ysep=7pt,
        font=\scriptsize
    },
    roadmaparrow/.style={->, line width=0.7pt},
    roadmaplabel/.style={
        align=center,
        fill=white,
        inner sep=1.5pt,
        font=\scriptsize
    }
]
\node[roadmapbox, text width=8.6cm] (dual) at (0,0) {%
    $\displaystyle nR_n(\mathbf h)
      =\sup_{\zeta\in\RR^d}\{-\zeta^\top H_n-M_n(\zeta)\}$\\[2pt]
    \textbf{Proposition~\ref{prop:rescale}}
};

\node[roadmapbox, text width=6.35cm] (approx) at (-3.6,-2.55) {%
    $\displaystyle \sup_{\zeta\in\mc Z_n}
      \lvert M_n(\zeta)-F_n(\zeta)\rvert
      \leq \widetilde O(n^{-1})$\\[2pt]
    \textbf{Proposition~\ref{prop:expanMn}}
};

\node[roadmapbox, text width=6.35cm] (localize) at (3.6,-2.55) {%
    $\displaystyle \mathop{\arg\sup}_{\zeta\in\RR^d}
      \{-\zeta^\top H_n-M_n(\zeta)\}\in\mc Z_n$\\[2pt]
    \textbf{Proposition~\ref{prop:bound_subgradient}}
};

\node[roadmapbox, text width=11.8cm] (conclusion) at (0,-5.35) {%
    $\displaystyle
    \begin{aligned}
    nR_n(\mathbf h)
      &=\sup_{\zeta\in\mc Z_n}
        \{-\zeta^\top H_n-F_n(\zeta)\}+\widetilde O(n^{-1})\\
      &=\xi_n^\top\mathcal J_n\xi_n
        +\frac{1}{\sqrt n}\mc K_n+\widetilde O(n^{-1})
    \end{aligned}$\\[2pt]
    \textbf{Proposition~\ref{prop:expanFn*}, Theorem 3.1}
};

\draw[roadmaparrow] (dual.south west) --
    node[roadmaplabel, midway, left=2pt] {Assumptions~2.5, 2.6\\\ref{cond:condA}}
    (approx.north);
\draw[roadmaparrow] (dual.south east) --
    node[roadmaplabel, midway, right=2pt] {Assumptions~2.4, 2.5\\\ref{cond:condA}}
    (localize.north);
\draw[roadmaparrow] (approx.south) -- ([xshift=-2.2cm]conclusion.north);
\draw[roadmaparrow] (localize.south) -- ([xshift=2.2cm]conclusion.north);
\node[roadmaplabel] at (0,-4.05) {Assumption~2.5\\\ref{cond:condA}};
\end{tikzpicture}
    \caption{Proof roadmap for Theorem~3.1.  Proposition~\ref{prop:tailbnd}
    shows that \ref{cond:condA} holds with probability $1-O(n^{-1})$.}
    \label{fig:proofmap}
\end{figure}

\begin{enumerate}[label=\textit{Step \arabic*:},leftmargin=0.72in]
    \item \textit{Dual representation.}
    Proposition~\ref{prop:prelim_rescale} gives
    \[
        nR_n(\mf h)
        =\sup_{\zeta\in\RR^{d_h}}
          \{-\zeta^\top H_n-M_n(\zeta)\},
        \qquad
        H_n=\frac1{\sqrt n}\sum_{i=1}^n\mf h(X_i).
    \]

    \item \textit{High-probability bounds.}
    Assumption~2.9 controls the empirical quantities entering the
    dual problem with probability at least $1-O(n^{-1})$.

    \item \textit{Localization of the dual optimizer.}
    On the preceding event, the supremum may be restricted to the Euclidean ball
    $\mc Z_n$ with the score-adaptive radius
    $r_n=8(\|H_n\|_2\vee n^{-1/2})/\lambda_{\min}(\mc V_n)$:
    \[
        nR_n(\mf h)
        =\max_{\zeta\in\mc Z_n}\{-\zeta^\top H_n-M_n(\zeta)\}.
    \]

    \item \textit{Expansion of the dual objective.}
    Uniformly over $\mc Z_n$,
    \[
        M_n(\zeta)=F_n(\zeta)+\eps_n^M(\zeta),
        \qquad
        \sup_{\zeta\in\mc Z_n}|\eps_n^M(\zeta)|=\widetilde{O}(n^{-1}).
    \]

    \item \textit{Expansion of the supremum.}
    With high probability,
    \[
        \sup_{\zeta\in\mc Z_n}\{-\zeta^\top H_n-F_n(\zeta)\}
        =\brak{\mc V_n,\xi_n^{\otimes2}}
         +\frac1{\sqrt n}\brak{\mc K_n,\xi_n^{\otimes3}}
         +\widetilde{O}(n^{-1}).
    \]
\end{enumerate}

\subsubsection{Step 1: dual form of the projection distance \texorpdfstring{$R_n(\mf h)$}{Lg}}\label{sec:step2}

The main result of this step is presented in Proposition~\ref{prop:rescale}, where we provide the dual form characterization of the scaled projected distance $n R_n(\mf h)$.

\begin{proposition}[Rescaling of the WP function] \label{prop:rescale}
Under Assumptions~2.4 and~2.5, we have
    \begin{align} \label{eq:dual}
        n R_n(\mf h) = \sup_{\zeta \in \RR^{d_h}}\{-\zeta^{\top} H_n - M_n(\zeta)\}, \tag{Dual}
   \end{align}
    where $H_n \Let \frac{1}{\sqrt{n}} \sum_{i=1}^n\mf h(X_i)$, and
    \begin{align*}
     M_{n}(\zeta) = \frac{1}{n} \sum_{i=1}^{n} \sup _{\Delta\in\RR^{d_X}}\left\{\zeta^{\top} \int_{0}^{1} \D \mf h\left(X_{i}+n^{-1 / 2} \Delta u\right) \Delta \diff u-\norm{\Delta}_{\Sigma}^2\right\} \qquad\forall \zeta\in \RR^{d_h}.
    \end{align*}   
    Further, $M_n(\cdot)$ is a proper convex function taking values in $[0, \infty]$.
\end{proposition}

The proof of Proposition~\ref{prop:rescale} relies on the following strong duality result that formulates the value $R_n(\mf h)$ as the optimal value of an optimization problem on $\RR^{d_h} \times \RR^{d_X}$. 

\begin{lemma}[Strong duality of $R_n(\mf h)$]
\label{prop: preliminary-thm:strong_duality}
Under Assumptions~2.4 and~2.5, we have
\begin{align}
    R_n(\mf h) = \sup_{\zeta \in \RR^{d_h}} \left\{-\frac{1}{n}\sum_{i=1}^{n}\sup_{ x \in \RR^{d_X}} \left\{\zeta^\top \mf h( x) - \norm{x - X_i}_{\Sigma}^2 \right\}\right\}. \label{eq: dual rwpi}
\end{align}
Moreover, $R_n(\mf h) < \infty$, and the supremum over $\zeta \in \RR^{d_h}$ in~\eqref{eq: dual rwpi} is always attainable.
\end{lemma}

\begin{proof}[Proof of Lemma~\ref{prop: preliminary-thm:strong_duality}] By \cite[Proposition~3]{blanchet2019robust}, we obtain~\eqref{eq: dual rwpi}. For the second claim, Assumption~2.4 gives
$\mf 0 \in \conv\{\mf h(x):x\in\RR^{d_X}\}$. Thus, by
\cite[Theorem 2.3]{Rockconvex}, there exists a discrete distribution
$\PP' = \sum_{i=1}^{s} w_i \delta_{X'_i} \in \mc P(\RR^{d_X})$ such that
$\mf 0 = \EE_{\PP'}[\mf h(X)]$, where $\PP'$ is supported on
$(X_i')_{i=1}^s$ with corresponding weights $(w_i)_{i=1}^s$. Therefore,
$R_n(\mf h)\leq \Wass_c(\QQ_n, \PP') < \infty$. By
\cite[Theorem 1]{isii1962sharpness}, the supremum over
$\zeta \in \RR^{d_h}$ in~\eqref{eq: dual rwpi} is attained.
\end{proof}    

We now prove Proposition~\ref{prop:rescale}.
\begin{proof}[Proof of Proposition~\ref{prop:rescale}]
    By Lemma~\ref{prop: preliminary-thm:strong_duality}, we have
    \begingroup\compactdisplaytrue
    \begin{align*}
        R_n(\mf h) &~= \sup_{\zeta \in \RR^{d_h}} \left\{-\frac{1}{n}\sum_{i=1}^{n}\sup_{\Delta \in \RR^{d_X}} \left\{\zeta^\top \mf h\left(\Delta + X_i\right) - \norm{\Delta}_{\Sigma}^2 \right\}\right\}\\
        &~= \sup_{\zeta \in \RR^{d_h}} \left\{-\frac{1}{n}\sum_{i=1}^n \zeta^{\top} \mf h(X_i) -\frac{1}{n}\sum_{i=1}^{n}\sup_{\Delta \in \RR^{d_X}} \left\{\zeta^\top \left(\mf h(\Delta + X_i) - \mf h(X_i)\right) - \norm{\Delta}_{\Sigma}^2 \right\}\right\}\\
        &~= \sup_{\zeta \in \RR^{d_h}} \left\{-\frac{1}{n}\sum_{i=1}^n \zeta^{\top} \mf h(X_i) -\frac{1}{n}\sum_{i=1}^{n}\sup_{\Delta \in \RR^{d_X}} \left\{\zeta^\top\int _{0}^{1} \D \mf h(X_i + \Delta u) \Delta\diff u - \norm{\Delta}_{\Sigma}^2 \right\}\right\}.
    \end{align*}
    \endgroup
    Multiplying both sides by $n$ and applying the change of variables $\zeta^{\prime} \leftarrow \sqrt{n}\zeta$ and $\Delta^{\prime} \leftarrow \sqrt{n} \Delta$, we get
    \begingroup\compactdisplaytrue
    \begin{align*}
        nR_n(\mf h) = \sup_{\zeta^{\prime} \in \RR^{d_h}} \left\{-\frac{1}{\sqrt{n}}\sum_{i=1}^n \zeta^{\prime\top} \mf h(X_i) -\frac{1}{n}\sum_{i=1}^{n}\sup_{\Delta^\prime \in \RR^{d_X}} \left\{\zeta^{\prime\top}\int _{0}^{1} \D \mf h(X_i + n^{-\half}\Delta^{\prime} u) \Delta^{\prime}\diff u - \norm{\Delta^{\prime}}_{\Sigma}^2 \right\}\right\}.
    \end{align*}
    \endgroup
    Renaming the variables completes the first claim of the proof.

    For the second claim, to establish the convexity of $M_n(\cdot)$, it suffices to note that
    \begin{align*}
        \sup _{\Delta\in\RR^{d_X}}\left\{\zeta^{\top} \int_{0}^{1} \D \mf h\left(X_{i}+n^{-1 / 2} \Delta u\right) \Delta \diff u-\norm{\Delta}_{\Sigma}^2\right\}
    \end{align*}
    is the supremum of linear functions of $\zeta$ and is therefore convex in $\zeta$. Since $M_n(\mf 0) = 0$, the effective domain of $M_n(\cdot)$ is nonempty. Also, for all $\zeta \in \RR^{d_h}$,
    \begin{align*}
        M_n(\zeta) \geq \frac{1}{n} \sum_{i=1}^{n} \left\{\zeta^{\top} \int_{0}^{1} \D \mf h\left(X_{i}+n^{-1 / 2} \mf 0 u\right) \mf 0 \diff u-\norm{\mf 0}_{\Sigma}^2\right\} = 0.
    \end{align*}
    Thus, $M_n(\cdot)$ is a proper convex function, which completes the proof.
\end{proof}

\subsubsection{Step 2: bounding the tail probabilities}
At this step, we establish the bounds in \ref{cond:condA} with failure probability $O(n^{-1})$, matching the error order in Theorem~3.1.
\begin{lemma}[Tail bound of supremum]\label{lem:prob_tailbound}
    Suppose that $(X_i, 1\leq i \le n)$ are independent copies of a random variable $X$ with law $\PP$ and $\EE|X|^{\rho} < \infty$ for some $\rho > 0$. Then
    \begin{align*}
        \PP\left(\sup_{1\leq i\leq n} |X_i| \geq n^{\frac{2}{\rho}}\right) \leq \frac{\EE|X|^{\rho}}{n}.
    \end{align*}
\end{lemma}
\begin{proof}[Proof of Lemma~\ref{lem:prob_tailbound}]
    By Markov's inequality, we find
    \begin{align*}
        \PP\left(\sup_{1\leq i\leq n} |X_i| \geq n^{\frac{2}{\rho}}\right) \leq \frac{\EE\left[\sup_{1\leq i\leq n} |X_i|^{\rho}\right]}{n^2} \leq \frac{ \EE\left[\sum_{i=1}^n |X_i|^{\rho}\right]}{n^2} = \frac{\EE|X|^{\rho}}{n}.
    \end{align*}
    This observation completes the proof.
\end{proof}

\begin{lemma}[Tail bound of average I]\label{lem:tail_avg}
    Suppose that $(X_i, 1 \leq i \leq n)$ are independent copies of a random variable $X$ with law $\PP$, with $\EE[X] = 0$ and $\EE[|X|^\rho]<\infty$ for $\rho \in [1,2]$. Then
    \begin{align*}
        \PP\left(\left|\frac{1}{n}\sum_{i=1}^n X_i\right| \geq n^{\frac{2 - \rho}{\rho}}\right) \leq \frac{2 \EE[|X|^\rho]}{n}.
    \end{align*}
\end{lemma}
\begin{proof}[Proof of Lemma~\ref{lem:tail_avg}]
    By \cite[Lemma 4.2]{cherapanamjeri2022optimal}, we get $\EE\left[\left|\sum_{i=1}^n X_i\right|^\rho\right] \leq 2n \EE\left[\left|X\right|^\rho\right]$ for $\rho \in [1,2]$. By Markov's inequality, we have
    \begingroup\compactdisplaytrue
    \begin{align*}
        \PP\left(\left|\frac{1}{n}\sum_{i=1}^n X_i\right| \geq g(n)\right) \leq \frac{\EE\left[\left|\frac{1}{n}\sum_{i=1}^n X_i\right|^\rho\right]}{g(n)^\rho} = \frac{\EE\left[\left|\sum_{i=1}^n X_i\right|^\rho\right]}{n^\rho g(n)^\rho} \leq  \frac{2 n \EE[|X|^\rho]}{n ^\rho g(n)^\rho} = \frac{2 \EE[|X|^\rho]}{n^{\rho-1} g(n)^\rho}.
    \end{align*}
    \endgroup
    Taking $g(n) = n^{(2 - \rho)/\rho}$ gives the desired result.
\end{proof}

\begin{lemma}[Tail bound of average II]\label{lem:tail_avg2}
    Suppose that $(X_i, i \in [n])$ are independent copies of a
    $d_h$-dimensional random vector $X$ with law $\PP$, mean zero, and
    covariance matrix $V=\Cov(X,X)\in\PD$.  If, for an integer $r\geq1$, the
    distribution of $\frac{1}{\sqrt{n}}\sum_{i=1}^n X_i$ admits an Edgeworth
    expansion through terms of order $n^{-r}$ with remainder $o(n^{-r})$, then
    for every $a>\sqrt{2r\lambda_{\max}(V)}$,
    \begin{align*}
        \PP\left(\norm{\frac{1}{\sqrt{n}}\sum_{i=1}^n X_i}_2 \geq a \sqrt{\log(n)}\right) = o\left(\frac{1}{n^r}\right).
    \end{align*}
\end{lemma}
\begin{proof}[Proof of Lemma~\ref{lem:tail_avg2}]
    By the Edgeworth expansion \cite[Corollary 20.4 and Lemma 7.2]{bhattacharya2010normal}, we have
    \begin{align}\label{eq:Bedgeworth}
        \sup_{B \in \mc B} \left| \PP\left(\frac{1}{\sqrt{n}}\sum_{i=1}^n X_i \in B\right) - \sum_{s=0}^{2r} n^{-\frac{s}{2}} \int_{v \in B} p_s(v)\phi_{V}(v)\diff v\right| = o\left(\frac{1}{n^r}\right),
    \end{align}
    where $\mc B$ is the set of Borel-measurable convex sets in $\RR^{d_h}$, $\phi_{V}(x)$ is the density of a normal vector with mean zero and covariance matrix $V \in \RR^{d_h\times d_h}$, and for each $s$, $p_s(x)$ is a polynomial of degree $3s$ with coefficients depending on cumulants of $X$. Apply~\eqref{eq:Bedgeworth} to the convex ball $B_n=\{x:\norm{x}_2<a\sqrt{\log(n)}\}$ and then take complements. Since $\int p_0\phi_V=1$ and $\int p_s\phi_V=0$ for $s\geq1$, we get
    \begin{align*}
        \PP\left(\norm{\frac{1}{\sqrt{n}}\sum_{i=1}^n X_i}_2 \geq a \sqrt{\log(n)}\right) =  \sum_{s=0}^{2r} n^{-\frac{s}{2}} \int_{\norm{v}_2 \geq a \sqrt{\log(n)}} p_s(v)\phi_{V}(v)\diff v + o\left(\frac{1}{n^r}\right).
    \end{align*}
    For each $s=0,\ldots,2r$, the corresponding integral on the right-hand side is bounded by a constant times
    \begin{align}\label{eq:1n-r}
        \int_{\norm{v}_2 \geq a\sqrt{\log(n)}} \norm{v}_2^{3s} \phi_{V}(v)\diff v
        \leq C_s\{\log(n)\}^{c_s}\exp\left\{-\frac{a^2\log(n)}{2\lambda_{\max}(V)}\right\}
        = o\left(\frac{1}{n^r}\right),
    \end{align}
    for constants $C_s,c_s<\infty$, where the last equality uses $a^2>2r\lambda_{\max}(V)$. The derivation of the Gaussian tail bound in \eqref{eq:1n-r} is deferred to Supplementary Section~\ref{sec:deri_ineq}. Therefore, we obtain
    \begin{align*}
        \PP\left(\norm{\frac{1}{\sqrt{n}}\sum_{i=1}^n X_i}_2 \geq a \sqrt{\log(n)}\right) = o\left(\frac{1}{n^r}\right).
    \end{align*}
    This observation completes the proof.
\end{proof}

\begin{lemma}[Matrix Bernstein inequality {\cite[Theorem 6.6.1 and Discussion 6.6.2]{tropp2015introduction}}]\label{lem:matberstein}
    Suppose that $(A_i, 1\leq i\leq n)$ are independent copies of a self-adjoint $d_h \times d_h$ random matrix $A$ with law $\mathbb{P}$, with
        $\EE[A] = 0$, $\EE[A^2]$ finite, and $\lambda_{\min}(A) \geq - L$ almost surely. Then for $B_n  =\sum_{i=1}^{n} A_i$ and $t \geq 0$,
    \begin{align*}
        \PP\left(\lambda_{\min}(B_n) \leq - t\right) \leq d_h \exp\left(\frac{-t^2 / 2}{\norm{\EE B_n^2}_2 + L t / 3}\right).
    \end{align*}
\end{lemma}

\begin{lemma}[Tail bound of the smallest eigenvalue]\label{lem:tail_eigen}
   Suppose that $(A_i, 1\leq i\leq n)$ are independent copies of a $d_h \times d_h$ symmetric random matrix $A$ such that $A\succeq0$ almost surely, $\EE[A]\in\PD$, and $\EE[A^2]$ is finite. When
    \begin{align*}
        n \geq \left(\frac{2}{9} \vee \frac{8}{\lambda_{\min}(\EE[A])^2}\right)\left(\norm{\EE[(A - \EE[A])^2]}_2 + \norm{\EE[A]}_2\right)\log(d_h n),
    \end{align*}
    we have $\lambda_{\min}\left(n^{-1}\sum_{i=1}^{n} A_i\right) \geq \half \lambda_{\min}(\EE[A])$ with probability at least $1 - n^{-1}$.
\end{lemma}
\begin{proof}[Proof of Lemma~\ref{lem:tail_eigen}]
    Let $B_n  = n^{-1} \sum_{i=1}^{n} \left(A_i - \EE[A]\right)$ and $L = n^{-1} \norm{\EE[A]}_2$. Thus, 
    \[
    \norm{\EE B_n^2}_2 = n^{-1}\norm{\EE[(A - \EE[A])^2]}_2.\]
    Define momentarily $\hat t = \sqrt{2 (\norm{\EE B_n^2}_2 + L) \log(d_h n)}$ and consider a concentration range $[-\hat t, \infty)$ for $\lambda_{\min}(B_n)$. When the sample size $n$ satisfies
    \begin{align*}
        n \geq \left(\frac{2}{9} \vee \frac{8}{\lambda_{\min}(\EE[A])^2}\right)\left(\norm{\EE[(A - \EE[A])^2]}_2 + \norm{\EE[A]}_2\right)\log(d_h n),
    \end{align*}
    we get $\hat t \leq \min\{3, \frac{1}{2}\lambda_{\min}(\EE[A])\}$. By Lemma~\ref{lem:matberstein}, since $A_i\succeq0$ almost surely,
    \begin{align*}
        \lambda_{\min}(A_i - \EE[A]) \geq - \norm{\EE[A]}_2 \quad \forall i \in [n],
    \end{align*}
    we obtain
    \begingroup\compactdisplaytrue
    \begin{align*}
        \PP \left(\lambda_{\min}(B_n) \leq -\frac{1}{2}\lambda_{\min}(\EE[A]) \right) &\leq
        \PP\left(\lambda_{\min}(B_n) \leq - \hat t\right) & \text{(because $\hat t \leq  \frac{1}{2}\lambda_{\min}(\EE[A])$)} \\
        &\leq d_h \exp\left(\frac{-{\hat t}^2 / 2}{\norm{\EE B_n^2}_2 + L \hat t /3}\right) & \text{(Lemma~\ref{lem:matberstein})}\\
        &\leq d_h \exp\left(\frac{-{\hat t}^2 / 2}{\norm{\EE B_n^2}_2 + L}\right) &\text{(because $\hat t \le 3$)}\\
        &= \frac{1}{n} &\text{(plug in $\hat t$ and $\EE B_n^2$)}.
    \end{align*} 
    \endgroup
    Thus, with probability at least $1- n^{-1}$, for any $\zeta \in \RR^{d_h}$ such that $\norm{\zeta}_2 = 1$, we have
    \begin{align*}
        \zeta^\top \left(\frac{1}{n} \sum_{i=1}^{n} A_i\right) \zeta
        & = \zeta^\top B_n \zeta +  \zeta^\top \EE[A] \zeta & \text{(because $B_n  = n^{-1} \sum_{i=1}^{n} A_i - \EE[A]$)}\\
        & \geq \lambda_{\min}(B_n) + \lambda_{\min}(\EE[A]) \\
        & \geq -\hat t + \lambda_{\min}(\EE[A])\\
        & \geq \half \lambda_{\min}(\EE[A]).
    \end{align*}
    Therefore, $\lambda_{\min}\left(n^{-1}\sum_{i=1}^{n} A_i\right) \geq \half \lambda_{\min}(\EE[A])$.
\end{proof}

\begin{proposition}[Joint probability tail bounds]\label{prop:tailbnd}
    Suppose that $\mf h: \RR^{d_X} \rightarrow \RR^{d_h}$ satisfies
    Assumptions~2.5--2.6 and that $(X_i,1\leq i\leq n)$ are independent copies
    of a $d_X$-dimensional random vector $X$ with law $\PP^\star$. Let
    $A=\D\mf h(X)\Sigma\D\mf h(X)^\top$. If Assumptions~2.7--2.9 hold and
    \begin{align*}
        n \geq \left(\frac{2}{9} \vee \frac{8}{\lambda_{\min}(\EE[A])^2}\right)\left(\norm{\EE[(A - \EE[A])^2]}_2 + \norm{\EE[A]}_2\right)\log(d_h n),
    \end{align*}
    define
    \begin{align*}
        W&\Let\Cov(\mf h(X),\mf h(X)),
        &a_H&\Let2\sqrt{1\vee\lambda_{\max}(W)},\\
        \sigma_{\min}&\Let\lambda_{\min}(\EE[A]),
        &\bar r_n&\Let\frac{16a_H}{\sigma_{\min}}\sqrt{\log(n)}.
    \end{align*}
    Then, with probability at least $1 - O(n^{-1})$, the following inequalities hold simultaneously:
    \begin{align*}
    \begin{split}
        &\text{(A1)}~\lambda_{\min}\left(\frac{1}{n} \sum_{i=1}^{n} \D \mf h(X_i) \Sigma \D \mf h(X_i)^\top \right) \geq \half \lambda_{\min}\left(\EE\left[\D \mf h(X) \Sigma \D \mf h(X)^\top\right]\right). \\
        &\text{(A2)}~\norm{\frac{1}{\sqrt{n}} \sum_{i=1}^n \mf h(X_i)}_2 \leq a_H\sqrt{\log(n)}.\\
        &\text{(A3)}~\text{(i)}~\sup_{1\leq i \leq n}\norm{\D \mf h(X_i)}_2 \leq \frac{\hat\delta\sqrt n}{\norm{\Sigma}_2\bar r_n}. \qquad
        \text{(ii)}~\sup_{1\leq i \leq n}\kappa_1(X_i) \leq \frac{\sqrt n}{\norm{\Sigma}_2\bar r_n}.\\
    \end{split}
    \end{align*}
    \begin{align*}
    \begin{split}
        &\text{(A4)}~\text{(i)}~\left|\frac{1}{n}\sum_{i=1}^n \norm{\D \mf h(X_i)}_2^2 \kappa_1(X_i)\right| \leq 1 + \EE\left[\norm{\D \mf h(X)}_2^2 \kappa_1(X)\right].\\
        &\quad\quad~\text{(ii)}~\left|\frac{1}{n}\sum_{i=1}^n \norm{\D \mf h(X_i)}_2^2 \kappa_1(X_i)^2\right| \leq 1 + \EE\left[\norm{\D \mf h(X)}_2^2 \kappa_1(X)^2\right].\\
        &\quad\quad~\text{(iii)}~\left|\frac{1}{n}\sum_{i=1}^n \norm{\D \mf h(X_i)}_2^3 \kappa_2(X_i)\right| \leq 1 + \EE\left[\norm{\D \mf h(X)}_2^3 \kappa_2(X)\right].\\
    \end{split}
    \end{align*}
    \begin{align}\label{cond:condA}
    \begin{split}
        &\text{(A5)}~\text{(i)}~\left|\frac{1}{n}\sum_{i=1}^n \norm{\D \mf h(X_i)}_2^3\kappa_1(X_i)\kappa_2(X_i) \right| \leq \sqrt{n} + \EE\left[\norm{\D \mf h(X)}_2^3\kappa_1(X)\kappa_2(X)\right].\\
        &\quad\quad~\text{(ii)}~\left|\frac{1}{n}\sum_{i=1}^n \norm{\D \mf h(X_i)}_2^2\kappa_1(X_i)^3 \right| \leq \sqrt{n} + \EE\left[\norm{\D \mf h(X)}_2^2\kappa_1(X)^3\right].\\
        &\text{(A6)}~\left|\frac{1}{n}\sum_{i=1}^n \norm{\D \mf h(X_i)}_2^3 \kappa_1(X_i)^2 \kappa_2(X_i)\right| \leq n + \EE\left[\norm{\D \mf h(X)}_2^3 \kappa_1(X)^2 \kappa_2(X)\right]. 
    \end{split}\tag{Condition A}
    \end{align} 
\end{proposition}
We introduce \ref{cond:condA} for use in the following proof.
\begin{proof}[Proof of Proposition~\ref{prop:tailbnd}]
We prove this proposition using Assumptions~2.7--2.9 together with the
preceding lemmas.  Each item below is established with probability at least
$1-O(n^{-1})$.
    \begin{itemize}
        \item By Assumption~2.9, $\EE\brac{\norm{\D \mf h(X)}_2^4} < \infty$. Using Lemma~\ref{lem:tail_eigen}, we get (A1).
        \item By Assumption~2.8, $\EE[\mf h(X)] = \mf 0$. By Assumption~2.9, $\EE\brac{\norm{\mf h(X)}_2^4} < \infty$. The Cram\'er regularity in Assumption~2.7 therefore yields an Edgeworth expansion of $\frac{1}{\sqrt{n}}\sum_{i=1}^n \mf h(X_i)$ up to order $n^{-1}$. Since $a_H^2>2\lambda_{\max}(W)$, Lemma~\ref{lem:tail_avg2} gives (A2).
        \item By Assumption~2.9, $\EE\brac{\norm{\D \mf h(X)}_2^8} < \infty$ and $\EE\brac{\kappa_1(X)^8} < \infty$. A union bound and Markov's inequality give
        \begin{align*}
            \PP\left(\max_{i\leq n}\kappa_1(X_i)>\frac{\sqrt n}{\norm{\Sigma}_2\bar r_n}\right)
            &\leq\frac{\EE[\kappa_1(X)^8]\norm{\Sigma}_2^8\bar r_n^8}{n^3},\\
            \PP\left(\max_{i\leq n}\norm{\D\mf h(X_i)}_2>\frac{\hat\delta\sqrt n}{\norm{\Sigma}_2\bar r_n}\right)
            &\leq\frac{\EE[\norm{\D\mf h(X)}_2^8]\norm{\Sigma}_2^8\bar r_n^8}{\hat\delta^8n^3}.
        \end{align*}
        Because $\bar r_n^8=O\{\log(n)^4\}$, both probabilities are $O(n^{-1})$, proving (A3).
        \item By Assumption~2.9,
        \begin{align*}
        \EE\left[\norm{\D \mf h(X)}_2^4 \kappa_1(X)^2\right]&<\infty,
        &\EE\left[\norm{\D \mf h(X)}_2^4 \kappa_1(X)^4\right]&<\infty,\\
        \EE\left[\norm{\D \mf h(X)}_2^6 \kappa_2(X)^2\right]&<\infty.
        \end{align*}
        Lemma~\ref{lem:tail_avg} with $\rho=2$ gives (A4).
        \item By Assumption~2.9, $\EE\left[\norm{\D \mf h(X)}_2^4\kappa_1(X)^{\frac{4}{3}}\kappa_2(X)^{\frac{4}{3}}\right] < \infty,
        \EE\left[\norm{\D \mf h(X)}_2^{\frac{8}{3}}\kappa_1(X)^4\right] < \infty$. Using Lemma~\ref{lem:tail_avg} with $\rho = \frac{4}{3}$, we get (A5).
        \item By Assumption~2.9, $\EE\left[\norm{\D \mf h(X)}_2^3 \kappa_1(X)^2 \kappa_2(X)\right] < \infty$. Using Lemma~\ref{lem:tail_avg} with $\rho = 1$, we get (A6).
    \end{itemize}
    Since each preceding result holds with probability at least $1-O(n^{-1})$, the Bonferroni inequality shows that they hold simultaneously with probability at least $1-O(n^{-1})$.
\end{proof}

\subsubsection{Step 3: bounding the optimal dual variable}\label{sec:step4}

In this section, we show that, under appropriate conditions, the supremum in the dual form~\eqref{eq:dual} is attained and can therefore be replaced by a maximum. Write
\begin{align*}
    s_n&\Let\norm{H_n}_2\vee n^{-1/2},
    &\lambda_n&\Let\lambda_{\min}(\mc V_n),
    &r_n&\Let\frac{8s_n}{\lambda_n},
\end{align*}
where $\mc V_n=n^{-1}\sum_{i=1}^n\D\mf h(X_i)\Sigma\D\mf h(X_i)^\top$. On \ref{cond:condA} and for $n\geq2$, $\lambda_n>0$ and $r_n\leq\bar r_n$. We localize to the score-adaptive ball
\begin{equation} \label{eq:Z}
\mc Z_n \Let \{ \zeta \in \RR^{d_h}: \| \zeta \|_2 \le r_n\}.
\end{equation}
Consequently, without any loss of optimality, we can add the constraint $\zeta \in \mc Z_n$ into the optimization problem~\eqref{eq:dual}. The following proposition asserts this result.

\begin{proposition}[Bound of the optimal dual variable]\label{prop:bound_subgradient}
If Assumptions~2.4 and~2.5 hold and \ref{cond:condA} is satisfied, then for
$n\geq2$ satisfying
\begin{align}\label{eq:largen_Mn}
    \frac{32a_HB_\Sigma}{\sigma_{\min}^2}
    \sqrt{\frac{\log(n)}{n}}<1,
\end{align}
where
\begin{align*}
    B_\Sigma\Let\norm{\Sigma}_2^2
    \left(1+\EE[\norm{\D\mf h(X)}_2^2\kappa_1(X)]\right),
\end{align*}
and $a_H,\sigma_{\min}$ are defined in Proposition~\ref{prop:tailbnd}, we have
\begin{equation} \label{eq:dual_Zn}
    n R_n(\mf h) = \max_{\zeta \in \mc Z_n}\{-\zeta^{\top} H_n - M_n(\zeta)\}.
\end{equation}
\end{proposition}

Unless stated otherwise in this section, $(X_i,1\leq i\leq n)$ are independent copies of a $d_X$-dimensional random vector $X$ with law $\PP^\star$, and $\EE$ denotes expectation under $\PP^\star$. The phrase ``under \ref{cond:condA}'' means that $(X_i,1\leq i\leq n)$ and $X$ satisfy the bounds listed there.

The proof of Proposition~\ref{prop:bound_subgradient} requires the following preliminary result.

\begin{lemma}[Lower bound of $M_n(\cdot)$]\label{lem:lowbd_Mn}
    Under Assumption~2.5, for any $\zeta \in \RR^{d_h}$, we have
    \begingroup\compactdisplaytrue
    \begin{align*}
        M_n(\zeta) \geq \frac{\left(\zeta^\top \left(\frac{1}{n}\sum_{i=1}^{n}\D \mf h(X_i) \Sigma \D \mf h(X_i)^\top\right) \zeta \right)^2}{4 \zeta^\top \left(\frac{1}{n}\sum_{i=1}^{n}\D \mf h(X_i) \Sigma \D \mf h(X_i)^\top\right) \zeta + \frac{2}{\sqrt{n}} \norm{\zeta}_2 \left(\zeta^\top \left( \frac{1}{n}\sum_{i=1}^{n} \kappa_1(X_i)\D \mf h(X_i) \Sigma^2 \D \mf h(X_i)^\top\right) \zeta\right)}.
    \end{align*}
    \endgroup
    Here the ratio is defined to be zero when its denominator vanishes (in
    particular, when $\zeta=\mf 0$).
\end{lemma}
\begin{proof}[Proof of Lemma~\ref{lem:lowbd_Mn}]
    If the displayed denominator vanishes, the result follows from
    $M_n(\zeta)\geq0$. Otherwise, $\zeta\neq\mf 0$, and the following
    construction is well defined.
    For any $c >0$, let $\Delta_c \Let c\Sigma^{\half}\D \mf h(X_i)^\top \zeta/\norm{\zeta}_2$. Then
\begingroup\compactdisplaytrue
\begin{align*}
    & \sup_{\Delta \in \RR^{d_X}} \left\{\zeta^\top \int_{0}^{1}\D \mf h(X_i+n^{-\half}\Sigma^{\half}\Delta u)\Sigma^{\half}\Delta \diff u - \norm{\Delta}_2^2 \right\}\\
    = &  \sup_{\Delta \in \RR^{d_X}} \left\{\zeta^\top \D \mf h(X_i)\Sigma^{\half}\Delta - \norm{\Delta}_2^2 + \zeta^\top \int_{0}^{1} \left(\D \mf h(X_i+n^{-\half}\Sigma^{\half}\Delta u) - \D \mf h(X_i) \right)\Sigma^{\half}\Delta \diff u  \right\}\\
    \geq &~ \zeta^\top \D \mf h(X_i)\Sigma^{\half}\Delta_c - \norm{\Delta_c}_2^2 + \zeta^\top \int_{0}^{1} \left(\D \mf h(X_i+n^{-\half}\Sigma^{\half}\Delta_c u) - \D \mf h(X_i) \right)\Sigma^{\half}\Delta_c \diff u \\
    \geq &~ \zeta^\top \D \mf h(X_i)\Sigma^{\half}\Delta_c - \norm{\Delta_c}_2^2 - \norm{\zeta}_2 \int_{0}^{1} \norm{\D \mf h(X_i+n^{-\half}\Sigma^{\half}\Delta_c u) - \D \mf h(X_i)}_2 \norm{\Sigma^\half \Delta_c}_2 \diff u\\
    \geq &~ \zeta^\top \D \mf h(X_i)\Sigma^{\half}\Delta_c - \norm{\Delta_c}_2^2 - \frac{1}{2\sqrt{n}}\norm{\zeta}_2 \kappa_1(X_i) \norm{\Sigma^{\half}\Delta_c}_2^2\\
    =&~ c \times \zeta^\top \D \mf h(X_i) \Sigma \D \mf h(X_i)^\top \zeta / \norm{\zeta}_2 - c^2 \times \zeta^\top \D \mf h(X_i) \Sigma \D \mf h(X_i)^\top \zeta / \norm{\zeta}_2^2 \\
    & - c^2 \times \frac{\norm{\zeta}_2}{2\sqrt{n}} \kappa_1(X_i) \zeta^\top \D \mf h(X_i) \Sigma^2 \D \mf h(X_i)^\top \zeta / \norm{\zeta}_2^2,
\end{align*}
\endgroup
where the first inequality follows by substituting the chosen vector $\Delta_c$,
the second follows from H\"older's inequality, and the third follows from
Assumption~2.5.

Since $M_n(\zeta)$ is the average of the left-hand side of the preceding inequality over $(X_i, 1\leq i \leq n)$, we obtain
\begingroup\compactdisplaytrue
\begin{align*}
    M_n(\zeta) \geq &~  c \times \zeta^\top \left( \frac{1}{n}\sum_{i=1}^{n}\D \mf h(X_i) \Sigma \D \mf h(X_i)^\top\right) \zeta / \norm{\zeta}_2 - c^2 \times \zeta^\top \left( \frac{1}{n}\sum_{i=1}^{n} \D \mf h(X_i) \Sigma \D \mf h(X_i)^\top\right) \zeta / \norm{\zeta}_2^2 \\
    & - c^2 \times \frac{\norm{\zeta}_2}{2\sqrt{n}} \zeta^\top\left(\frac{1}{n}\sum_{i=1}^{n} \kappa_1(X_i)  \D \mf h(X_i) \Sigma^2 \D \mf h(X_i)^\top \right)\zeta / \norm{\zeta}_2^2.
\end{align*}
\endgroup
Maximizing the right-hand side of the inequality with respect to $c$ gives
\begingroup\compactdisplaytrue
\begin{align*}
    M_n(\zeta) \geq \frac{\left(\zeta^\top \left(\frac{1}{n}\sum_{i=1}^{n}\D \mf h(X_i) \Sigma \D \mf h(X_i)^\top\right) \zeta \right)^2}{4 \zeta^\top \left(\frac{1}{n}\sum_{i=1}^{n}\D \mf h(X_i) \Sigma \D \mf h(X_i)^\top\right) \zeta + \frac{2}{\sqrt{n}} \norm{\zeta}_2 \left(\zeta^\top \left( \frac{1}{n}\sum_{i=1}^{n} \kappa_1(X_i)\D \mf h(X_i) \Sigma^2 \D \mf h(X_i)^\top\right) \zeta\right)}. 
\end{align*}
\endgroup
This observation completes the proof.
\end{proof}

\begin{proof}[Proof of Proposition~\ref{prop:bound_subgradient}]
    Under Assumptions~2.4 and~2.5, we can invoke Proposition~\ref{prop:rescale} and Lemma~\ref{prop: preliminary-thm:strong_duality} to conclude that the optimization problem
    \begin{align*}
        \max_{\zeta \in \RR^{d_h}}\left\{-\zeta^\top H_n - M_n(\zeta)\right\}
    \end{align*}
    always has at least one finite optimizer. 
    
    Define
    \begin{align*}
        \widehat B_n
        \Let\norm{\frac1n\sum_{i=1}^n
        \kappa_1(X_i)\D\mf h(X_i)\Sigma^2\D\mf h(X_i)^\top}_2.
    \end{align*}
    \ref{cond:condA}(A4)(i) gives $\widehat B_n\leq B_\Sigma$. Moreover, (A1)--(A2) give
    \begin{align*}
        \lambda_n&\geq\frac12\sigma_{\min},
        &s_n&\leq a_H\sqrt{\log(n)},
        &r_n&\leq\bar r_n.
    \end{align*}
    Consequently, condition~\eqref{eq:largen_Mn} implies
    \begin{align*}
        \frac{8\widehat B_ns_n}{\lambda_n^2\sqrt n}
        \leq\frac{32a_HB_\Sigma}{\sigma_{\min}^2}
        \sqrt{\frac{\log(n)}{n}}<1.
    \end{align*}
    Now take $\zeta\notin\mc Z_n$ and write $r=\norm{\zeta}_2>r_n$. Lemma~\ref{lem:lowbd_Mn}, monotonicity in $r$, and the definitions of $r_n$ and $s_n$ yield
    \begin{align*}
        \frac{M_n(\zeta)}r
        &\geq\frac{\lambda_n^2r}{4\lambda_n+2\widehat B_nr/\sqrt n}
        \geq\frac{\lambda_n^2r_n}{4\lambda_n+2\widehat B_nr_n/\sqrt n}\\
        &=\frac{2s_n}{1+4\widehat B_ns_n/(\lambda_n^2\sqrt n)}
        >s_n\geq\norm{H_n}_2.
    \end{align*}
    Hence $r^{-1}\{-\zeta^\top H_n-M_n(\zeta)\}<0$ outside $\mc Z_n$.
    Every optimizer $\zeta_n\opt$ of the optimization problem $\max_{\zeta \in \RR^{d_h}}\left\{-\zeta^\top H_n - M_n(\zeta)\right\}$ satisfies
    \begin{align*}
        -(\zeta_n\opt)^\top H_n - M_n(\zeta_n\opt) \geq - \mf 0^\top H_n - M_n(\mf 0) = 0.
    \end{align*}
    Therefore, $\norm{\zeta_n\opt}_2 \leq r_n$, which proves the claim.
\end{proof}

\subsubsection{Step 4: expanding \texorpdfstring{$M_n(\zeta)$}{Lg}}\label{sec:step3}
In this section, we show that when \ref{cond:condA} holds, $\mc Z_n$ defined in~\eqref{eq:Z} lies inside the effective domain of $M_n(\cdot)$, that is, $M_n(\zeta)$ is finite for any $\zeta \in \mc Z_n$. Moreover, we also show that $M_n(\zeta)$ can be approximated by a polynomial function of $\zeta$ on the set $\mc Z_n$.

\begin{lemma}[$\Delta_{n,i}$ inside $M_n$]\label{lem:DelinMn}
    If Assumption~2.5 and \ref{cond:condA} hold, then for $\zeta \in \mc Z_n$ and every $i \in [n]$, the optimization problem
\begin{subequations}
\begin{align}\label{eq:opt_Delta}
    \max_{\Delta \in \RR^{d_X}}\left\{\zeta^{\top} \int_{0}^{1} \D \mf h\left(X_{i}+n^{-1 / 2} \Delta u\right) \Delta \diff u-\norm{\Delta}_{\Sigma}^2\right\}
\end{align}
has a non-empty solution set; hence $M_n(\zeta) < \infty$ for $\zeta \in \mc Z_n$. Let $\Delta_{n,i}(\zeta)$ be any solution to the maximization problem. We have:
\begin{enumerate}[label=(\roman*)]
    \item $\Delta_{n,i}(\zeta)$ satisfies the first-order optimality condition: 
    \begin{align}
    2 \Delta_{n,i}(\zeta) &=  \Sigma\D \mf h(X_i + n^{-\half}\Delta_{n,i}(\zeta))^\top \zeta.\label{eq:opt_Delta2}
\end{align}
    \item The norm of $\Delta_{n,i}(\zeta)$ is bounded by
    \begin{align}
        \norm{\Delta_{n,i}(\zeta)}_{2} &\leq \norm{\D \mf h(X_i)}_2 \norm{\Sigma}_2\norm{\zeta}_2. \label{eq:opt_Delta3} 
    \end{align}
\end{enumerate}

\end{subequations}    
\end{lemma}

\begin{proof}[Proof of Lemma~\ref{lem:DelinMn}]
    To prove the first claim, by Assumption~2.5, we have, for any $i \in [n]$,
    \begin{align*}
        & \zeta^{\top} \int_{0}^{1} \D \mf h\left(X_{i}+n^{-1 / 2} \Delta u\right) \Delta \diff u-\norm{\Delta}_{\Sigma}^2 \\
        \leq & \zeta^{\top} \int_{0}^{1} \left(\D \mf h\left(X_{i}+n^{-1 / 2} \Delta u\right) - \D \mf h\left(X_{i}\right)\right) \Delta \diff u + \zeta^{\top}\D \mf h\left(X_{i}\right)\Delta - \norm{\Delta}_{\Sigma}^2 \\
        \leq & \frac{\norm{\zeta}_2 \kappa_1(X_i)}{2\sqrt{n}} \norm{\Delta}_2^2 +
        \zeta^{\top}\D \mf h\left(X_{i}\right)\Delta - \norm{\Delta}_{\Sigma}^2\\
        \leq & \Big( \frac{\norm{\zeta}_2 \kappa_1(X_i)}{2\sqrt{n}} - \frac{1}{\norm{\Sigma}_2} \Big) \norm{\Delta}_2^2 +
        \zeta^{\top}\D \mf h\left(X_{i}\right)\Delta .
    \end{align*}
    Under \ref{cond:condA}(A3)(ii), $\norm{\zeta}_2\leq r_n\leq\bar r_n$ gives
    \begin{align*}
        \frac{\norm{\zeta}_2 \kappa_1(X_i)}{2\sqrt{n}} - \frac{1}{\norm{\Sigma}_2}
        \leq \frac{1}{2\norm{\Sigma}_2}-\frac{1}{\norm{\Sigma}_2}
        =-\frac{1}{2\norm{\Sigma}_2}<0.
    \end{align*}
    As a result, as $\|\Delta\|_2$ goes to infinity, the objective value in~\eqref{eq:opt_Delta} goes to $-\infty$. Since the objective function is continuous with respect to $\Delta$, the maximum is attained.

    To prove the second claim, since 
    \begin{align*}
        \int_{0}^{1} \D \mf h\left(X_{i}+n^{-1 / 2} \Delta u\right) \Delta \diff u = \sqrt{n} \left(\mf h\left(X_{i}+n^{-1 / 2} \Delta\right) - \mf h\left(X_{i}\right)\right),
    \end{align*}
    then by setting the gradient of the objective in~\eqref{eq:opt_Delta} to be zero, we get
    \begin{align*}
         \mf 0 = \D \mf h(X_i + n^{-\half}\Delta_{n,i}(\zeta))^\top \zeta - 2 \Sigma^{-1} \Delta_{n,i}(\zeta).
    \end{align*}
    This leads to~\eqref{eq:opt_Delta2}.

    To prove the third claim, by Assumption~2.5, we have
    \begin{align*}
        \norm{2 \Delta_{n,i}(\zeta)}_2 &= \norm{\Sigma \D \mf h(X_i + n^{-\half}\Delta_{n,i}(\zeta))^\top\zeta}_2 \\
        & \leq \norm{\Sigma \left(\D \mf h(X_i + n^{-\half}\Delta_{n,i}(\zeta)) - \D \mf h(X_i)\right)^\top\zeta}_2  + \norm{\Sigma \D \mf h(X_i)^\top \zeta}_2\\
        & \leq \frac{\norm{\Sigma}_2 \norm{\zeta}_2 \kappa_1(X_i)}{\sqrt{n}} \norm{\Delta_{n,i}(\zeta)}_2+ \norm{\D \mf h(X_i)}_2 \norm{\Sigma}_2\norm{\zeta}_2\\
        & \leq \norm{\Delta_{n,i}(\zeta)}_2 +  \norm{\D \mf h(X_i)}_2 \norm{\Sigma}_2\norm{\zeta}_2,
    \end{align*}
    where the last inequality follows from $\norm{\Sigma}_2\norm{\zeta}_2\kappa_1(X_i)/\sqrt n\leq1$. Rearranging gives~\eqref{eq:opt_Delta3}.
\end{proof}

\begin{proposition}[Asymptotic expansion of $M_n$]\label{prop:expanMn}
    If Assumptions~2.5,~2.6 and \ref{cond:condA} hold, then for $\zeta \in \mc Z_n$, we have
        \begin{subequations}
    \begin{align}\label{eq:expanMn}
        M_n(\zeta) = \underbrace{\frac{1}{4} \zeta^\top \mc V_n \zeta}_{\text{quadratic term}} + \underbrace{\frac{1}{\sqrt{n}} L_n(\zeta)}_{\text{cubic term}}+ \underbrace{\eps_n^{M}(\zeta)}_{\text{error term}},
    \end{align}
    where $L_n(\zeta)$ is a cubic form of $\zeta = (\zeta^{(\beta)})_{\beta \in [d_h]}$ defined by
    \begin{align}\label{eq:expanL}
        L_{n}(\zeta) = \sum_{\beta, \gamma, \omega \in [d_h]} \frac{1}{8n}\sum_{i=1}^{n} (\D h^{\beta}(X_i)\Sigma\D^2 h^{\gamma}(X_i)\Sigma\D h^{\omega}(X_i)^\top) \zeta^{(\beta)}\zeta^{(\gamma)}\zeta^{(\omega)}.
    \end{align}  
    \end{subequations}
    Moreover, the higher-order error term $\eps_n^{M}(\zeta)$ satisfies
    \begin{align*}
        \left|\eps_n^{M}(\zeta)\right| \leq C  \frac{\log(n)^6}{n} \qquad \forall \zeta \in \mc Z_n,
    \end{align*}    
    where the constant $C$ is independent of $n$. 
\end{proposition}

\begin{remark}
    In~\eqref{eq:expanMn}, the quadratic term is obtained by replacing $\D \mf h\left(X_{i}+n^{-1 / 2} \Delta u\right)$ with $\D \mf h\left(X_{i}\right)$ in the definition of $M_n(\zeta)$ for every $i \in [n]$. The cubic term results from the replacement gap, whose scale is $n^{-1 / 2} \Delta \approx \frac{1}{2} n^{-1 / 2} \Sigma \D \mf h(X_i) \zeta$, as in~\eqref{eq:opt_Delta2}.
\end{remark}

\begin{proof}[Proof of Proposition~\ref{prop:expanMn}]
If Assumption~2.5 and \ref{cond:condA} hold, then by Lemma~\ref{lem:DelinMn}, for any $\zeta \in \mc Z_n$, we have
\begin{subequations}
    \begin{align}
    & M_n(\zeta) = \frac{1}{n}\sum_{i=1}^n\left\{ \zeta^\top \int_{0}^{1} \D \mf h\left(X_{i}+n^{-1 / 2} \Delta_{n,i}(\zeta) u\right)\Delta_{n,i}(\zeta) \diff u - \norm{\Delta_{n,i}(\zeta)}_{\Sigma}^2\right\}, \label{eq:Mntmp}\\
    & \norm{\Delta_{n,i}(\zeta)}_{2} \leq \norm{\D \mf h(X_i)}_2 \norm{\Sigma}_2\norm{\zeta}_2, \label{eq:Deltaleq}
\end{align}
\end{subequations}
where $\Delta_{n,i}(\zeta)$ is an optimal solution to~\eqref{eq:opt_Delta}. By \eqref{eq:Mntmp},
\begingroup\compactdisplaytrue
\begin{align}
    & M_n(\zeta) - \frac{1}{4} \zeta^\top \mc V_n \zeta\notag\\
    =~& \frac{1}{n}\sum_{i=1}^n\left\{ \zeta^\top \int_{0}^{1} \D \mf h\left(X_{i}+n^{-\frac{1}{2}} \Delta_{n,i}(\zeta) u\right)\Delta_{n,i}(\zeta) \diff u - \norm{\Delta_{n,i}(\zeta)}_{\Sigma}^2 - \frac{1}{4} \zeta^\top \D \mf h(X_i)\Sigma\D \mf h(X_i)^\top \zeta \right\}  \notag\\
    \overset{\text{by}~\eqref{eq:opt_Delta2}}{=}~& \frac{1}{n}\sum_{i=1}^n\zeta^\top\left\{ \frac{1}{2}\int_{0}^{1} \D \mf h\left(X_{i}+n^{-\frac{1}{2}} \Delta_{n,i}(\zeta) u\right)\Sigma\D \mf h\left(X_{i}+n^{-\frac{1}{2}}\Delta_{n,i}(\zeta) \right)^\top \diff u \right. \notag\\
    & -
    \left.\frac{1}{4}\D \mf h\left(X_{i}+n^{-\frac{1}{2}} \Delta_{n,i}(\zeta)\right) \Sigma\D \mf h\left(X_{i}+n^{-\frac{1}{2}} \Delta_{n,i}(\zeta)\right)^\top - \frac{1}{4}\D \mf h(X_i)\Sigma\D \mf h(X_i)^\top \right\}\zeta \label{eq:Mn_expansion}.
\end{align}
\endgroup
Consider the scalar-valued function $F_{n,i}: \RR^{d_X} \rightarrow \RR$ defined as
\begingroup\compactdisplaytrue
\begin{align}\label{eq:Fni}
    F_{n,i}(v) \Let  \frac{1}{2}\int_{0}^{1} \zeta^{\top}\D \mf h \left(X_{i}+ v u\right) \Sigma~\D \mf h\left(X_{i}+ v \right)^\top \zeta \diff u - \frac{1}{4}\zeta^{\top}\D \mf h \left(X_{i}+ v\right)\Sigma \D \mf h\left(X_{i}+v\right)^\top \zeta.
\end{align}
\endgroup
Then~\eqref{eq:Mn_expansion} becomes
\begin{align}\label{eq: 11}
    ~\eqref{eq:Mn_expansion}=~~~& \frac{1}{n}\sum_{i=1}^n \left(F_{n,i}\left(n^{-\frac{1}{2}} \Delta_{n,i}(\zeta)\right)  - F_{n,i}(\mf 0) \right) \notag\\
    =~~~& \frac{1}{n^\frac{3}{2}}\sum_{i=1}^n \int_{0}^{1} \D F_{n,i}\left(n^{-\frac{1}{2}}\Delta_{n,i}(\zeta) s\right) \Delta_{n,i}(\zeta) \diff s \notag\\
    \overset{\text{by~\eqref{eq:opt_Delta2}}}{=}~& \frac{1}{2n^\frac{3}{2}}\sum_{i=1}^n \int_{0}^{1} \D F_{n,i}\left(n^{-\frac{1}{2}} \Delta_{n,i}(\zeta) s\right)\Sigma \D \mf h\left(X_{i}+n^{-\frac{1}{2}} \Delta_{n,i}(\zeta)\right)^\top\zeta \diff s ,
\end{align}
where the second equality follows from the fundamental theorem of calculus. Writing $\zeta=(\zeta^{(1)},\ldots,\zeta^{(d_h)})^\top$ and defining $D^{(2)}:\RR^{d_X} \rightarrow \RR^{d_X\times d_X}$ by
\begin{align*}
    D^{(2)}(x) \Let \sum_{\beta = 1}^{d_h} \zeta^{(\beta)} \D^2 h^{\beta}(x),
\end{align*}
we can compute
\begingroup\compactdisplaytrue
\begin{align*}
    \D F_{n,i}(v) 
    =~& \frac{1}{2} \left(\int_{0}^{1} \zeta^{\top}\D\mf h\left(X_{i}+v\right)\Sigma D^{(2)}\left(X_{i}+vu \right) u + \zeta^{\top}\D\mf h\left(X_{i}+v u\right)\Sigma D^{(2)}\left(X_{i}+v \right) \diff u \right) \notag\\
    &- \frac{1}{2} \zeta^{\top}\D \mf h \left(X_{i}+v\right)\Sigma D^{(2)}\left(X_{i}+v \right). 
\end{align*}
\endgroup
As a consequence, we obtain
\begin{align*}
    \eqref{eq: 11} &= \frac{1}{2n^{\frac{3}{2}}}\sum_{i=1}^n \left(\D F_{n,i}(\mf 0) \Sigma \D \mf h \left(X_{i}\right)^\top\zeta\right) + \eps_n^{M}(\zeta) \\
    &=  \frac{L_{n}(\zeta)}{\sqrt{n}}  + \eps_n^{M}(\zeta) \\
    &= \sum_{\beta, \gamma, \omega \in [d_h]} \frac{1}{8n^\frac{3}{2}}\sum_{i=1}^{n} (\D h^{\beta}(X_i)\Sigma \D^2 h^{\gamma}(X_i)\Sigma\D h^{\omega}(X_i)^\top) \zeta^{(\beta)}\zeta^{(\gamma)}\zeta^{(\omega)} + \eps_n^{M}(\zeta).
\end{align*}

We now consider the error term $\eps_n^{M}(\zeta)$ and obtain
\begingroup\compactdisplaytrue
\begin{align}\label{eq:defTn}
    \eps_n^{M}(\zeta) = \frac{1}{2n^\frac{3}{2}}\sum_{i=1}^n \left(\int_{0}^{1} \D F_{n,i}(n^{-\frac{1}{2}} \Delta_{n,i}(\zeta) s) \Sigma \D \mf h\left(X_{i}+n^{-\frac{1}{2}} \Delta_{n,i}(\zeta)\right)^\top - \D F_{n,i}(\mf 0)\Sigma \D \mf h\left(X_{i}\right)^\top \diff s\right)\zeta. 
\end{align}
\endgroup
First, we compute the term inside the integral:
\begin{align*}
    & \D F_{n,i}(n^{-\frac{1}{2}} \Delta_{n,i}(\zeta) s) \Sigma \D \mf h\left(X_{i}+n^{-\frac{1}{2}} \Delta_{n,i}(\zeta)\right)^\top - \D F_{n,i}(\mf 0)\Sigma \D \mf h\left(X_{i}\right)^\top \\
    =& \left(\D F_{n,i}(n^{-\frac{1}{2}} \Delta_{n,i}(\zeta) s) - \D F_{n,i}(\mf 0)\right)\Sigma\left(\D \mf h\left(X_{i}+n^{-\frac{1}{2}} \Delta_{n,i}(\zeta)\right) - \D \mf h\left(X_{i}\right)\right)^\top \\
    &+ \D F_{n,i}(\mf 0) \Sigma\left(\D \mf h\left(X_{i}+n^{-\frac{1}{2}} \Delta_{n,i}(\zeta)\right) - \D \mf h\left(X_{i}\right)\right)^\top \\
    &+ \left(\D F_{n,i}(n^{-\frac{1}{2}} \Delta_{n,i}(\zeta) s) - \D F_{n,i}(\mf 0)\right) \Sigma \D \mf h\left(X_{i}\right)^\top.
\end{align*}
    In addition, \eqref{eq:Deltaleq} and \ref{cond:condA}(A3)(i) give
    \begin{align*}
        n^{-1/2}\norm{\Delta_{n,i}(\zeta)}_2
        \leq n^{-1/2}\norm{\D\mf h(X_i)}_2\norm{\Sigma}_2\norm{\zeta}_2
        \leq\hat\delta\frac{\norm{\zeta}_2}{\bar r_n}
        \leq\hat\delta.
    \end{align*}
    Plugging the definition of $F_{n,i}$ in~\eqref{eq:Fni}, we therefore have, for $s \in [0,1]$,
\begin{subequations}
    \begin{align}
    &\norm{\D F_{n,i}(\mf 0)}_2 \leq \frac{\norm{\Sigma}_2}{4}\kappa_1(X_i)\norm{\D \mf h(X_i)}_2\norm{\zeta}_2^2, \label{eq:DF}\\
    &\norm{\D F_{n,i}(n^{-\frac{1}{2}}\Delta_{n,i}(\zeta)s) - \D F_{n,i}(\mf 0)}_2 \notag\\ 
    \leq & \left(
    \kappa_1(X_i) \kappa_2(X_i) \norm{n^{-\frac{1}{2}}\Delta_{n,i}(\zeta)}_2^2 + \kappa_2(X_i)\norm{\D \mf h(X_i)}_2 \norm{n^{-\frac{1}{2}}\Delta_{n,i}(\zeta)}_2 \right. \notag \\ 
    &\left. + \kappa_1(X_i)^2 \norm{n^{-\frac{1}{2}}\Delta_{n,i}(\zeta)}_2\right) \frac{5\norm{\Sigma}_2\norm{\zeta}_2^2}{4},  \label{eq:DDF}
    \end{align}
due to Assumptions~2.5 and~2.6 (the detailed algebra is deferred to Supplementary Section~\ref{sec:deri_ineq}).

Thus, by \eqref{eq:Deltaleq},    
    \begin{align}
    \norm{\D F_{n,i}(n^{-\frac{1}{2}}\Delta_{n,i}(\zeta)s) - \D F_{n,i}(\mf 0)}_2
    \leq & C \left( n^{-1} \kappa_1(X_i) \kappa_2(X_i) \norm{\D \mf h(X_i)}_2^2 \norm{\zeta}_2^4 \right. \notag\\
    & + n^{-\half} \kappa_2(X_i) \norm{\D \mf h(X_i)}_2^2 \norm{\zeta}_2^3  \notag\\
    &\left. + n^{-\half} \kappa_1(X_i)^2 \norm{\D \mf h(X_i)}_2 \norm{\zeta}_2^3 \right), \label{eq:DFbd}
    \end{align}    
where $C$ is a constant dependent on $\norm{\Sigma}_2$. We also have
\begingroup\compactdisplaytrue
\begin{align}\label{eq:DFbd2}
    \norm{\D \mf h\left(X_{i}+n^{-\frac{1}{2}} \Delta_{n,i}(\zeta)\right) - \D \mf h\left(X_{i}\right)}_2 &\leq \kappa_1(X_i) \norm{n^{-\frac{1}{2}}\Delta_{n,i}(\zeta)}_2 \notag\\
    &\leq n^{-\half} \norm{\Sigma}_2 \kappa_1(X_i) \norm{\D \mf h(X_i)}_2 \norm{\zeta}_2\qquad \text{by}~\eqref{eq:Deltaleq}.
\end{align}
\endgroup
Combining \eqref{eq:DF}, \eqref{eq:DFbd}, and \eqref{eq:DFbd2} and applying
the resulting bounds to \eqref{eq:defTn} gives
\begin{align*}
    \left|\eps_n^{M}(\zeta)\right|
    \leq &~ \frac{C'}{2n^\frac{3}{2}}\sum_{i=1}^n \left(n^{-\half} \times \left( \norm{\D \mf h(X_i)}_2^2 \kappa_1(X_i)^2 \norm{\zeta}_2^4 + \norm{\D \mf h(X_i)}_2^3 \kappa_2(X_i) \norm{\zeta}_2^4 \right)\right. \notag\\
    & \left. + n^{-1} \times \left(\norm{\D \mf h(X_i)}_2^3\kappa_1(X_i)\kappa_2(X_i)\norm{\zeta}_2^5 + \norm{\D \mf h(X_i)}_2^2\kappa_1(X_i)^3 \norm{\zeta}_2^5\right) \right.\notag\\
    & \left. + n^{-\frac{3}{2}} \times \norm{\D \mf h(X_i)}_2^3 \kappa_1(X_i)^2 \kappa_2(X_i)\norm{\zeta}_2^6 \right).
\end{align*}
\end{subequations}
Further, under parts (A4)--(A6) of \ref{cond:condA}, we find
\begin{align*}
    \left|\eps_n^{M}(\zeta)\right| &\leq C'\left(n^{-1} \times (\norm{\zeta}_2^4 + \norm{\zeta}_2^5 + \norm{\zeta}_2^6) + n^{-\frac{3}{2}} \times \norm{\zeta}_2^5 + n^{-2} \times \norm{\zeta}_2^6 \right),\\
    &\leq C''  n^{-1} \left(\norm{\zeta}_2^4 + \norm{\zeta}_2^5 + \norm{\zeta}_2^6\right),
\end{align*}
where the constants $C, C', C''$ are independent of $n$. For $\zeta \in \mc Z_n$, $\norm{\zeta}_2\leq r_n\leq\bar r_n=O\{\sqrt{\log(n)}\}$, so the displayed bound is $O\{n^{-1}\log(n)^3\}$ and, in particular, is bounded by $C\log(n)^6/n$ after enlarging $C$. This proves the desired result.
\end{proof}

\subsubsection{Step 5: expanding \texorpdfstring{$\max_{\zeta \in \mc Z_n} \left\{-\zeta^\top H_n - F_n(\zeta)\right\}$}{Lg}}
Recalling~\eqref{eq:expanMn}, we define $F_n(\zeta)$ from the expansion of $M_n$ as
\begin{align}\label{eq:def_Fn}
        M_n(\zeta) = \underbrace{\frac{1}{4} \zeta^\top \mc V_n \zeta + \frac{1}{\sqrt{n}} L_n(\zeta)}_{\Let F_n(\zeta)} + \eps_n^M(\zeta).
    \end{align}
Consider now the modification of~\eqref{eq:largen_Mn} obtained by replacing $M_n$ with $F_n$:
\begin{align*}
    \sup_{\zeta \in \mc Z_n} \left\{-\zeta^\top H_n - F_n(\zeta)\right\}.
\end{align*}

The main result of this step is the following asymptotic expansion of the above quantity.

\begin{proposition}[Expansion of $F_n^*$]\label{prop:expanFn*}
Under Assumption~2.5 and \ref{cond:condA}, if~\eqref{eq:largen_Mn} holds, then
\begin{align}\label{eq:Fn*_tmp}
    \max_{\zeta \in \mc Z_n} \left\{- \zeta^\top H_n - F_n(\zeta)\right\}  =\langle \mc V_n, \xi_n^{\otimes 2} \rangle  + \frac{1}{\sqrt{n}} \langle \mc K_n, \xi_n^{\otimes 3} \rangle  + \eps^{F}_n,
\end{align}
where $\mc V_n,~\xi_n$, and $\mc K_n$ are defined immediately after (6) in
Theorem~3.1, and
\begin{align*}
    \left| \eps^{F}_n \right| \leq C \frac{\log(n)^3}{n},
\end{align*}
where $C$ is a constant independent of $n$.
\end{proposition}

The proof of Proposition~\ref{prop:expanFn*} follows these steps: 
\begin{enumerate}
    \item We first present a preparatory result that bounds $L_n(\zeta)$ in Lemma~\ref{lem:Ln}.
    \item We show in Lemma~\ref{prop:lowerbnd_Fn} that, with high probability, the optimizer on the left-hand side of \eqref{eq:Fn*_tmp} is attained in $\inte(\mc Z_n)$.
    \item When the optimizer is in $\inte(\mc Z_n)$, we compute the asymptotic expansion of the optimizer in Lemma~\ref{lem:expan_subgradient}, where Lemma~\ref{lem:Ln} is applied.
    \item Combining these steps leads to the complete proof.
\end{enumerate}

\begin{lemma}[Bounding $L_n(\zeta)$]\label{lem:Ln}
Under Assumption~2.5 and \ref{cond:condA}, the following bounds hold for any $\zeta \in \RR^{d_h}$:
    \begingroup\compactdisplaytrue
    \begin{align*}
        &|L_{n}(\zeta)| \leq \frac{\norm{\Sigma}_2^2 \norm{\zeta}_2^3}{8}\left(1 + \EE\left[\norm{\D \mf h(X)}_2^2\kappa_1(X)\right]\right),\\
        &\norm{\D L_n(\zeta)}_2 \leq \frac{3\norm{\Sigma}_2^2 \norm{\zeta}_2^2}{8}\left(1 + \EE\left[\norm{\D \mf h(X)}_2^2\kappa_1(X)\right]\right).
    \end{align*}
    \endgroup
\end{lemma}
\begin{proof}[Proof of Lemma~\ref{lem:Ln}]
    Recall that by \eqref{eq:expanL}, $L_n$ is a cubic form of $\zeta$:
    \begingroup\compactdisplaytrue
    \begin{align*}
        L_{n}(\zeta) &= \sum_{\beta, \gamma, \omega \in [d_h]} \frac{1}{8n}\sum_{i=1}^{n} (\D h^{\beta}(X_i)\Sigma\D^2 h^{\gamma}(X_i)\Sigma\D h^{\omega}(X_i)^\top) \zeta^{(\beta)}\zeta^{(\gamma)}\zeta^{(\omega)}\\
        &= \frac{1}{8n}\sum_{i=1}^{n} \left(\left(\sum_{\beta \in [d_h]}\zeta^{(\beta)}\D h^{\beta}(X_i)\right)\Sigma \left(\sum_{\gamma \in [d_h]} \zeta^{(\gamma)}\D^2 h^{\gamma}(X_i)\right)\Sigma \left(\sum_{\omega \in [d_h]}\zeta^{(\omega)}\D h^{\omega}(X_i)^\top\right)\right).
    \end{align*}
    \endgroup
    Then, by Assumption~2.5, we get
    \begin{align*}
        |L_n(\zeta)| \leq \frac{1}{8n}\sum_{i=1}^{n} \norm{\Sigma}_2^2 \norm{\zeta}_2^3 \norm{\D \mf h(X_i)}_2^2 \kappa_1(X_i).
    \end{align*}
    Thus, by part (A4)(i) of \ref{cond:condA}, we get
    \begin{align*}
        |L_n(\zeta)| \leq \frac{\norm{\Sigma}_2^2 \norm{\zeta}_2^3}{8} \left(1 + \EE\left[\norm{\D \mf h(X)}_2^2\kappa_1(X)\right]\right),
    \end{align*}
    which is the first inequality in the statement of the lemma.
    
    As for the gradient $\D L_n(\zeta)$, we have, for $\eta^\top = (\eta^{(1)}, \ldots, \eta^{(d_h)})$,
    \begin{align*}
        \eta^\top \D L_n(\zeta) = & \sum_{\beta, \gamma, \omega \in [d_h]} \frac{1}{8n}\sum_{i=1}^{n} (\D h^{\beta}(X_i)\Sigma\D^2 h^{\gamma}(X_i)\Sigma\D h^{\omega}(X_i)^\top) \eta^{(\beta)}\zeta^{(\gamma)}\zeta^{(\omega)}\\
        & + \sum_{\beta, \gamma, \omega \in [d_h]} \frac{1}{8n}\sum_{i=1}^{n} (\D h^{\beta}(X_i)\Sigma\D^2 h^{\gamma}(X_i)\Sigma\D h^{\omega}(X_i)^\top) \zeta^{(\beta)}\eta^{(\gamma)}\zeta^{(\omega)}\\
        & + \sum_{\beta, \gamma, \omega \in [d_h]} \frac{1}{8n}\sum_{i=1}^{n} (\D h^{\beta}(X_i)\Sigma\D^2 h^{\gamma}(X_i)\Sigma\D h^{\omega}(X_i)^\top) \zeta^{(\beta)}\zeta^{(\gamma)}\eta^{(\omega)}.
    \end{align*}
    Then by Assumption~2.5, we get
    \begin{align*}
        |\eta^\top \D L_n(\zeta)| \leq \frac{3}{8n}\sum_{i=1}^{n} \norm{\Sigma}_2^2 \norm{\zeta}_2^2 \norm{\eta}_2 \norm{\D \mf h(X_i)}_2^2 \kappa_1(X_i),
    \end{align*}
    thus, by part (A4)(i) of \ref{cond:condA},
    \begingroup\compactdisplaytrue
    \begin{align*}
        \norm{\D L_n(\zeta)}_2 \leq \frac{3}{8n}\sum_{i=1}^{n} \norm{\Sigma}_2^2 \norm{\zeta}_2^2 \norm{\D \mf h(X_i)}_2^2 \kappa_1(X_i) \leq \frac{3\norm{\Sigma}_2^2 \norm{\zeta}_2^2}{8}  \left(1 + \EE\left[\norm{\D \mf h(X)}_2^2\kappa_1(X)\right]\right),
    \end{align*}
    \endgroup
    which completes the proof.
\end{proof}

\begin{lemma}[Tail bound of the solution in $F_n^*(-H_n)$]\label{prop:lowerbnd_Fn}
Under Assumption~2.5 and \ref{cond:condA}, define $\mc Z_n$ as in~\eqref{eq:Z} and define $\zeta_n^{\dagger}$ as an optimizer:
\begin{align*}
    \zeta_n^{\dagger} \in \arg\max_{\zeta \in \mc Z_n } \left\{- \zeta^\top H_n - F_n(\zeta)\right\}.
\end{align*}
If~\eqref{eq:largen_Mn} holds, then $\zeta_n^{\dagger} \in \inte(\mc Z_n)$.
\end{lemma}

\begin{proof}[Proof of Lemma~\ref{prop:lowerbnd_Fn}]
    This proof follows the same pattern as the proof of
    Proposition~\ref{prop:bound_subgradient}.
    For $\zeta$ satisfying $\norm{\zeta}_2=r_n$, Lemma~\ref{lem:Ln} and \ref{cond:condA}(A4)(i) give
    \begin{align*}
        & \frac{1}{\norm{\zeta}_2} \left(- \zeta^\top H_n - F_n(\zeta) \right)\\
        \leq&\norm{H_n}_2-\frac14\lambda_nr_n+\frac{B_\Sigma r_n^2}{8\sqrt n}\\
        \leq&s_n-2s_n+\frac{8B_\Sigma s_n^2}{\lambda_n^2\sqrt n}\\
        =&-s_n\left(1-\frac{8B_\Sigma s_n}{\lambda_n^2\sqrt n}\right)<0,
    \end{align*}
    where the last inequality follows from~\eqref{eq:largen_Mn} and (A1)--(A2). Thus $- \zeta^\top H_n - F_n(\zeta) < 0$ on the boundary. Note that for any optimizer $\zeta_n^{\dagger}$, we have
    \begin{align*}
        -\zeta_n^{\dagger \top} H_n - F_n(\zeta_n^{\dagger}) \geq - \mf 0^\top H_n - F_n(\mf 0) = 0.
    \end{align*} 
    Therefore, $\norm{\zeta_n^{\dagger}}_2<r_n$, implying that $\zeta_n^{\dagger} \in \inte(\mc Z_n)$.
\end{proof}

\begin{lemma}[Expansion of the solution in $F_n^*(-H_n)$]\label{lem:expan_subgradient}
    Under Assumption~2.5 and \ref{cond:condA}, let $\zeta_n^{\dagger}$ be the
    optimizer in Lemma~\ref{prop:lowerbnd_Fn}.  If
    $\zeta_n^{\dagger} \in \inte(\mc Z_n)$, then
\begin{align}\label{eq:zetan_dagger}
        \zeta_n^{\dagger} = - 2 \mc V_n^{-1} H_n - \frac{8}{\sqrt{n}} \mc V_n^{-1} \D L_n(- \mc V_n^{-1} H_n) + \eps_n^{\dagger},
    \end{align}
    with 
    \begin{align*}
        \norm{\eps_n^{\dagger}}_2 \leq C\left(\frac{\log(n)^{\frac{5}{2}}}{n} + \frac{\log(n)^4}{n^{\frac{3}{2}}}\right),
    \end{align*}
    where $\D L_n(\cdot)$ is the gradient function of the function $L_n(\cdot)$ defined in \eqref{eq:expanL}, and $C$ is a constant independent of the sample size $n$.
\end{lemma}
\begin{proof}[Proof of Lemma~\ref{lem:expan_subgradient}]
    First, note that $F_n(\cdot)$ is continuously differentiable with respect to its argument because $F_n$ is a third-order polynomial. When $\zeta_n^{\dagger} \in \inte(\mc Z_n)$, $\zeta_n^{\dagger}$ satisfies the first-order optimality condition:
    \begin{align*}
        - H_n - \D F_n(\zeta_n^{\dagger}) = \mf 0.
    \end{align*}
    By the definition of $F_n(\cdot)$ in \eqref{eq:def_Fn}, we have
    \begin{align*}
        - H_n - \frac{1}{2} \mc V_n \zeta_n^{\dagger} - \frac{1}{\sqrt{n}} \D L_n(\zeta_n^{\dagger}) = \mf 0 \quad \implies \quad \frac{1}{2} \mc V_n  \zeta_n^{\dagger}  &= - H_n - \frac{1}{\sqrt{n}} \D L_n(\zeta_n^{\dagger}).
    \end{align*}
    By part (A1) of \ref{cond:condA}, $\lambda_{\min}(\mc V_n) > 0$, so $\mc V_n$ is invertible. Consequently, we have
    \begin{align*}
        \zeta_n^{\dagger} & = - 2 \mc V_n^{-1} H_n - \frac{2}{\sqrt{n}} \mc V_n^{-1} \D L_n(\zeta_n^{\dagger})\\
        &= - 2 \mc V_n^{-1} H_n - \frac{2}{\sqrt{n}} \mc V_n^{-1} \D L_n\left(- 2 \mc V_n^{-1} H_n - \frac{2}{\sqrt{n}} \mc V_n^{-1} \D L_n(\zeta_n^{\dagger})\right)\\
        &= - 2 \mc V_n^{-1} H_n - \frac{8}{\sqrt{n}} \mc V_n^{-1} \D L_n\left(- \mc V_n^{-1} H_n - \frac{1}{\sqrt{n}} \mc V_n^{-1} \D L_n(\zeta_n^{\dagger})\right).
    \end{align*}
    In the last equality above, we have exploited that $L_n$ is a cubic form; hence, $\D L_n$ is a quadratic form. 
        
    By part (A2) of \ref{cond:condA}, we have $\norm{H_n}_2 \leq a_H\sqrt{\log(n)}$, and by part (A1)
    \begin{align*}
        \norm{\mc V_n^{-1}}_2 \leq 2 \lambda_{\min}\left(\EE\left[\D \mf h(X)\Sigma\D \mf h(X)^\top\right]\right)^{-1} \Let 2 \sigma_{\min}^{-1}.
    \end{align*}
    By Lemma~\ref{lem:Ln}, we have
    \begin{align*}
        \norm{\D L_n(\zeta_n^{\dagger})}_2 \leq \frac{3\norm{\Sigma}_2^2\norm{\zeta_n^{\dagger}}_2^2}{8}\left(1 + \EE\left[\norm{\D \mf h(X)}_2^2\kappa_1(X)\right]\right).
    \end{align*}
    Therefore, we obtain that
    \begin{align*}
        \norm{\eps_n^{\dagger}}_2 &= \norm{\frac{8}{\sqrt{n}} \mc V_n^{-1} \D L_n\left(- \mc V_n^{-1} H_n - \frac{1}{\sqrt{n}} \mc V_n^{-1} \D L_n(\zeta_n^{\dagger})\right) - \frac{8}{\sqrt{n}} \mc V_n^{-1} \D L_n\left(- \mc V_n^{-1} H_n\right)}_2 \\
        & \leq \frac{C}{\sqrt{n}}\left(\frac{\norm{\zeta_n^{\dagger} }_2^2\norm{H_n}_2}{\sqrt{n}} + \frac{\norm{\zeta_n^{\dagger}}_2^4}{n}\right) = C'\left(\frac{\log(n)^{\frac{5}{2}}}{n} + \frac{\log(n)^4}{n^{\frac{3}{2}}}\right),
    \end{align*}
    where $C$ and $C'$ are constants independent of $n$. 
\end{proof}

We are now ready to prove Proposition~\ref{prop:expanFn*}.

\begin{proof}[Proof of Proposition~\ref{prop:expanFn*}]
    By Lemma~\ref{prop:lowerbnd_Fn}, when \eqref{eq:largen_Mn} holds, we have $\zeta_n^{\dagger} \in \inte(\mc Z_n)$ and $\|\zeta_n^{\dagger}\|_2<r_n\leq\bar r_n=O\{\sqrt{\log(n)}\}$. We have
    \begin{align*}
        \max_{\zeta \in \mc Z_n} \left\{- \zeta^\top H_n - F_n(\zeta)\right\} = - (\zeta_n^{\dagger})^ \top H_n -  F_n(\zeta_n^{\dagger}) = \underbrace{- (\zeta_n^{\dagger})^ \top H_n}_{(D)} - \underbrace{\frac{1}{4} \zeta_n^{\dagger\top} \mc V_n \zeta_n^{\dagger}}_{(E)} - \underbrace{\frac{1}{\sqrt{n}} L_n(\zeta_n^{\dagger})}_{(F)},
    \end{align*}
    where the second equality follows from the definition of $F_n(\cdot)$ in
    \eqref{eq:def_Fn}. Next, we substitute the expansion~\eqref{eq:zetan_dagger}
    from Lemma~\ref{lem:expan_subgradient} and compute each term, noting that
    part (A2) of \ref{cond:condA} gives
    $\norm{H_n}_2 \leq a_H\sqrt{\log(n)}$.
    \begin{itemize}
    \item For part (D), by plugging in \eqref{eq:zetan_dagger}, we get
    \begin{align*}
        - (\zeta_n^{\dagger})^ \top H_n = & 2 H_n^\top \mc V_n^{-1} H_n + \frac{8}{\sqrt{n}} H_n^\top \mc V_n^{-1} \D L_n(- \mc V_n^{-1} H_n) - H_n^\top \eps_n^{\dagger} \\
        = & 2 \langle \mc V_n, \xi_n^{\otimes 2} \rangle + \frac{3}{\sqrt{n}} \langle \mc K_n, \xi_n^{\otimes 3} \rangle - H_n^\top \eps_n^{\dagger}.  
        \end{align*} 
        By Lemma~\ref{lem:expan_subgradient}, we have an upper bound of $\norm{\eps_n^{\dagger}}_2$. Thus
        \begin{align*}
            \left|H_n^\top \eps_n^{\dagger}\right| \leq C_1\left(\frac{\log(n)^{3}}{n} + \frac{\log(n)^{\frac{9}{2}}}{n^{\frac{3}{2}}}\right),
        \end{align*}
        where $C_1$ is a constant independent of the sample size $n$. 

    \item For part (E), by plugging in \eqref{eq:zetan_dagger}, we get
    \begingroup\compactdisplaytrue
    \begin{align*}
        & \frac{1}{4} \zeta_n^{\dagger\top} \mc V_n \zeta_n^{\dagger} \\
        =& H_n^\top \mc V_n^{-1} H_n + \frac{8}{\sqrt{n}}H_n^\top \mc V_n^{-1} \D L_n(-\mc V_n^{-1}H_n) - H_n^\top \eps^{\dagger}_n + \frac{16}{n} \D L_n(- \mc V_n^{-1} H_n)^\top \mc V_n^{-1} \D L_n(- \mc V_n^{-1} H_n) \\
        & - \frac{4}{\sqrt{n}} \D L_n(- \mc V_n^{-1} H_n)^\top \eps^{\dagger}_n + \frac14(\eps^{\dagger}_n)^\top\mc V_n\eps^{\dagger}_n\\
        =& \langle \mc V_n, \xi_n^{\otimes 2} \rangle  + \frac{3}{\sqrt{n}} \langle \mc K_n, \xi_n^{\otimes 3} \rangle
        - H_n^\top \eps^{\dagger}_n + \frac{16}{n} \D L_n(- \mc V_n^{-1} H_n)^\top \mc V_n^{-1} \D L_n(- \mc V_n^{-1} H_n) \\
        & - \frac{4}{\sqrt{n}} \D L_n(- \mc V_n^{-1} H_n)^\top \eps^{\dagger}_n + \frac14(\eps^{\dagger}_n)^\top\mc V_n\eps^{\dagger}_n.
    \end{align*}
    \endgroup

    By Lemma~\ref{lem:Ln}, we have
    \begin{align*}
            & \left|\frac{16}{n} \D L_n(- \mc V_n^{-1} H_n)^\top \mc V_n^{-1} \D L_n(- \mc V_n^{-1} H_n)\right| \leq C_2 n^{-1} \log(n)^2,\\
            & \left|\frac{4}{\sqrt{n}} \D L_n(- \mc V_n^{-1} H_n)^\top  \eps^{\dagger}_n\right| \leq C_3 n^{-\half} \log(n) \norm{\eps^{\dagger}_n}_2.
        \end{align*}
    \item For part (F), we have
    \begin{align*}
        & \frac{1}{\sqrt{n}} L_n(\zeta_n^{\dagger}) = \frac{1}{\sqrt{n}} L_n\left(- 2 \mc V_n^{-1} H_n - \frac{8}{\sqrt{n}} \mc V_n^{-1} \D L_n(- \mc V_n^{-1} H_n) + \eps^{\dagger}_n\right) \qquad \text{by~\eqref{eq:zetan_dagger}}\\
        & \frac{1}{\sqrt{n}} L_n\left(- 2 \mc V_n^{-1} H_n\right) = - \frac{1}{\sqrt{n}}  \langle \mc K_n, \xi_n^{\otimes 3} \rangle \qquad \text{by $L_n$: \eqref{eq:expanL}, $\mc K_n$: Theorem~3.1}.
    \end{align*}
    By Lemma~\ref{lem:Ln}, we have
    \begin{align*}
        \norm{- \frac{8}{\sqrt{n}} \mc V_n^{-1} \D L_n(- \mc V_n^{-1} H_n) + \eps^{\dagger}_n}_2 \leq C_4\left(n^{-\frac{1}{2}} \log(n) + \norm{\eps_n^{\dagger}}_2\right) \Let b_n.
    \end{align*}
    Thus, we have
    \begin{align*}
        \left|\frac{1}{\sqrt{n}} L_n(\zeta_n^{\dagger}) - \frac{1}{\sqrt{n}} L_n\left(- 2 \mc V_n^{-1} H_n\right)\right|
        \leq C_5 n^{-\half}\left(\log(n)^2 b_n + \log(n) b_n^2 + b_n^3\right).
    \end{align*}
    \end{itemize}
    Combining these inequalities gives
    \begin{align*}
        \max_{\zeta \in \mc Z_n} \left\{- \zeta^\top H_n - F_n(\zeta)\right\} = \langle \mc V_n, \xi_n^{\otimes 2} \rangle  + \frac{1}{\sqrt{n}} \langle \mc K_n, \xi_n^{\otimes 3} \rangle + \eps^{F}_n,
    \end{align*}
    and we have
    \begingroup\compactdisplaytrue
    \begin{align*}
        |\eps^{F}_n| \leq & 2 \left|H_n^\top \eps^{\dagger}_n\right| + \left|\frac{16}{n} \D L_n(- \mc V_n^{-1} H_n)^\top \mc V_n^{-1} \D L_n(- \mc V_n^{-1} H_n)\right| + \left|\frac{4}{\sqrt{n}} \D L_n(- \mc V_n^{-1} H_n)^\top \eps^{\dagger}_n\right| + \frac14\norm{\mc V_n}_2\norm{\eps^{\dagger}_n}_2^2 \\
        & + \left|\frac{1}{\sqrt{n}} L_n(\zeta_n^{\dagger}) - \frac{1}{\sqrt{n}} L_n\left(- 2 \mc V_n^{-1} H_n\right)\right| \\
        \leq & C_6\left(\frac{\log(n)^{3}}{n} + \frac{\log(n)^{\frac{9}{2}}}{n^{\frac{3}{2}}} + \frac{\log(n)^2}{n} + \frac{\log(n)}{\sqrt{n}}\left(\frac{\log(n)^{3}}{n} + \frac{\log(n)^{\frac{9}{2}}}{n^{\frac{3}{2}}}\right) + \left(\frac{\log(n)^{3}}{n} + \frac{\log(n)^{\frac{9}{2}}}{n^{\frac{3}{2}}}\right)^2 \right. \\
        & \left.+ \frac{\log(n)^2 b_n + \log(n) b_n^2 + b_n^3}{\sqrt{n}}\right)\\ 
        = & C \frac{\log(n)^3}{n},
    \end{align*}
    \endgroup
    where all $C, C_1, \ldots, C_6$ are constants independent of $n$. This completes the proof.
\end{proof}

\subsubsection{Combining steps}\label{sec:comb_step}

Combining the previous five steps completes the proof of the asymptotic expansion in Theorem~3.1.

\begin{proof}[Proof of Theorem~3.1]
Applying Proposition~\ref{prop:rescale}, under Assumptions~2.4 and~2.5, we have
\begin{align}\label{eq:step5supMn}
    n R_n(\mf h) = \sup_{\zeta \in \RR^{d_h}}\{-\zeta^{\top} H_n - M_n(\zeta)\}.
\end{align}
Applying Proposition~\ref{prop:bound_subgradient}, under Assumptions~2.4 and~2.5, there is a deterministic integer $N_1$, such that when $n \geq N_1$, \eqref{eq:largen_Mn} holds, and thus
\begin{align}\label{eq:step5supMn2}
    \eqref{eq:step5supMn} = \sup_{\zeta \in \mc Z_n}\{-\zeta^{\top} H_n - M_n(\zeta)\}.
\end{align}
Applying Proposition~\ref{prop:expanMn}, if Assumptions~2.5,~2.6 and \ref{cond:condA} hold, we obtain
\begin{align}\label{eq:step5supMn3}
    \eqref{eq:step5supMn2} &= \sup_{\zeta \in \mc Z_n} \left\{-\zeta^{\top} H_n - F_n(\zeta) - \eps^M_n(\zeta) \right\} = \sup_{\zeta \in \mc Z_n} \left\{-\zeta^{\top} H_n - F_n(\zeta)\right\} + \eps^{\ddagger}_n,
\end{align}
where for $\zeta \in \mc Z_n$,
\begin{align*}
    & F_n(\zeta) = \frac{1}{4n} \sum_{i=1}^n \zeta^\top \D \mf h(X_i)\Sigma\D \mf h(X_i)^\top \zeta + \frac{1}{\sqrt{n}} L_n(\zeta) ,\\
    & \left|\eps_n^{M}(\zeta)\right| \leq C  \frac{\log(n)^6}{n},
\end{align*}
and $C$ is a constant independent of $n$. As a result, $|\eps^{\ddagger}_n| \leq C  \frac{\log(n)^6}{n}$.

Applying Proposition~\ref{prop:expanFn*}, under \ref{cond:condA} and the same condition~\eqref{eq:largen_Mn}, we obtain
\begin{align}
    \eqref{eq:step5supMn3} = \brak{\mc V_n, \xi_n^{\otimes 2}}  + \frac{1}{\sqrt{n}} \brak{ \mc K_n, \xi_n^{\otimes 3} }  + \eps^{F}_n + \eps^{\ddagger}_n,
\end{align}
where the error term $\eps^{F}_n$ satisfies
\begin{align*}
     \left| \eps^{F}_n \right| \leq C' \frac{\log(n)^3}{n},
  \end{align*}
where $C'$ is a constant independent of $n$. As a result, we obtain~(6).

Finally, Proposition~\ref{prop:tailbnd} gives a deterministic integer $N_2$
such that, when $n\geq N_2$, \ref{cond:condA} holds with probability at least
$1-O(n^{-1})$.  Taking $N=\max\{N_1,N_2\}$ concludes the proof.
\end{proof}

\subsection{Other proofs in Section~3}

\subsubsection{Signed-root reduction for Theorem~3.2}
\label{app:signed_root_edgeworth}

Recall that
\[
    V=\EE[\D\mf h(X)\Sigma\D\mf h(X)^\top],
    \qquad
    W=\EE[\mf h(X)^{\otimes2}],
\]
and define
$\tilde\xi_n=n^{-1/2}\sum_{i=1}^nV^{-1}\mf h(X_i)$.  The expansion in
Theorem~3.1 can be rewritten as
\begin{align}\label{eq:nRn_expan_gamma}
    nR_n(\mf h)
    &=\underbrace{
      \brak{V,\tilde\xi_n^{\otimes2}}
      -\brak{\mc V_n-V,\tilde\xi_n^{\otimes2}}
      +\frac1{\sqrt n}\brak{\EE[\mc K_n],\tilde\xi_n^{\otimes3}}
      }_{(A)}
      +\widetilde{O}_p(n^{-1})
      =\norm{\gamma}_2^2+\widetilde{O}_p(n^{-1}).
\end{align}
Here $(A)$ is expressed as the squared norm of an asymptotically normal vector
$\gamma=(\gamma_j)_{j\in[d_h]}\in\RR^{d_h}$, commonly called the signed-root
statistic; see, e.g.,~\cite[Section~3.2]{hall1990methodology}.  Specifically, for
every $\eta\in\RR^{d_h}$,
\begin{align}\label{eq:defgamma_n}
    \langle\gamma,\eta\rangle
    &=\left\langle V^{\half},\tilde\xi_n\otimes\eta\right\rangle
      -\frac12\left\langle\mc V_n-V,
        \tilde\xi_n\otimes(V^{-\half}\eta)\right\rangle \notag\\
    &\quad+\frac1{2\sqrt n}\left\langle\EE[\mc K_n],
        \tilde\xi_n\otimes\tilde\xi_n\otimes(V^{-\half}\eta)\right\rangle.
\end{align}
Writing
$m_k^{\alpha_1\ldots\alpha_k}
=\EE[\gamma_{\alpha_1}\cdots\gamma_{\alpha_k}]$, its first three moments satisfy
\begin{align}\label{eq:gamma_moment}
    m_1^j
    &=\frac{\mu_{1,n}^j}{\sqrt n},
    &m_2^{jk}
    &=\left(V^{-\half}WV^{-\half}\right)_{jk}+O(n^{-1}),
    &m_3^{jk\ell}
    &=\frac{\mu_{3,n}^{jk\ell}}{\sqrt n},
\end{align}
for $j,k,\ell\in[d_h]$, where $\mu_{1,n}^j$ and
$\mu_{3,n}^{jk\ell}$ are bounded deterministic coefficients.
Appendix~\ref{sec:gamma_moment} gives a deterministic expectation
representation for these coefficients.

Applying the multivariate delta method and the Edgeworth expansion
\cite[Theorem~2(b)]{chandra1979valid} gives
\begingroup\compactdisplaytrue
\begin{align}\label{eq:edgeworth_gamma}
    \PP\opt(\norm{\gamma}_2^2\leq z)
    &=\int_{\norm{v}_2^2\leq z}\tilde\phi(v)\diff v \notag\\
    &\quad+\underbrace{\frac1{\sqrt n}
      \int_{\norm{v}_2^2\leq z}
      \left(
        \sum_{j\in[d_h]}\mu_{1,n}^j\mathrm h_1^j(v)
        +\sum_{j,k,\ell\in[d_h]}\frac{\mu_{3,n}^{jk\ell}}6
          \mathrm h_3^{jk\ell}(v)
      \right)\tilde\phi(v)\diff v}_{(B)=0}
      +O(n^{-1}),
\end{align}
\endgroup
where $\tilde\phi$ is the density of $V^{-\half}H$ for $H\sim \mc N(\mf0,W)$,
and $\mathrm h_1^j$ and $\mathrm h_3^{jk\ell}$ are the first- and third-order
multivariate Hermite polynomials; see~\eqref{eq:hermit_poly}.  Both polynomials
are odd, whereas the integration region and $\tilde\phi$ are centrally symmetric.
Consequently, $(B)=0$, so the null CDF error is of order $n^{-1}$ rather than
$n^{-1/2}$.

\subsubsection{Proof of Theorem~3.2}
\begin{lemma}[Tail bounds of average III]\label{lem:tail_avg3}
    Under Assumptions~2.7, 2.9, and~2.10, let $Y_{VK}$ be the centered joint
    collection of one-observation summands that generate
    $\mc V_n-V$ and $\mc K_n-\EE[\mc K_n]$, represented in an orthonormal
    basis of its covariance support after deterministic linear redundancies
    have been removed.  Fix any constant
    $a_{VK}>\sqrt{2\lambda_{\max}\{\Cov(Y_{VK},Y_{VK})\}}$.  With probability
    at least $1-o(n^{-1})$, we have
    \begin{subequations}
    \begin{align}\label{eq:iq1lem4.17}
         \norm{\mc V_n - V}_F \leq a_{VK} \sqrt{\frac{\log(n)}{n}}, \quad \norm{\mc K_n - \EE[\mc K_n]}_F \leq a_{VK} \sqrt{\frac{\log(n)}{n}}.
    \end{align}
    As a consequence, there exists a deterministic sequence $\bar \delta_n = \widetilde{O}(n^{-1})$, such that
    \begin{align}\label{eq:eq2lem4.17}
        \brak{\mc V_n, \xi_n^{\otimes 2}}  + \frac{1}{\sqrt{n}} \brak{\mc K_n, \xi_n^{\otimes 3}}
        = \norm{\gamma}_2^2 + \bar \eps_n,
    \end{align}
    with $\gamma$ defined by \eqref{eq:defgamma_n} and $|\bar \eps_n| \leq \bar \delta_n$, where $\mc V_n, \mc K_n$ are defined immediately after (6) in Theorem~3.1, and $V, \tilde \xi_n$ are defined in the discussion preceding \eqref{eq:nRn_expan_gamma}.
    \end{subequations}
\end{lemma}
\begin{proof}[Proof of Lemma~\ref{lem:tail_avg3}]
    In this proof,  the constants $C_1, \ldots, C_9$ are independent of the sample size $n$. 
    
    In the covariance-support representation just defined,
    $\Cov(Y_{VK},Y_{VK})$ is positive definite, and the Euclidean norm of
    $Y_{VK}$ agrees with the norm of the full centered collection.  Moreover,
    Assumptions~2.9 and~2.10 give $Y_{VK}$ a finite fourth moment.
    The Cram\'er regularity in Assumption~2.7 supplies the required nonlattice
    condition. Hence these empirical averages admit joint Edgeworth expansions
    up to order $n^{-1}$.
 
    Applying Lemma~\ref{lem:tail_avg2} to $Y_{VK}$ shows that, with probability at least $1-o(n^{-1})$, \eqref{eq:iq1lem4.17} holds, and hence
    \[
        \norm{\mc V_n - V}_2 \leq a_{VK} \sqrt{\frac{\log(n)}{n}}.
    \]
    As a result,
    \begin{align*}
        & \norm{\mc V_n^{-1} - V^{-1}}_2 \leq C_1 \sqrt{\frac{\log(n)}{n}},\\
        & \norm{\mc V_n^{-1} - \left(V^{-1} - V^{-1} (\mc V_n - V) V^{-1}\right)}_2 \leq C_2 \frac{\log(n)}{n}.
    \end{align*}

    To prove the second claim, recall that $\tilde \xi_n = \frac{1}{\sqrt{n}}\sum_{i=1}^n V^{-1} \mf h(X_i)$. By Proposition~\ref{prop:tailbnd},
    \[
        \norm{\frac{1}{\sqrt{n}}\sum_{i=1}^n \mf h(X_i)}_2 \leq a_H\sqrt{\log(n)},
    \]
    and thus
    \[
        \norm{\xi_n - \tilde \xi_n}_2 \leq C_3 \frac{\log(n)}{\sqrt{n}}.
    \]

    Finally, we have
    \begin{align*}
        &\left|\left\langle \mc V_n, \xi_n^{\otimes 2} \right\rangle - \left\langle V, \tilde \xi_n^{\otimes 2}\right\rangle + \left\langle \mc V_n - V, \tilde \xi_n^{\otimes 2} \right\rangle\right|\\
        =& \left|\left\langle \mc V_n^{-1}, \left(\frac{1}{\sqrt{n}}\sum_{i=1}^n \mf h(X_i)\right)^{\otimes 2}\right\rangle - \left\langle V, \tilde \xi_n^{\otimes 2}\right\rangle + \left\langle \mc V_n - V, \tilde \xi_n^{\otimes 2} \right\rangle\right|\\
        =& \left|\left\langle\mc V_n^{-1} - \left(V^{-1} - V^{-1} (\mc V_n - V) V^{-1}\right), \left(\frac{1}{\sqrt{n}}\sum_{i=1}^n \mf h(X_i)\right)^{\otimes 2}\right\rangle\right| \leq C_4 \frac{\log(n)^2}{n},
    \end{align*}
    and
    \begingroup\compactdisplaytrue
    \begin{align*}
         &\left|\frac{1}{\sqrt{n}} \left\langle \mc K_n, \xi_n^{\otimes 3} \right\rangle - \frac{1}{\sqrt{n}} \left\langle \EE[\mc K_n], \tilde \xi_n^{\otimes 3} \right\rangle\right|\\
         \leq & \left|\frac{1}{\sqrt{n}} \left\langle \EE[\mc K_n], \xi_n^{\otimes 3} - \tilde \xi_n^{\otimes 3}\right\rangle\right| + \left|\frac{1}{\sqrt{n}} \left\langle \mc K_n - \EE[\mc K_n], \xi_n^{\otimes 3} \right\rangle\right| \\
         \leq & \left|\frac{1}{\sqrt{n}} \left\langle \EE[\mc K_n], \xi_n^{\otimes 3} - \tilde \xi_n^{\otimes 3}\right\rangle\right| + \left|\frac{1}{\sqrt{n}} \left\langle \mc K_n - \EE[\mc K_n], \tilde \xi_n^{\otimes 3} \right\rangle\right| + \left|\frac{1}{\sqrt{n}} \left\langle \mc K_n - \EE[\mc K_n], \xi_n^{\otimes 3} - \tilde \xi_n^{\otimes 3} \right\rangle\right| \\
         \leq & C_5\frac{\log(n)^{\frac{7}{2}}}{n}.
    \end{align*}
    \endgroup
    Therefore, we obtain 
    \[
    \left\langle \mc V_n, \xi_n^{\otimes 2} \right\rangle  + \frac{1}{\sqrt{n}} \left\langle \mc K_n, \xi_n^{\otimes 3} \right\rangle = 
    \left\langle V, \tilde \xi_n^{\otimes 2}\right\rangle - \left\langle \mc V_n - V, \tilde \xi_n^{\otimes 2} \right\rangle + \frac{1}{\sqrt{n}} \left\langle \EE[\mc K_n], \tilde \xi_n^{\otimes 3} \right\rangle + \eps'_n, 
    \]
    with $|\eps'_n| \leq C_6 \log(n)^{\frac{7}{2}} n^{-1}$. Finally, note that
    \begin{align}\label{eq:Vn-VEKn}
        \left|\left \langle \mc V_n - V, \tilde \xi_n^{\otimes 2} \right\rangle \right| \leq C_7 \frac{\log(n)}{\sqrt{n}},\quad \left|\frac{1}{\sqrt{n}} \left\langle \EE[\mc K_n], \tilde \xi_n^{\otimes 3} \right\rangle\right| \leq C_8 \frac{\log(n)^{\frac{3}{2}}}{\sqrt{n}},
    \end{align}
    we can write 
    \begin{align*}
        \left\langle V, \tilde \xi_n^{\otimes 2}\right\rangle - \left\langle \mc V_n - V, \tilde \xi_n^{\otimes 2} \right\rangle + \frac{1}{\sqrt{n}} \left\langle \EE[\mc K_n], \tilde \xi_n^{\otimes 3} \right\rangle = \norm{\gamma}_2^2 + \bar \eps_n^{(0)},
    \end{align*}
    for the vector $\gamma$ defined by \eqref{eq:defgamma_n}, where
    $|\bar \eps_n^{(0)}|=O\{\log(n)^3n^{-1}\}$. Combining this relation
    with the preceding display and absorbing $\eps'_n+\bar\eps_n^{(0)}$
    into $\bar\eps_n$ gives~\eqref{eq:eq2lem4.17}, with
    $|\bar\eps_n|\leq C_9\log(n)^{7/2}n^{-1}$.
\end{proof}

\begin{proof}[Proof of Theorem~3.2]
    By Theorem~3.1 and Lemma~\ref{lem:tail_avg3}, with probability at least $1 - O(n^{-1})$,
    \begin{align*}
        n R_n(\mf h) &= \left\langle \mc V_n, \xi_n^{\otimes 2} \right\rangle  + \frac{1}{\sqrt{n}} \left\langle \mc K_n, \xi_n^{\otimes 3} \right\rangle + \eps_n = \underbrace{\norm{\gamma}_2^2}_{\text{(A')}} + \underbrace{\bar \eps_n + \eps_n}_{\text{(B')}},
    \end{align*}
    where $\gamma$ is defined in \eqref{eq:defgamma_n}, and $|\bar \eps_n + \eps_n| \leq \delta_n$ for a deterministic sequence $\delta_n$ of order $\widetilde{O}(n^{-1})$.
    
    Now consider the Edgeworth expansion of part (A'). By Assumption~2.9, we have $\EE\brac{\norm{\mf h(X)}_2^4} < \infty$ and $\EE\brac{\norm{\D \mf h(X)}_2^8}<\infty$. The Cram\'er regularity in Assumption~2.7 therefore permits an Edgeworth expansion for $\gamma$ with error of order $n^{-1}$. Thus, \eqref{eq:edgeworth_gamma} is valid.
    
    To deal with part (B'), write $Q_n=\norm{\gamma}_2^2$ and
    $r_n=\bar\eps_n+\eps_n$, and let $\mc E_n$ be the event on which
    $nR_n(\mf h)=Q_n+r_n$ and $|r_n|\leq\delta_n$. The preceding probability
    bound gives $\PP\opt(\mc E_n^c)=O(n^{-1})$, and, for every $z$,
    \begin{align}\label{eq:deltamethod}
        \PP\opt(Q_n\leq z-\delta_n)-\PP\opt(\mc E_n^c)
        \leq \PP\opt(nR_n(\mf h)\leq z)
        \leq \PP\opt(Q_n\leq z+\delta_n)+\PP\opt(\mc E_n^c).
    \end{align}
    Fix $z_0>0$. For all sufficiently large $n$, $\delta_n\leq z_0/2$.
    The distribution $\Psi$ has a continuous density on $(0,\infty)$ that is
    bounded on $[z_0/2,\infty)$; hence, for some constant $C_{z_0}$,
    \[
        \sup_{z\geq z_0}|\Psi(z\pm\delta_n)-\Psi(z)|
        \leq C_{z_0}\delta_n=\widetilde{O}(n^{-1}).
    \]
    Combining this bound, \eqref{eq:deltamethod}, and the expansion
    \eqref{eq:edgeworth_gamma} yields
    \begin{align}\label{eq:edgeworth_nRn_2}
        \sup_{z\geq z_0}\left|
        \PP\opt(nR_n(\mf h)\leq z)-\Psi(z)
        \right|=\widetilde{O}(n^{-1}),
    \end{align}
    which completes the proof.
\end{proof}

\subsubsection{Proof of Theorem~3.3}\label{sec:pf_localalter_expan}\label{app:local_alternative_conditions}
The stochastic expansion under local alternatives uses the following uniform
version of the null regularity conditions.

\begin{assumption}[Uniform regularity under local alternatives]\label{a:alter_edgeworth}
Assume that
\begin{enumerate}[label=(\roman*)]
    \item $(\PP_n\opt)_{n\geq1}$ is contained in a class $\mc P_c$ that is
    compact in total variation, and every element of $\mc P_c$ satisfies
    Assumption~2.7.

    \item $\EE_n[\mf h(X)]=O(n^{-1/2})$ and
    $\EE_n[\mf h(X)^{\otimes2}]\in\PD$ for every $n$.  Moreover, for some
    $\sigma_0>0$,
    \begin{align*}
        \lambda_{\min}\!\left(\EE_n[\D\mf h(X)\Sigma\D\mf h(X)^\top]\right)
        \wedge
        \lambda_{\min}\!\left(\EE_n[\mf h(X)^{\otimes2}]\right)
        \geq \sigma_0
        \qquad\forall n\geq1.
    \end{align*}

    \item The required moments are bounded uniformly:
    \[
        \sup_{n\geq1}\EE_n\!\left[
        \left(\norm{\mf h(X)}_2+\norm{\D\mf h(X)}_2
        +\kappa_1(X)+\kappa_2(X)\right)^8\right]<\infty.
    \]
\end{enumerate}
\end{assumption}

The compactness condition yields the uniform Cram\'{e}r bound used below.  For
example, total-variation convergence follows from Scheff\'{e}'s lemma
\cite{scheffe1947useful} when the laws have densities $f_n$ converging almost
everywhere to a density $f_\infty$.  The eigenvalue bounds prevent $\xi_n$ and
$\mc V_n^{-1}$ from becoming unstable, while the uniform moments control the
remainder independently of $n$.  These conditions hold, for instance, when
$\PP_n\opt$ and the relevant moments converge to a null law satisfying
Assumptions~2.8 and~2.9, with
$\EE_n[\mf h(X)]=O(n^{-1/2})$.

\begin{lemma}[Uniform Cram\'er condition]\label{lem:cramerII}
    Assume that $Y_n$ and $Y_{\infty}$ have laws $\PP_{n}\opt$ and $\PP\opt$, respectively, and that every $Y_n$ and $Y_{\infty}$ satisfies the Cram\'er condition. If $\PP_{n}\opt$ converges to $\PP\opt$ in total variation as $n\to\infty$, then $(\PP_{n}\opt)_{n \geq 1}$ satisfies the uniform Cram\'er condition: for every $b >0$,
    \begin{align*}
        \limsup_{n\rightarrow\infty}\sup_{t: \norm{t}_2 \geq b}\abs{\EE_{n}\brac{\exp\pare{\mathrm{i}t^\top Y_{n}}}} < 1.
    \end{align*}
\end{lemma}
\begin{proof}[Proof of Lemma~\ref{lem:cramerII}]
    In the following, we use $\EE_n$ and $\EE_{\infty}$ to denote the expectation under $\PP\opt_{n}$ and $\PP\opt$, respectively. We prove this by contradiction.
    Suppose that there is a sequence of $(n_k)_{k \geq 1}$ and $b_0 > 0$, such that 
    \[
    \lim_{k\rightarrow\infty}\sup_{t: \norm{t}_2 \geq b_0}\abs{\EE_{n_k}\brac{\exp\pare{\mathrm{i}t^\top Y_{n_k}}}} = 1.
    \]
    Because $\PP_{n}\opt$ converges to $\PP_{\infty}\opt$ in total variation, uniformly over $t$,
    \[
    \lim_{k\rightarrow\infty}\abs{\EE_{n_k}\brac{\exp\pare{\mathrm{i}t^\top Y_{n_k}}} - \EE_{\infty}\brac{\exp\pare{\mathrm{i}t^\top Y_{\infty}}}} \leq \lim_{k\rightarrow\infty} \text{TV}\pare{\PP_{n_k}\opt, \PP_{\infty}\opt} = 0,
    \]
    where $\text{TV}(\PP,\QQ)$ denotes the total variation distance between $\PP$ and $\QQ$. Consequently,
    $\sup_{t: \norm{t}_2 \geq b_0}\abs{\EE_{\infty}\brac{\exp\pare{\mathrm{i}t^\top Y_{\infty}}}} = 1$,
    contradicting the Cram\'er condition for $Y_{\infty}$ under $\PP_{\infty}\opt$.
\end{proof}

\begin{proof}[Proof of Theorem~3.3]
    Note that the proof of Theorem~3.1 relies only on Assumptions~2.4--2.6 once
    \ref{cond:condA} is available. Proposition~\ref{prop:tailbnd} establishes
    that condition with probability $1-O(n^{-1})$ under the null. Therefore, in
    the setting of Theorem~3.3, it remains to establish the same probability
    bound under Assumption~\ref{a:alter_edgeworth}.

    Items 2 and 3 of Assumption~\ref{a:alter_edgeworth} provide the uniform
    moment bounds needed to repeat the argument establishing
    \ref{cond:condA} in Proposition~\ref{prop:tailbnd}, except for
    part (A2).

    For part (A2), we establish an Edgeworth expansion for
    $\frac{1}{\sqrt{n}}\sum_{i=1}^n (\mf h(X_i) - \EE_{n}[\mf h(X)])$
    under the local alternative $\PP_{n}\opt$:
    \begin{itemize}[leftmargin=0.24in]
        \item Let $X(n)$ have law $\PP_{n}\opt$. By item 1 of Assumption~\ref{a:alter_edgeworth}, $\mf h(X(n))$ satisfies the Cram\'er condition under every $\PP_n\opt$. Since $(\PP_{n}\opt)_{n \geq 1} \subset \mc P_c$ and $\mc P_c$ is compact in total variation, a subsequence argument using Lemma~\ref{lem:cramerII} gives the required uniform Cram\'er condition.
        \item By item 2 of Assumption~\ref{a:alter_edgeworth} and
        $\EE_n[\mf h(X)]=O(n^{-1/2})$, for all sufficiently large $n$,
        \[
            \lambda_{\min}\!\left\{\Cov_n(\mf h(X),\mf h(X))\right\}
            \geq
            \lambda_{\min}\!\left\{\EE_n[\mf h(X)^{\otimes2}]\right\}
            -\norm{\EE_n[\mf h(X)]}_2^2
            \geq\frac{\sigma_0}{2}.
        \]
        \item By item 3 of Assumption~\ref{a:alter_edgeworth}, $\sup_{n \geq 1}\EE_{n}{\norm{\mf h(X)}_2^4} < \infty$, and
        \[
        \EE_{n}\brac{\norm{\mf h(X)}_2^4; \norm{\mf h(X)}_2 \geq \eps n^{\half}}
        \leq \frac{1}{\eps^4n^2}
        \sup_{m\geq1}\EE_m\brac{\norm{\mf h(X)}_2^{8}}
        \longrightarrow0.
        \]
    \end{itemize}

    Therefore, \cite[Theorem 20.6]{bhattacharya2010normal} implies that
    $\frac{1}{\sqrt{n}}\sum_{i=1}^n (\mf h(X_i) - \EE_{n}[\mf h(X)])$ admits an Edgeworth expansion under the local alternative $\PP_{n}\opt$. Letting $\Lambda_{H,\mathrm{loc}}=\sup_{m\geq1}\lambda_{\max}\{\Cov_m(\mf h(X),\mf h(X))\}<\infty$ and fixing any $a_{H,\mathrm{loc}}>\sqrt{2\Lambda_{H,\mathrm{loc}}}$, Lemma~\ref{lem:tail_avg2} gives, with probability at least $1-o(n^{-1})$,
    \[
        \norm{\frac{1}{\sqrt{n}}\sum_{i=1}^n \pare{\mf h(X_i) - \EE_{n}\brac{\mf h(X)}}}_2 \leq a_{H,\mathrm{loc}} \sqrt{\log(n)}.
    \]
    By item 2 of Assumption~\ref{a:alter_edgeworth}, $\norm{\EE_{n}\brac{\mf h(X)}}_2 = O\pare{n^{-\half}}$. Thus, we get
    \[
        \norm{\frac{1}{\sqrt{n}}\sum_{i=1}^n \mf h(X_i)}_2 \leq (a_{H,\mathrm{loc}}+1) \sqrt{\log(n)},
    \]
    when $n$ is sufficiently large. The change of constant in front of $\sqrt{\log(n)}$ does not affect the proof of Theorem~3.1. Therefore, we obtain the desired result by following the same proof of Theorem~3.1.
\end{proof}

\subsection{Proofs in Section~3.3}
\begin{lemma}[Gradient of $z_{1-\varrho}$]\label{lem:gradientZn}
    Let $z_{W,V} = F_{W,V}^{-1}(1-\varrho)$ for $W,V \in \RR^{d_h\times d_h}$, where
    \begin{align*}
        F_{W,V}(z) \Let \int_{\left\{v \in \RR^{d_h}:~v^\top W^{\half} V^{-1} W^{\half} v \leq z\right\}} \phi(v) \diff v,
    \end{align*}
    then, for $W,V \in \PD$, the function $z_{W,V}$ is twice continuously differentiable with respect to $(W,V)$, and
    \begingroup\compactdisplaytrue
    \begin{align}
         - \frac{\diff z_{W,V}}{z_{W,V}} 
        = \frac{\left\langle L_{W,V}, \diff \left( V^\half W^{-1} V^{\half} \right) \right\rangle}{\left\langle L_{W,V}, V^\half W^{-1} V^{\half}\right \rangle} 
        \Let  \mc L_{W,V}(\diff W, \diff V), \label{eq:mcLn}
    \end{align}
    where
    \begin{align}
        &\diff \left( V^\half W^{-1} V^{\half} \right) \notag\\
        =&\int_{0}^\infty e^{-s V^\half} \left(\diff V\right)  e^{-s V^\half} W^{-1} V^{\half} \diff s -V^\half W^{-1} \left(\diff W\right) W^{-1} V^{\half}
        +\int_{0}^\infty V^\half W^{-1} e^{-s V^\half} \left(\diff V\right)  e^{-s V^\half}\diff s,\notag\\
        &L_{W,V} \notag\\
        =&\int_{\norm{v}_2^2 \leq z_{W,V}} (2\pi)^{-\frac{d_h}{2}} \exp\left(-\frac{1}{2} v^\top V^\half W^{-1} V^{\half} v\right) \sqrt{\det \left(V^\half W^{-1} V^{\half}\right)} \left( -\half v v^\top + \half V^{-\half} W V^{-\half} \right)\diff v. \label{eq:LWV}
    \end{align}
    \endgroup
\end{lemma}
\begin{proof}[Proof of Lemma~\ref{lem:gradientZn}]
    Under the change of variables $v \leftarrow z^{-\half} V^{-\half}W^{\half}v$, we can write $F_{W,V}$ as
    \begin{align*}
        F_{W,V}(z) &=  \int_{\norm{v}_2^2 \leq 1} (2\pi)^{-\frac{d_h}{2}} \exp\left(-\frac{z}{2} v^\top V^\half W^{-1} V^{\half} v\right) \sqrt{\det \left(V^\half W^{-1} V^{\half}\right)} z^{d_h/2} \diff v \\
        &\Let \int_{\norm{v}_2^2 \leq 1} f_{W,V,z}(v) \diff v.
    \end{align*}
    On compact subsets of $\PD\times\PD\times(0,\infty)$,
    $f_{W,V,z}(v)$ is twice continuously differentiable with respect to
    $(W,V,z)$, and its derivatives through order two are bounded by an
    integrable envelope. Applying \cite[Theorem 2.27]{folland1999real} twice
    shows that $F_{W,V}(z)$ is twice continuously differentiable with respect
    to $(W,V,z)$. Further, $F'_{W,V}(z) \neq 0$ for $z > 0$.

    By definition, $F_{W,V}(z_{W,V})-(1-\varrho)=0$. The implicit function
    theorem \cite[Theorem 3.2.1]{krantz2002implicit} therefore shows that
    $z_{W,V}$ is twice continuously differentiable with respect to $(W,V)$, and
    \begin{align*}
        \mf 0_{d_h \times d_h} &= F'_{W,V}(z_{W,V}) \frac{\diff z_{W,V}}{\diff W} + \int_{\norm{v}_2^2 \leq 1} \frac{\diff f_{W,V,z_{W,V}}(v)}{\diff W} \diff v,\\
        \mf 0_{d_h \times d_h} &= F'_{W,V}(z_{W,V}) \frac{\diff z_{W,V}}{\diff V} + \int_{\norm{v}_2^2 \leq 1} \frac{\diff f_{W,V,z_{W,V}}(v)}{\diff V} \diff v,
    \end{align*}
    where $\mf 0_{d_h \times d_h}$ is the $d_h \times d_h$ zero matrix. As a result, we have 
    \begingroup\compactdisplaytrue
    \begin{align*}
        \diff z_{W,V} =& - \frac{1}{F'_{W,V}(z_{W,V})}\int_{\norm{v}_2^2 \leq z_{W,V}} (2\pi)^{-\frac{d_h}{2}} \exp\left(-\frac{1}{2} v^\top V^\half W^{-1} V^{\half} v\right) \sqrt{\det \left(V^\half W^{-1} V^{\half}\right)} \\
        & \brak{ -\half v v^\top + \half V^{-\half} W V^{-\half}, \diff \left( V^\half W^{-1} V^{\half} \right)} \diff v\\
        =&  - \frac{1}{F'_{W,V}(z_{W,V})} 
        \left\langle L_{W,V}, \diff \left( V^\half W^{-1} V^{\half} \right) \right \rangle,
    \end{align*}
    \endgroup
    where $L_{W,V}$ is defined in \eqref{eq:LWV} and we apply the differentiation equations: 
    \begin{align*}
        &\diff \exp\left(-\half v^\top X v\right) = \exp\left(-\half v^\top X v\right)\left\langle -\half v v^\top, \diff X \right \rangle, \\
        &\diff \sqrt{\det X} = \sqrt{\det X} \left\langle \half X^{-\top}, \diff X \right \rangle.
    \end{align*}

    As for $F'_{W,V}(z_{W,V})$, we have
    \begin{align*}
        F'_{W,V}(z_{W,V}) =& \int_{\norm{v}_2^2 \leq z_{W,V}} (2\pi)^{-\frac{d_h}{2}} \exp\left(-\frac{1}{2} v^\top V^\half W^{-1} V^{\half} v\right) \sqrt{\det \left(V^\half W^{-1} V^{\half}\right)}  \\
        & \left(-\frac{1}{2z_{W,V}} v^\top V^\half W^{-1} V^{\half} v + \frac{d_h}{2z_{W,V}}\right) \diff v\\
        =& \frac{1}{z_{W,V}} \left\langle L_{W,V}, V^\half W^{-1} V^{\half} \right \rangle.
    \end{align*}
    
    As for $\diff \left( V^\half W^{-1} V^{\half} \right)$, we have
    \begin{align*}
        \diff \left( V^\half W^{-1} V^{\half} \right) &= \left(\diff V^\half\right) W^{-1} V^{\half} +  V^\half \left(\diff  W^{-1}\right) V^{\half} +  V^\half W^{-1} \diff \left(V^{\half} \right),
    \end{align*}
    where $\diff  W^{-1}$ and $\diff V^{\half}$ are
    \begin{align*}
        \diff  W^{-1} = - W^{-1} \left(\diff W\right) W^{-1}, \qquad 
        \diff V^{\half} = \int_{0}^\infty e^{-s V^\half} \left(\diff V\right)  e^{-s V^\half}\diff s.
    \end{align*}
    The first identity follows from the Neumann series, while the second follows
    from \cite[Theorem 1.1]{del2018taylor} for a symmetric infinitesimal matrix
    $\diff V$.  It is straightforward to verify that $\diff z_{W,V}/\diff W$
    and $\diff z_{W,V}/\diff V$ are continuously differentiable in $(W,V)$;
    for the second-order derivative of $V^{\half}$, see
    \cite[Theorem 1.1]{del2018taylor}.  This completes the proof.
\end{proof}
\subsubsection{Quantile expansion for the size theorem}
\label{app:size_quantile_expansion}

\begin{proposition}[Second-order expansion of $\hat z_{1-\varrho}$]\label{prop:expanZn}
    Under Assumptions~2.7--2.9, there is a deterministic sequence $\tilde \delta_n = \widetilde{O}\left(n^{-1}\right)$ and an integer $N$ such that, when $n \geq N$, we have
    \begin{align*}
        \hat z_{1-\varrho} = z_{1-\varrho}\left(1- \mc L_{W,V}\pare{\mc W_n - W, \mc V_n - V}\right) + \tilde \eps_n,
    \end{align*}
    where $\PP\opt(|\tilde \eps_n| \leq \tilde \delta_n) = 1 - O(n^{-1})$, $W= \EE\left[\mf h(X)^{\otimes 2}\right], V = \EE\left[\D \mf h(X)\Sigma\D\mf h(X)^\top\right]$, and $\mc L_{W,V}(\cdot, \cdot)$ is the linear functional defined in~\eqref{eq:mcLn}.
\end{proposition}

\begin{proof}[Proof of Proposition~\ref{prop:expanZn}]
    Assumption~2.9 gives finite eighth moments for both $\mf h(X)$ and
    $\D\mf h(X)$. Let $Y_{WV}$ be the centered joint collection of
    one-observation summands that generate $\mc W_n-W$ and $\mc V_n-V$,
    represented in an orthonormal basis of its covariance support after
    deterministic linear redundancies have been removed, and fix any
    $a_{WV}>\sqrt{2\lambda_{\max}\{\Cov(Y_{WV},Y_{WV})\}}$. Assumption~2.7
    supplies the Cram\'er condition, so Lemma~\ref{lem:tail_avg2} yields the
    following bounds with probability at least $1-o(n^{-1})$:
    \begin{align}\label{eq:bndVnWn}
    \begin{split}
        &\norm{\frac{1}{\sqrt{n}}\sum_{i=1}^n \left(\mf h(X_i)^{\otimes 2} - \EE\left[\mf h(X)^{\otimes 2}\right]\right) }_F \leq a_{WV} \sqrt{\log(n)},\\
        &\norm{\frac{1}{\sqrt{n}}\sum_{i=1}^n \left(\D \mf h(X_i)\Sigma\D \mf h(X_i)^\top  - \EE\left[\D \mf h(X)\Sigma \D \mf h(X)^\top\right]\right) }_F \leq a_{WV} \sqrt{\log(n)}.
    \end{split}
    \end{align}
    Consequently, $\norm{\mc W_n-W}_2$ and $\norm{\mc V_n-V}_2$ are of order $\sqrt{\log(n)/n}$, where $W=\EE[\mf h(X)^{\otimes2}]$ and $V=\EE[\D\mf h(X)\Sigma\D\mf h(X)^\top]$. On this event, $\mc W_n$ and $\mc V_n$ are invertible for all sufficiently large $n$, and hence $\hat z_{1-\varrho}$ is well-defined. By Lemma~\ref{lem:gradientZn}, for $\hat z_{1-\varrho} = F_{\mc W_n, \mc V_n}^{-1}(1-\varrho)$ and $z_{1-\varrho} = F_{W, V}^{-1}(1-\varrho)$, we have
    \begin{align*}
        &\left|\hat z_{1-\varrho} - z_{1-\varrho}\left(1- \mc L_{W,V}(\mc W_n - W, \mc V_n - V)\right)\right|\\
        = &\left|\hat z_{1-\varrho} - z_{1-\varrho} - \brak{\frac{\diff z_{W,V}}{\diff W}, \mc W_n - W} - \brak{\frac{\diff z_{W,V}}{\diff V}, \mc V_n - V}\right|\\
        \leq &C \left(\norm{\mc W_n - W}_2^2 + \norm{\mc V_n - V}_2^2\right)\\
        \leq &2Ca_{WV}^2\log(n)n^{-1},
    \end{align*}
    where $\mc L_{W,V}(\cdot, \cdot)$ is defined in~\eqref{eq:mcLn}, and $C$ is a constant independent of $n$.
\end{proof}

\subsubsection{Proof of Theorem~3.4}

\begin{proof}[Proof of Theorem~3.4]
This proof follows the pattern of the proof of Theorem~3.2.  With probability
at least $1-O(n^{-1})$ under $\PP\opt$, we have
\[n R_n(\mf h) = \norm{\gamma}_2^2 + \bar \eps_n + \eps_n.\]
Also, by Proposition~\ref{prop:expanZn}, with probability at least $1 - O(n^{-1})$,
\[\hat z_{1-\varrho} = z_{1-\varrho}\left(1- \mc L_{W,V}(\mc W_n - W, \mc V_n - V)\right)+ \tilde \eps_n.\]
Here, we have
\begin{align}\label{eq:3eps}
    |\bar \eps_n| + |\eps_n| + |\tilde \eps_n| \leq \delta_n   
\end{align}
for a deterministic sequence $\delta_n = \widetilde{O}(n^{-1})$. As a consequence,
\begin{align}\label{eq:nRn_expand1}
    n R_n(\mf h) - \hat z_{1-\varrho}
    =\left(\norm{\gamma}_2^2 - z_{1-\varrho}\left(1- \mc L_{W,V}(\mc W_n - W, \mc V_n - V)\right)\right) + \bar \eps_n + \eps_n - \tilde \eps_n.
\end{align}
Let $a_n=\mc L_{W,V}(\mc W_n-W,\mc V_n-V)$ and define
$\tilde\gamma=(1+a_n/2)\gamma$.  Since $a_n=\widetilde{O}_p(n^{-1/2})$,
\begin{align}\label{eq:nRn-hatz}
    nR_n(\mf h)-\hat z_{1-\varrho}
    =(1-a_n)\left(\norm{\tilde\gamma}_2^2-z_{1-\varrho}\right)
      +\widetilde{O}_p(n^{-1}).
\end{align}
The remaining argument makes this reduction uniform on the required high-probability event
and applies the Edgeworth expansion to the corresponding signed root.

Under Assumption~2.9, $\EE\left[\norm{\mf h(X)}_2^8\right] < \infty, \EE\left[\norm{\D \mf h(X)}_2^8\right] < \infty$. Using Lemma~\ref{lem:tail_avg2}, by the same reasoning as in the proof of Proposition~\ref{prop:expanZn}, we have $\norm{\mc W_n - W}_2, \norm{\mc V_n - V}_2$ are of order $\sqrt{\log(n)} n^{-\half}$ with probability at least $1 - o(n^{-1})$.

Therefore, there are constants $C_1, C_2$ independent of $n$, such that
\begin{align}
    \abs{\mc L_{W,V}(\mc W_n - W, \mc V_n - V)} \leq C_1 \sqrt{\frac{\log(n)}{n}}, \label{eq:bndLWV}
\end{align}
and thus
\[\left|\left(1- \mc L_{W,V}(\mc W_n - W, \mc V_n - V)\right)^{-\half} - \left(1 + \half \mc L_{W,V}(\mc W_n - W, \mc V_n - V)\right) \right| \leq C_2 \frac{\log(n)}{n}.\]

Let $\Gamma_\gamma=V^{-\half}WV^{-\half}$ and fix any
$a_\gamma>\sqrt{2\lambda_{\max}(\Gamma_\gamma)}$. Since $\gamma$ has an
Edgeworth expansion up to order $n^{-1}$ (see the proof of Theorem~3.2), the
argument of Lemma~\ref{lem:tail_avg2} gives, with probability at least
$1-o(n^{-1})$, $\norm{\gamma}_2\leq a_\gamma\sqrt{\log(n)}$.

Thus, we have
\begingroup\compactdisplaytrue
\begin{align}
    \eqref{eq:nRn_expand1} =& \left(1- \mc L_{W,V}(\mc W_n - W, \mc V_n - V)\right)\left(\norm{\left(1- \mc L_{W,V}(\mc W_n - W, \mc V_n - V)\right)^{-\half}\gamma}_2^2 - z_{1-\varrho}\right) \notag\\
    &+ \bar \eps_n + \eps_n - \tilde \eps_n \notag\\
    =& \left(1- \mc L_{W,V}(\mc W_n - W, \mc V_n - V)\right)\left(\norm{\left(1 + \half \mc L_{W,V}(\mc W_n - W, \mc V_n - V)\right)\gamma}_2^2 - z_{1-\varrho} + \eps'_n\right) \notag\\
    &+ \bar \eps_n + \eps_n - \tilde \eps_n,\label{eq:nRnexpand2}
\end{align}
\endgroup
with the newly added error $\eps'_n$ satisfying $|\eps'_n| \leq C_3 \log(n)^3 n^{-1}$ for a constant $C_3$.

We further expand the term $\left(1 + \half \mc L_{W,V}(\mc W_n - W, \mc V_n - V)\right)\gamma$ to remove the higher order terms (higher than $n^{-\half}$). Concretely, there is a vector $\gamma^{\dagger}$ defined by
\begingroup\compactdisplaytrue
\begin{align}\label{eq:defgammada_n}
    \langle \gamma^{\dagger}, \eta \rangle =& \left\langle V^{\half}, \tilde \xi_n \otimes \eta\right\rangle + \half \mc L_{W,V}(\mc W_n - W, \mc V_n - V) \left\langle V^{\half}, \tilde \xi_n \otimes \eta\right\rangle - \frac{1}{2} \left\langle \mc V_n - V, \tilde \xi_n \otimes \left(V^{-\half} \eta\right) \right\rangle \notag\\
    & + \frac{1}{2\sqrt{n}} \left\langle \EE[\mc K_n], \tilde \xi_n \otimes \tilde \xi_n \otimes \left(V^{-\half} \eta\right) \right\rangle \quad \forall \eta \in \RR^{d_h}.
\end{align}
\endgroup
Further, by \eqref{eq:Vn-VEKn} and \eqref{eq:bndLWV}, 
\begin{align}\label{eq:eps4}
    \norm{\left(1 + \half \mc L_{W,V}(\mc W_n - W, \mc V_n - V)\right)\gamma}_2^2 =  \norm{\gamma^{\dagger}}_2^2 + \eps^{\dagger}_n,~~ |\eps^{\dagger}_n| \leq C_4 \frac{\log(n)^4}{n}.
\end{align}
Therefore, we obtain
\begin{align}\label{eq:nRnexpand3}
    \eqref{eq:nRnexpand2} = \left(1- \mc L_{W,V}(\mc W_n - W, \mc V_n - V)\right)\left(\norm{\gamma^{\dagger}}_2^2 - z_{1-\varrho} + \eps^{\dagger}_n\right) + \bar \eps_n + \eps_n - \tilde \eps_n.
\end{align}
When $n$ is sufficiently large, $|\mc L_{W,V}(\mc W_n - W, \mc V_n - V)| \leq \half$; thus $1- \mc L_{W,V}(\mc W_n - W, \mc V_n - V) > 0$.

Finally, following the delta method again as in \eqref{eq:deltamethod}, we have
\begingroup\compactdisplaytrue
\begin{align*}
    &\PP\opt\left(n R_n(\mf h) - \hat z_{1-\varrho} \leq 0\right) & \\
    =& \PP\opt\left(\left(1- \mc L_{W,V}(\mc W_n - W, \mc V_n - V)\right)\left(\norm{\gamma^{\dagger}}_2^2 - z_{1-\varrho} + \eps^{\dagger}_n\right) + \bar \eps_n + \eps_n - \tilde \eps_n \leq 0\right)  &\text{by \eqref{eq:nRnexpand3}}\\
    \leq&  \PP\opt\left(\left(1- \mc L_{W,V}(\mc W_n - W, \mc V_n - V)\right)\left(\norm{\gamma^{\dagger}}_2^2 - z_{1-\varrho} + \eps^{\dagger}_n\right) \leq 0\right) + \widetilde{O}\left(\frac{1}{n}\right)  &\text{by \eqref{eq:3eps}}\\
    \leq& \PP\opt\left(\norm{\gamma^{\dagger}}_2^2 - z_{1-\varrho} + \eps^{\dagger}_n \leq 0\right) + \widetilde{O}\left(\frac{1}{n}\right) &\text{$|\mc L_{W,V}(\cdot,\cdot) | \leq \half$}\\
    \leq& \PP\opt\left(\norm{\gamma^{\dagger}}_2^2 - z_{1-\varrho} \leq 0\right) + \widetilde{O}\left(\frac{1}{n}\right)  &\text{by \eqref{eq:eps4}}.
\end{align*}
\endgroup
The remaining steps and the reverse inequality are similar to those in \eqref{eq:deltamethod}.

For the Edgeworth expansion of $\gamma^{\dagger}$, note that its moment expansion has the same form as that of $\gamma$ in~\eqref{eq:gamma_moment}. Assumption~2.7 supplies the required Cram\'er condition, and the moment conditions above validate the expansion. Consequently,
\[
    \PP\opt\left(n R_n(\mf h) \leq \hat z_{1-\varrho} \right) = \Psi(z_{1-\varrho}) + \widetilde{O}\left(\frac{1}{n}\right)
    = 1-\varrho + \widetilde{O}\left(\frac{1}{n}\right).
\]
Therefore, the proof is complete.
\end{proof}

\subsubsection{Proof of Theorem~\MainPowerExpansionNumber}\label{pf:thm4_5}\label{app:local_power_conditions}

The second-order power expansion requires the following strengthening of
Assumption~\ref{a:alter_edgeworth}.

\begin{assumption}[Additional Edgeworth regularity]\label{a:alter_edgeworthII}
Assume that Assumption~\ref{a:alter_edgeworth} holds and, in addition,
\begin{enumerate}[label=(\roman*)]
    \item for some $\delta'>0$,
    \begingroup\compactdisplaytrue
    \[
        \sup_{n\geq1}\EE_n\!\left[
        \left(\norm{\mf h(X)}_2+\norm{\D\mf h(X)}_2\right)^{8+\delta'}
        \right]<\infty,
        \qquad
        \sup_{n\geq1}\EE_n\!\left[
        \left(\norm{\D\mf h(X)}_2^2\kappa_1(X)\right)^{4+\delta'}
        \right]<\infty;
    \]
    \endgroup

    \item there is a limiting law $\PP\opt\in\mc P_0$ such that
    \begin{enumerate}[label=(\alph*)]
        \item $\PP\opt$ satisfies Assumptions~2.7--2.10;
        \item $\PP_n\opt$ converges to $\PP\opt$ in total variation.
    \end{enumerate}
\end{enumerate}
\end{assumption}

The first condition validates an order-$n^{-1}$ Edgeworth expansion for the
joint fluctuations of $\mc W_n$, $\mc V_n$, and $\mc K_n$; the formulation in
\cite[Assumption~(c)]{chandra1980valid} takes $\delta'=1$.  We do not require
their moments to possess expansions in powers of $n^{-1/2}$, so the correction
coefficients may depend explicitly on $n$.  The scalar location-shift case in
Proposition~\ref{prop:expan_rwpi} gives those coefficients explicitly.

For each $i \in [n]$, define $Y_i \Let f(X_i)$ by
\begin{align*}
Y_i=\big(& (h^\alpha)_{\alpha\in[d_h]},
          (h^\beta h^\gamma)_{1\leq\beta\leq\gamma\leq d_h},\\
        & (\D h^\beta\Sigma\D h^{\gamma\top})_{1\leq\beta\leq\gamma\leq d_h},
          (\D h^\beta\Sigma\D^2h^\gamma\Sigma\D h^{\omega\top})_{
          1\leq\beta\leq\omega\leq d_h,\ \gamma\in[d_h]}\big)(X_i).
\end{align*}
For each $n$, the vectors $(Y_i)_{i\in[n]}$ are i.i.d. under $\PP_n\opt$; we
write $Y$ for a generic copy.

\begin{proof}[Proof of Theorem~\MainPowerExpansionNumber]
    First, Assumption~\ref{a:alter_edgeworthII} and Lemma~\ref{lem:cramerII} validate the Edgeworth expansion of
    \[
        \pare{\sqrt{n}\pare{\mc W_n - W_{n}},
        \sqrt{n}\pare{\mc V_n - V_{n}}}.
    \]
    Proposition~\ref{prop:expanZn} therefore continues to hold under the local alternatives $\pare{\PP_{n}\opt}_{n\geq 1}$.

    Next, following the same steps as in the proof of Theorem~3.4, it remains to show a valid Edgeworth expansion of $\gamma^{\dagger}$ \eqref{eq:defgammada_n}. Under local alternatives $\pare{\PP_{n}\opt}_{n\geq 1}$, the definition of $\gamma^{\dagger}$ is different. We therefore introduce $\gamma^{\ddagger}$, defined by
    \begingroup\compactdisplaytrue
    \begin{align*}
        \langle \gamma^{\ddagger}, \eta \rangle =& \left\langle V_{n}^{\half}, \bar \xi_n \otimes \eta\right\rangle + \half \mc L_{W_{n},V_{n}}(\mc W_n - W_{n}, \mc V_n - V_{n}) \left\langle V_{n}^{\half}, \bar \xi_n \otimes \eta\right\rangle \notag\\ 
        & - \frac{1}{2} \left\langle \mc V_n - V_{n}, \bar \xi_n \otimes \left(V_{n}^{-\half} \eta\right) \right\rangle + \frac{1}{2\sqrt{n}} \left\langle \EE_{n}[\mc K_n], \bar \xi_n \otimes \bar \xi_n \otimes \left(V_{n}^{-\half} \eta\right) \right\rangle \quad \forall \eta \in \RR^{d_h},
    \end{align*}
    \endgroup
    where 
    \begin{align*}
        \bar \xi_n 
        &= \frac{1}{\sqrt{n}}\sum_{i=1}^n  V_{n}^{-1} \mf h(X_i)= \sqrt{n} V_{n}^{-1} \pare{\frac{1}{n}\sum_{i=1}^n  \mf h(X_i) - \EE_{n}\brac{\mf h(X)}} + \sqrt{n} V_{n}^{-1}\EE_{n}\brac{\mf h(X)}.
    \end{align*}

    Lemma~\ref{lem:cramerII} and
    \cite[Theorem 2(b)]{bhattacharya1978validity} imply that
    $\gamma^{\ddagger}$ admits an Edgeworth expansion up to order $n^{-1}$.

    Finally, the computation of the explicit Edgeworth-expansion formula is deferred to Lemma~\ref{lem:alter_cumulantII} in Supplementary Section~\ref{sec:expan_cum_compute}.
\end{proof}

\subsubsection{Proof of Proposition~\ref{prop:expan_rwpi}}\label{sec:pf_prop_5_7}
\paragraph{Regularity for location shifts}\label{app:location_shift_regularity}
With $\theta_n=\tau_0/\sqrt n$, assume that $\PP_0\opt$ has a density $f_0$
on $\RR^{d_X}$ satisfying $f_0(x-\theta_n)\to f_0(x)$ for almost every $x$.
The model in Assumption~\MainLocationShiftNumber{} is equivalently
\[
    \PP_n\opt(S)=\PP_0\opt(S-\theta_n),
    \qquad S\in\mc B(\RR^{d_X}),
\]
where $S-\theta_n=\{x:x+\theta_n\in S\}$.  Thus $\PP_n\opt$ has density
$f_n(x)=f_0(x-\theta_n)$, and Scheff\'{e}'s lemma gives
$\norm{f_n-f_0}_{L^1}\to0$.

\begin{proposition}[Power expansion of WP-based test under location shift]\label{prop:expan_rwpi}
    Suppose Assumption~\MainLocationShiftNumber{} and the conditions of
    Theorem~\MainPowerExpansionNumber{} hold with $d_h=1$.
    Then
    \begin{align}\label{eq:power_expan_one_dim}
        \PP_n\opt\left(nR_n(h)>\hat z_{1-\varrho}\right)
        ={}& \underbrace{\int_{|v+\tau|^2>\chi^2_{1;1-\varrho}}\phi(v)\diff v}_{\text{first-order}}
        +\underbrace{\frac{1}{\sqrt n}E_2(\varrho,\tau_0)}_{\text{second-order}}
        +\widetilde{O}\left(\frac1n\right),
    \end{align}
    where $\phi$ is the density of $\mc N(0,1)$,
    $\chi^2_{1;1-\varrho}$ is the $(1-\varrho)$-quantile of the chi-square
    distribution with one degree of freedom, and
    \begin{align*}
        \tau & = \frac{\EE_0[\D h(X)]\tau_0}{\sqrt{\alpha_2}},\\
        E_2(\varrho,\tau_0)
        & = \int_{|v+\tau|^2>\chi^2_{1;1-\varrho}}
        \left\{k_1v+\frac{k_2}{2}(v^2-1)+\frac{k_3}{6}(v^3-3v)\right\}
        \phi(v)\diff v.
    \end{align*}
    The coefficients are
    \begin{subequations}\label{eq:k11k21}
    \begin{align}
        k_1={}&\frac{\tau_0^\top\EE_0[\D^2h(X)]\tau_0}{2\sqrt{\alpha_2}}
        -\frac{\tau_0^\top\EE_0[\D h(X)^\top]\EE_0[h(X)\D h(X)]\tau_0}{\alpha_2^{3/2}}
        -\frac{\alpha_3}{2\alpha_2^{3/2}}
        +\frac{\tilde\alpha_3\sqrt{\alpha_2}}{2\tilde\alpha_2^2}\notag\\
        &\quad+\frac{\tilde\alpha_3}{2\sqrt{\alpha_2}\tilde\alpha_2^2}
        \tau_0^\top\EE_0[\D h(X)^\top]\EE_0[\D h(X)]\tau_0,\\
        k_2={}&\left(-\frac{\alpha_3}{\alpha_2^2}
        +\frac{2\tilde\alpha_3}{\tilde\alpha_2^2}\right)\EE_0[\D h(X)]\tau_0,\\
        k_3={}&-\frac{2\alpha_3}{\alpha_2^{3/2}}
        +\frac{3\tilde\alpha_3\sqrt{\alpha_2}}{\tilde\alpha_2^2}.
    \end{align}
    \end{subequations}
    Here
    \begin{align*}
        \alpha_j&=\EE_0[h(X)^j],\quad j=2,3,\\
        \tilde\alpha_2&=\EE_0[\D h(X)\Sigma\D h(X)^\top],\\
        \tilde\alpha_3&=\EE_0[\D h(X)\Sigma\D^2h(X)\Sigma\D h(X)^\top],
    \end{align*}
    with $\alpha_2,\tilde\alpha_2>0$.
\end{proposition}

\noindent\textbf{Notation.} In this section, we set
\begin{align*}
\alpha_{j,n}&\Let\EE_{\PP_n\opt}[h(X)^j],\\
\tilde\alpha_{2,n}&\Let\EE_n[\D h(X)\Sigma\D h(X)^\top],\\
\tilde\alpha_{3,n}&\Let\EE_n[\D h(X)\Sigma\D^2h(X)\Sigma\D h(X)^\top],\\[2pt]
A_1&\Let n^{-1}\sum_{i=1}^n h(X_i)-\alpha_{1,n},
&A_2&\Let n^{-1}\sum_{i=1}^n h(X_i)^2-\alpha_{2,n},\\[2pt]
\tilde A_2&\Let n^{-1}\sum_{i=1}^n\D h(X_i)\Sigma\D h(X_i)^\top-\tilde\alpha_{2,n},\\
\tilde A_3&\Let n^{-1}\sum_{i=1}^n\D h(X_i)\Sigma\D^2h(X_i)\Sigma\D h(X_i)^\top-\tilde\alpha_{3,n}.
\end{align*}
\begin{proof}[Proof of Proposition~\ref{prop:expan_rwpi}]
Using Theorem~3.3 and the identity, valid when $d_h=1$,
\[
\hat z_{1-\varrho} = \frac{\mc W_n}{\mc V_n} \chi^2_{1;1-\varrho}
\]
for $\mc W_n,\mc V_n\in\RR$, we obtain the asymptotic expansion
    \begin{align*}
            n R_n(h) - \hat z_{1-\varrho} = \frac{\mc W_n}{\mc V_n} \left(\bar \gamma^2 - \chi^2_{1;1-\varrho}\right) + \eps_n,
        \end{align*}
        where $|\eps_n| \leq \delta_n$ for a deterministic sequence $\delta_n = \widetilde{O}\pare{n^{-1}}$, $\bar \gamma = \sqrt{n}\pare{R_1 + R_2}$ for $R_j = O_p(n^{-j/2}),$ $j = 1,2$. Specifically,
        \begin{subequations}
        \begin{align*}
            R_1 &= \frac{1}{\sqrt{\alpha_{2,n}}} (A_1 + \alpha_{1,n}), \\
            R_2 &= - \frac{1}{2 \alpha_{2,n}^{\frac{3}{2}}} A_2 (A_1 + \alpha_{1,n}) + \frac{\tilde \alpha_{3,n}}{2 \sqrt{\alpha_{2,n}} \tilde \alpha_{2,n}^2}(A_1 + \alpha_{1,n})^2.
        \end{align*}
        \end{subequations}
        
        Then, using Theorem~\MainPowerExpansionNumber, $n R_n(h) - \hat z_{1-\varrho}$ has a valid Edgeworth expansion whose coefficients depend on the cumulant expansion of $\bar \gamma$:
        \begin{align*}
            \PP_{n}\opt\left(n R_n(h) > \hat z_{1-\varrho}\right) = \int_{|v + \tau|^2 > \chi^2_{1;1-\varrho}} \phi(v) \diff v + \frac{1}{\sqrt{n}}\int_{|v + \tau|^2 > \chi^2_{1;1-\varrho}} q(v) \phi(v) \diff v + \widetilde{O}\left(\frac{1}{n}\right),
        \end{align*}
        where $\phi(\cdot)$ is the density function of $\mc N(0,1)$, and
        \begin{align*}
            \tau = & \frac{\EE_0\left[\D h(X)\right] \tau_0 }{\sqrt{\alpha_{2}}} \in \RR,\\
            q(x) = & k_{1} x + \frac{k_{2}}{2}\pare{x^2 - 1} + \frac{k_{3}}{6}\pare{x^3 - 3x}.
        \end{align*}         
        The constants $k_1,k_2,k_3$ are the first-order coefficients in the
        expansions of the first three cumulants of $\bar\gamma$.  Their
        derivation is deferred to Supplementary Section~\ref{app:prop_5_7}.
\end{proof}

\subsubsection{Proof of Proposition~\MainPowerComparisonNumber}
The next lemma supplies the triangular-array EL expansion used in the
comparison.  It is stated separately because applying a fixed-law expansion
pointwise for each $n$ would not by itself control its remainder along the
sequence of transformed laws.

\begin{lemma}[Uniform scalar EL expansion under location shifts]
\label{lem:uniform_scalar_el}
Suppose Assumption~\MainLocationShiftNumber{} and the conditions of Theorem~\MainPowerExpansionNumber{} hold with $d_h=1$,
and let $Z_{n,i}=h(X_i)$ for $X_i\sim\PP_n\opt$.  Then
\begin{align}
    \PP_n\opt\left(nR_n^{EL}(h)>\chi^2_{1;1-\varrho}\right)
    ={}&\int_{|v+\tau|^2>\chi^2_{1;1-\varrho}}\phi(v)\diff v
      +\frac{1}{\sqrt n}F_2(\varrho,\tau_0)+o(n^{-1/2}),
    \label{eq:uniform_el_power}
\end{align}
uniformly along the location-shift sequence, where
\begin{align*}
    \tau&=\frac{\EE_0[\D h(X)]\tau_0}{\sqrt{\alpha_2}},\\
    F_2(\varrho,\tau_0)
    &=\int_{|v+\tau|^2>\chi^2_{1;1-\varrho}}
      \left\{\tilde k_1v+\frac{\tilde k_2}{2}(v^2-1)\right\}
      \phi(v)\diff v,
\end{align*}
and
\begin{subequations}\label{eq:tk11tk21}
\begin{align}
    \tilde k_1={}&
      \frac{\tau_0^\top\EE_0[\D^2h(X)]\tau_0}{2\sqrt{\alpha_2}}
      -\frac{\tau_0^\top\EE_0[\D h(X)^\top]
      \EE_0[h(X)\D h(X)]\tau_0}{\alpha_2^{3/2}}
      -\frac{\alpha_3}{6\alpha_2^{3/2}}\notag\\
      &\quad+\frac{\alpha_3}{3\alpha_2^{5/2}}
      \tau_0^\top\EE_0[\D h(X)^\top]\EE_0[\D h(X)]\tau_0,\\
    \tilde k_2={}&
      \left(-\frac{\alpha_3}{\alpha_2^2}
      +\frac{4\alpha_3}{3\alpha_2^2}\right)
      \EE_0[\D h(X)]\tau_0.
\end{align}
\end{subequations}
\end{lemma}

\begin{proof}
Write $\widehat m_{r,n}=n^{-1}\sum_{i=1}^nZ_{n,i}^r$ and
$\bar Z_n=\widehat m_{1,n}$.  Because $\EE_0[Z_0]=0$ and
$\alpha_2>0$, both $\PP_0\opt(Z_0<0)$ and $\PP_0\opt(Z_0>0)$ are
positive.  Total-variation convergence of $\PP_n\opt$ to $\PP_0\opt$
therefore gives a constant $p>0$ such that, for all sufficiently large $n$,
each of these two probabilities is at least $p$.  Consequently,
\[
    \PP_n\opt\left(0\notin
      (\min_{i\leq n}Z_{n,i},\max_{i\leq n}Z_{n,i})\right)
    \leq 2(1-p)^n=o(n^{-1/2}).
\]
Thus the ordinary scalar EL multiplier exists with probability
$1-o(n^{-1/2})$.

On this event, let $\lambda_n$ solve
\[
    \frac1n\sum_{i=1}^n\frac{Z_{n,i}}{1+\lambda_nZ_{n,i}}=0.
\]
The uniform moments and nondegeneracy in the conditions of Theorem~\MainPowerExpansionNumber{} give
\[
    \bar Z_n=O_p(n^{-1/2}),\qquad
    \widehat m_{2,n}^{-1}=O_p(1),\qquad
    \widehat m_{r,n}=O_p(1),\quad r=3,4,
\]
where all stochastic orders here and below are uniform along the array.  They
also imply $\max_{i\leq n}|Z_{n,i}|=o_p(\sqrt n)$.  Expanding the multiplier
equation uniformly on $|\lambda_n|\max_i|Z_{n,i}|=o_p(1)$ gives
\begin{align*}
    \lambda_n
    =\frac{\bar Z_n}{\widehat m_{2,n}}
     +\frac{\widehat m_{3,n}\bar Z_n^2}{\widehat m_{2,n}^3}
     +O_p(n^{-3/2}).
\end{align*}
Substitution into the scalar EL dual yields
\begin{align}
    nR_n^{EL}(h)
    &=2\sum_{i=1}^n\log(1+\lambda_nZ_{n,i})\notag\\
    &=\frac{n\bar Z_n^2}{\widehat m_{2,n}}
      +\frac{2n\widehat m_{3,n}\bar Z_n^3}
      {3\widehat m_{2,n}^3}+O_p(n^{-1}).
    \label{eq:el_dual_uniform}
\end{align}
Since $\widehat m_{3,n}-\alpha_{3,n}=O_p(n^{-1/2})$, this can be written as
\begin{align}
    nR_n^{EL}(h)=\Gamma_{n,EL}^2+r_n,
    \qquad
    \Gamma_{n,EL}
    =\frac{\sqrt n\bar Z_n}{\sqrt{\widehat m_{2,n}}}
     +\frac{\sqrt n\alpha_{3,n}\bar Z_n^2}
      {3\widehat m_{2,n}^{5/2}}.
    \label{eq:el_signed_root_uniform}
\end{align}
The truncation and empirical-moment bounds used in the proof of Theorem~\MainPowerExpansionNumber
strengthen the displayed stochastic remainders as follows: for a deterministic
$\delta_n=\widetilde{O}(n^{-1})$,
\begin{equation}
    \PP_n\opt(|r_n|\leq\delta_n)=1-o(n^{-1/2}).
    \label{eq:el_remainder_uniform}
\end{equation}
Indeed, on the same truncation event the Taylor remainder in the multiplier
equation is bounded by a polynomial in the uniformly bounded empirical
moments times $n^{-3/2}$, and the omitted fourth-order term in the EL dual is
$\widetilde{O}(n^{-1})$.  The complement has probability $o(n^{-1/2})$ by
Assumption~\ref{a:alter_edgeworthII}.

It remains to justify the distributional expansion uniformly in $n$.  The
vector $(Z_{n,i},Z_{n,i}^2)$ is a subvector of the empirical-moment vector
used in the proof of Theorem~\MainPowerExpansionNumber. Hence the uniform moment and Cram\'{e}r
conditions in Assumption~\ref{a:alter_edgeworthII} and
Lemma~\ref{lem:cramerII} apply; if this vector has a redundant coordinate, we
use the maximal nondegenerate subvector reduction described in
Section~\ref{sec:valid_Edgeworth}.  The triangular-array Edgeworth theorem
used in the proof of Theorem~3.3, followed by the smooth-transformation
argument of \cite[Theorem~2(b)]{bhattacharya1978validity}, therefore applies
to \eqref{eq:el_signed_root_uniform}.  The transformation is uniformly smooth
on a neighborhood of the population moments because
$\inf_n\EE_n[Z_n^2]>0$.

To make the coefficient calculation explicit, set
$U_n=\sqrt n(\bar Z_n-\alpha_{1,n})$,
$T_n=\sqrt n(\widehat m_{2,n}-\alpha_{2,n})$, and
$a_n=\sqrt n\alpha_{1,n}$.  A second-order Taylor expansion of
\eqref{eq:el_signed_root_uniform} gives
\begin{align*}
    \Gamma_{n,EL}
    ={}&\frac{a_n+U_n}{\sqrt{\alpha_{2,n}}}
    +\frac1{\sqrt n}\left\{
      -\frac{(a_n+U_n)T_n}{2\alpha_{2,n}^{3/2}}
      +\frac{\alpha_{3,n}(a_n+U_n)^2}
      {3\alpha_{2,n}^{5/2}}\right\}
    +O_p(n^{-1}).
\end{align*}
Here
$\operatorname{Cov}_n(U_n,T_n)=\alpha_{3,n}-\alpha_{1,n}\alpha_{2,n}$
and
$\kappa_{3,n}(U_n)=n^{-1/2}
(\alpha_{3,n}-3\alpha_{1,n}\alpha_{2,n}+2\alpha_{1,n}^3)$.
Combining these identities with the translated-moment expansions in
\eqref{eq:expan_alpha} gives the first three cumulants of
$\Gamma_{n,EL}-\tau$ as
\begin{align*}
    \kappa_1(\Gamma_{n,EL}-\tau)
      &=\frac{\tilde k_1}{\sqrt n}+O(n^{-1}),\\
    \kappa_2(\Gamma_{n,EL}-\tau)
      &=1+\frac{\tilde k_2}{\sqrt n}+O(n^{-1}),\\
    \kappa_3(\Gamma_{n,EL}-\tau)&=O(n^{-1}),
\end{align*}
with $\tilde k_1,\tilde k_2$ in \eqref{eq:tk11tk21}.  In particular, the
order-$n^{-1/2}$ skewness term cancels in the EL signed root.  In CDF form,
the resulting uniform Edgeworth expansion is
\begin{align*}
    \sup_{x\in\RR}\left|\PP_n\opt(\Gamma_{n,EL}-\tau\leq x)
    -\int_{-\infty}^x\left[1+\frac1{\sqrt n}
      \left\{\tilde k_1v+\frac{\tilde k_2}{2}(v^2-1)\right\}\right]
      \phi(v)\diff v\right|=o(n^{-1/2}).
\end{align*}
Finally, \eqref{eq:el_remainder_uniform} and the preceding expansion imply
that replacing $nR_n^{EL}(h)$ by $\Gamma_{n,EL}^2$ changes the rejection
probability by at most $o(n^{-1/2})+\widetilde{O}(n^{-1})$: off the exceptional
event, a disagreement can occur only within $\delta_n$ of either fixed
boundary.  Integrating over the two rejection tails therefore proves
\eqref{eq:uniform_el_power}.  The coefficients coincide with the scalar
mean formula in \cite[Theorem~3.1]{chen1994comparing}, which provides an
independent fixed-law check of the calculation.
\end{proof}

\begin{proof}[Proof of Proposition~\MainPowerComparisonNumber]
Since the power expansion of Hotelling's $T^2$ test is a special instance of
Proposition~\ref{prop:expan_rwpi} obtained by setting $\tilde\alpha_3=0$, it
remains to derive the EL comparison.  Set $Z_{n,i}=h(X_i)$.  The EL constraint
depends on $X_i$ only through $Z_{n,i}$, and hence
\[
    R_n^{EL}(h;X_1,\ldots,X_n)
    =R_n^{EL}(\operatorname{id};Z_{n,1},\ldots,Z_{n,n}).
\]
Lemma~\ref{lem:uniform_scalar_el} gives its uniform triangular-array power
expansion and the coefficients in \eqref{eq:tk11tk21}.

Comparing the power expansion of the WP-based test in Proposition~\ref{prop:expan_rwpi} with that of the EL test, we examine the difference in the $n^{-1/2}$ term. Specifically, we have
\begin{align*}
    E_2(\varrho, \tau_0) - F_2(\varrho, \tau_0)
    &=\sqrt{\alpha_{2}} \left(\frac{\tilde \alpha_{3} }{\tilde \alpha_{2}^2}  - \frac{2 \alpha_{3} }{3\alpha_{2}^2}\right) I_2(\varrho, \tau),
\end{align*}
where $\tau=\EE_0[\D h(X)]\tau_0/\sqrt{\alpha_2}$ is the scalar standardized
drift from Proposition~\ref{prop:expan_rwpi}.
For $w_1 \Let \sqrt{\chi^2_{1;1-\varrho}}-\tau$ and
$w_2 \Let \sqrt{\chi^2_{1;1-\varrho}}+\tau$, define
$I_2:(0,1)\times\RR\to\RR$ by
\begin{equation}
\begin{aligned}
    I_2(\varrho, \tau) &\Let \int_{|v + \tau|^2 > \chi^2_{1;1-\varrho}} \left(\frac{1+\tau^2}{2} v + \tau (v^2 - 1) + \frac{1}{2}(v^3 - 3v) \right) \phi(v) \diff v \\
    &= \frac{\chi^2_{1;1-\varrho}}{2}(\phi(w_1) - \phi(w_2)).
\end{aligned}
\end{equation} 

The power gap between the WP test and Hotelling's $T^2$
test has the same formula as above, with $0$ in place of
$2\alpha_3/(3\alpha_2^2)$.
\end{proof}

\subsubsection{Proof sketch of Theorem~\MainHigherExpansionNumber}\label{sec:pf_thm18}\label{app:higher_order_details}
The fourth-order tensor in Theorem~\MainHigherExpansionNumber{} is
\begin{equation*}
\resizebox{0.98\linewidth}{!}{$\displaystyle
\mc L_n \Let \left(
-\frac1n\sum_{i=1}^n
 h^{\omega^1}_{\alpha\alpha'}h^{\omega^2}_{\beta'\gamma'}
 h^{\omega^3}_{\beta}h^{\omega^4}_{\omega'}
 \Sigma_{\omega'\gamma'}\Sigma_{\beta\alpha}\Sigma_{\beta'\alpha'}
-\frac1{3n}\sum_{i=1}^n
 h^{\omega^1}_{\alpha\alpha'\alpha''}h^{\omega^2}_{\beta}
 h^{\omega^3}_{\beta'}h^{\omega^4}_{\beta''}
 \Sigma_{\beta\alpha}\Sigma_{\beta'\alpha'}\Sigma_{\beta''\alpha''}
+S_{\gamma\omega^1\omega^2}(\mc V_n^{-1})_{\gamma\theta}
 S_{\theta\omega^3\omega^4}
\right)_{\omega^1,\omega^2,\omega^3,\omega^4\in[d_h]}$}.
\end{equation*}
Here all summands are evaluated at $X_i$, and
\begin{align*}
A_{\gamma\omega^1\omega^2}
&\Let n^{-1}\sum_{i=1}^n
 h^\gamma_{\alpha\alpha'}h^{\omega^1}_\beta h^{\omega^2}_{\beta'}
 \Sigma_{\beta\alpha}\Sigma_{\beta'\alpha'},\\
S_{\gamma\omega^1\omega^2}
&\Let \frac12 A_{\gamma\omega^1\omega^2}
       +A_{\omega^1\gamma\omega^2}.
\end{align*}
The tensor $A$ is symmetric in its last two indices.  Thus $S$ is precisely the
coefficient obtained by differentiating the cubic form: if
$P_{3,n}(\zeta)=\frac18A_{\alpha\beta\gamma}
\zeta^{(\alpha)}\zeta^{(\beta)}\zeta^{(\gamma)}$, then
\begin{equation*}
    \left\{\D P_{3,n}(-2\xi)\right\}_\gamma
    =S_{\gamma\omega^1\omega^2}
      \xi^{(\omega^1)}\xi^{(\omega^2)}.
\end{equation*}
Repeated indices are summed using the Einstein convention.  Only the full
symmetrization of $\mc L_n$ contributes to
$\brak{\mc L_n,\xi_n^{\otimes4}}$.

\begin{example}[Quadratic moment function II]
    Suppose $\mf h(x)=\norm{x}_2^2-1$, $\Sigma=I_{d_X}$,
    $\EE[\norm{X}_2^2]=1$, and the support of $\PP\opt$ is compact.  The
    projection problem~(4) has the exact expansion
    \begin{equation}
        nR_n(\mf h)=n\left(\sqrt{1+\frac{\Delta_n}{\sqrt n}}-1\right)^2
        =n\left(\frac{\Delta_n}{2\sqrt n}-\frac{\Delta_n^2}{8n}
        +\frac{\Delta_n^3}{16n^{3/2}}+O_p(n^{-2})\right)^2,\label{eq:19aII}
    \end{equation}
    where $\Delta_n=n^{-1/2}\sum_{i=1}^n(\norm{X_i}_2^2-1)$.  Theorem~\MainHigherExpansionNumber
    instead gives
    \begin{equation}\label{eq:approx_rwpi_quadII}
        \text{RHS of~(12)}
        =\frac{\Delta_n^2}{4}-\frac{\Delta_n^3}{8\sqrt n}
        +\frac{5\Delta_n^4}{64n}+O_p(n^{-3/2}),
    \end{equation}
    and the two expressions agree through order $n^{-1}$.
\end{example}

\begin{remark}[Proof adjustments for Theorem~\MainHigherExpansionNumber]\label{rmk:diff_oper_pf}
This proof follows the steps used for Theorem~3.1 and described in Section~\ref{sec:proofroadmap}. Specifically, Steps~2, 4, and~5 require the following adjustments:
\begin{enumerate}[label=(\arabic*)]
    \item For \textit{Step 2}, Assumption~\MainSmoothnessNumber{} includes the global
    Jacobian-Lipschitz condition of Assumption~2.5 and gives sub-exponential
    tails for its envelope $\kappa_1(X)$, for $\mf h(X)$, and for the local
    derivatives up to order four. By choosing the constants in the logarithmic
    envelopes sufficiently large, a union bound for the sub-exponential tails
    gives an event of probability at least $1-O\pare{n^{-3/2}}$ on which these
    quantities are bounded by deterministic logarithmic envelopes. On this event, the tail
    bounds (A3)--(A6) in
    \ref{cond:condA} of Proposition~\ref{prop:tailbnd} follow with corresponding
    logarithmic bounds, which are absorbed by the $\widetilde{O}(\cdot)$ notation.
    The checks of (A1) and (A2) follow from the same concentration arguments,
    using Lemma~\ref{lem:tail_avg2} and Lemma~\ref{lem:matberstein}.
    
    \item Additionally, on the same high-probability event, the third-order
    derivatives of $\mf h$ and the envelopes $\kappa_1(\cdot),\kappa_2(\cdot)$
    are bounded by deterministic polylogarithmic quantities; here one may take
    $\kappa_2(x)=\sup_{\norm{u-x}_2\leq1}\norm{\D^3\mf h(u)}_{\mathrm{op}}$.
    These factors are again absorbed into the $\widetilde{O}(\cdot)$ remainders.
    
    \item As for \textit{Step 4}, we will instead show (in the sequel) that 
    \[
    M_n(\zeta) = F^{\dagger}_n(\zeta) + \eps^{M\dagger}_n(\zeta) \quad \forall \zeta \in \mc Z_n
    \]
    for a quartic form $F^{\dagger}_n(\zeta)$. Above, $\mc Z_n$ is defined as in~\eqref{eq:Z}.

    \item As for \textit{Step 5}, we will instead provide the expansion: with high probability, 
    \[
    \sup_{\zeta \in \mc Z_n}\left\{-\zeta^{\top} H_n - F^{\dagger}_n(\zeta)\right\} = \brak{\mc V_n, \xi_n^{\otimes 2}} + \frac{1}{\sqrt{n}} \brak{\mc K_n, \xi_n^{\otimes 3}} + \frac{1}{n}\brak{\mc L_n, \xi_n^{\otimes 4}} + \widetilde{O}\left(n^{-\frac{3}{2}}\right).
    \]
\end{enumerate}

With the above adjustments, the final step will be the same as stated in Section~\ref{sec:comb_step} with $F_n^{\dagger}$ (see Proposition~\ref{prop:expanMn2}) in place of $F_n$.
\end{remark}

We now prove items (3) and (4) in the following two sections.

\subsubsection{Proof of Theorem~\MainHigherExpansionNumber{} -- part I}
This section proves item (3) listed in Section~\ref{sec:pf_thm18}. We fix a
vector $\zeta \in \mc Z_n$ and introduce the following notation.

\textbf{Notation.} In this section only, let
\begin{itemize}
    \item $\Delta_i = \pare{\Delta^{\alpha}_i, \alpha \in [d_X]} \Let \Delta_{n,i}(\zeta)$ for every $i \in [n]$.
    \item $\bar h \Let \zeta^\top \mf h,~~ \pare{\bar h_{\beta}, \beta \in [d_X]} \Let \D (\zeta^\top \mf h),~~ \pare{\bar h_{\beta\gamma}, \beta,\gamma \in [d_X]} \Let \D^2 (\zeta^\top \mf h)$.
\end{itemize}

We formalize the result in the following proposition.
\begin{proposition}[Expanding $M_n(\cdot)$, cf. Proposition~\ref{prop:expanMn}]\label{prop:expanMn2}
    Under Assumption~\MainSmoothnessNumber{}, there exists a deterministic
    integer $N$ such that, on the logarithmic-envelope event of
    Remark~\ref{rmk:diff_oper_pf} and whenever \ref{cond:condA} holds, for every
    $n \geq N$,
    \begin{subequations}
    \begin{align*}
        M_n(\zeta) = F^{\dagger}_n(\zeta) + \eps^{M^\dagger}_n(\zeta),
    \end{align*}
    where $\sup_{\zeta \in \mc Z_n} |\eps^{M^\dagger}_n(\zeta)| \leq \delta_n^{\dagger}$ for a deterministic sequence $\delta_n^{\dagger} = \widetilde{O}\pare{n^{-\frac{3}{2}}}$, and
    \begingroup\compactdisplaytrue
    \begin{align*}
    F^{\dagger}_n(\zeta) =&  \frac{1}{n}\sum_{i=1}^n \left(\underbrace{\frac{1}{4}  \bar h_{\alpha} \bar h_{\beta} \Sigma_{\beta\alpha}}_{\text{quadratic term}} + \underbrace{\frac{1}{\sqrt{n}} \frac{\bar h_{\alpha \alpha'} \bar h_{\beta}\bar h_{\beta'} \Sigma_{\beta \alpha}\Sigma_{\beta'\alpha'}}{8}}_{\text{cubic term}} + \underbrace{\frac{1}{n} \frac{\bar h_{\alpha \alpha'} \bar h_{\beta'\gamma'} \bar h_{\beta} \bar h_{\omega'} \Sigma_{\omega'\gamma'}\Sigma_{\beta \alpha}\Sigma_{\beta'\alpha'}}{16}}_{\text{quartic term}}\right.\\
    & \left. + \underbrace{\frac{1}{n}\frac{\bar h_{\alpha \alpha' \alpha''}\bar h_{\beta}\bar h_{\beta'}\bar h_{\beta''}\Sigma_{\beta \alpha}\Sigma_{\beta'\alpha'}\Sigma_{\beta'' \alpha''}}{48}}_{\text{quartic term}}\right)(X_i).
    \end{align*} 
    \endgroup
    \end{subequations}
\end{proposition}

To prove this proposition, we need the following expansion of $\Delta_i$ for $i \in [n]$.
\begin{lemma}[Expanding $\Delta_{i}$, cf.~Lemma~\ref{lem:DelinMn}]\label{lem:DelinMn3}
    Under Assumption~\MainSmoothnessNumber{}, there exists a deterministic
    integer $N$ such that, on the logarithmic-envelope event of
    Remark~\ref{rmk:diff_oper_pf} and whenever \ref{cond:condA} holds, the
    following expansion holds uniformly for every $n\geq N$, $i\in[n]$,
    $\zeta\in\mc Z_n$, and $\alpha\in[d_X]$:
    \begingroup\compactdisplaytrue
    \begin{align*}
        \Delta^{\alpha}_i =& \Sigma_{\beta\alpha}\left(\frac{1}{2}\bar h_{\beta} + \frac{1}{\sqrt{n}}\frac{\bar h_{\beta \gamma} \bar h_{\omega} \Sigma_{\omega \gamma}}{4} + \frac{1}{n}\frac{2 \bar h_{\beta \gamma} \bar h_{\omega \omega'} \bar h_{\gamma'}\Sigma_{\omega \gamma} \Sigma_{\gamma'\omega'} + \bar h_{\beta\gamma \omega} \bar h_{\gamma'} \bar h_{\omega'}\Sigma_{\gamma' \gamma}\Sigma_{\omega' \omega}}{16}\right)(X_i) \\
        &+ \widetilde{O}\left(n^{-\frac{3}{2}}\right).
    \end{align*}
    \endgroup
\end{lemma}
\begin{proof}[Proof of Lemma~\ref{lem:DelinMn3}]
    Using Lemma~\ref{lem:DelinMn}, equation~\eqref{eq:opt_Delta2} reads
    \begin{align}\label{eq:tmp522}
        \Delta_i^{\alpha} = \frac{1}{2} \Sigma_{\beta \alpha}(\zeta^\top \D \mf h)^{\beta}\pare{X_i + n^{-\half}\Delta_i} .
    \end{align} 
	On the logarithmic-envelope event in Remark~\ref{rmk:diff_oper_pf}, we have
    \[
        \max_{i\leq n}\norm{\D \mf h(X_i)}_2 \leq C\log(n)
    \]
    for a constant $C$. Then by equation~\eqref{eq:opt_Delta3} and $\zeta \in \mc Z_n$, we get
    \begin{align*}
        \norm{\Delta_{i}}_{2} \leq \norm{\D \mf h(X_i)}_2 \norm{\Sigma}_2\norm{\zeta}_2 \leq C \norm{\Sigma}_2 \log(n)^2.
    \end{align*}
    As a result, we apply the Taylor expansion of $\D \mf h(\cdot)$ around $X_i$ in \eqref{eq:tmp522} and get 
    \begingroup\compactdisplaytrue
    \begin{align*}
        &\Delta_i^{\alpha} \\
        =& \frac{1}{2} \Sigma_{\beta \alpha}\left(\bar h_{\beta} + \frac{1}{\sqrt{n}}\bar h_{\beta \gamma} \Delta_i^{\gamma} + \frac{1}{2n} \bar h_{\beta \gamma \omega} \Delta_i^{\gamma}\Delta_i^{\omega}\right)(X_i)  + \widetilde{O}\left(n^{-\frac{3}{2}}\right)\\ =&\Sigma_{\beta\alpha}\left(\frac{1}{2}\bar h_{\beta} + \frac{1}{\sqrt{n}}\frac{\bar h_{\beta \gamma} \bar h_{\omega} \Sigma_{\omega \gamma}}{4} + \frac{1}{n}\frac{2 \bar h_{\beta \gamma} \bar h_{\omega \omega'} \bar h_{\gamma'}\Sigma_{\omega \gamma} \Sigma_{\gamma'\omega'} + \bar h_{\beta\gamma \omega} \bar h_{\gamma'} \bar h_{\omega'}\Sigma_{\gamma' \gamma}\Sigma_{\omega' \omega}}{16}\right)(X_i) + \widetilde{O}\left(n^{-\frac{3}{2}}\right).
    \end{align*}
    \endgroup
    This completes the proof.
\end{proof}

\begin{proof}[Proof of Proposition~\ref{prop:expanMn2}]
    By Lemma~\ref{lem:DelinMn} and \ref{lem:DelinMn3}, to compute the asymptotic expansion of $M_n(\zeta)$ for $\zeta \in \mc Z_n$, we only need to plug in the expansion of $\Delta_i$ as shown in Lemma~\ref{lem:DelinMn3} into 
    \begin{align*}
        \sqrt{n}\left(\bar h\left(X_i + n^{-\half} \Delta_{i}\right) - \bar h\left(X_i\right)\right) - \Delta_{i}^\top \Sigma^{-1} \Delta_{i}
    \end{align*}
    for each $i \in [n]$. The concrete algebra is deferred to Supplementary Section~\ref{app:algebra_5.21}.
\end{proof}

\subsubsection{Proof of Theorem~\MainHigherExpansionNumber{} -- part II}
We formalize item (4) in the following proposition.
\begin{proposition}[Expansion of $F_n^{\dagger *}$, cf. Proposition~\ref{prop:expanFn*}]\label{prop:exp_Fndagger}
    Under Assumption~\MainSmoothnessNumber{}, there exists a deterministic
    integer $N$ such that, on the logarithmic-envelope event of
    Remark~\ref{rmk:diff_oper_pf} and whenever \ref{cond:condA} holds, for every
    $n \geq N$,
    \[
    \sup_{\zeta \in \mc Z_n}\left\{-\zeta^{\top} H_n - F^{\dagger}_n(\zeta)\right\} = \brak{\mc V_n, \xi_n^{\otimes 2}} + \frac{1}{\sqrt{n}} \brak{\mc K_n, \xi_n^{\otimes 3}} + \frac{1}{n}\brak{\mc L_n, \xi_n^{\otimes 4}} + \eps_n^{F^\dagger},
    \]
    where $\abs{\eps_n^{F^\dagger}} \leq \delta_n^{\ddagger}$ for a deterministic sequence $\delta_n^{\ddagger} = \widetilde{O}\pare{n^{-\frac{3}{2}}}$.
\end{proposition}

\begin{lemma}[Tail bound of the solution in $F_n^{\dagger *}(-H_n)$, cf. Lemma~\ref{prop:lowerbnd_Fn}]\label{prop:lowerbnd_FnII}
Under Assumption~\MainSmoothnessNumber{}, define $\mc Z_n$ as in~\eqref{eq:Z}
and $\zeta_n^{\ddagger}$ as any optimizer of
\begin{align*}
    \max_{\zeta \in \mc Z_n } \left\{- \zeta^\top H_n - F^{\dagger}_n(\zeta)\right\}.
\end{align*}
There exists a deterministic integer $N$ such that, on the
logarithmic-envelope event of Remark~\ref{rmk:diff_oper_pf} and whenever
\ref{cond:condA} holds, $\zeta_n^{\ddagger}\in\inte(\mc Z_n)$ for
every $n\geq N$.
\end{lemma}
\begin{proof}[Proof of Lemma~\ref{prop:lowerbnd_FnII}]
    Recall the definition of $F_n$ in \eqref{eq:def_Fn}. On the
    logarithmic-envelope event of Remark~\ref{rmk:diff_oper_pf}, the first three
    derivatives of $\mf h$ are bounded by deterministic logarithmic envelopes.
    Hence, for a deterministic polylogarithmic sequence $b_n$,
    \[
        \abs{F_n(\zeta) - F^{\dagger}_n(\zeta)}
        \leq \frac{Cb_n}{n}\norm{\zeta}_2^4
        \qquad \forall\zeta\in\mc Z_n.
    \]

    On the boundary $\norm{\zeta}_2=r_n$, the proof of
    Lemma~\ref{prop:lowerbnd_Fn} gives, for all sufficiently large $n$,
    \[
        -\zeta^\top H_n-F_n(\zeta)
        \leq-\frac{\sigma_{\min}}{32}r_n^2.
    \]
    Moreover, $r_n^2=O\{\log(n)\}$, so
    $Cb_nr_n^4/n=o(r_n^2)$. Therefore
    $-\zeta^\top H_n-F_n^\dagger(\zeta)<0$ throughout the boundary for
    all sufficiently large $n$. Since the objective equals zero at
    $\zeta=\mathbf0$, every optimizer $\zeta_n^{\ddagger}$ is in
    $\inte(\mc Z_n)$.
\end{proof}
\begin{proof}[Proof of Proposition~\ref{prop:exp_Fndagger}]
    Lemma~\ref{prop:lowerbnd_FnII} shows that the optimal solution $\zeta_n^{\ddagger}$ of
    \[
    \sup_{\zeta \in \mc Z_n}\left\{-\zeta^{\top} H_n - F^{\dagger}_n(\zeta)\right\}
    \]
    is achieved in $\inte(\mc Z_n)$. Thus, $\zeta_n^{\ddagger}$ satisfies $-H_n - \D F_n^\dagger(\zeta_n^{\ddagger}) = \mf 0$. Using this relation, we obtain
    \[
    -\zeta_n^{\ddagger\top} H_n - F^{\dagger}_n(\zeta_n^{\ddagger}) = \brak{\mc V_n, \xi_n^{\otimes 2}} + \frac{1}{\sqrt{n}} \brak{\mc K_n, \xi_n^{\otimes 3}} + \frac{1}{n}\brak{\mc L_n, \xi_n^{\otimes 4}} + \eps_n^{F^\dagger},
    \]
    where $\abs{\eps_n^{F^\dagger}} \leq \delta_n^{\ddagger}$ for a deterministic sequence $\delta_n^{\ddagger} = \widetilde{O}\pare{n^{-\frac{3}{2}}}$. The concrete algebra is deferred to Supplementary Section~\ref{sec:pf_prop725}.
\end{proof}

\subsubsection{Proof of Theorem~\MainHigherExpansionNumber{} -- combining steps}\label{sec:comb_stepII}
\begin{proof}[Proof of Theorem~\MainHigherExpansionNumber]
Applying Proposition~\ref{prop:rescale}, under Assumptions~2.4 and~2.5, we have
\begin{align}\label{eq:step5supMn_II}
    n R_n(\mf h) = \sup_{\zeta \in \RR^{d_h}}\{-\zeta^{\top} H_n - M_n(\zeta)\}.
\end{align}
Applying Proposition~\ref{prop:bound_subgradient}, under Assumptions~2.4
and~2.5 (the latter is included in Assumption~\MainSmoothnessNumber), there is a deterministic
    integer $N_1$ such that, when $n \geq N_1$ and \ref{cond:condA}
holds,
\begin{align}\label{eq:step5supMn2II}
    \eqref{eq:step5supMn_II} = \sup_{\zeta \in \mc Z_n}\{-\zeta^{\top} H_n - M_n(\zeta)\}.
\end{align}
On the logarithmic-envelope event of Remark~\ref{rmk:diff_oper_pf}, applying
Proposition~\ref{prop:expanMn2} under Assumption~\MainSmoothnessNumber{} and
\ref{cond:condA} gives a deterministic integer $N_2$ such that, when
$n \geq N_2$,
\begin{align}\label{eq:step5supMn3II}
    \eqref{eq:step5supMn2II} &= \sup_{\zeta \in \mc Z_n} \left\{-\zeta^{\top} H_n -  F^{\dagger}_n(\zeta) - \eps^{M^\dagger}_n(\zeta) \right\} = \sup_{\zeta \in \mc Z_n} \left\{-\zeta^{\top} H_n - F^{\dagger}_n(\zeta)\right\} + \eps^{F^\ddagger}_n,
\end{align}
where $\abs{\eps^{F^\ddagger}_n} \leq \delta_n$ for a deterministic sequence $\delta_n = \widetilde{O}\pare{n^{-\frac{3}{2}}}$.

On the same logarithmic-envelope event, Proposition~\ref{prop:exp_Fndagger},
under \ref{cond:condA}, gives a deterministic
integer $N_3$ such that, when $n \geq N_3$,
\begin{align}
    \eqref{eq:step5supMn3II} = \brak{\mc V_n, \xi_n^{\otimes 2}} + \frac{1}{\sqrt{n}} \brak{\mc K_n, \xi_n^{\otimes 3}} + \frac{1}{n}\brak{\mc L_n, \xi_n^{\otimes 4}} + \eps_n^{F^\dagger}   + \eps^{F^\ddagger}_n,
\end{align}
where $\abs{\eps^{F^\dagger}_n} \leq \delta'_n$ for a deterministic sequence $\delta'_n = \widetilde{O}\pare{n^{-\frac{3}{2}}}$. Therefore, we obtain~(12).

Finally, by the reasoning in item (1) of Section~\ref{sec:pf_thm18}, there is
a deterministic integer $N_4$ such that, when $n\geq N_4$, the intersection
of the event in \ref{cond:condA} and
the logarithmic-envelope event has probability at least $1-O(n^{-3/2})$,
which concludes the proof.
\end{proof}

\subsubsection{Proof of Theorem~\MainBartlettCoverageNumber}\label{sec:pf_thm6_5}\label{app:bartlett_coefficients}
\noindent\textbf{Bartlett coefficient definitions.}
For $k=1,2,3$, define
\begin{align*}
    C_k=-2u_k^{-1}\sum_{r=k}^3B_r,
    \qquad
    u_k=2^k\Gamma(k+\tfrac12)/\Gamma(\tfrac12),
\end{align*}
where
\begin{align*}
    B_0={}&-\frac{k_{11}^2+k_{22}}2+\frac{4k_{11}k_{31}+k_{42}}8
       -\frac{5k_{31}^2}{24},
    &B_1={}&\frac{k_{11}^2+k_{22}}2-\frac{4k_{11}k_{31}+k_{42}}4
       +\frac{5k_{31}^2}{8},\\
    B_2={}&\frac{4k_{11}k_{31}+k_{42}}8-\frac{5k_{31}^2}{8},
    &B_3={}&\frac{5k_{31}^2}{24}.
\end{align*}
The cumulant coefficients are
\begingroup\compactdisplaytrue
\begin{align*}
    k_{11}={}&-\frac{\alpha_3}{2\alpha_2^{3/2}}
       +\frac{\tilde\alpha_3\sqrt{\alpha_2}}{2\tilde\alpha_2^2},
    &k_{31}={}&-\frac{2\alpha_3}{\alpha_2^{3/2}}
       +\frac{3\tilde\alpha_3\sqrt{\alpha_2}}{\tilde\alpha_2^2},\\
    k_{22}={}&-\frac{3\alpha_3\tilde\alpha_3}{2\alpha_2\tilde\alpha_2^2}
       +\frac{7\alpha_3^2}{4\alpha_2^3}
       +\frac{-6\tilde\alpha_3\EE[(h\D h\Sigma\D h^\top)(X)]
       +3\tilde\alpha_4\alpha_2}{\tilde\alpha_2^3}\notag\\
       &+\frac{3\EE[(h\D h\Sigma\D^2h\Sigma\D h^\top)(X)]}{\tilde\alpha_2^2}
       -\frac{\alpha_2\tilde\alpha_3^2}{\tilde\alpha_2^4},\\
    k_{42}={}&-\frac{2\alpha_4}{\alpha_2^2}+\frac{12\alpha_3^2}{\alpha_2^3}
       -\frac{18\alpha_3\tilde\alpha_3}{\alpha_2\tilde\alpha_2^2}
       +\frac{12\EE[(h\D h\Sigma\D^2h\Sigma\D h^\top)(X)]}{\tilde\alpha_2^2}\notag\\
       &+\frac{12\alpha_2\tilde\alpha_4
       -24\tilde\alpha_3\EE[(h\D h\Sigma\D h^\top)(X)]}{\tilde\alpha_2^3}
       +\frac{9\alpha_2\tilde\alpha_3^2}{\tilde\alpha_2^4}.
\end{align*}
\endgroup
Here
\begin{align*}
    \alpha_j={}&\EE[h(X)^j],\quad j=1,2,3,4,
    &\tilde\alpha_2={}&\EE[(h_\alpha h_{\alpha'}\Sigma_{\alpha\alpha'})(X)],\\
    \tilde\alpha_3={}&\EE[(h_{\alpha'\beta'}h_\alpha h_\beta
       \Sigma_{\alpha\alpha'}\Sigma_{\beta'\beta})(X)],
\end{align*}
and
\begingroup\compactdisplaytrue
\begin{align*}
    \tilde\alpha_4={}&-\EE[(h_{\alpha\alpha'}h_{\beta'\gamma'}h_\beta h_{\omega'}
       \Sigma_{\omega'\gamma'}\Sigma_{\beta\alpha}\Sigma_{\beta'\alpha'})(X)]
       -\frac13\EE[(h_{\alpha\alpha'\alpha''}h_\beta h_{\beta'}h_{\beta''}
       \Sigma_{\beta\alpha}\Sigma_{\beta'\alpha'}\Sigma_{\beta''\alpha''})(X)]\notag\\
       &+\frac{9}{4\tilde\alpha_2}
       \left\{\EE[(h_{\alpha\alpha'}h_\beta h_{\beta'}
       \Sigma_{\beta\alpha}\Sigma_{\beta'\alpha'})(X)]\right\}^2.
\end{align*}
\endgroup
Here $k_{\beta\gamma}$ denotes the $\gamma$th-order term in the expansion of
the $\beta$th cumulant of the signed-root statistic associated with
(12).  Under no location shift, $k_{11}$ and $k_{31}$
coincide with $k_1$ and $k_3$ in~\eqref{eq:k11k21}, while $k_2=0$.

\medskip
\noindent\textbf{Notation for the proof.} In addition to the quantities above, set
\begin{itemize}
        \item $A_1 \Let \frac{1}{n}\sum_{i=1}^n h(X_i),~A_2 \Let \frac{1}{n}\sum_{i=1}^n h(X_i)^2 - \alpha_{2}$.
        \item $\tilde A_2 \Let \frac{1}{n}\sum_{i=1}^n \D h(X_i) \Sigma \D h(X_i)^\top - \tilde \alpha_{2},~\tilde A_3 \Let \frac{1}{n}\sum_{i=1}^n \D h(X_i) \Sigma \D^2 h(X_i) \Sigma \D h(X_i)^\top - \tilde \alpha_{3}$.
\end{itemize}
\begin{proof}[Proof of Theorem~\MainBartlettCoverageNumber]
    Using Theorem~\MainHigherExpansionNumber{} and noting that when $d_h = 1$,
    \[
    \hat z_{1-\varrho} = \frac{A_2 + \alpha_2}{\tilde A_2 + \tilde \alpha_2} \chi^2_{1;1-\varrho},
    \]
    for the $(1-\varrho)$-quantile $\chi^2_{1;1-\varrho}$ of the chi-square distribution with one degree of freedom, we get
    \begin{align}\label{eq:nRnh_higher}
        n R_n(h) - \hat z_{1-\varrho} = \frac{A_2 + \alpha_2}{\tilde A_2 + \tilde \alpha_2} \left(\hat \gamma^2 - \chi^2_{1;1-\varrho} \right) + \hat \eps_n,
    \end{align}
    where we recall that $\abs{\hat \eps_n} \leq \hat \delta_n$ with probability at least $1 - O\pare{n^{-\frac{3}{2}}}$, for a deterministic sequence $\hat \delta_n = \widetilde{O}\left(n^{-\frac{3}{2}}\right)$, $\hat \gamma = \sqrt{n}\pare{R_1 + R_2 + R_3}$, and $R_j = O_p(n^{-j / 2}), j=1,2,3$. Concretely,
    \begin{align*}
        & R_1 = \frac{A_1}{\alpha_2^{\half}},~R_2 = -\frac{A_1 A_2}{2 \alpha_2^{\frac{3}{2}}}  + \frac{\tilde \alpha_3 A_1^2}{2\alpha_2^{\half}\tilde \alpha_2^2}, \\
        & R_3 = \frac{3 A_1 A_2^2}{8\alpha_2^{\frac{5}{2}}} - \frac{\tilde \alpha_3 A_1^2 A_2}{4\alpha_2^{\frac{3}{2}}\tilde \alpha_2^2}  - \frac{\tilde \alpha_3 A_1^2 \tilde A_2}{\alpha_2^{\frac{1}{2}}\tilde \alpha_2^{3}} 
        + \frac{A_1^2 \tilde A_3}{2\alpha_2^{\half}\tilde \alpha_2^2} + \left( - \frac{\tilde \alpha_3^2}{8 \alpha_{2}^{\half}\tilde \alpha_2^4} + \frac{\tilde \alpha_4}{2\alpha_2^{\half}\tilde \alpha_2^3}\right)A_1^3,
    \end{align*} 
    where
    \begingroup\compactdisplaytrue
    \begin{align*}
        \tilde \alpha_4 \Let & -\EE\left[\left( h_{\alpha \alpha'}  h_{\beta'\gamma'}  h_{\beta} h_{\omega'} \Sigma_{\omega'\gamma'}\Sigma_{\beta \alpha}\Sigma_{\beta'\alpha'}\right)(X)\right] - \frac{1}{3} \EE\left[\left(h_{\alpha \alpha' \alpha''} h_{\beta} h_{\beta'}  h_{\beta''}\Sigma_{\beta \alpha}\Sigma_{\beta'\alpha'}\Sigma_{\beta'' \alpha''}\right)(X)\right]\\
        & + \frac{9}{4\tilde \alpha_2} \EE \left[\left( h_{\alpha \alpha'}  h_{\beta}  h_{\beta'} \Sigma_{\beta \alpha}\Sigma_{\beta'\alpha'}\right)(X)\right]^2.
    \end{align*}
    \endgroup
    
    Applying the Edgeworth expansion to $\hat \gamma$ \cite[Theorem 2]{bhattacharya1978validity} and using the same reasoning as in \eqref{eq:deltamethod}, we get
    \begin{align*}
        \PP^\star\left(n R_n(h) \leq \hat z_{1-\varrho}\right) = \int_{-\sqrt{\chi^2_{1;1-\varrho}}}^{\sqrt{\chi^2_{1;1-\varrho}}}\phi(v) \diff v + \underbrace{\frac{1}{n} \int_{-\sqrt{\chi^2_{1;1-\varrho}}}^{\sqrt{\chi^2_{1;1-\varrho}}} p(v)\phi(v) \diff v}_{(F)} + \widetilde{O}\left(n^{-\frac{3}{2}}\right),
    \end{align*}
    where $\phi$ is the density function of a standard normal distribution $\mc N(0,1)$, and the polynomial $p$ is
    \begin{align*}
        p(v) = \frac{k_{11}^2 + k_{22}}{2}(v^2 - 1) + \frac{4k_{11}k_{31}+k_{42}}{24}(v^4 - 6v^2 + 3) + \frac{k_{31}^2}{72}(v^6 - 15 v^4 + 45 v^2 - 15).
    \end{align*}
    The derivation of the constants $k$'s (the cumulants expansion of $\hat \gamma$) is collected in Section~\ref{app:thm6_5}. 

    Then we rewrite part $(F)$ in Theorem~\MainBartlettCoverageNumber{} as
    \begin{align*}
        (F) = \frac{1}{n}\sum_{r = 0}^{3} B_{r} G_{2r+1}(\chi^2_{1;1-\varrho})
        = \frac{1}{n}  g_{1}\pare{\chi^2_{1;1-\varrho}} \sum_{i = 1}^{3} C_{i} \pare{\chi^2_{1;1-\varrho}}^i,
    \end{align*}
    where $G_r(\cdot)$, $g_r(\cdot)$, and $\chi^2_{r;1-\varrho}$ are the CDF, the density, and the $(1-\varrho)$-quantile of a chi-square distribution with $r$ degrees of freedom. The proof is complete.
\end{proof}

\subsubsection{Proof of Proposition~\MainBartlettOneNumber}
\begin{proof}[Proof of Proposition~\MainBartlettOneNumber]
    Using Theorem~\MainBartlettCoverageNumber, we get
    \begin{align*}
        \PP^\star\left(n R_n(h) \leq \hat z_{1-\varrho}\right) = G_1\pare{\chi^2_{1;1-\varrho}} +
        \frac{1}{n}  g_{1}\pare{\chi^2_{1;1-\varrho}} q\pare{\chi^2_{1;1-\varrho}} + \widetilde{O}\left(n^{-\frac{3}{2}}\right).
    \end{align*}
    Let $c = \frac{1}{n} q\pare{\chi^2_{1;1-\varrho}} / \chi^2_{1;1-\varrho}$. Then
    \begin{align*}
        &\PP^\star\left(n R_n(h) \leq (1 - c)\hat z_{1-\varrho}\right) \\
        =& G_1\pare{\pare{1-c}\chi^2_{1;1-\varrho}} +
        \frac{1}{n}  g_{1}\pare{(1 - c)\chi^2_{1;1-\varrho}} q\pare{(1 - c)\chi^2_{1;1-\varrho}} + \widetilde{O}\left(n^{-\frac{3}{2}}\right)\\
        =& G_1\pare{\chi^2_{1;1-\varrho} - \frac{1}{n} q\pare{\chi^2_{1;1-\varrho}}} +
        \frac{1}{n}  g_{1}\pare{\chi^2_{1;1-\varrho}} q\pare{\chi^2_{1;1-\varrho}} + \widetilde{O}\left(n^{-\frac{3}{2}}\right)\\
        =& G_1\pare{\chi^2_{1;1-\varrho}} -  \frac{1}{n}  g_{1}\pare{\chi^2_{1;1-\varrho}} q\pare{\chi^2_{1;1-\varrho}} +
        \frac{1}{n}  g_{1}\pare{\chi^2_{1;1-\varrho}} q\pare{\chi^2_{1;1-\varrho}} + \widetilde{O}\left(n^{-\frac{3}{2}}\right)\\
        =& G_1\pare{\chi^2_{1;1-\varrho}} + \widetilde{O}\left(n^{-\frac{3}{2}}\right)\\
        =& 1-\varrho +\widetilde{O}\left(n^{-\frac{3}{2}}\right),
    \end{align*}
    where the second equality uses $c = O\pare{n^{-1}}$ and the third uses the
    Taylor expansion of $G_1(\cdot)$ around $\chi^2_{1;1-\varrho}$.
\end{proof}
\subsubsection{Proof of Proposition~\MainBartlettTwoNumber}
\begin{proof}[Proof of Proposition~\MainBartlettTwoNumber]
    Note that 
    \begin{align*}
        \PP^\star\left(n R_n(h) \leq \hat z_{1-\varrho}\right) =& \PP^\star\left(\frac{\chi^2_{1;1-\varrho}}{\hat z_{1-\varrho}}n R_n(h) \leq \chi^2_{1;1-\varrho} \right)\\
        =& \PP^\star\left(\frac{\mc V_n}{\mc W_n} n R_n(h) \leq \chi^2_{1;1-\varrho} \right).
    \end{align*}
    Let $S = \frac{\mc V_n}{\mc W_n} n R_n(h)$. Comparing with \eqref{eq:nRnh_higher}, where $A_2 + \alpha_2 = \mc W_n$ and $\tilde A_2 + \tilde \alpha_2 = \mc V_n$, we get 
    \[
    S = \hat \gamma^2 + \frac{\mc V_n}{\mc W_n} \hat \eps_n.
    \]

    Recall the bounds derived in \eqref{eq:bndVnWn}. Under Assumption~\MainSmoothnessNumber, $h(X)$ and $\D h(X)$ have exponential moments, so $h(X)^2$ and $\D h(X)\Sigma\D h(X)^\top$ have finite moments of all orders. Using the same truncation and concentration argument as in Lemma~\ref{lem:tail_avg2}, we obtain
    $\abs{\mc W_n - W}, \abs{\mc V_n - V}
    =\widetilde{O}\pare{\sqrt{\log(n)/n}}$
    with probability at least $1 - O\pare{n^{-\frac{3}{2}}}$, where $W= \EE\left[h(X)^2\right]$ and $V = \EE\left[\D h(X)\Sigma\D h(X)^\top\right]$. Therefore, $\abs{\frac{\mc V_n}{\mc W_n} \hat \eps_n} \leq  \delta_n$ with probability at least $1 - O\pare{n^{-\frac{3}{2}}}$, for a deterministic sequence $\delta_n = \widetilde{O}\left(n^{-\frac{3}{2}}\right)$.

    Then, the same reasoning as in \eqref{eq:deltamethod} shows that $S$ and
    $\hat \gamma^2$ share the same Edgeworth expansion up to order $n^{-1}$.
    Finally, applying \cite[Theorem 1]{cox1987approximations} gives
    \begin{align*}
        \PP\opt\pare{\pare{1 + \frac{1}{n} \sum_{k=1}^3 C_k S^{k-1}} S \leq \chi^2_{1;1-\varrho}} = 1-\varrho + \widetilde{O}\left(n^{-\frac{3}{2}}\right),
    \end{align*}
    or equivalently,
    \begin{align*}
        \PP\opt\pare{\pare{1 + \frac{1}{n} \sum_{k=1}^3 C_k S^{k-1}} nR_n(h) \leq \hat z_{1-\varrho}} = 1-\varrho + \widetilde{O}\left(n^{-\frac{3}{2}}\right).
    \end{align*}
    This completes the proof.
\end{proof}

\subsection{Proofs in Section~4}
\label{sec:pf_computation_algorithm}

Throughout this subsection, the significance level is written $\varrho$, while
$\alpha$ denotes the transport multiplier. The notation of
Section~4 is in force: the dual objective $F_\delta$ and value
$\Phi_n(\delta)$ of~(17), the localized value
$\Phi^{\mathrm{loc}}_n(\delta)$ of~(25), the matrix
$A(x)=\D\mf h(x)\Sigma^{\half}$, the curvature constant $K_{\mathrm{curv}}$
of~(21), the eigenvalues $\ell_n,U_n$
of~(23), the box
$\mc K_\delta=[\underline\alpha_\delta,\overline\alpha_\delta]\times\mathbb B_2$
of~(24).  We work on the
observable event $\ell_n>0$.  The following elementary consequences of
$0<\delta\leq\delta_0$ are used repeatedly; each follows from
$K_{\mathrm{curv}}\sqrt{\delta/\ell_n}\leq(3+2\sqrt{\varkappa_n})^{-1}\leq1/5$:
\begin{align}
\begin{split}
    \text{(B1)}\quad
    &\underline\alpha_\delta\geq K_{\mathrm{curv}}\pare{1+\sqrt{\varkappa_n}}\geq 2K_{\mathrm{curv}},
    \qquad
    \underline\alpha_\delta\geq\frac{2}{5}\sqrt{\frac{\ell_n}{\delta}},
    \qquad
    \overline\alpha_\delta-\underline\alpha_\delta
    \leq\sqrt{\frac{U_n}{\delta}}\ ;\\
    \text{(B2)}\quad
    &2\alpha-K_{\mathrm{curv}}\geq\frac35\sqrt{\frac{\ell_n}{\delta}}
    \ \text{ and }\
    2\alpha-K_{\mathrm{curv}}\geq K_{\mathrm{curv}}(1+2\sqrt{\varkappa_n}),
    \qquad\forall\,\alpha\geq\underline\alpha_\delta\ ;\\
    \text{(B3)}\quad
    &\frac{K_{\mathrm{curv}}\sqrt{U_n}}{2\alpha-K_{\mathrm{curv}}}\leq\frac{\sqrt{\ell_n}}{2},
    \qquad
    \frac{2\alpha+K_{\mathrm{curv}}}{2\alpha-K_{\mathrm{curv}}}\leq\frac53,
    \qquad\forall\,\alpha\geq\underline\alpha_\delta\ ;\\
    \text{(B4)}\quad
    &2\alpha+K_{\mathrm{curv}}\leq\frac75\sqrt{\frac{U_n}{\delta}},
    \qquad\forall\,\alpha\leq\overline\alpha_\delta\ .
\end{split}
\label{eq:comp_box_facts}
\end{align}
Indeed, (B1) restates $\underline\alpha_\delta\geq\underline\alpha_{\delta_0}
=K_{\mathrm{curv}}(1+\sqrt{\varkappa_n})$ and $K_{\mathrm{curv}}\leq\frac15\sqrt{\ell_n/\delta}$; (B2) follows from
$2\underline\alpha_\delta-K_{\mathrm{curv}}=\sqrt{\ell_n/\delta}-2K_{\mathrm{curv}}$ and (B1); the first part of (B3) uses
$2\alpha-K_{\mathrm{curv}}\geq K_{\mathrm{curv}}(1+2\sqrt{\varkappa_n})\geq2K_{\mathrm{curv}}\sqrt{\varkappa_n}$, the second uses
$\alpha\geq2K_{\mathrm{curv}}$; and (B4) uses $2\overline\alpha_\delta+K_{\mathrm{curv}}=\sqrt{U_n/\delta}+2K_{\mathrm{curv}}$ and
$K_{\mathrm{curv}}\leq\frac15\sqrt{\ell_n/\delta}\leq\frac15\sqrt{U_n/\delta}$.  In particular,
$\sqrt{\ell_n}-K_{\mathrm{curv}}\sqrt{\delta}\geq\frac45\sqrt{\ell_n}>0$ for $\delta\leq\delta_0$, and
$\underline\alpha_{\delta_0}=K_{\mathrm{curv}}(1+\sqrt{\varkappa_n})\geq2K_{\mathrm{curv}}$: the condition
$\delta\leq\delta_0$ places the entire box strictly inside the regime $2\alpha>K_{\mathrm{curv}}$ in which
the inner problems are strongly convex, with quantitative slack.  When $K_{\mathrm{curv}}=0$ these facts
hold trivially with $\delta_0=+\infty$.

\subsubsection{Proof of Proposition~4.1}
\label{sec:pf_wp_dro_reformulation}

\begin{proof}
The growth condition needed for the DRO reformulation is already contained in
Assumption~2.5.  Indeed, fix any $x_0\in\RR^{d_X}$.  Applying that assumption with
base point $x_0$ and integrating the Jacobian along the segment from $x_0$ to $x$ gives
\begin{align}
    \norm{\mf h(x)}_2
    &\leq \norm{\mf h(x_0)}_2
    +\norm{\D\mf h(x_0)}_2\norm{x-x_0}_2
    +\frac{\kappa_1(x_0)}{2}\norm{x-x_0}_2^2
    \leq C_{x_0}\left(1+\norm{x-x_0}_{\Sigma}^{2}\right)
    \label{eq:comp_quadratic_growth}
\end{align}
for some finite $C_{x_0}$.  Thus no separate computational growth assumption is required,
but the moment image need not be closed at this quadratic growth rate.  We therefore use its
closure below.

\textit{Step 1: the finite-dimensional moment image.}
For $\PP\in\mc B_\delta(\QQ_n)$, couple $\PP$ with $\QQ_n$ and apply the pointwise
inequality $\norm{x-x_0}_{\Sigma}^2\leq
2\norm{x-X_i}_{\Sigma}^2+2\norm{X_i-x_0}_{\Sigma}^2$ to obtain
\[
    \EE_{\PP}\!\left[\norm{X-x_0}_{\Sigma}^2\right]
    \leq2\delta+\frac{2}{n}\sum_{i=1}^n\norm{X_i-x_0}_{\Sigma}^2.
\]
Together with~\eqref{eq:comp_quadratic_growth}, this ensures that every moment below is finite.
Moreover, $\mc B_\delta(\QQ_n)$ is convex because couplings can be mixed. Define its moment image by
\begin{align}
    \mc M_\delta
    \Let
    \left\{\EE_{\PP}[\mf h(X)]:\PP\in\mc B_\delta(\QQ_n)\right\}
    \subset\RR^{d_h}.
    \label{eq:comp_moment_image}
\end{align}
This set is nonempty and convex, although it need not be closed. By definition,
\begin{align}
    \inf_{\PP\in\mc B_\delta(\QQ_n)}
    \norm{\EE_{\PP}[\mf h(X)]}_2
    =\inf_{m\in\mc M_\delta}\norm{m}_2.
    \label{eq:comp_moment_distance}
\end{align}

\textit{Step 2: finite-dimensional distance duality.}
Let $d_\delta\Let\inf_{m\in\mc M_\delta}\norm{m}_2$ and let
$m^\star$ be the Euclidean projection of the origin onto the nonempty closed
convex set $\overline{\mc M}_\delta$.  For every $\xi\in\mathbb B_2$, a
sequence in $\mc M_\delta$ converging to $m^\star$ gives
$\inf_{m\in\mc M_\delta}\xi^\top m\leq\xi^\top m^\star\leq d_\delta$.
If $d_\delta>0$, the projection inequality
$(m-m^\star)^\top m^\star\geq0$ for
$m\in\overline{\mc M}_\delta$ shows that
$\xi^\star=m^\star/d_\delta$ attains the reverse bound; if $d_\delta=0$,
$\xi=\mf0$ does so.  Consequently,
\begin{align}
    \inf_{m\in\mc M_\delta}\norm{m}_2
    =\sup_{\xi\in\mathbb B_2}\inf_{m\in\mc M_\delta}\xi^\top m
    =
    \sup_{\xi\in\mathbb B_2}\,
    \inf_{\PP\in\mc B_\delta(\QQ_n)}
    \EE_{\PP}\!\left[\xi^\top\mf h(X)\right].
    \label{eq:comp_distance_duality}
\end{align}
This is the only minimax step required; it remains valid even when
$\mc M_\delta$ is not closed.  Also,
$\max_{\xi\in\mathbb B_2}\EE_{\PP}[\xi^\top\mf h(X)]
=\norm{\EE_{\PP}[\mf h(X)]}_2$ for each $\PP$.

\textit{Step 3: Wasserstein DRO duality.}
Fix $\xi\in\mathbb B_2$ and set $f_\xi(x)\Let\xi^\top\mf h(x)$.  The function
$-f_\xi$ is continuous and $\QQ_n$-integrable, and
\eqref{eq:comp_quadratic_growth} gives it a finite (possibly nonzero) quadratic
growth rate in the sense of \cite{gao2023distributionally}.  Applying the
Wasserstein DRO strong duality \cite[Theorem~1]{gao2023distributionally} to $-f_\xi$ over
$\mc B_\delta(\QQ_n)$ with radius $\delta>0$, and changing signs, gives
\begin{align}
    \inf_{\PP\in\mc B_\delta(\QQ_n)}
    \EE_{\PP}[f_\xi(X)]
    =
    \sup_{\alpha\geq0}
    \left\{
        -\alpha\delta +
        \frac1n\sum_{i=1}^n
            \inf_{x\in\RR^{d_X}}
            \left[
                f_\xi(x)+\alpha\norm{x-X_i}_{\Sigma}^{2}
            \right]
    \right\}
    =\sup_{\alpha\geq0}F_\delta(\alpha,\xi).
    \label{eq:dro_duality_wp_threshold}
\end{align}
For multipliers below the effective quadratic growth rate, an inner infimum
on the right may equal $-\infty$; the extended-valued formulation above
includes such multipliers harmlessly.  Taking the supremum over
$\xi\in\mathbb B_2$ and combining \eqref{eq:dro_duality_wp_threshold} with
\eqref{eq:comp_distance_duality}, \eqref{eq:comp_moment_distance}, and the
dual representation of the Euclidean norm proves all equalities in
(19) and their nonnegativity.

\textit{Step 4: equivalence of the threshold comparison.}
By part~(i) and $\Phi_n(\delta)\geq0$, it suffices to show
\begin{align}
    R_n(\mf h)\leq\delta
    \quad\Longleftrightarrow\quad
    \inf_{\PP\in\mc B_\delta(\QQ_n)}
    \norm{\EE_{\PP}[\mf h(X)]}_2=0.
    \label{eq:comp_zero_equivalence}
\end{align}

($\Rightarrow$)  Suppose $R_n(\mf h)\leq\delta$.  Choose $\PP_k$ with
$\EE_{\PP_k}[\mf h(X)]=\mf0$ and write
$\rho_k\Let\Wass_c(\QQ_n,\PP_k)\leq\delta+1/k$.  If $\rho_k\leq\delta$,
put $t_k=1$ and $\widetilde\PP_k=\PP_k$.  Otherwise put
$t_k\Let\delta/\rho_k$ and
$\widetilde\PP_k\Let t_k\PP_k+(1-t_k)\QQ_n$.  Mixing a coupling for
$(\QQ_n,\PP_k)$ with the diagonal coupling shows that
$\widetilde\PP_k\in\mc B_\delta(\QQ_n)$, while
\[
    \EE_{\widetilde\PP_k}[\mf h(X)]
    =(1-t_k)\overline{\mf h}_n\longrightarrow\mf0,
    \qquad
    \overline{\mf h}_n\Let\frac1n\sum_{i=1}^n\mf h(X_i),
\]
because $t_k\to1$ whenever the second case occurs.  Hence the infimum in
\eqref{eq:comp_zero_equivalence} is zero.

($\Leftarrow$)  Assumption~2.4 and finite dimensionality
provide fixed points $y_1,\ldots,y_J\in\RR^{d_X}$ and $r>0$ such that
\begin{align}
    r\mathbb B_2\subseteq
    \conv\{\mf h(y_1),\ldots,\mf h(y_J)\}.
    \label{eq:comp_finite_correction_set}
\end{align}
Choose $\PP_k\in\mc B_\delta(\QQ_n)$ with
$m_k\Let\EE_{\PP_k}[\mf h(X)]\to\mf0$.  If some $m_k=\mf0$, then that
$\PP_k$ is feasible for~(4) and $R_n(\mf h)\leq\delta$.  Otherwise,
let $u_k\Let m_k/\norm{m_k}_2$ and choose a probability measure $S_k$
supported on $\{y_1,\ldots,y_J\}$ with
$\EE_{S_k}[\mf h(X)]=-r u_k$.  Set
\[
    s_k\Let\frac{\norm{m_k}_2}{r+\norm{m_k}_2},
    \qquad
    \PP_k^0\Let(1-s_k)\PP_k+s_kS_k.
\]
Then $\EE_{\PP_k^0}[\mf h(X)]=\mf0$.  Moreover,
$C_0\Let\max_{i\in[n],j\in[J]}\norm{X_i-y_j}_\Sigma^2<\infty$ and the
product coupling give $\Wass_c(\QQ_n,S_k)\leq C_0$, so convexity of the
transport cost yields
\[
    R_n(\mf h)
    \leq\Wass_c(\QQ_n,\PP_k^0)
    \leq(1-s_k)\delta+s_kC_0
    \leq\delta+s_kC_0\longrightarrow\delta.
\]
Thus $R_n(\mf h)\leq\delta$.  Taking
complements in \eqref{eq:comp_zero_equivalence} proves
(20).
\end{proof}

\subsubsection{Proof of Proposition~4.2}
\label{sec:pf_wp_localization}

\begin{proof}
Fix $0<\delta\leq\delta_0$.  We first record two facts used throughout.

\emph{Fact 1 (Lipschitz Jacobian).}  For all $x,y\in\RR^{d_X}$,
\begin{align}
    \norm{A(x)-A(y)}_2
    \leq K_{\mathrm{curv}}\norm{x-y}_{\Sigma}.
    \label{eq:comp_A_lipschitz}
\end{align}
Indeed, for unit vectors $u\in\RR^{d_h}$, $v\in\RR^{d_X}$, writing
$x-y=\Sigma^{1/2}w$ with $\norm{w}_2=\norm{x-y}_\Sigma$, the fundamental
theorem of calculus gives
\[
    \begin{aligned}
    \left|u^\top\left(A(x)-A(y)\right)v\right|
    &\leq\int_0^1
    \left|w^\top\Sigma^{1/2}\D^2\!\left(u^\top\mf h\right)(y+t(x-y))
    \Sigma^{1/2}v\right|\,\diff t\\
    &\leq K_{\mathrm{curv}}\norm{w}_2,
    \end{aligned}
\]
by (22).

\emph{Fact 2 (spectral bounds).}  For every $u\in\RR^{d_h}$,
$\frac1n\sum_{i=1}^n\norm{A(X_i)^\top u}_2^2=u^\top\mc V_nu$, so by
(23),
\begin{align}
    \sqrt{\ell_n}\norm{u}_2
    \leq
    \left(\frac1n\sum_{i=1}^n\norm{A(X_i)^\top u}_2^2\right)^{1/2}
    \leq
    \sqrt{U_n}\norm{u}_2.
    \label{eq:comp_spectral_sandwich}
\end{align}

\textit{Step 1: strong convexity of the inner problems (part (ii)).}
Fix $\alpha\geq\underline\alpha_\delta$, $\xi\in\mathbb B_2$, and $i\in[n]$.  The
inner objective is twice continuously differentiable with Hessian
$H_i(x)=\D^2(\xi^\top\mf h)(x)+2\alpha\Sigma^{-1}$.  For any
$v=\Sigma^{1/2}w\in\RR^{d_X}$, (22) and
$\norm{\xi}_2\leq1$ give
$\abs{v^\top\D^2(\xi^\top\mf h)(x)v}\leq K_{\mathrm{curv}}\norm{w}_2^2
=K_{\mathrm{curv}}\norm{v}_{\Sigma}^2$, whence
\begin{align}
    \lambda(\alpha)\,\Sigma^{-1}
    \preccurlyeq H_i(x)
    \preccurlyeq(2\alpha+K_{\mathrm{curv}})\,\Sigma^{-1}
    \qquad\text{for all }x\in\RR^{d_X}.
    \label{eq:comp_inner_hessian_bounds}
\end{align}
By (B2) of~\eqref{eq:comp_box_facts},
$\lambda(\alpha)\geq2\underline\alpha_\delta-K_{\mathrm{curv}}=\sqrt{\ell_n/\delta}-2K_{\mathrm{curv}}
\geq\frac35\sqrt{\ell_n/\delta}>0$, which is~(28),
so the inner objective is $\lambda(\alpha)$-strongly convex and
$(2\alpha+K_{\mathrm{curv}})$-smooth with respect to $\norm{\cdot}_{\Sigma}$; strong
convexity implies coercivity and hence a unique minimizer $g_i(\alpha,\xi)$.
The condition-ratio bound on $\mc K_\delta$ is (B3).  For the displacement bound, the
first-order condition for
$g_i\Let g_i(\alpha,\xi)$ gives
$\Sigma^{-\half}(g_i-X_i)=-\frac{1}{2\alpha}A(g_i)^\top\xi$, hence, by
\eqref{eq:comp_A_lipschitz} and $\norm{\xi}_2\leq1$,
\[
    2\alpha\norm{g_i-X_i}_\Sigma
    =\norm{A(g_i)^\top\xi}_2
    \leq\norm{A(X_i)}_2+K_{\mathrm{curv}}\norm{g_i-X_i}_\Sigma,
\]
so rearranging and applying (B2) yields
\begin{align}
    \norm{g_i(\alpha,\xi)-X_i}_\Sigma
    \leq\frac{\norm{A(X_i)}_2}{\lambda(\alpha)}
    \leq\frac{5}{3}\,\norm{A(X_i)}_2\sqrt{\frac{\delta}{\ell_n}}.
    \label{eq:comp_per_i_displacement}
\end{align}
This proves part~(ii).

\textit{Step 2: existence and localization of maximizers.}
Evaluating each inner infimum in~(17) at $x=X_i$ gives,
with $\bar{\mf h}_n\Let\frac1n\sum_{i=1}^n\mf h(X_i)$,
\begin{align}
    F_\delta(\alpha,\xi)
    \leq
    -\alpha\delta+\xi^\top\bar{\mf h}_n
    \leq
    -\alpha\delta+\norm{\xi}_2\norm{\bar{\mf h}_n}_2,
    \label{eq:comp_G_linear_upper}
\end{align}
while $F_\delta(\alpha,\mf 0)=-\alpha\delta\leq0$ and $F_\delta(0,\mf 0)=0$,
so $\Phi_n(\delta)\geq0$ and the supremum in~(17) is
unchanged if we restrict to the compact set
$[0,\norm{\bar{\mf h}_n}_2/\delta]\times\mathbb B_2$, on which the upper
semicontinuous function $F_\delta$ (a pointwise infimum of affine functions
of $(\alpha,\xi)$, minus $\alpha\delta$) attains its maximum.  This proves
attainment.

Now suppose $\Phi_n(\delta)>0$ and let $(\alpha^\star,\xi^\star)$ be any
maximizer.  Since $F_\delta(\alpha,\mf 0)\leq0$, we have
$\xi^\star\neq\mf 0$.  By the homogeneity~(18), if
$\norm{\xi^\star}_2<1$ then $t\Let1/\norm{\xi^\star}_2>1$ gives the feasible
point $(t\alpha^\star,t\xi^\star)$ with strictly larger value
$t\Phi_n(\delta)>\Phi_n(\delta)$, a contradiction; hence
$\norm{\xi^\star}_2=1$.  Next, Assumption~2.4 provides
$r_{\mf h}>0$ with $\{u:\norm{u}_2\leq r_{\mf h}\}\subseteq
\conv\{\mf h(x):x\in\RR^{d_X}\}$.  Since the infimum of a linear functional
over a set coincides with its infimum over the convex hull,
\begin{align}
    F_\delta(0,\xi^\star)
    =\inf_{x\in\RR^{d_X}}(\xi^\star)^\top\mf h(x)
    \leq(\xi^\star)^\top\pare{-r_{\mf h}\xi^\star}
    =-r_{\mf h}<0<\Phi_n(\delta),
    \label{eq:comp_G_at_zero}
\end{align}
so $\alpha^\star>0$.

We next bracket $\alpha^\star$.  More generally, fix $\xi\neq\mf 0$ and
suppose $\alpha^\star>0$ maximizes
$G_\xi(\alpha)\Let F_\delta(\alpha,\xi)$ over $[0,\infty)$.  Write
\[
    \alpha_c\Let\frac12\norm{\xi}_2\underline\alpha_\delta.
\]
By (B1),
$2\alpha_c-K_{\mathrm{curv}}\norm{\xi}_2
=\norm{\xi}_2(\underline\alpha_\delta-K_{\mathrm{curv}})>0$.
Thus, for every $\alpha\geq\alpha_c$, each inner objective is strongly
convex and has a unique minimizer $g_i(\alpha,\xi)$.  Define
\[
    D_\xi(\alpha)
    \Let
    \left(\frac1n\sum_{i=1}^n
    \norm{g_i(\alpha,\xi)-X_i}_\Sigma^2\right)^{1/2}.
\]
The classical envelope formula and the inner first-order conditions give
\begin{align}
    G_\xi'(\alpha)
    &=-\delta+D_\xi(\alpha)^2,
    \label{eq:comp_profile_derivative}\\
    2\alpha D_\xi(\alpha)
    &=
    \left(\frac1n\sum_{i=1}^n
    \norm{A(g_i(\alpha,\xi))^\top\xi}_2^2\right)^{1/2}
    \notag\\
    &\geq
    \left(\frac1n\sum_{i=1}^n
    \norm{A(X_i)^\top\xi}_2^2\right)^{1/2}
    -\left(\frac1n\sum_{i=1}^n
    \norm{\left(A(g_i(\alpha,\xi))-A(X_i)\right)^\top\xi}_2^2
    \right)^{1/2}
    \notag\\
    &\geq
    \sqrt{\ell_n}\norm{\xi}_2
    -K_{\mathrm{curv}}\norm{\xi}_2D_\xi(\alpha).
    \label{eq:comp_profile_displacement_lower}
\end{align}
Here the last two inequalities use Minkowski,
\eqref{eq:comp_spectral_sandwich}, and
\eqref{eq:comp_A_lipschitz}.  Rearranging gives
\[
    D_\xi(\alpha)
    \geq
    \frac{\sqrt{\ell_n}\norm{\xi}_2}
    {2\alpha+K_{\mathrm{curv}}\norm{\xi}_2}.
\]
Consequently, whenever
$\alpha\in[\alpha_c,\norm{\xi}_2\underline\alpha_\delta)$,
the definition of $\underline\alpha_\delta$ gives
$2\alpha+K_{\mathrm{curv}}\norm{\xi}_2
<\norm{\xi}_2\sqrt{\ell_n/\delta}$, and hence
$D_\xi(\alpha)>\sqrt\delta$ and $G_\xi'(\alpha)>0$.

Because $G_\xi$ is concave on $[0,\infty)$ and differentiable at
$\alpha_c$, its supporting-line inequality gives, for $0\leq\beta<\alpha_c$,
\[
    G_\xi(\beta)
    \leq G_\xi(\alpha_c)+G_\xi'(\alpha_c)(\beta-\alpha_c)
    <G_\xi(\alpha_c).
\]
Together with strict increase on
$[\alpha_c,\norm{\xi}_2\underline\alpha_\delta)$, this excludes every
maximizer below $\norm{\xi}_2\underline\alpha_\delta$ and proves
\begin{align}
    \alpha^\star
    \geq
    \norm{\xi}_2\,\underline\alpha_\delta .
    \label{eq:wp_alpha_lower_bound}
\end{align}
In particular, $2\alpha^\star>K_{\mathrm{curv}}\norm{\xi}_2$, so $G_\xi$ is differentiable
at the positive unconstrained maximizer $\alpha^\star$ and
\begin{align}
    D_\xi(\alpha^\star)^2=\delta.
    \label{eq:comp_alpha_envelope}
\end{align}
Writing $g_i^\star\Let g_i(\alpha^\star,\xi)$, the inner first-order
condition gives
\begin{align}
    2\alpha^\star\norm{g_i^\star-X_i}_{\Sigma}
    =\norm{A(g_i^\star)^\top\xi}_2.
    \label{eq:comp_inner_foc_energy}
\end{align}
Minkowski, \eqref{eq:comp_spectral_sandwich},
\eqref{eq:comp_A_lipschitz}, and
\eqref{eq:comp_alpha_envelope} now yield the two-sided comparison
\begin{align}
    \pare{\sqrt{\ell_n}-K_{\mathrm{curv}}\sqrt\delta}\norm{\xi}_2
    \leq
    2\alpha^\star\sqrt{\delta}
    \leq
    \pare{\sqrt{U_n}+K_{\mathrm{curv}}\sqrt{\delta}}\norm{\xi}_2,
    \label{eq:comp_rms_gradient_comparison}
\end{align}
so $\alpha^\star\leq\norm{\xi}_2\,\overline\alpha_\delta$.  Applying the
two-sided bracket to the joint
maximizer $(\alpha^\star,\xi^\star)$, which has $\norm{\xi^\star}_2=1$ and
$\alpha^\star>0$ maximizing $G_{\xi^\star}$, yields
$(\alpha^\star,\xi^\star)\in\mc K_\delta$.

\textit{Step 3: the localization identity (part (i)).}
The maximum in~(25) is attained because $\mc K_\delta$ is
compact and $F_\delta$ upper semicontinuous, and
$\Phi^{\mathrm{loc}}_n(\delta)\leq\Phi_n(\delta)$ trivially.  If
$\Phi_n(\delta)>0$, Step~2 shows that some (indeed every) maximizer of
$\Phi_n(\delta)$ lies in $\mc K_\delta$, whence
$\Phi^{\mathrm{loc}}_n(\delta)=\Phi_n(\delta)>0$
and~(26) holds.  If $\Phi_n(\delta)=0$, then
$\Phi^{\mathrm{loc}}_n(\delta)\leq0$, so
$\max\{\Phi^{\mathrm{loc}}_n(\delta),0\}=0=\Phi_n(\delta)$.  Combining
(26) with
Proposition~4.1(ii) yields
(27).

For the sharper boundary statement used after the proposition, suppose
$R_n(\mf h)<\delta$.  By the definition of $R_n(\mf h)$, there is a
distribution $\PP_0$ satisfying $\EE_{\PP_0}[\mf h(X)]=\mf0$ and a coupling
with $\QQ_n$ whose transportation cost $r$ is strictly smaller than $\delta$.
Weak duality then gives, for every $(\alpha,\xi)\in\mc K_\delta$,
\[
    F_\delta(\alpha,\xi)
    \leq-\alpha\delta
    +\xi^\top\EE_{\PP_0}[\mf h(X)]+\alpha r
    =-\alpha(\delta-r)
    \leq-\underline\alpha_\delta(\delta-r)<0.
\]
Hence $\Phi_n^{\mathrm{loc}}(\delta)<0$.  Together with
(27), this proves
$\{\Phi_n^{\mathrm{loc}}(\delta)=0\}\subseteq
\{R_n(\mf h)=\delta\}$.

\textit{Step 4: strong concavity in $\xi$ (part (iii)).}
Fix $\alpha\geq\underline\alpha_\delta$ and write $g_i=g_i(\alpha,\xi)$ for
the unique minimizers from Step~1, and
\[
    D_g
    \Let
    \left(
        \frac{1}{n}\sum_{i=1}^n
        \norm{g_i-X_i}_{\Sigma}^2
    \right)^{1/2}.
\]
The first-order condition $2\alpha\norm{g_i-X_i}_\Sigma
=\norm{A(g_i)^\top\xi}_2$ and the triangle inequality in the empirical $L^2$
norm give, using \eqref{eq:comp_spectral_sandwich},
\eqref{eq:comp_A_lipschitz}, and $\norm{\xi}_2\leq1$,
\begin{align}
    2\alpha D_g
    =
    \left(
        \frac{1}{n}\sum_{i=1}^n
        \norm{A(g_i)^\top\xi}_2^2
    \right)^{1/2}
    \leq
    \sqrt{U_n}+K_{\mathrm{curv}}D_g,
    \qquad\text{whence}\qquad
    D_g\leq\frac{\sqrt{U_n}}{\lambda(\alpha)}.
    \label{eq:comp_g_rms_displacement}
\end{align}
Transferring the empirical nondegeneracy from the atoms $X_i$ to the
minimizers $g_i$, for every $u\in\RR^{d_h}$,
\begin{align}
    \left(
        \frac{1}{n}\sum_{i=1}^n\norm{A(g_i)^\top u}_2^2
    \right)^{1/2}
    \geq
    \sqrt{\ell_n}\norm{u}_2-K_{\mathrm{curv}}\norm{u}_2D_g
    \geq
    \left(
        \sqrt{\ell_n}
        -\frac{K_{\mathrm{curv}}\sqrt{U_n}}{\lambda(\alpha)}
    \right)\norm{u}_2,
    \label{eq:comp_nondegeneracy_at_g}
\end{align}
and the coefficient is positive: by (B2),
$\lambda(\alpha)\geq K_{\mathrm{curv}}(1+2\sqrt{\varkappa_n})
\geq K_{\mathrm{curv}}(1+\sqrt{\varkappa_n})$, so (when $K_{\mathrm{curv}}>0$)
\begin{align}
    \frac{K_{\mathrm{curv}}\sqrt{U_n}}{\lambda(\alpha)}
    \leq
    \frac{\sqrt{U_n}}{1+\sqrt{\varkappa_n}}
    =
    \sqrt{\ell_n}\,\frac{\sqrt{\varkappa_n}}{1+\sqrt{\varkappa_n}}
    \leq
    \sqrt{\ell_n}\left(1-\frac{1}{1+\sqrt{\varkappa_n}}\right),
    \label{eq:comp_margin_bound}
\end{align}
whence
$\sqrt{\ell_n}-K_{\mathrm{curv}}\sqrt{U_n}/\lambda(\alpha)
\geq\sqrt{\ell_n}/(1+\sqrt{\varkappa_n})>0$; for $K_{\mathrm{curv}}=0$ the coefficient
equals $\sqrt{\ell_n}$ and the same lower bound holds trivially.  This also
proves the second inequality in (29) after
squaring and dividing by $2\alpha+K_{\mathrm{curv}}$.

It remains to bound the Hessian of $F_\delta(\alpha,\cdot)$.  Since
$2\alpha>K_{\mathrm{curv}}\norm{\xi}_2$ holds on an open neighborhood of $\mathbb B_2$, the
minimizers $g_i(\alpha,\xi)$ are unique there, and by
$\mf h\in C^3$ (Assumption~2.5), $H_i(g_i)\succ0$, and the implicit
function theorem, $\xi\mapsto g_i(\alpha,\xi)$ is continuously
differentiable with
$\D_\xi g_i=-H_i(g_i)^{-1}\D\mf h(g_i)^\top$.  Danskin's envelope theorem
for the (unique-minimizer) infimum gives
$\D_\xi F_\delta(\alpha,\xi)=\frac1n\sum_{i=1}^n\mf h(g_i(\alpha,\xi))^\top$,
and differentiating once more,
\begin{align}
    \D_\xi^2F_\delta(\alpha,\xi)
    =
    -\frac{1}{n}\sum_{i=1}^n
    \D\mf h(g_i)H_i(g_i)^{-1}\D\mf h(g_i)^\top.
    \label{eq:comp_F_xi_hessian}
\end{align}
The upper bound in \eqref{eq:comp_inner_hessian_bounds} implies
$H_i(g_i)^{-1}\succcurlyeq(2\alpha+K_{\mathrm{curv}})^{-1}\Sigma$, so for every
$u\in\RR^{d_h}$,
\[
    u^\top\D\mf h(g_i)H_i(g_i)^{-1}\D\mf h(g_i)^\top u
    \geq
    \frac{\norm{A(g_i)^\top u}_2^2}{2\alpha+K_{\mathrm{curv}}}.
\]
Averaging over $i$ and applying \eqref{eq:comp_nondegeneracy_at_g} yields
\[
    \D_\xi^2F_\delta(\alpha,\xi)
    \preccurlyeq
    -\mu(\alpha)I_{d_h}
    \qquad\text{for every }\xi\in\mathbb B_2.
\]
Since $\mathbb B_2$ is convex, integrating this bound along line segments in
$\mathbb B_2$ shows that $F_\delta(\alpha,\cdot)$ is
$\mu(\alpha)$-strongly concave on $\mathbb B_2$, proving part~(iii).
\end{proof}

\subsubsection{Localization and margin lemmas}
\label{sec:pf_comp_margin}

\begin{lemma}[Upper bound for the optimal transport multiplier]
\label{lem:alpha_upper}
Suppose Assumption~2.5 holds and
(21) is satisfied.  Fix $\delta>0$ and
$\xi\in\mathbb B_2$.  Every maximizer of
$\alpha\mapsto F_\delta(\alpha,\xi)$ over $[0,\infty)$ satisfies
\begin{align}
    \alpha^\star
    \leq
    \frac{\norm{\xi}_2}{2}
    \pare{\sqrt{\frac{U_n}{\delta}}+K_{\mathrm{curv}}}
    \leq
    \overline\alpha_\delta .
    \label{eq:comp_alpha_upper}
\end{align}
\end{lemma}

\begin{proof}
If $\alpha^\star=0$ the claim is trivial, so assume $\alpha^\star>0$.  The
claim is immediate if
$2\alpha^\star\leq K_{\mathrm{curv}}\norm{\xi}_2$.  Otherwise each inner objective is
strongly convex in a neighborhood of $\alpha^\star$, so its unique minimizer
$g_i^\star$ obeys the classical envelope formula.  Since $\alpha^\star$ is
a positive unconstrained maximizer,
\[
    \frac1n\sum_{i=1}^n\norm{g_i^\star-X_i}_\Sigma^2=\delta.
\]
The inner first-order conditions, Minkowski,
\eqref{eq:comp_spectral_sandwich}, and \eqref{eq:comp_A_lipschitz} give
\[
    2\alpha^\star\sqrt\delta
    =\left(\frac1n\sum_{i=1}^n
    \norm{A(g_i^\star)^\top\xi}_2^2\right)^{1/2}
    \leq\pare{\sqrt{U_n}+K_{\mathrm{curv}}\sqrt\delta}\norm{\xi}_2,
\]
which is~\eqref{eq:comp_alpha_upper}.  Thus the bound holds on both sides of
the strong-convexity threshold, without requiring subquadratic coercivity.
\end{proof}

For the uniform curvature bounds, stopping rules, and complexity ledger below, define
\begin{align}
\begin{split}
    \mu_\delta
    &\Let\frac{\ell_n}{8}\sqrt{\frac{\delta}{U_n}},
    \qquad
    \mathsf L_\delta
    \Let\frac{15U_n}{4}\sqrt{\frac{\delta}{\ell_n}},
    \qquad
    \frac{\mathsf L_\delta}{\mu_\delta}=30\varkappa_n^{3/2},\\
    \Lambda_\delta&\Let3\varkappa_n\delta,
    \qquad
    \mathsf s_n\Let\max_{i\in[n]}\norm{A(X_i)}_2,
    \qquad
    \bar{\mathsf m}_n^2\Let\frac1n\sum_{i=1}^n\norm{\mf h(X_i)}_2^2,\\
    \mathsf M_\delta&\Let\bar{\mathsf m}_n
    +\frac{20}{9}\,d_hU_n\sqrt{\frac{\delta}{\ell_n}}.
\end{split}
    \label{eq:comp_constants}
\end{align}
Thus $\mu_\delta$ and $\mathsf L_\delta$ are box-uniform strong-concavity and smoothness
moduli in $\xi$, $\Lambda_\delta$ bounds $|\partial_\alpha F_\delta|$ on $\mc K_\delta$,
$\mathsf M_\delta$ bounds $\norm{\D_\xi F_\delta}_2$, and the data scales $\mathsf s_n$ and
$\bar{\mathsf m}_n$ enter the complexity bound only through logarithms.  These claims are
established in the next lemma.

\begin{lemma}[Joint curvature structure of the dual objective]
\label{lem:joint_concavity}
Suppose Assumption~2.5 holds,
(21) is satisfied, and $\ell_n>0$.  Let
$0<\delta\leq\delta_0$.  Write $z=(\alpha,\xi)$ and let
$\Omega\Let\{(\alpha,\xi):2\alpha>K_{\mathrm{curv}}\norm{\xi}_2\}$, an open convex cone
containing $\mc K_\delta$.
\begin{enumerate}[label=\textup{(\roman*)}]
    \item \textup{(Concavity and homogeneity.)}  $F_\delta$ is jointly
    concave and upper semicontinuous on $[0,\infty)\times\RR^{d_h}$, and
    positively homogeneous of degree one:
    $F_\delta(tz)=tF_\delta(z)$ for all $t>0$.
    \item \textup{(Smoothness.)}  On $\Omega$, each inner problem has a
    unique minimizer $g_i(z)$ which is $C^1$ in $z$, and
    $F_\delta\in C^2(\Omega)$ with
    \begin{align}
        \partial_\alpha F_\delta(z)
        =-\delta+\frac1n\sum_{i=1}^n\norm{g_i(z)-X_i}_\Sigma^2,
        \qquad
        \D_\xi F_\delta(z)
        =\frac1n\sum_{i=1}^n\mf h\pare{g_i(z)}.
        \label{eq:comp_gradients}
    \end{align}
    \item \textup{(Failure of joint strong concavity.)}  For every
    $z\in\Omega$,
    \begin{align}
        \D^2F_\delta(z)\,[z,z]=0 .
        \label{eq:comp_euler_degeneracy}
    \end{align}
    Consequently $F_\delta$ is not strongly concave on any convex subset of
    $\Omega$ with nonempty interior.  In particular, if
    $\underline\alpha_\delta<\overline\alpha_\delta$, then $\mc K_\delta$
    has nonempty interior and joint strong concavity fails on $\mc K_\delta$.
    \item \textup{(Uniform partial curvature on the box.)}  For every
    $(\alpha,\xi)\in\mc K_\delta$,
    \begin{align}
        -\mathsf L_\delta I_{d_h}
        \preccurlyeq
        \D_\xi^2F_\delta(\alpha,\xi)
        \preccurlyeq
        -\mu_\delta I_{d_h},
        \qquad
        \abs{\partial_\alpha F_\delta(\alpha,\xi)}\leq\Lambda_\delta,
        \qquad
        \norm{\D_\xi F_\delta(\alpha,\xi)}_2\leq\mathsf M_\delta .
        \label{eq:comp_box_curvature}
    \end{align}
    \item \textup{(Profile.)}  $G_\delta(\alpha)=
    \max_{\xi\in\mathbb B_2}F_\delta(\alpha,\xi)$ is concave and
    $\Lambda_\delta$-Lipschitz on
    $[\underline\alpha_\delta,\overline\alpha_\delta]$, and for each such
    $\alpha$ the maximizer $\xi^\star(\alpha)$ is unique with the partial
    quadratic growth
    \begin{align}
        G_\delta(\alpha)-F_\delta(\alpha,\xi)
        \geq\frac{\mu_\delta}{2}\norm{\xi-\xi^\star(\alpha)}_2^2,
        \qquad\forall\,\xi\in\mathbb B_2 .
        \label{eq:comp_partial_QG}
    \end{align}
\end{enumerate}
\end{lemma}

\begin{proof}
\textit{(i)}  For fixed $x$ and $i$, the map
$(\alpha,\xi)\mapsto\xi^\top\mf h(x)+\alpha\norm{x-X_i}_\Sigma^2$ is affine;
an infimum of affine functions is concave and upper semicontinuous, and so
is the average, minus the linear term $\alpha\delta$.  Homogeneity holds
because $\inf_x\{t\xi^\top\mf h(x)+t\alpha\norm{x-X_i}_\Sigma^2\}
=t\inf_x\{\xi^\top\mf h(x)+\alpha\norm{x-X_i}_\Sigma^2\}$ for $t>0$, and
$-(t\alpha)\delta=t(-\alpha\delta)$.

\textit{(ii)}  Fix $z_0=(\alpha_0,\xi_0)\in\Omega$.  On a neighborhood of
$z_0$ contained in $\Omega$, the inner Hessian satisfies
$\D^2\pare{\xi^\top\mf h}(x)+2\alpha\Sigma^{-1}
\succcurlyeq(2\alpha-K_{\mathrm{curv}}\norm{\xi}_2)\Sigma^{-1}\succ0$
by~(22) (applied to the unit vector
$\xi/\norm{\xi}_2$ and scaled), so the inner objective is strongly convex,
coercive, and has a unique minimizer $g_i(z)$.  The first-order condition
$\Psi_i(x;z)\Let\D\mf h(x)^\top\xi+2\alpha\Sigma^{-1}(x-X_i)=\mf 0$ has a
continuously differentiable defining map ($\mf h\in C^3$ by
Assumption~2.5) whose $x$-Jacobian at $g_i(z)$ is the inner Hessian,
which is nonsingular on $\Omega$; the implicit function theorem
\cite[Theorem~3.2.1]{krantz2002implicit} then makes $g_i$ locally $C^1$.
Since the minimizer is unique, Danskin's theorem yields that
$\phi_i(z)\Let\min_x\{\xi^\top\mf h(x)+\alpha\norm{x-X_i}_\Sigma^2\}$ is
differentiable with
$\nabla\phi_i(z)=\pare{\norm{g_i(z)-X_i}_\Sigma^2,\ \mf h(g_i(z))}$, which is
$C^1$ because $g_i$ is; hence $\phi_i\in C^2(\Omega)$,
$F_\delta\in C^2(\Omega)$, and~\eqref{eq:comp_gradients} holds.

\textit{(iii)}  Each $\phi_i$ is positively homogeneous of degree one on the
cone $\Omega$ by the computation in (i).  Differentiating
$\phi_i(tz)=t\phi_i(z)$ in $z$ gives $\nabla\phi_i(tz)=\nabla\phi_i(z)$ for
$t>0$; differentiating this identity in $t$ at $t=1$ gives
$\D^2\phi_i(z)\,z=\mf 0$.  Since the term $-\alpha\delta$ is linear, the same
holds for $F_\delta$, proving~\eqref{eq:comp_euler_degeneracy}.  If
$F_\delta$ were $\mu$-strongly concave on a convex set
$\mc C\subseteq\Omega$ with interior point $\bar z$, then
$s\mapsto F_\delta(\bar z+s\bar z)$, defined for $|s|$ small, would satisfy
$\frac{\diff^2}{\diff s^2}F_\delta(\bar z+s\bar z)\big|_{s=0}
=\D^2F_\delta(\bar z)[\bar z,\bar z]\leq-\mu\norm{\bar z}_2^2<0$,
contradicting~\eqref{eq:comp_euler_degeneracy}.

\textit{(iv)}  The upper bound
$\D_\xi^2F_\delta\preccurlyeq-\mu(\alpha)I$ with
$\mu(\alpha)$ as in~(29) is
Proposition~4.2(iii).  On the box, (B3) gives
$\sqrt{\ell_n}-K_{\mathrm{curv}}\sqrt{U_n}/(2\alpha-K_{\mathrm{curv}})\geq\sqrt{\ell_n}/2$ and (B4) gives
$2\alpha+K_{\mathrm{curv}}\leq\frac75\sqrt{U_n/\delta}$, whence
$\mu(\alpha)\geq\frac{\ell_n/4}{(7/5)\sqrt{U_n/\delta}}
=\frac{5\ell_n}{28}\sqrt{\delta/U_n}\geq\mu_\delta$.  For the smoothness bound,
\eqref{eq:comp_F_xi_hessian} and
$H_i(g_i)^{-1}\preccurlyeq(2\alpha-K_{\mathrm{curv}})^{-1}\Sigma$
(from~\eqref{eq:comp_inner_hessian_bounds}) give
\[
    \mf 0\succcurlyeq
    \D^2_\xi F_\delta
    \succcurlyeq
    -\frac{1}{2\alpha-K_{\mathrm{curv}}}\cdot\frac1n\sum_{i=1}^nA(g_i)A(g_i)^\top ,
\]
and, by the triangle inequality in the product $L^2$ space,
\eqref{eq:comp_A_lipschitz}, \eqref{eq:comp_g_rms_displacement}, and (B3),
\[
\begin{aligned}
    \lambda_{\max}\pare{\frac1n\sum_{i=1}^nA(g_i)A(g_i)^\top}
    &\leq
    \pare{\sqrt{U_n}+K_{\mathrm{curv}}\,
    \frac{\sqrt{U_n}}{2\alpha-K_{\mathrm{curv}}}}^{\!2}\\
    &\leq\pare{\sqrt{U_n}+\frac{\sqrt{\ell_n}}{2}}^{\!2}
    \leq\frac94\,U_n ,
\end{aligned}
\]
so that, using (B2),
$\norm{\D^2_\xi F_\delta}_2\leq\frac{9U_n/4}{(3/5)\sqrt{\ell_n/\delta}}
=\mathsf L_\delta$.  For the $\alpha$-derivative,
\eqref{eq:comp_gradients} and~\eqref{eq:comp_g_rms_displacement} give
$-\delta\leq\partial_\alpha F_\delta\leq-\delta+U_n/(2\alpha-K_{\mathrm{curv}})^2
\leq-\delta+\frac{25}{9}\varkappa_n\delta$, and since $\varkappa_n\geq1$,
$\abs{\partial_\alpha F_\delta}\leq\frac{25}{9}\varkappa_n\delta
\leq\Lambda_\delta$.  Finally, by~\eqref{eq:comp_per_i_displacement} and
\eqref{eq:comp_A_lipschitz}, the segment from $X_i$ to $g_i$ satisfies
$\sup_{y\in[X_i,g_i]}\norm{A(y)}_2
\leq\norm{A(X_i)}_2(1+\tfrac53K_{\mathrm{curv}}\sqrt{\delta/\ell_n})
\leq\tfrac43\norm{A(X_i)}_2$, whence
\begin{align}
    \norm{\mf h(g_i)-\mf h(X_i)}_2
    \leq\frac43\norm{A(X_i)}_2\norm{g_i-X_i}_\Sigma
    \leq\frac{20}{9}\norm{A(X_i)}_2^2\sqrt{\frac{\delta}{\ell_n}},
    \label{eq:comp_h_displacement}
\end{align}
using $\norm{\mf h(y')-\mf h(y)}_2=\|\int_0^1A(y+t(y'-y))
\Sigma^{-\half}(y'-y)\diff t\|_2\leq\sup_{[y,y']}\norm{A}_2
\norm{y'-y}_\Sigma$.  Averaging over $i$, and using
$\frac1n\sum_i\norm{A(X_i)}_2^2\leq\frac1n\sum_i\norm{A(X_i)}_F^2
=\Tr(\mc V_n)\leq d_hU_n$ together with
$\frac1n\sum_i\norm{\mf h(X_i)}_2\leq\bar{\mathsf m}_n$, gives
$\norm{\D_\xi F_\delta}_2\leq\mathsf M_\delta$.

\textit{(v)}  Concavity of $G_\delta$ follows from joint concavity of
$F_\delta$ by partial maximization over the convex set $\mathbb B_2$.
Uniqueness of $\xi^\star(\alpha)$ and~\eqref{eq:comp_partial_QG} follow from
$\mu_\delta$-strong concavity in $\xi$ on the convex set $\mathbb B_2$
(part (iv)); the Lipschitz bound follows because, for
$\alpha,\alpha'$ in the box,
$\abs{G_\delta(\alpha)-G_\delta(\alpha')}\leq
\sup_{\xi\in\mathbb B_2}\abs{F_\delta(\alpha,\xi)-F_\delta(\alpha',\xi)}
\leq\Lambda_\delta\abs{\alpha-\alpha'}$ by part (iv).
\end{proof}

\begin{remark}[Degeneracy along homogeneous rays]
\label{rem:comp_degeneracy}
Identity~\eqref{eq:comp_euler_degeneracy} is Euler's relation for the
$1$-homogeneous function $F_\delta+\alpha\delta$ and shows that its Hessian is
singular along every homogeneous ray in the original $(\alpha,\xi)$
parameterization.  Consequently, whenever an optimizer can be contracted
along a nontrivial feasible homogeneous ray, a quadratic-growth (PL-type)
constant on $\mc K_\delta$ is necessarily controlled by the optimal value.
Quantitatively, suppose that $z^\star=(\alpha^\star,\xi^\star)$ maximizes
$F_\delta$ over $\mc K_\delta$ with value $\Phi_n(\delta)>0$ and that
$\alpha^\star>\underline\alpha_\delta$.  Then, along the feasible segment
$t\mapsto tz^\star$, $t\in[\underline\alpha_\delta/\alpha^\star,1]$, one has
$F_\delta(z^\star)-F_\delta(tz^\star)=(1-t)\Phi_n(\delta)$ while
$\norm{z^\star-tz^\star}_2^2=(1-t)^2\norm{z^\star}_2^2$, so any quadratic
growth constant is at most
$\Phi_n(\delta)\big/\pare{(1-\underline\alpha_\delta/\alpha^\star)
\norm{z^\star}_2^2}$.  Along sequences for which
$1-\underline\alpha_\delta/\alpha^\star$ is bounded away from zero, this
bound vanishes as the decision boundary $\Phi_n(\delta)\downarrow0$ is
approached and is $O\pare{\Phi_n(\delta)\,\delta/\ell_n}$, since
$\norm{z^\star}_2^2\geq\underline\alpha_\delta^2\gtrsim \ell_n/\delta$.
A one-dimensional witness of the same
phenomenon in the $\alpha$ direction alone is obtained as follows: for $\mf h(x)=x$,
$\Sigma=I_1$, one computes
$G_\delta(\alpha)=-\alpha\delta+\abs{\bar X_n}-1/(4\alpha)$ on the relevant
range, whose curvature $\abs{G_\delta''(\alpha)}=1/(2\alpha^3)
\asymp\delta^{3/2}$ vanishes at the rate $n^{-3/2}$ when
$\delta\asymp1/n$.  This is why Theorem~\MainComputableTestNumber
treats $\alpha$ by a one-dimensional search on the concave,
$\Lambda_\delta$-Lipschitz profile $G_\delta$ rather than by joint
gradient ascent.
\end{remark}

\begin{remark}[Euler identity as an implementation check]
\label{rem:comp_euler_check}
On $\Omega$ the gradients~\eqref{eq:comp_gradients} satisfy
$\alpha\,\partial_\alpha F_\delta(\alpha,\xi)
+\xi^\top\D_\xi F_\delta(\alpha,\xi)=F_\delta(\alpha,\xi)$, an exact
identity that provides a costless correctness test for implementations of
the inner solver.
\end{remark}

\begin{lemma}[Two-sided margin for the dual value]
\label{lem:margin}
Suppose Assumptions~2.4 and 2.5 hold,
(21) is
satisfied, and $\ell_n>0$.
\begin{enumerate}[label=\textup{(\roman*)}]
    \item The map $\delta\mapsto\Phi_n(\delta)$ is convex, nonincreasing,
    and nonnegative on $(0,\infty)$, with $\Phi_n(\delta)=0$ if and only if
    $R_n(\mf h)\leq\delta$.
    \item \textup{(Lower bound.)}  For all
    $0<\delta<\delta'\leq\min\{R_n(\mf h),\delta_0\}$,
    \begin{align}
        \Phi_n(\delta)
        \geq
        \underline\alpha_{\delta'}\pare{\delta'-\delta}
        =
        \frac12\pare{\sqrt{\frac{\ell_n}{\delta'}}-K_{\mathrm{curv}}}\pare{\delta'-\delta}.
        \label{eq:comp_margin_lower}
    \end{align}
    \item \textup{(Upper bound.)}  For all $0<\delta<R_n(\mf h)$,
    \begin{align}
        \Phi_n(\delta)
        \leq
        \overline\alpha_\delta\pare{R_n(\mf h)-\delta}
        =
        \frac12\pare{\sqrt{\frac{U_n}{\delta}}+K_{\mathrm{curv}}}\pare{R_n(\mf h)-\delta}.
        \label{eq:comp_margin_upper}
    \end{align}
\end{enumerate}
\end{lemma}

\begin{proof}
Write $F_0(\alpha,\xi)\Let\frac1n\sum_i\inf_x\{\xi^\top\mf h(x)
+\alpha\norm{x-X_i}_\Sigma^2\}$, so that
$F_\delta(\alpha,\xi)=F_0(\alpha,\xi)-\alpha\delta$ and
$\Phi_n(s)=\sup_{\alpha\geq0,\xi\in\mathbb B_2}\{F_0(\alpha,\xi)-\alpha s\}$
is a supremum of affine functions of $s$, hence convex; it is nonincreasing
because the affine pieces have nonpositive slopes $-\alpha$.  Nonnegativity
and the sign characterization are
Proposition~4.1.  This proves (i).

(ii)  Fix $\delta'\in(\delta,\ R_n(\mf h)\wedge\delta_0]$ and take any
$\delta''\in(\delta,\delta')$; then $R_n(\mf h)>\delta''$, so
$\Phi_n(\delta'')>0$, and by
Proposition~4.2(i) (applicable since
$\delta''<\delta_0$) there is a maximizer $(\alpha'',\xi'')$ of
$\Phi_n(\delta'')$ with $\alpha''\geq\underline\alpha_{\delta''}$.  Using
this point as a feasible candidate at radius $\delta$,
\begingroup\compactdisplaytrue
\[
    \Phi_n(\delta)
    \geq F_0(\alpha'',\xi'')-\alpha''\delta
    =\brac{F_0(\alpha'',\xi'')-\alpha''\delta''}+\alpha''(\delta''-\delta)
    =\Phi_n(\delta'')+\alpha''(\delta''-\delta)
    \geq\underline\alpha_{\delta''}\pare{\delta''-\delta},
\]
\endgroup
where the last step uses $\Phi_n(\delta'')\geq0$.  Since
$\delta''\mapsto\underline\alpha_{\delta''}(\delta''-\delta)$ is continuous,
letting $\delta''\uparrow\delta'$ gives~\eqref{eq:comp_margin_lower}.

(iii)  Since $0<\delta<R_n(\mf h)$, we have $\Phi_n(\delta)>0$ by
Proposition~4.1(ii); let
$(\alpha^\star,\xi^\star)$ be a maximizer (existence by
Proposition~4.2(i) when $\delta\leq\delta_0$; for
$\delta>\delta_0$ existence still holds by the compactification argument in
Step~2 of its proof, and Lemma~\ref{lem:alpha_upper} applies for every
$\delta>0$).  Then, with $\rho\Let R_n(\mf h)$, using the candidate
$(\alpha^\star,\xi^\star)$ at radius $\rho$,
\[
    0=\Phi_n(\rho)
    \geq F_0(\alpha^\star,\xi^\star)-\alpha^\star\rho
    =\Phi_n(\delta)-\alpha^\star(\rho-\delta),
\]
where $\Phi_n(\rho)=0$ by part (i).  Hence
$\Phi_n(\delta)\leq\alpha^\star(\rho-\delta)
\leq\overline\alpha_\delta(\rho-\delta)$ by
Lemma~\ref{lem:alpha_upper}.
\end{proof}

\begin{corollary}[Decision band of an $\eps$-accurate dual solve]
\label{cor:comp_decision_band}
In the setting of Lemma~\ref{lem:margin}, let $0<\delta\leq\delta_0/2$ and
$0\leq\eps\leq\frac18\sqrt{\ell_n\delta}$.  If
$R_n(\mf h)>\delta$ and $\Phi_n(\delta)\leq2\eps$, then
\begin{align}
    R_n(\mf h)-\delta
    \leq
    5\,\eps\,\sqrt{\frac{2\delta}{\ell_n}} .
    \label{eq:comp_band}
\end{align}
Conversely, by~\eqref{eq:comp_margin_upper}, whenever
$0<R_n(\mf h)-\delta\leq
2\eps\big/\overline\alpha_\delta$ one has $\Phi_n(\delta)\leq2\eps$; hence
the width~\eqref{eq:comp_band} is sharp up to the factor
$5\sqrt{2\varkappa_n}$.
\end{corollary}

\begin{proof}
Suppose first $R_n(\mf h)\geq2\delta$.  Then
$\delta'\Let2\delta\leq\min\{R_n(\mf h),\delta_0\}$
in~\eqref{eq:comp_margin_lower} gives, using
$K_{\mathrm{curv}}\leq\frac15\sqrt{\ell_n/(2\delta)}$ (valid since $2\delta\leq\delta_0$),
\[
    \Phi_n(\delta)
    \geq\frac12\pare{\sqrt{\frac{\ell_n}{2\delta}}-K_{\mathrm{curv}}}\delta
    \geq\frac25\sqrt{\frac{\ell_n}{2\delta}}\,\delta
    =\frac{\sqrt2}{5}\sqrt{\ell_n\delta}
    >\frac{1}{4}\sqrt{\ell_n\delta}
    \geq2\eps,
\]
contradicting $\Phi_n(\delta)\leq2\eps$.  Hence $R_n(\mf h)<2\delta$, and
$\delta'\Let R_n(\mf h)\leq2\delta\leq\delta_0$
in~\eqref{eq:comp_margin_lower} gives
\[
    2\eps\geq\Phi_n(\delta)
    \geq\frac12\pare{\sqrt{\frac{\ell_n}{R_n(\mf h)}}-K_{\mathrm{curv}}}
    \pare{R_n(\mf h)-\delta}
    \geq\frac25\sqrt{\frac{\ell_n}{2\delta}}\pare{R_n(\mf h)-\delta},
\]
using again $K_{\mathrm{curv}}\leq\frac15\sqrt{\ell_n/R_n(\mf h)}$ (valid since
$R_n(\mf h)\leq\delta_0$).  Rearranging yields~\eqref{eq:comp_band}.
\end{proof}

\subsubsection{Algorithmic building blocks and pseudocode}
\label{sec:pf_comp_blocks}

\begin{lemma}[Inner solver and error propagation]
\label{lem:comp_inner}
In the setting of Proposition~4.2, fix
$(\alpha,\xi)\in\mc K_\delta$ and $i\in[n]$, and write
$\psi_i(x)\Let\xi^\top\mf h(x)+\alpha\norm{x-X_i}_\Sigma^2$,
$m\Let2\alpha-K_{\mathrm{curv}}$, $M\Let2\alpha+K_{\mathrm{curv}}$, $g_i\Let g_i(\alpha,\xi)$.
Consider the preconditioned gradient iteration
$x^{(k+1)}=x^{(k)}-\frac1M\Sigma\nabla\psi_i(x^{(k)})$ started at
$x^{(0)}=X_i$.  Then:
\begin{enumerate}[label=\textup{(\roman*)}]
    \item $\psi_i(x^{(k)})-\psi_i(g_i)\leq
    \pare{1-\frac mM}^k\pare{\psi_i(X_i)-\psi_i(g_i)}
    \leq\pare{\frac25}^{k}\,\frac{2\,\mathsf s_n^2\sqrt{U_n\delta}}{\ell_n}$.
    \item If $\tilde g_i$ satisfies
    $\psi_i(\tilde g_i)-\psi_i(g_i)\leq\eps_x$ with
    $0\leq\eps_x\leq\mathsf s_n^2\sqrt{\delta/\ell_n}$, then
    \begin{align}
        \norm{\tilde g_i-g_i}_\Sigma\leq2\sqrt{\eps_x}
        \pare{\frac{\delta}{\ell_n}}^{\!1/4},
        \qquad
        \norm{\mf h(\tilde g_i)-\mf h(g_i)}_2
        \leq 4\,\mathsf s_n\sqrt{\eps_x}\pare{\frac{\delta}{\ell_n}}^{\!1/4}
        \Let b(\eps_x),
        \label{eq:comp_inner_bias}
    \end{align}
    and $0\leq\psi_i(\tilde g_i)-\psi_i(g_i)\leq\eps_x$ transfers to the
    dual objective: replacing every $g_i$ by $\tilde g_i$ produces
    $\tilde F_\delta$ with $F_\delta\leq\tilde F_\delta\leq
    F_\delta+\eps_x$ and a $\xi$-gradient surrogate
    $\frac1n\sum_i\mf h(\tilde g_i)$ within $b(\eps_x)$ of
    $\D_\xi F_\delta$ in Euclidean norm.
\end{enumerate}
\end{lemma}

\begin{proof}
By Proposition~4.2(ii), $\psi_i$ is $m$-strongly
convex and $M$-smooth with respect to $\norm\cdot_\Sigma$, with
$M/m\leq5/3$.  \textit{(i)}  The descent step gives, by $M$-smoothness in
the $\Sigma$-geometry,
$\psi_i(x^{(k+1)})\leq\psi_i(x^{(k)})
-\frac{1}{2M}\nabla\psi_i(x^{(k)})^\top\Sigma\nabla\psi_i(x^{(k)})$, while
$m$-strong convexity gives the Polyak--{\L}ojasiewicz bound
$\psi_i(x)-\psi_i(g_i)\leq\frac{1}{2m}\nabla\psi_i(x)^\top\Sigma
\nabla\psi_i(x)$; combining, the optimality gap contracts by the factor
$1-m/M\leq2/5$ per step.  The initial gap is bounded by smoothness and
\eqref{eq:comp_per_i_displacement}:
$\psi_i(X_i)-\psi_i(g_i)\leq\frac M2\norm{X_i-g_i}_\Sigma^2
\leq\frac{7}{10}\sqrt{U_n/\delta}\cdot\frac{25}{9}
\norm{A(X_i)}_2^2\frac{\delta}{\ell_n}
\leq\frac{2\mathsf s_n^2\sqrt{U_n\delta}}{\ell_n}$, using (B2) and (B4).

\textit{(ii)}  Strong convexity gives
$\frac m2\norm{\tilde g_i-g_i}_\Sigma^2\leq\eps_x$, so
$\norm{\tilde g_i-g_i}_\Sigma\leq\sqrt{2\eps_x/m}
\leq\sqrt{\frac{10}{3}\eps_x\sqrt{\delta/\ell_n}}
\leq2\sqrt{\eps_x}(\delta/\ell_n)^{1/4}$ by (B2).  Under the standing bound
    $\eps_x\leq\mathsf s_n^2\sqrt{\delta/\ell_n}$, this is at most
$2\mathsf s_n\sqrt{\delta/\ell_n}$, so every point $y$ on the segment from
$g_i$ to $\tilde g_i$ satisfies, by
\eqref{eq:comp_per_i_displacement} and~\eqref{eq:comp_A_lipschitz},
\begingroup\compactdisplaytrue
\[
    \norm{A(y)}_2
    \leq\norm{A(X_i)}_2+K_{\mathrm{curv}}\pare{\tfrac53\norm{A(X_i)}_2
    \sqrt{\tfrac{\delta}{\ell_n}}+2\mathsf s_n\sqrt{\tfrac{\delta}{\ell_n}}}
    \leq\norm{A(X_i)}_2+\tfrac{11}{3}\,\mathsf s_n
    K_{\mathrm{curv}}\sqrt{\tfrac{\delta}{\ell_n}}
    \leq2\,\mathsf s_n,
\]
\endgroup
using $K_{\mathrm{curv}}\sqrt{\delta/\ell_n}\leq1/5$.  The fundamental theorem of calculus
along that segment then gives
$\norm{\mf h(\tilde g_i)-\mf h(g_i)}_2\leq2\mathsf s_n
\norm{\tilde g_i-g_i}_\Sigma\leq b(\eps_x)$.  The value statement is
immediate from $\psi_i(g_i)\leq\psi_i(\tilde g_i)\leq\psi_i(g_i)+\eps_x$
termwise, and the gradient statement follows by averaging
\eqref{eq:comp_inner_bias} over $i$ and comparing
with~\eqref{eq:comp_gradients}.
\end{proof}

\begin{lemma}[Inexact projected gradient ascent in $\xi$]
\label{lem:comp_pga}
In the setting of Proposition~4.2, fix
$\alpha\in[\underline\alpha_\delta,\overline\alpha_\delta]$ and a target
$\eps_G\in(0,\mathsf M_\delta+\mathsf L_\delta]$, and set
\[
    r\Let\frac{\eps_G}{\mathsf M_\delta+\mathsf L_\delta},
    \qquad
    T\Let\Big\lceil30\,\varkappa_n^{3/2}\ln\frac{2}{r}\Big\rceil,
    \qquad
    \eps_x\Let
    \min\cure{\eps_G,\ \mathsf s_n^2\sqrt{\tfrac{\delta}{\ell_n}},\
    \pare{\frac{\mu_\delta\,r}{16\,\mathsf s_n}}^{\!2}
    \sqrt{\tfrac{\ell_n}{\delta}}}.
\]
Run, from $\xi^{(0)}=\mf 0$, the iteration
\[
    \xi^{(t+1)}
    =\Pi_{\mathbb B_2}\!\pare{\xi^{(t)}
    +\frac{1}{\mathsf L_\delta}\cdot\frac1n\sum_{i=1}^n
    \mf h\pare{\tilde g_i^{(t)}}},
    \qquad t=0,\ldots,T-1,
\]
where each $\tilde g_i^{(t)}$ is an $\eps_x$-accurate inner solution at
$(\alpha,\xi^{(t)})$ in the sense of Lemma~\ref{lem:comp_inner}(ii).  Then
\[
    G_\delta(\alpha)-F_\delta\pare{\alpha,\xi^{(T)}}\leq\eps_G,
    \qquad\text{and}\qquad
    \abs{\tilde F_\delta\pare{\alpha,\xi^{(T)}}-G_\delta(\alpha)}
    \leq\eps_G,
\]
where $\tilde F_\delta(\alpha,\xi^{(T)})$ is the value computed from the
$\eps_x$-accurate inner solutions.
\end{lemma}

\begin{proof}
Write $\xi^\star\Let\xi^\star(\alpha)$, $\mu\Let\mu_\delta$,
$\mathsf L\Let\mathsf L_\delta$, and let
$e_t\Let\frac1n\sum_i\mf h(\tilde g_i^{(t)})-\D_\xi F_\delta(\alpha,\xi^{(t)})$,
so $\norm{e_t}_2\leq b(\eps_x)\leq\mu r/4$ by
Lemma~\ref{lem:comp_inner}(ii) and the choice of $\eps_x$.  Since
$\xi^\star$ maximizes the concave $F_\delta(\alpha,\cdot)$ over
$\mathbb B_2$, it is a fixed point:
$\xi^\star=\Pi_{\mathbb B_2}(\xi^\star
+\frac{1}{\mathsf L}\D_\xi F_\delta(\alpha,\xi^\star))$.  By nonexpansiveness
of $\Pi_{\mathbb B_2}$,
\[
    \norm{\xi^{(t+1)}-\xi^\star}_2
    \leq
    \norm{\pare{\xi^{(t)}+\tfrac1{\mathsf L}\D_\xi F_\delta(\alpha,\xi^{(t)})}
    -\pare{\xi^\star+\tfrac1{\mathsf L}\D_\xi F_\delta(\alpha,\xi^\star)}}_2
    +\frac{\norm{e_t}_2}{\mathsf L}.
\]
By Lemma~\ref{lem:joint_concavity}(ii,iv), $F_\delta(\alpha,\cdot)$ is
$C^2$ on a neighborhood of the segment $[\xi^{(t)},\xi^\star]
\subseteq\mathbb B_2$ with
$\mu I\preccurlyeq-\D_\xi^2F_\delta\preccurlyeq\mathsf LI$ there, so with
$\bar H\Let-\int_0^1\D^2_\xi F_\delta(\alpha,\xi^\star
+t(\xi^{(t)}-\xi^\star))\diff t$ we get
$\D_\xi F_\delta(\alpha,\xi^{(t)})-\D_\xi F_\delta(\alpha,\xi^\star)
=-\bar H(\xi^{(t)}-\xi^\star)$ and
\[
    \norm{\pare{I-\tfrac1{\mathsf L}\bar H}\pare{\xi^{(t)}-\xi^\star}}_2
    \leq\pare{1-\frac{\mu}{\mathsf L}}\norm{\xi^{(t)}-\xi^\star}_2 .
\]
Setting $q\Let1-\mu/\mathsf L=1-\frac{1}{30\varkappa_n^{3/2}}$ and iterating,
\[
    \norm{\xi^{(T)}-\xi^\star}_2
    \leq q^T\norm{\xi^{(0)}-\xi^\star}_2
    +\frac{b(\eps_x)}{\mathsf L}\sum_{j\geq0}q^j
    \leq q^T+\frac{b(\eps_x)}{\mu}
    \leq\frac r2+\frac r4<r,
\]
using $\norm{\xi^{(0)}-\xi^\star}_2\leq1$ and
$q^T\leq e^{-T\mu/\mathsf L}\leq r/2$ by the choice of $T$.  By
$\mathsf L$-smoothness and $\norm{\D_\xi F_\delta(\alpha,\xi^\star)}_2
\leq\mathsf M_\delta$ (Lemma~\ref{lem:joint_concavity}(iv)),
\[
    G_\delta(\alpha)-F_\delta(\alpha,\xi^{(T)})
    \leq\mathsf M_\delta\norm{\xi^{(T)}-\xi^\star}_2
    +\frac{\mathsf L}{2}\norm{\xi^{(T)}-\xi^\star}_2^2
    \leq\pare{\mathsf M_\delta+\mathsf L_\delta}\,r=\eps_G,
\]
using $r\leq1$.  Finally, Lemma~\ref{lem:comp_inner}(ii) gives
$F_\delta(\alpha,\xi^{(T)})\leq\tilde F_\delta(\alpha,\xi^{(T)})\leq
F_\delta(\alpha,\xi^{(T)})+\eps_x\leq G_\delta(\alpha)+\eps_G$, and
$\tilde F_\delta(\alpha,\xi^{(T)})\geq F_\delta(\alpha,\xi^{(T)})\geq
G_\delta(\alpha)-\eps_G$.
\end{proof}

\begin{lemma}[Trisection with inexact values]
\label{lem:comp_trisection}
Let $f$ be concave on $[a,b]$ and $\Lambda$-Lipschitz, and suppose an oracle
returns $\tilde f$ with $\abs{\tilde f-f}\leq\eps_G$ at any query point.
Run $k$ rounds of trisection: given the current bracket, query the two
interior points $m_1,m_2$ trisecting it, discard the right third if
$\tilde f(m_1)\geq\tilde f(m_2)$ and the left third otherwise.  Let
$\hat\alpha$ be the queried point with the largest $\tilde f$-value over all
rounds.  Then
\begingroup\compactdisplaytrue
\begin{align}
    f(\hat\alpha)
    \geq\max_{[a,b]}f-(2k+2)\,\eps_G-\Lambda\,(b-a)\pare{\frac23}^{k},
    \qquad
    \abs{\tilde f(\hat\alpha)-\max_{[a,b]}f}
    \leq(2k+3)\,\eps_G+\Lambda(b-a)\pare{\frac23}^{k}.
    \label{eq:comp_trisection}
\end{align}
\endgroup
\end{lemma}

\begin{proof}
We claim each discard lowers the bracket maximum by at most $2\eps_G$.
Consider a round with bracket $[u,v]$, $w\Let v-u$, interior points
$m_1=u+w/3<m_2=v-w/3$, and suppose $\tilde f(m_1)\geq\tilde f(m_2)$, so
$(m_2,v]$ is discarded and $f(m_1)\geq f(m_2)-2\eps_G$.  Let
$x^\star\in\argmax_{[u,v]}f$.  If $x^\star\leq m_2$ nothing is lost.
Otherwise $x^\star\in(m_2,v]$ and
$m_2=\lambda m_1+(1-\lambda)x^\star$ with
$\lambda=\frac{x^\star-m_2}{x^\star-m_1}\leq\frac12$, because
$x^\star-m_2\leq w/3$ and $x^\star-m_1\geq w/3$; concavity gives
$f(m_2)\geq\lambda f(m_1)+(1-\lambda)f(x^\star)
\geq\lambda\pare{f(m_2)-2\eps_G}+(1-\lambda)f(x^\star)$, hence
$f(m_2)\geq f(x^\star)-2\eps_G\frac{\lambda}{1-\lambda}
\geq f(x^\star)-2\eps_G$, and $m_2$ remains in the new bracket.  The
discarded-left case is symmetric.  By induction, after $k$ rounds the
bracket $B_k$ has width $(b-a)(2/3)^k$ and
$\max_{B_k}f\geq\max_{[a,b]}f-2k\eps_G$.  Both round-$k$ query points lie in
$B_k$, so one of them, say $m$, satisfies
$f(m)\geq\max_{B_k}f-\Lambda(b-a)(2/3)^k$ by the Lipschitz property.  Then
$f(\hat\alpha)\geq\tilde f(\hat\alpha)-\eps_G\geq\tilde f(m)-\eps_G
\geq f(m)-2\eps_G$, which proves the first claim; the second follows from
$\abs{\tilde f(\hat\alpha)-f(\hat\alpha)}\leq\eps_G$ and
$f(\hat\alpha)\leq\max_{[a,b]}f$.
\end{proof}

The certified nested scheme is summarized in
Algorithm~\ref{alg:wp_nested}; its parameter choices implement
Lemmas~\ref{lem:comp_inner}--\ref{lem:comp_trisection}, and its guarantees
are those of Theorem~\MainComputableTestNumber.

\begin{algorithm}[!t]
    \caption{Certified nested scheme for the computable WP test
    $\hat{\mathsf T}_n$}
    \label{alg:wp_nested}
    \begin{algorithmic}[1]
        \REQUIRE data $(X_i)_{i=1}^n$, level $\varrho$, accuracy $\eps_{\mathrm{alg},n}$,
        curvature bound $K_{\mathrm{curv}}$ from~(21).
        \STATE Compute $\mc V_n$, $\ell_n$, $U_n$, $\varkappa_n$,
        $\hat z\Let\hat z_{1-\varrho}$ from~(9),
        $\delta\Let\hat z/n$, and the constants of~(24)
        and~\eqref{eq:comp_constants}.
        \IF{the observable condition $\mathfrak E_n$
        of~(30) fails}
            \STATE \textbf{return} ``fail to reject $\mc H_0$''.
        \ENDIF
        \STATE Set
        $k^\star\Let\max\cure{1,\big\lceil\log_{3/2}
        \pare{6\Lambda_\delta\sqrt{U_n/\delta}\big/\eps_{\mathrm{alg},n}}\big\rceil}$,
        $\eps_G\Let\eps_{\mathrm{alg},n}/(3(k^\star+2))$,
        $r\Let\eps_G/(\mathsf M_\delta+\mathsf L_\delta)$,
        $T\Let\big\lceil30\varkappa_n^{3/2}\ln(2/r)\big\rceil$,
        $\eps_x\Let\min\cure{\eps_G,\,\mathsf s_n^2\sqrt{\delta/\ell_n},\,
        \pare{\mu_\delta r/(16\mathsf s_n)}^{2}\sqrt{\ell_n/\delta}}$, and
        $k_{\mathrm{in}}\Let\big\lceil\log_{5/2}
        \pare{2\mathsf s_n^2\sqrt{U_n\delta}/(\ell_n\eps_x)}\big\rceil$.
        \STATE $[a,b]\leftarrow[\underline\alpha_\delta,
        \overline\alpha_\delta]$;\quad $\hat\Phi_n\leftarrow-\infty$.
        \FOR{$k=1,\ldots,k^\star$}
            \STATE $m_1\leftarrow a+(b-a)/3$;\quad
            $m_2\leftarrow b-(b-a)/3$.
            \FOR{$\alpha\in\{m_1,m_2\}$}
                \STATE $\xi\leftarrow\mf 0$.
                \FOR{$t=1,\ldots,T$}
                    \STATE For each $i\in[n]$ in parallel: starting from
                    $x\leftarrow X_i$, run $k_{\mathrm{in}}$ steps
                    $x\leftarrow x-\frac{1}{2\alpha+K_{\mathrm{curv}}}
                    \Sigma\nabla\psi_i(x)$ for
                    $\psi_i(x)=\xi^\top\mf h(x)
                    +\alpha\norm{x-X_i}_\Sigma^2$, and set
                    $\tilde g_i\leftarrow x$.
                    \STATE $\xi\leftarrow\Pi_{\mathbb B_2}
                    \pare{\xi+\frac{1}{\mathsf L_\delta}\cdot
                    \frac1n\sum_{i=1}^n\mf h(\tilde g_i)}$.
                \ENDFOR
                \STATE Recompute $(\tilde g_i)_{i\in[n]}$ at the current
                $(\alpha,\xi)$ as above, and set
                $\tilde G(\alpha)\leftarrow-\alpha\delta
                +\frac1n\sum_{i=1}^n\psi_i(\tilde g_i)$.
                \STATE $\hat\Phi_n\leftarrow
                \max\cure{\hat\Phi_n,\tilde G(\alpha)}$.
            \ENDFOR
            \IF{$\tilde G(m_1)\geq\tilde G(m_2)$}
                \STATE $b\leftarrow m_2$
            \ELSE
                \STATE $a\leftarrow m_1$
            \ENDIF
        \ENDFOR
        \ENSURE reject $\mc H_0$ if and only if $\hat\Phi_n>\eps_{\mathrm{alg},n}$.
    \end{algorithmic}
\end{algorithm}

\subsubsection{Statistical inheritance lemmas}
\label{sec:pf_comp_inherit}

The next lemma isolates the two inputs from the coverage analysis of the
paper that the inheritance argument needs; its proof consists of precise
pointers into the corresponding existing proofs.

\begin{lemma}[Edgeworth-expansion inputs from the coverage proofs]
\label{lem:comp_edgeworth_input}
The following hold.
\begin{enumerate}[label=\textup{(\roman*)}]
    \item Under Assumptions~2.4--2.10, there
    exist an event $\mc G_n$ with $\PP^\star(\mc G_n^c)=O(1/n)$, a constant
    $C<\infty$, and a deterministic sequence $\delta_n=\widetilde{O}(1/n)$ such
    that, on $\mc G_n$,
    \begin{align}
        nR_n(\mf h)-\hat z_{1-\varrho}
        =\pare{1-\mc L_n^\circ}
        \pare{\norm{\gamma^\dagger}_2^2-z_{1-\varrho}+\eps^\dagger_n}
        +\eps_n^\Sigma,
        \label{eq:comp_coverage_representation}
    \end{align}
    where $\gamma^\dagger$ is the random vector defined
    in~\eqref{eq:defgammada_n},
    $\mc L_n^\circ\Let\mc L_{W,V}(\mc W_n-W,\mc V_n-V)$ satisfies
    $\abs{\mc L_n^\circ}\leq\half$ for all $n$ large,
    $\abs{\eps^\dagger_n}\leq C\log^4(n)/n$, and
    $\abs{\eps_n^\Sigma}\leq\delta_n$.  Moreover, there are distribution
    functions $\Psi_n^\dagger$ with
    $\sup_z\abs{\PP^\star(\norm{\gamma^\dagger}_2^2\leq z)
    -\Psi_n^\dagger(z)}=\widetilde{O}(1/n)$, whose densities are bounded on the
    interval $[z_{1-\varrho}/2,\,2z_{1-\varrho}]$ by a constant
    $C_\Psi'<\infty$ uniformly in $n$.
    \item \textup{($d_h=1$.)}  Under the assumptions of
    Theorem~\MainHigherExpansionNumber, there exist an event of probability
    $1-O(n^{-3/2})$, constants $0<c_1\leq c_2<\infty$, and a deterministic
    sequence $\hat\delta_n=\widetilde{O}(n^{-3/2})$ such that, on that event,
    \begin{align}
        nR_n(h)-\hat z_{1-\varrho}
        =\frac{\mc W_n}{\mc V_n}
        \pare{\hat\gamma^2-\chi^2_{1;1-\varrho}}+\hat\eps_n,
        \qquad
        c_1\leq\frac{\mc W_n}{\mc V_n}\leq c_2,
        \qquad
        \abs{\hat\eps_n}\leq\hat\delta_n,
        \label{eq:comp_bartlett_representation}
    \end{align}
    where $\hat\gamma$ is the signed-root statistic constructed in the proof
    of Theorem~\MainBartlettCoverageNumber, which admits an Edgeworth expansion
    valid to order $n^{-3/2}$ whose distribution functions have densities
    uniformly bounded on compact subsets of $\RR$.
\end{enumerate}
\end{lemma}

\begin{proof}
\textit{(i)}  The representation is established in the proof of
Theorem~3.4: it is
display~\eqref{eq:nRnexpand3} there, with
$\eps_n^\Sigma\Let\bar\eps_n+\eps_n-\tilde\eps_n$ (in the notation of that
proof) controlled
by~\eqref{eq:3eps}, $\eps_n^\dagger$ controlled by~\eqref{eq:eps4},
$\gamma^\dagger$ defined in~\eqref{eq:defgammada_n}, and
$\abs{\mc L_n^\circ}\leq\half$ for $n$ large by~\eqref{eq:bndLWV}; the event
$\mc G_n$ is the intersection of the probability-$1-O(1/n)$ events on which
those displays hold.  The validity of the Edgeworth expansion of
$\gamma^\dagger$ up to order $n^{-1}$ is established in the closing
paragraph of that proof (via the moment expansion~\eqref{eq:gamma_moment}
and Lemma~\ref{lem:cramer}); the leading term of $\Psi_n^\dagger$ is the
distribution function $\Psi$ of Theorem~3.2, which is
continuously differentiable on $(0,\infty)$, and the $n^{-1/2}$ Edgeworth
correction is a polynomial multiple of a Gaussian density, so the densities
are uniformly bounded on the stated compact interval.

\textit{(ii)}  The representation is display~\eqref{eq:nRnh_higher} in the
proof of Theorem~\MainBartlettCoverageNumber
(Section~\ref{sec:pf_thm6_5}), together with the bound
$\abs{\hat\eps_n}\leq\hat\delta_n$ stated there, and
$\hat z_{1-\varrho}=(\mc W_n/\mc V_n)\chi^2_{1;1-\varrho}$ for $d_h=1$.  The
two-sided bound on $\mc W_n/\mc V_n$ at the stated probability level follows
from the concentration bounds recorded in the proof of
Proposition~\MainBartlettTwoNumber{} (under Assumption~\MainSmoothnessNumber, both
$\mc W_n$ and $\mc V_n$ deviate from the positive constants $W$ and $V$ by
$\widetilde{O}(\sqrt{\log(n)/n})$ with probability $1-O(n^{-3/2})$).  The
Edgeworth expansion of $\hat\gamma$ valid to order $n^{-3/2}$ is
established in the same proof via \cite[Theorem~2]{bhattacharya1978validity}
(with Theorem~\MainHigherExpansionNumber{} supplying the expansion of
$nR_n(h)$); its terms are a Gaussian distribution function plus polynomial
multiples of the Gaussian density, whose densities are bounded on compacts.
\end{proof}

\begin{lemma}[Crossing-band probability]
\label{lem:comp_band_prob}
The following hold.
\begin{enumerate}[label=\textup{(\roman*)}]
    \item Under Assumptions~2.4--2.10 (the
    hypotheses of Theorem~3.4),
    there is a constant $C_\Psi<\infty$ such that for every deterministic
    sequence $w_n\in[0,z_{1-\varrho}/4]$,
    \begin{align}
        \PP^\star\pare{0<nR_n(\mf h)-\hat z_{1-\varrho}\leq w_n}
        \leq
        C_\Psi\,w_n+\widetilde{O}\pare{\frac1n}.
        \label{eq:comp_band_prob}
    \end{align}
    \item \textup{($d_h=1$, Bartlett radius.)}  Under the assumptions of
    Theorem~\MainHigherExpansionNumber, with $C=O(1/n)$ as in
    Proposition~\MainBartlettOneNumber, there is a constant
    $C_\Psi^\flat<\infty$ such that for every deterministic sequence
    $w_n\in[0,\chi^2_{1;1-\varrho}/4]$,
    \begin{align}
        \PP^\star\pare{0<nR_n(h)-(1-C)\hat z_{1-\varrho}\leq w_n}
        \leq
        C_\Psi^\flat\,w_n+\widetilde{O}\pare{n^{-3/2}}.
        \label{eq:comp_band_prob_bartlett}
    \end{align}
\end{enumerate}
\end{lemma}

\begin{proof}
\textit{(i)}  Let $\mc G_n$, $\gamma^\dagger$, $\mc L_n^\circ$,
$\eps_n^\dagger$, $\eps_n^\Sigma$, $\delta_n$, and $\Psi_n^\dagger$ be as in
Lemma~\ref{lem:comp_edgeworth_input}(i).  On $\mc G_n$, the event
$\cure{0<nR_n(\mf h)-\hat z_{1-\varrho}\leq w_n}$
and~\eqref{eq:comp_coverage_representation} force
\[
    -\delta_n
    <\pare{1-\mc L_n^\circ}
    \pare{\norm{\gamma^\dagger}_2^2-z_{1-\varrho}+\eps^\dagger_n}
    \leq w_n+\delta_n,
\]
and since $\half\leq1-\mc L_n^\circ\leq\frac32$ this implies
$\norm{\gamma^\dagger}_2^2-z_{1-\varrho}\in
\big(-2\delta_n-\abs{\eps^\dagger_n},\,2(w_n+\delta_n)
+\abs{\eps^\dagger_n}\big]$, an interval whose endpoints lie in
$[z_{1-\varrho}/2,\,2z_{1-\varrho}]-z_{1-\varrho}$ for $n$ large because
$w_n\leq z_{1-\varrho}/4$.  Hence, by the Edgeworth approximation and the
density bound of Lemma~\ref{lem:comp_edgeworth_input}(i),
\begingroup\compactdisplaytrue
\[
    \PP^\star\pare{0<nR_n(\mf h)-\hat z_{1-\varrho}\leq w_n}
    \leq
    C_\Psi'\pare{2w_n+4\delta_n+2C\frac{\log^4 n}{n}}
    +\widetilde{O}\pare{\frac1n}+\PP^\star(\mc G_n^c)
    =C_\Psi\,w_n+\widetilde{O}\pare{\frac1n}.
\]
\endgroup

\textit{(ii)}  Work on the event of
Lemma~\ref{lem:comp_edgeworth_input}(ii).  Since
$\hat z_{1-\varrho}=(\mc W_n/\mc V_n)\chi^2_{1;1-\varrho}$,
representation~\eqref{eq:comp_bartlett_representation} gives
\[
    nR_n(h)-(1-C)\hat z_{1-\varrho}
    =\frac{\mc W_n}{\mc V_n}
    \pare{\hat\gamma^2-(1-C)\chi^2_{1;1-\varrho}}+\hat\eps_n .
\]
On the band, using $c_1\leq\mc W_n/\mc V_n\leq c_2$,
\[
    \hat\gamma^2-(1-C)\chi^2_{1;1-\varrho}
    \in\Big(-\frac{\hat\delta_n}{c_1},\ \frac{w_n+\hat\delta_n}{c_1}\Big].
\]
Since $C=O(1/n)$ and $w_n\leq\chi^2_{1;1-\varrho}/4$, the resulting
interval $(a_n,b_n]$ for $\hat\gamma^2$ satisfies
$a_n\geq\chi^2_{1;1-\varrho}/2$ for $n$ large and
$b_n-a_n\leq(w_n+2\hat\delta_n)/c_1$.  Writing
$\PP^\star(a_n<\hat\gamma^2\leq b_n)
=\PP^\star(\sqrt{a_n}<\hat\gamma\leq\sqrt{b_n})
+\PP^\star(-\sqrt{b_n}\leq\hat\gamma<-\sqrt{a_n})$ and applying the
Edgeworth expansion of $\hat\gamma$ (valid to order $n^{-3/2}$, with
densities bounded on compacts) on each interval of length
$\sqrt{b_n}-\sqrt{a_n}\leq(b_n-a_n)/(2\sqrt{a_n})$ yields, enlarging the
generic constant $C'$ if necessary and absorbing the probability of the
complement of the working event into the $\widetilde{O}(n^{-3/2})$ term,
\begingroup\compactdisplaytrue
\[
    \PP^\star\pare{0<nR_n(h)-(1-C)\hat z_{1-\varrho}\leq w_n}
    \leq
    C'\,\frac{w_n+2\hat\delta_n}{c_1\sqrt{2\chi^2_{1;1-\varrho}}}
    +\widetilde{O}\pare{n^{-3/2}}
    =C_\Psi^\flat\,w_n+\widetilde{O}\pare{n^{-3/2}} .
\]
\endgroup
\end{proof}

\subsubsection{Proof of Theorem~\MainComputableTestNumber}
\label{sec:pf_wp_computable_test}

\begin{proof}
Part (i) and the second inclusion in part (ii) are deterministic statements on
$\mathfrak E_n$.  The first inclusion in part (ii) holds globally because the
computable test fails to reject off $\mathfrak E_n$.  Part (iii) combines a
deterministic work bound with a high-probability event controlling the empirical
condition number and data scales.

The total work count used in part (iii), measured in first-order oracle calls for $\mf h$, is
\begin{equation}
    \underbrace{2k^\star}_{\substack{\text{trisection}\\ \text{queries}}}
    \times
    \underbrace{(T+1)}_{\substack{\text{PGA steps and value}\\ \text{evaluation per query}}}
    \times
    \underbrace{n}_{\substack{\text{inner problems}\\ \text{per step}}}
    \times
    \underbrace{(k_{\mathrm{in}}+1)}_{\substack{\text{GD steps and output evaluation}\\
    \text{per inner problem}}}.
    \label{eq:comp_ledger}
\end{equation}

\textit{(i)}  On $\mathfrak E_n$ we have $\delta\leq\delta_0$, so
Lemma~\ref{lem:joint_concavity}(v) makes $G_\delta$ concave and
$\Lambda_\delta$-Lipschitz on
$[\underline\alpha_\delta,\overline\alpha_\delta]$, an interval of length at
most $\sqrt{U_n/\delta}$ by (B1); each query is $\eps_G$-accurate by
Lemma~\ref{lem:comp_pga}, whose admissibility condition
$\eps_G\leq\mathsf M_\delta+\mathsf L_\delta$ holds on $\mathfrak E_n$
since $\eps_G\leq\eps_{\mathrm{alg},n}\leq\tfrac18\sqrt{\ell_n\delta}
\leq\mathsf L_\delta$; and
$\max_{[\underline\alpha_\delta,\overline\alpha_\delta]}G_\delta
=\Phi^{\mathrm{loc}}_n(\delta)$.
Lemma~\ref{lem:comp_trisection} with $k=k^\star$ then gives
\[
\begin{aligned}
    \Big|\hat\Phi_n-\Phi^{\mathrm{loc}}_n(\delta)\Big|
    &\leq(2k^\star+3)\,\eps_G
    +\Lambda_\delta\sqrt{\frac{U_n}{\delta}}\pare{\frac23}^{k^\star}\\
    &\leq\frac{2k^\star+3}{3(k^\star+2)}\,\eps_{\mathrm{alg},n}
    +\frac{\eps_{\mathrm{alg},n}}{6}\\
    &\leq\frac23\eps_{\mathrm{alg},n}+\frac16\eps_{\mathrm{alg},n}
    \leq\eps_{\mathrm{alg},n},
\end{aligned}
\]
where $(2/3)^{k^\star}\leq\eps_{\mathrm{alg},n}/(6\Lambda_\delta\sqrt{U_n/\delta})$ by the
choice of $k^\star$.

\textit{(ii)}  First inclusion: off $\mathfrak E_n$ the computable test fails to
reject.  On $\mathfrak E_n$, if $\hat\Phi_n>\eps_{\mathrm{alg},n}$, then by (i)
$\Phi^{\mathrm{loc}}_n(\delta)>0$, so
Proposition~4.2(i)
(display~(27)) gives $R_n(\mf h)>\delta$, i.e.
$nR_n(\mf h)>\hat z_{1-\varrho}$.  Second inclusion: on $\mathfrak E_n$, suppose
$nR_n(\mf h)>\hat z_{1-\varrho}+w_n$, i.e.
$R_n(\mf h)-\delta>w_n/n=5\eps_{\mathrm{alg},n}\sqrt{2\delta/\ell_n}$.  By the
contrapositive of Corollary~\ref{cor:comp_decision_band} (applicable on
$\mathfrak E_n$), $\Phi_n(\delta)>2\eps_{\mathrm{alg},n}$; by
Proposition~4.2(i),
$\Phi^{\mathrm{loc}}_n(\delta)=\Phi_n(\delta)$, and by (i),
$\hat\Phi_n\geq\Phi_n(\delta)-\eps_{\mathrm{alg},n}>\eps_{\mathrm{alg},n}$: the test rejects.

\textit{(iii)}  The count~\eqref{eq:comp_ledger} is immediate from the
construction: $k^\star$ rounds of two queries each; each query is one run of
Lemma~\ref{lem:comp_pga} with $T$ gradient steps plus one terminal value
evaluation; each step solves $n$ inner problems by
Lemma~\ref{lem:comp_inner}(i), each to accuracy $\eps_x$ in
$k_{\mathrm{in}}$ preconditioned gradient steps, each step one oracle call.
The expressions for $k^\star,T,k_{\mathrm{in}}$ restate the choices in
Lemmas~\ref{lem:comp_pga} and~\ref{lem:comp_inner} (the bound on
$k_{\mathrm{in}}$ uses the initial-gap bound of
Lemma~\ref{lem:comp_inner}(i) and $(1-m/M)\leq2/5$).  First impose the
additional bounds $\mathsf s_n\vee\bar{\mathsf m}_n\leq n^{1/4}$,
$\varkappa_n\leq\bar\varkappa$, and $\hat z_{1-\varrho}\leq n$.  On
$\mathfrak E_n$, every quantity inside the three logarithms is then bounded
above and below by fixed
powers of $n$ times fixed constants: $U_n\leq\mathsf s_n^2\leq n^{1/2}$
and $\delta=\hat z_{1-\varrho}/n\leq1$ give
$\Lambda_\delta\sqrt{U_n/\delta}=3\varkappa_n\sqrt{U_n\delta}
\leq3\bar\varkappa\,n^{1/4}$ and
$\mathsf M_\delta+\mathsf L_\delta\leq Cd_h\,n^{1/2}$; the third condition
in~(30) and $\eps_{\mathrm{alg},n}^{-1}=O(n^c)$ give
$\sqrt{\ell_n\delta}\geq8\eps_{\mathrm{alg},n}\geq c_\eps n^{-c}$,
whence $\mu_\delta^{-1}=8\ell_n^{-1}\sqrt{U_n/\delta}
=8\sqrt{\varkappa_n}\big/\sqrt{\ell_n\delta}
\leq C_\eps\sqrt{\bar\varkappa}\,n^{c}$; and $\eps_x^{-1}$ is polynomial in
$n,\eps_{\mathrm{alg},n}^{-1},\mathsf s_n,\mu_\delta^{-1}$.  Hence
$k^\star,\,T/\varkappa_n^{3/2},\,k_{\mathrm{in}}=O((1+c)\log n)$ and the
product of the three factors is $O((1+c)^3\bar\varkappa^{3/2}\,n\log^3 n)$.

It remains to show that these empirical bounds hold simultaneously with high
probability.  Let $\sigma_{\min}\Let\lambda_{\min}(V)>0$ and define
\begingroup\compactdisplaytrue
\[
    \mc A_n\Let
    \cure{\ell_n\geq\tfrac{\sigma_{\min}}2}
    \cap\cure{\Tr(\mc V_n)\leq\Tr(V)+1}
    \cap\cure{\Tr(\mc W_n)\leq\Tr(W)+1}
    \cap\cure{\lambda_{\min}(\mc W_n)\geq\tfrac{\lambda_{\min}(W)}2}.
\]
\endgroup
The first and last events have probability $1-O(1/n)$ by
Lemma~\ref{lem:tail_eigen} applied to
$A=\D\mf h(X)\Sigma\D\mf h(X)^\top$ (as in Proposition~\ref{prop:tailbnd},
part (A1)) and to $A=\mf h(X)^{\otimes2}$ (both admissible since
Assumption~2.9 gives finite fourth moments of
$\norm{\D\mf h(X)}_2$ and $\norm{\mf h(X)}_2$).  The two trace events have
probability $1-O(1/n)$ by Chebyshev's inequality applied to the scalar
averages $\Tr(\mc V_n)$ and $\Tr(\mc W_n)$, whose variances are $O(1/n)$
under Assumption~2.9.  On $\mc A_n$,
$U_n\leq\Tr(V)+1$ and $\varkappa_n\leq2(\Tr(V)+1)/\sigma_{\min}
\Let\bar\varkappa$, so
$\delta_0\geq\underline\delta_0\Let\frac{\sigma_{\min}/2}
{K_{\mathrm{curv}}^2(3+2\sqrt{\bar\varkappa})^2}>0$.  Moreover, since
$\hat z_{1-\varrho}$ is the $(1-\varrho)$-quantile of
$v^\top\mc W_n^{\half}\mc V_n^{-1}\mc W_n^{\half}v$ with $v$ standard
normal,
\[
    \underline z
    \Let\frac{\lambda_{\min}(W)\,\chi^2_{d_h;1-\varrho}}{2(\Tr(V)+1)}
    \leq\hat z_{1-\varrho}\leq
    \frac{2(\Tr(W)+1)\,\chi^2_{d_h;1-\varrho}}{\sigma_{\min}}
    \Let\bar z
    \qquad\text{on }\mc A_n .
\]
Define
\[
    \mathfrak C_n\Let\mc A_n\cap\cure{\mathsf s_n\leq n^{1/4}}.
\]
Because $\mathsf s_n=\max_{i\in[n]}\norm{\D\mf h(X_i)\Sigma^{1/2}}_2$,
Assumption~2.9 and a union--Markov bound give
\[
    \PP^\star\pare{\mathsf s_n>n^{1/4}}
    \leq n\,\PP^\star\pare{\norm{\D\mf h(X)\Sigma^{1/2}}_2>n^{1/4}}
    \leq\frac{\EE\norm{\D\mf h(X)\Sigma^{1/2}}_2^8}{n}
    =O(n^{-1}).
\]
Consequently, $\PP^\star(\mathfrak C_n^c)=O(n^{-1})$.  On $\mathfrak C_n$,
$\bar{\mathsf m}_n^2=\Tr(\mc W_n)\leq\Tr(W)+1$,
$\varkappa_n\leq\bar\varkappa$, and $\hat z_{1-\varrho}\leq\bar z$.
Thus, for all sufficiently large $n$, the empirical bounds used above hold on
$\mathfrak C_n$, proving the asserted work bound on
$\mathfrak E_n\cap\mathfrak C_n$.

For the accuracy in part~(iv), take $n\geq N_0$, where $N_0$ depends only on
\[
(\sigma_{\min},\Tr V,\Tr W,\lambda_{\min}(W),K_{\mathrm{curv}},\varrho).
\]
On $\mc A_n$,
$\delta=\hat z_{1-\varrho}/n\leq\bar z/n\leq\underline\delta_0/2
\leq\delta_0/2$ and
$\eps_{\mathrm{alg},n}=n^{-3/2}\log^{-2}n
\leq\frac18\sqrt{\ell_n\hat z_{1-\varrho}/n}$
(the right side is at least
$\frac18\sqrt{\sigma_{\min}\underline z/2}\,n^{-1/2}$).  Hence
$\mc A_n\subseteq\mathfrak E_n$, and therefore
$\PP^\star(\mathfrak E_n\cap\mathfrak C_n)
=\PP^\star(\mathfrak C_n)=1-O(n^{-1})$.

\textit{(iv)}  Now take $\eps_{\mathrm{alg},n}=1/(n^{3/2}\log^2n)$.  Write
$R\Let\{nR_n(\mf h)>\hat z_{1-\varrho}\}$ and
$\hat R\Let\{\hat{\mathsf T}_n\text{ rejects}\}$.  By part (ii),
$\hat R\subseteq R$ and
$R\setminus\hat R\subseteq
\{0<nR_n(\mf h)-\hat z_{1-\varrho}\leq w_n\}\cup\mathfrak E_n^c$ with
$w_n$ defined as in Theorem~\MainComputableTestNumber.  The construction in
part~(iii) gives $\mc A_n\subseteq\mathfrak E_n$ and
$\PP^\star(\mc A_n^c)=O(n^{-1})$.  Furthermore, on $\mc A_n$,
$w_n\leq5\eps_{\mathrm{alg},n}\sqrt{4n\bar z/\sigma_{\min}}
=\bar w_n\Let10\sqrt{\bar z/\sigma_{\min}}\;n^{-1}\log^{-2}n$, a
deterministic sequence with $\bar w_n=o(1/n)$; in
particular $\bar w_n\leq z_{1-\varrho}/4$ eventually.
Lemma~\ref{lem:comp_band_prob}(i) then gives
\[
    \PP^\star\pare{R\setminus\hat R}
    \leq\PP^\star\pare{0<nR_n(\mf h)-\hat z_{1-\varrho}\leq\bar w_n}
    +\PP^\star(\mc A_n^c)
    \leq C_\Psi\bar w_n+\widetilde{O}\pare{\tfrac1n}
    =\widetilde{O}\pare{\tfrac1n},
\]
where the residual event is $\mc A_n^c$ (not merely $\mathfrak E_n^c$)
because the bound $w_n\leq\bar w_n$ was established on $\mc A_n$; recall
$\mc A_n\subseteq\mathfrak E_n$ and $\PP^\star(\mc A_n^c)=O(1/n)$.
Combining with $\hat R\subseteq R$ and
Theorem~3.4 proves the size and
coverage claims; the confidence-region statement is the same computation
applied to each tested parameter value.

For the local-power claim, the argument is identical under
$\PP_n^\star$: the inclusion structure of part (ii) is deterministic; the
event bounds defining $\mc A_n$ hold with probability $1-O(1/n)$ uniformly
in $n$ under Assumption~\ref{a:alter_edgeworthII}(i)--(ii) (uniform fourth
moments and uniformly nondegenerate $V_n,W_n$); and the band probability
under $\PP_n^\star$ is $O(\bar w_n)+\widetilde{O}(1/n)$ because the proof of
Theorem~\MainPowerExpansionNumber{} establishes the same
representation~\eqref{eq:nRn-hatz} with $\gamma^\ddagger$ in place of
$\gamma^\dagger$ and a valid Edgeworth expansion for it, whose leading
distribution function has a density bounded near $z_{1-\varrho}$ uniformly
in $n$ (its covariance matrices $V_n^{-\half}W_nV_n^{-\half}$ are uniformly
nondegenerate).
\end{proof}

\subsubsection{The Bartlett-corrected variant}
\label{sec:pf_comp_bartlett}

\begin{corollary}[Bartlett-corrected computable test, $d_h=1$]
\label{cor:comp_coverage}
Under the assumptions of Theorem~\MainHigherExpansionNumber{} together with
Assumption~2.5 and~(21), the computable test of
Theorem~\MainComputableTestNumber, run with the threshold radius
$\delta=(1-C)\hat z_{1-\varrho}/n$ from Proposition~\MainBartlettOneNumber
(so that the deterministic guarantees of
Theorem~\MainComputableTestNumber(i)--(ii) apply with
$(1-C)\hat z_{1-\varrho}$ in place of $\hat z_{1-\varrho}$) and with
accuracy $\eps_{\mathrm{alg},n}=1/(n^{2}\log^2n)$, satisfies
\[
    \PP^\star\pare{\hat{\mathsf T}_n\text{ fails to reject}}
    =1-\varrho+\widetilde{O}\pare{n^{-3/2}} .
\]
Moreover, the corrected algorithm uses $O(n\log^3n)$ first-order oracle calls
with probability $1-O(n^{-3/2})$.
\end{corollary}

\begin{proof}
For $d_h=1$ the constants entering $\mathfrak E_n$ are scalar sample means,
and under Assumption~\MainSmoothnessNumber{} all polynomial moments of
$\abs{h(X)}$ and $\norm{\D h(X)}_2$ are finite; a fourth-moment Chebyshev
bound for each scalar average in the event $\mc A_n$ constructed in the proof of
Theorem~\MainComputableTestNumber(iii) (variance version, with the
elementary estimate $\EE[(\bar Y_n-\EE Y)^4]=O(n^{-2})$ for i.i.d.\ $Y_i$
with $\EE Y^4<\infty$) upgrades $\PP^\star(\mc A_n^c)$ to $O(n^{-2})$, and
the maxima bounds $\mathsf s_n\vee\bar{\mathsf m}_n\leq n^{1/4}$ hold with
probability $1-O(n^{-3/2})$ by Markov's inequality with tenth moments.
With $\eps_{\mathrm{alg},n}=n^{-2}\log^{-2}n$ the band width becomes
$\bar w_n=O(n^{-3/2}\log^{-2}n)$.
Lemma~\ref{lem:comp_band_prob}(ii) bounds the probability of the crossing
band at the radius $(1-C)\hat z_{1-\varrho}/n$ by
$C_\Psi^\flat\bar w_n+\widetilde{O}(n^{-3/2})=\widetilde{O}(n^{-3/2})$.  The chain of
inclusions is then identical to the proof of
Theorem~\MainComputableTestNumber(iv) (with the residual event
$\mc A_n^c$ there, as before), with every $\widetilde{O}(1/n)$
replaced by $\widetilde{O}(n^{-3/2})$, and Proposition~\MainBartlettOneNumber
supplies the oracle accuracy $1-\varrho+\widetilde{O}(n^{-3/2})$.
Finally, since $C=O(n^{-1})$, the corrected radius obeys the same fixed
upper and lower quantile bounds used in the proof of
Theorem~\MainComputableTestNumber(iii) for all sufficiently large $n$.
Combining the $O(n^{-2})$ bound for $\mc A_n^c$ with the
$O(n^{-3/2})$ maxima bound above shows that the corresponding complexity
event and $\mathfrak E_n$ hold jointly with probability $1-O(n^{-3/2})$.
The deterministic work bound in that proof therefore gives the asserted
$O(n\log^3n)$ oracle count.
\end{proof}

\subsubsection{The stochastic-gradient variant}
\label{sec:pf_comp_sgd}

\begin{proposition}[Biased projected stochastic gradient ascent at fixed
$\alpha$]
\label{prop:comp_sgd}
In the setting of Proposition~4.2, fix
$\alpha\in[\underline\alpha_\delta,\overline\alpha_\delta]$ and
$\eps_x\in(0,\mathsf s_n^2\sqrt{\delta/\ell_n}]$.  For an integer $T\geq1$,
choose any $\xi^{(1)}\in\mathbb B_2$ and run, for $t=1,\ldots,T$, the iteration
\[
    \xi^{(t+1)}
    =\Pi_{\mathbb B_2}\pare{\xi^{(t)}+\eta_t\,
    \mf h\pare{\tilde g_{I_t}\pare{\alpha,\xi^{(t)}}}},
    \qquad
    \eta_t=\frac{2}{\mu_\delta(t+1)},
    \qquad
    I_t\sim\mathrm{Unif}[n]\text{ i.i.d.},
\]
where $\tilde g_i(\alpha,\xi)$ is the deterministic output of
$k_{\mathrm{in}}$ inner gradient steps as in Lemma~\ref{lem:comp_inner},
achieving accuracy $\eps_x$.  Then the conditional mean of the update
direction is $\D_\xi F_\delta(\alpha,\xi^{(t)})+e_t$ with
$\norm{e_t}_2\leq b(\eps_x)$ as in~\eqref{eq:comp_inner_bias}, its
conditional second moment is at most
$\bar\sigma_\delta^2\Let2\bar{\mathsf m}_n^2+128\,\mathsf s_n^4\delta/\ell_n$,
and the weighted average
$\bar\xi^{(T)}\Let\frac{2}{T(T+1)}\sum_{t=1}^{T}t\,\xi^{(t)}$ satisfies
\begin{align}
    \EE\brac{G_\delta(\alpha)-F_\delta\pare{\alpha,\bar\xi^{(T)}}}
    \leq
    \frac{2\bar\sigma_\delta^2}{\mu_\delta\,(T+1)}
    +2\,b(\eps_x).
    \label{eq:comp_sgd_rate}
\end{align}
Hence a bias of order $\eps$ is absorbed ($b(\eps_x)\leq\eps/4$ requires
only $\eps_x\leq\eps^2\sqrt{\ell_n/\delta}/(256\,\mathsf s_n^2)$, i.e.
$k_{\mathrm{in}}=O(\log(1/\eps))$ inner steps), and reaching
$\EE[\,\text{gap}\,]\leq\eps$ requires
\begin{align}
    T
    =\Big\lceil\frac{4\bar\sigma_\delta^2}{\mu_\delta\,\eps}\Big\rceil
    =O\pare{\frac{\bar\sigma_\delta^2}{\ell_n}
    \sqrt{\frac{U_n}{\delta}}\cdot\frac1\eps}
    \label{eq:comp_sgd_iters}
\end{align}
stochastic iterations, i.e.\ $O(T\log(1/\eps_x))$ oracle calls in total.
With $\delta=\hat z_{1-\varrho}/n$ this is
$\widetilde{O}(\sqrt n/\eps)$: at $\eps=1/(n\log^2n)$ the stochastic scheme costs
$\widetilde{O}(n^{3/2})$, and at the coverage-preserving accuracy
$\eps=1/(n^{3/2}\log^2n)$ it costs $\widetilde{O}(n^{2})$.  The stochastic
variant is therefore dominated by the deterministic scheme of
Theorem~\MainComputableTestNumber{} in this problem, where the relevant
accuracy scales as a negative power of $n$.
\end{proposition}

\begin{proof}
Conditionally on $\mc F_t\Let\sigma(\xi^{(1)},I_1,\ldots,I_{t-1})$, the
index $I_t$ is uniform, so
$\EE[\mf h(\tilde g_{I_t}(\alpha,\xi^{(t)}))\mid\mc F_t]
=\frac1n\sum_i\mf h(\tilde g_i(\alpha,\xi^{(t)}))
=\D_\xi F_\delta(\alpha,\xi^{(t)})+e_t$ with $\norm{e_t}_2\leq b(\eps_x)$
by Lemma~\ref{lem:comp_inner}(ii) (the inner solver output is a
deterministic function of $(\alpha,\xi^{(t)},i)$, so no additional
randomness enters).  For the second moment, by
\eqref{eq:comp_h_displacement}, \eqref{eq:comp_inner_bias}, and
$\eps_x\leq\mathsf s_n^2\sqrt{\delta/\ell_n}$,
\[
    \norm{\mf h(\tilde g_i)}_2
    \leq\norm{\mf h(X_i)}_2
    +\frac{20}{9}\mathsf s_n^2\sqrt{\frac{\delta}{\ell_n}}
    +4\mathsf s_n\sqrt{\eps_x}\pare{\frac{\delta}{\ell_n}}^{1/4}
    \leq\norm{\mf h(X_i)}_2+8\,\mathsf s_n^2\sqrt{\frac{\delta}{\ell_n}},
\]
so $\frac1n\sum_i\norm{\mf h(\tilde g_i)}_2^2\leq
2\bar{\mathsf m}_n^2+128\,\mathsf s_n^4\delta/\ell_n=\bar\sigma_\delta^2$.

For~\eqref{eq:comp_sgd_rate}, write $\xi^\star=\xi^\star(\alpha)$,
$\mu=\mu_\delta$, $\mathsf g_t\Let\mf h(\tilde g_{I_t}(\alpha,\xi^{(t)}))$,
$d_t\Let\EE\norm{\xi^{(t)}-\xi^\star}_2^2$, and
$G_t\Let\EE\brac{G_\delta(\alpha)-F_\delta(\alpha,\xi^{(t)})}\geq0$.
Nonexpansiveness of the projection (with
$\xi^\star\in\mathbb B_2$) gives
$\norm{\xi^{(t+1)}-\xi^\star}_2^2
\leq\norm{\xi^{(t)}-\xi^\star}_2^2
+2\eta_t\langle\mathsf g_t,\,\xi^{(t)}-\xi^\star\rangle
+\eta_t^2\norm{\mathsf g_t}_2^2$.
Taking conditional expectations and using
$\mu$-strong concavity in the form
$\langle\D_\xi F_\delta(\alpha,\xi^{(t)}),\,\xi^{(t)}-\xi^\star\rangle
\leq-\pare{G_\delta(\alpha)-F_\delta(\alpha,\xi^{(t)})}
-\frac\mu2\norm{\xi^{(t)}-\xi^\star}_2^2$
together with
$\langle e_t,\xi^{(t)}-\xi^\star\rangle\leq2b(\eps_x)$ (diameter of
$\mathbb B_2$ is $2$),
\[
    d_{t+1}\leq(1-\mu\eta_t)\,d_t-2\eta_tG_t
    +2\eta_t\cdot2b(\eps_x)+\eta_t^2\bar\sigma_\delta^2 .
\]
With $\eta_t=\frac{2}{\mu(t+1)}$, so that
$\frac{1-\mu\eta_t}{2\eta_t}=\frac{\mu(t-1)}{4}$, this rearranges to
\begin{align*}
G_t\leq\frac{\mu}{4}\brac{(t-1)d_t-(t+1)d_{t+1}}
+\frac{\bar\sigma_\delta^2}{\mu(t+1)}+2b(\eps_x).
\end{align*}
Multiplying by $t$ and summing over $t=1,\ldots,T$, the first sum
telescopes to $-\frac{\mu}{4}T(T+1)\,d_{T+1}\leq0$ (the term
$t(t+1)d_{t+1}$ produced at step $t$ cancels against the term
$(t+1)t\,d_{t+1}$ at step $t+1$), and
$\sum_{t\leq T}\frac{t}{t+1}\leq T$, whence
$\sum_{t=1}^Tt\,G_t\leq\frac{T\bar\sigma_\delta^2}{\mu}
+2b(\eps_x)\frac{T(T+1)}{2}$.  Dividing by $T(T+1)/2$ and using concavity
of $F_\delta(\alpha,\cdot)$ (Jensen, with the weights $2t/(T(T+1))$)
gives~\eqref{eq:comp_sgd_rate}.  The iteration
count~\eqref{eq:comp_sgd_iters} follows by making each right-hand term at
most $\eps/2$, and the scaling claims follow from
$\mu_\delta=\frac{\ell_n}8\sqrt{\delta/U_n}$ with $\delta=\hat z_{1-\varrho}/n$.
\end{proof}

\begin{remark}[Why not stochastic gradients]
\label{rem:comp_sgd}
A natural first idea is to replace the full gradient $\frac1n\sum_i\mf h(g_i(\alpha,\xi))$ by
a single-sample estimate and run projected stochastic gradient ascent in $\xi$.  The honest
rate of this variant (Proposition~\ref{prop:comp_sgd}) is $\widetilde{O}(\sqrt n/\eps)$ inner
solves to reach accuracy $\eps$, because the strong-concavity modulus on the box is
$\mu_\delta\asymp\sqrt\delta\asymp n^{-1/2}$: at $\eps=1/(n\log^2n)$ the cost is already
$\widetilde{O}(n^{3/2})$, and at the coverage-preserving accuracy $\eps=1/(n^{3/2}\log^2n)$ it is
$\widetilde{O}(n^{2})$.  Moreover, by the accuracy-to-decision-band exchange rate in Theorem~\MainComputableTestNumber, the
accuracy $1/(n\log^2n)$ would in any case be statistically insufficient: it leaves a
misclassification band of width $\tilde\Theta(n^{-1/2})$ on the $nR_n$ scale, degrading the
size accuracy from $\widetilde{O}(n^{-1})$ to $O(n^{-1/2})$.  The deterministic nested scheme,
whose loops all converge linearly thanks to Proposition~4.2, achieves
the stronger accuracy at total cost $\widetilde{O}(n)$, and is therefore the scheme we certify.
\end{remark}

\subsection{{Derivation of inequalities}}\label{sec:deri_ineq}
\begin{proof}[Derivation of \eqref{eq:1n-r}]
Put $t_n=a\sqrt{\log(n)}$ and $\lambda=\lambda_{\max}(V)$. Since
$v^\top V^{-1}v\geq\norm{v}_2^2/\lambda$, polar coordinates and the
standard Gaussian tail integral bound give
    \begin{align*}
        \int_{\norm{v}_2 \geq t_n}\norm{v}_2^{3s}\phi_V(v)\diff v
        &\leq C_0\int_{t_n}^{\infty}
        u^{3s+d_h-1}\exp\left(-\frac{u^2}{2\lambda}\right)\diff u\\
        &\leq C_s\{\log(n)\}^{c_s}
        \exp\left\{-\frac{a^2\log(n)}{2\lambda}\right\}
        =o(n^{-r}),
    \end{align*}
    where $C_0,C_s,c_s$ are constants independent of $n$, and the last
    equality follows from $a^2>2r\lambda$.
\end{proof}

\begin{proof}[Derivation of Inequality~\eqref{eq:DF}]
Recall that
\begingroup\compactdisplaytrue
\begin{align*}
    \D F_{n,i}(v) 
    =~& \frac{1}{2} \left(\int_{0}^{1} \zeta^{\top}\D\mf h\left(X_{i}+v\right) \Sigma D^{(2)}\left(X_{i}+vu \right) u + \zeta^{\top}\D\mf h\left(X_{i}+v u\right) \Sigma D^{(2)}\left(X_{i}+v \right) \diff u \right) \\
    &- \frac{1}{2} \zeta^{\top}\D \mf h \left(X_{i}+v\right)\Sigma D^{(2)}\left(X_{i}+v \right),
\end{align*}
\endgroup
where $D^{(2)}(v) \Let \sum_{\beta \in [d_h]} \zeta^{(\beta)} \D^2 h^{\beta}(v)$. Thus, we find \begin{align*}
    \D F_{n,i}(\mf 0) = \frac{1}{4} \zeta^\top \D \mf h(X_i) \Sigma D^{(2)}(X_i).
\end{align*}
By Assumption~2.5, $\norm{D^{(2)}(x)}_2 \leq \kappa_1(x)\norm{\zeta}_2$; thus
\begin{align*}
    \norm{\D F_{n,i}(\mf 0)}_2 = \norm{\frac{1}{4} \zeta^\top \D \mf h(X_i) \Sigma D^{(2)}(X_i)}_2 \leq \frac{\norm{\Sigma}_2}{4} \norm{\D \mf h(X_i)}_2\kappa_1(X_i)\norm{\zeta}_2^2.
\end{align*}
\end{proof}

\begin{proof}[Derivation of Inequality~\eqref{eq:DDF}]
Note that the three terms in $\D F_{n,i}(v)$ are similar, so we only need to bound 
\[
    \zeta^{\top}\D \mf h \left(X_{i}+v\right) \Sigma D^{(2)}\left(X_{i}+v \right) - \zeta^{\top}\D \mf h \left(X_{i}\right)\Sigma D^{(2)}\left(X_{i} \right).
\]
We write 
\begingroup\compactdisplaytrue
\begin{align*}
    & \zeta^{\top}\D \mf h \left(X_{i}+v\right) \Sigma D^{(2)}\left(X_{i}+v \right) -  \zeta^{\top}\D \mf h \left(X_{i}\right) \Sigma D^{(2)}\left(X_{i} \right)\\
    =~&  \zeta^{\top}\left(\D \mf h \left(X_{i}+v\right) - \D \mf h \left(X_{i}\right)\right) \Sigma \left(D^{(2)}\left(X_{i}+v \right) -D^{(2)}\left(X_{i} \right)\right) \\
    &+ \zeta^{\top} \D \mf h \left(X_{i}\right) \Sigma \left(D^{(2)}\left(X_{i}+v \right) -D^{(2)}\left(X_{i} \right)\right) + \zeta^{\top}\left(\D \mf h \left(X_{i}+v\right) - \D \mf h \left(X_{i}\right)\right)\Sigma D^{(2)}\left(X_{i} \right).
\end{align*}
\endgroup
By Assumptions~2.5 and~2.6, when $\norm{v}_2 \leq \hat \delta$, we have
\begin{align*}
    & \norm{\zeta^{\top}\left(\D \mf h \left(X_{i}+v\right) - \D \mf h \left(X_{i}\right)\right) \Sigma \left(D^{(2)}\left(X_{i}+v \right) -D^{(2)}\left(X_{i} \right)\right)}_2 \\
    \overset{\text{by A.}~2.6}{\leq} & \kappa_2 (X_i) \norm{\zeta}_2 \norm{v}_2 \norm{\zeta^{\top}\left(\D \mf h \left(X_{i}+v\right) - \D \mf h \left(X_{i}\right)\right) \Sigma }_2\\
    \overset{\text{by A.}~2.5}{\leq} & \norm{\Sigma}_2 \kappa_1(X_i)\kappa_2 (X_i) \norm{\zeta}_2^2 \norm{v}_2^2.
\end{align*}
and 
\begin{align*}
    \norm{\zeta^{\top} \D \mf h \left(X_{i}\right) \Sigma \left(D^{(2)}\left(X_{i}+v \right) -D^{(2)}\left(X_{i} \right)\right)}_2 \overset{\text{by A.}~2.6}{\leq} \norm{\Sigma}_2\norm{\D \mf h(X_i)}_2\kappa_2(X_i)\norm{\zeta}_2^2\norm{v}_2.
\end{align*}

By Assumption~2.5, we have
\begingroup\compactdisplaytrue
\begin{align*}
    \norm{\zeta^{\top}\left(\D \mf h \left(X_{i}+v\right) - \D \mf h \left(X_{i}\right)\right)\Sigma D^{(2)}\left(X_{i} \right)}_2
    \overset{\text{by A.}~2.5}{\leq} & \norm{\zeta^{\top}\left(\D \mf h \left(X_{i}+v\right) - \D \mf h \left(X_{i}\right)\right)\Sigma}_2 \kappa_1(X_i) \norm{\zeta}_2 \\
    \overset{\text{by A.}~2.5}{\leq} & \norm{\Sigma}_2 \kappa_1(X_i)^2 \norm{\zeta}_2^2 \norm{v}_2.
\end{align*}
\endgroup
Since 
\begin{align*}
     \norm{n^{-\half}\Delta_{n,i}(\zeta)}_2 \leq & n^{-\half} \norm{\D \mf h(X_i)}_2 \norm{\Sigma}_2 \norm{\zeta}_2\\
    \leq & \hat\delta\frac{\norm{\zeta}_2}{\bar r_n}
    \leq \hat\delta,
\end{align*}
where the first inequality is due to \eqref{eq:Deltaleq}, the second is due to \ref{cond:condA}(A3)(i), and the last uses $\zeta\in\mc Z_n$ and $r_n\leq\bar r_n$. Thus, for $s \in [0,1]$,
\begin{align*}
& \norm{\D F_{n,i}(n^{-\frac{1}{2}}\Delta_{n,i}(\zeta)s) - \D F_{n,i}(\mf 0)}_2 \notag\\
\leq &~ \frac{5\norm{\Sigma}_2\norm{\zeta}_2^2}{4} \left(
\kappa_1(X_i) \kappa_2(X_i) \norm{n^{-\frac{1}{2}}\Delta_{n,i}(\zeta)}_2^2 + \kappa_2(X_i)\norm{\D \mf h(X_i)}_2 \norm{n^{-\frac{1}{2}}\Delta_{n,i}(\zeta)}_2 \right.\notag\\ 
& \left.+ \kappa_1(X_i)^2 \norm{n^{-\frac{1}{2}}\Delta_{n,i}(\zeta)}_2\right).
\end{align*}
This completes the proof.
\end{proof}

\subsection{Selection of {\texorpdfstring{$N$}{TEXT}} in Theorem~3.1}\label{app:ndpair}
Let
\begin{align*}
    B_\Sigma&=\norm{\Sigma}_2^2
    \left(1+\EE[\norm{\D\mf h(X)}_2^2\kappa_1(X)]\right),
    &a_H&=2\sqrt{1\vee\lambda_{\max}(W)}.
\end{align*}
The integer $N$ in Theorem~3.1 may be chosen so that, for every $n\geq N$,
\begingroup\compactdisplaytrue
\begin{align*}
    &n \geq \left(\frac{2}{9} \vee \frac{8}{\lambda_{\min}(\EE[A])^2}\right)\left(\norm{\EE[(A - \EE[A])^2]}_2 + \norm{\EE[A]}_2\right)\log(d_h n),&\text{(Prop.~\ref{prop:tailbnd})}\\
    &\frac{32a_HB_\Sigma}{\sigma_{\min}^2}
    \sqrt{\frac{\log(n)}{n}}<1.
    &\text{(Props.~\ref{prop:bound_subgradient} and~\ref{prop:expanFn*})}
\end{align*}
\endgroup
Here $N\geq3$, $W=\Cov(\mf h(X),\mf h(X))$,
$A=\D\mf h(X)\Sigma\D\mf h(X)^\top$, and
$\sigma_{\min}=\lambda_{\min}(\EE[A])$. The pointwise bounds in
\ref{cond:condA}(A3) jointly have failure probability
$O\{\log(n)^4/n^3\}$ and therefore impose no additional deterministic
lower bound on $n$. In particular, the localization requirement is now of
the polynomial form $n/\log(n)>(32a_HB_\Sigma/\sigma_{\min}^2)^2$, rather
than an exponential requirement on $n$.
    
\subsection{Moment expansions in Equation~\eqref{eq:gamma_moment}}\label{sec:gamma_moment}
Write $\gamma=G_n+n^{-1/2}Q_n$, where, for every
$\eta\in\RR^{d_h}$,
\begin{align*}
    \brak{G_n,\eta}
    &\Let \brak{V^{\half},\tilde\xi_n\otimes\eta},\\
    \brak{Q_n,\eta}
    &\Let-\frac{\sqrt n}{2}
       \brak{\mc V_n-V,\tilde\xi_n\otimes(V^{-\half}\eta)}
       +\frac12\brak{\EE[\mc K_n],
       \tilde\xi_n^{\otimes2}\otimes(V^{-\half}\eta)}.
\end{align*}
Because $\EE[G_n]=0$ under the null, the deterministic first-moment
coefficient is
\begin{align*}
    \mu_{1,n}^j
    &=\EE[Q_{n,j}]\\
    &=-\frac12\brak{\EE[\sqrt n(\mc V_n-V)\tilde\xi_n],
        V^{-\half}e_j}
      +\frac12\brak{\EE[\mc K_n],
        \EE[\tilde\xi_n^{\otimes2}]\otimes(V^{-\half}e_j)}.
\end{align*}
Similarly, the deterministic third-moment coefficient is
\begin{align*}
    \mu_{3,n}^{jk\ell}
    &\Let\sqrt n\,\EE[\gamma_j\gamma_k\gamma_\ell]\\
    &=\sqrt n\,\EE[G_{n,j}G_{n,k}G_{n,\ell}]
      +\EE[Q_{n,j}G_{n,k}G_{n,\ell}]\\
    &\quad+\EE[G_{n,j}Q_{n,k}G_{n,\ell}]
      +\EE[G_{n,j}G_{n,k}Q_{n,\ell}]
      +O(n^{-1/2}).
\end{align*}
The first term in the last display simplifies exactly to
\[
    \sqrt n\,\EE[G_{n,j}G_{n,k}G_{n,\ell}]
    =\EE[(V^{-\half}\mf h(X))_j
         (V^{-\half}\mf h(X))_k
         (V^{-\half}\mf h(X))_\ell].
\]
Here $e_j$ is the $j$th coordinate vector in $\RR^{d_h}$, and all
expectations are under $\PP\opt$. The moment assumptions used in the
Edgeworth expansion make $\mu_{1,n}$ and $\mu_{3,n}$ uniformly bounded;
when their limits exist, the limits may be used as the corresponding
fixed Edgeworth coefficients.

\section{Computation details in the proofs}
\subsection{Expansion of {\texorpdfstring{$M_n(\zeta)$}{TEXT}} in the proof of Proposition~\ref{prop:expanMn2}}\label{app:algebra_5.21}
    In the following computation, we suppress the dependence of $\bar h$ and its derivatives on $X_i$ for simplicity. For $\zeta\in \mc Z_n$, the same reasoning as in the proof of Lemma~\ref{lem:DelinMn3} gives $\norm{\Delta_{i}}_2 \leq C \norm{\Sigma}_2 \log(n)^2$ for a constant $C$; hence
    \begingroup\compactdisplaytrue
    \begin{align*}
        \sqrt{n}\left(\bar h\left(X_i + n^{-\half} \Delta_{i}\right) - \bar h\left(X_i\right)\right) &= \underbrace{\bar h_{\alpha} \Delta_i^{\alpha}}_{(A)} + \underbrace{\frac{1}{2\sqrt{n}} \bar h_{\alpha \alpha'}\Delta_i^{\alpha} \Delta_i^{\alpha'}}_{(B)} + \underbrace{\frac{1}{6n} \bar h_{\alpha \alpha' \alpha''}\Delta_i^{\alpha} \Delta_i^{\alpha'} \Delta_i^{\alpha''}}_{(C)} + \widetilde{O}\pare{n^{-\frac{3}{2}}}.
    \end{align*}
    \endgroup
    By Lemma~\ref{lem:DelinMn3}, we get
    \begin{itemize}
        \item For part $(A)$,
        \begin{align*}
            &\bar h_{\alpha} \Delta_i^{\alpha} \\ 
            =& \bar h_{\alpha}\Sigma_{\beta\alpha} \left(\frac{1}{2}\bar h_{\beta} + \frac{1}{\sqrt{n}}\frac{\bar h_{\beta \gamma} \bar h_{\omega} \Sigma_{\omega \gamma}}{4} + \frac{1}{n}\frac{2 \bar h_{\beta \gamma} \bar h_{\omega \omega'} \bar h_{\gamma'}\Sigma_{\omega \gamma} \Sigma_{\gamma'\omega'} + \bar h_{\beta\gamma \omega} \bar h_{\gamma'} \bar h_{\omega'}\Sigma_{\gamma' \gamma}\Sigma_{\omega' \omega}}{16}\right) \\
            &+ \widetilde{O}\left(n^{-\frac{3}{2}}\right).
        \end{align*}

        \item For part $(B)$,
        \begin{align*}
            &\frac{1}{2\sqrt{n}} \bar h_{\alpha \alpha'}\Delta_i^{\alpha} \Delta_i^{\alpha'} \\
            =& \frac{1}{2 \sqrt{n}} \bar h_{\alpha \alpha'} \Sigma_{\beta \alpha}\Sigma_{\beta'\alpha'} \left(\frac{1}{4}\bar h_{\beta}\bar h_{\beta'} + \frac{1}{\sqrt{n}}\frac{\bar h_{\beta} \bar h_{\omega'} \Sigma_{\omega'\gamma'}\bar h_{\beta'\gamma'} + \bar h_{\beta'}\bar h_{\omega} \Sigma_{\omega \gamma} \bar h_{\beta \gamma}}{8}\right)  + \widetilde{O}\left(n^{-\frac{3}{2}}\right).
        \end{align*}

        \item For part $(C)$,
        \begin{align*}
            \frac{1}{6n} \bar h_{\alpha \alpha' \alpha''}\Delta_i^{\alpha} \Delta_i^{\alpha'} \Delta_i^{\alpha''} 
            = \frac{1}{6n} \bar h_{\alpha \alpha' \alpha''} \Sigma_{\beta \alpha}\Sigma_{\beta'\alpha'}\Sigma_{\beta'' \alpha''} \left(\frac{1}{8}\bar h_{\beta}\bar h_{\beta'}\bar h_{\beta''}\right) + \widetilde{O}\left(n^{-\frac{3}{2}}\right).
        \end{align*}
    \end{itemize}

    Combining these expansions gives
    \begingroup\compactdisplaytrue
    \begin{align*}
        & \sqrt{n}\left(\bar h\left(X_i + n^{-\half} \Delta_{i}\right) - \bar h\left(X_i\right)\right) \\
        =&  \frac{1}{2} \bar h_{\alpha} \bar h_{\beta} \Sigma_{\beta\alpha} + \frac{1}{\sqrt{n}} \frac{3\bar h_{\alpha \alpha'} \bar h_{\beta}\bar h_{\beta'} \Sigma_{\beta \alpha}\Sigma_{\beta'\alpha'}}{8} \\
        &+ \frac{1}{n} \left(\frac{2 \bar h_{\beta \gamma} \bar h_{\omega \omega'} \bar h_{\gamma'}\Sigma_{\omega \gamma} \Sigma_{\gamma'\omega'} + \bar h_{\beta\gamma \omega} \bar h_{\gamma'} \bar h_{\omega'}\Sigma_{\gamma' \gamma}\Sigma_{\omega' \omega}}{16}\right) \Sigma_{\beta\alpha} \bar h_{\alpha} \\
        & + \frac{1}{n} \frac{\bar h_{\alpha \alpha'} \bar h_{\beta'\gamma'} \bar h_{\beta} \bar h_{\omega'} \Sigma_{\omega'\gamma'}\Sigma_{\beta \alpha}\Sigma_{\beta'\alpha'}}{8} + \frac{1}{n}\left(\frac{\bar h_{\alpha \alpha' \alpha''}\bar h_{\beta}\bar h_{\beta'}\bar h_{\beta''}}{48}\right)\Sigma_{\beta \alpha}\Sigma_{\beta'\alpha'}\Sigma_{\beta'' \alpha''} + \widetilde{O}\left(n^{-\frac{3}{2}}\right)
    \end{align*}
    \endgroup

    On the other hand, we have
    \begin{align*}
        \Delta_{i}^\top \Sigma^{-1} \Delta_{i} = & \frac{1}{4} \D \bar h\left(X_i +n^{-\half} \Delta_{i}\right) \Sigma \D \bar h\left(X_i + n^{-\half} \Delta_{i}\right)^\top\\
        = & \frac{1}{4} \bar h_{\alpha} \bar h_{\alpha'} \Sigma_{\alpha \alpha'} + \underbrace{
        \frac{1}{\sqrt{n}} \frac{\left(\bar h_{\alpha'} \bar h_{\alpha \beta} \Delta_i^{\beta} + \bar h_{\alpha} \bar h_{\alpha' \beta'} \Delta_i^{\beta'} \right)\Sigma_{\alpha \alpha'}}{4}}_{(D)} \\
        & + \underbrace{\frac{1}{n} \frac{\bar h_{\alpha \beta} \bar h_{\alpha'\beta'} \Delta_i^{\beta} \Delta_i^{\beta'} \Sigma_{\alpha \alpha'}}{4}}_{(E)} + \underbrace{\frac{1}{n} \frac{\left(\bar h_{\alpha'} \bar h_{\alpha\beta\gamma} \Delta_i^{\beta} \Delta_i^{\gamma} + \bar h_{\alpha} \bar h_{\alpha'\beta'\gamma'} \Delta_i^{\beta'} \Delta_i^{\gamma'}\right)\Sigma_{\alpha \alpha'}}{8}}_{(F)}\\
        & + \widetilde{O}\left(n^{-\frac{3}{2}}\right).
    \end{align*}

    \begin{itemize}
        \item For part $(D)$,
        \begin{align*}
            &\frac{1}{\sqrt{n}} \frac{\left(\bar h_{\alpha'} \bar h_{\alpha \beta} \Delta_i^{\beta} + \bar h_{\alpha} \bar h_{\alpha' \beta'} \Delta_i^{\beta'} \right)\Sigma_{\alpha \alpha'}}{4} \\
            =& \frac{1}{\sqrt{n}} \frac{\bar h_{\alpha'} \bar h_{\alpha \beta} \bar h_{\omega} \Sigma_{\omega \beta}\Sigma_{\alpha \alpha'}}{4} + \frac{1}{n} \frac{\bar h_{\alpha'} \bar h_{\alpha \beta} \bar h_{\omega\gamma} \bar h_{\gamma'} \Sigma_{\gamma \gamma'} \Sigma_{\omega \beta}\Sigma_{\alpha \alpha'}}{8}\\
            & + \widetilde{O}\left(n^{-\frac{3}{2}}\right).
        \end{align*}
        \item For parts $(E)$ and $(F)$,
        \begin{align*}
            & \frac{1}{n} \frac{\bar h_{\alpha \beta} \bar h_{\alpha'\beta'} \Delta_i^{\beta} \Delta_i^{\beta'} \Sigma_{\alpha \alpha'}}{4} + \frac{1}{n} \frac{\left(\bar h_{\alpha'} \bar h_{\alpha\beta\gamma} \Delta_i^{\beta} \Delta_i^{\gamma} + \bar h_{\alpha} \bar h_{\alpha'\beta'\gamma'} \Delta_i^{\beta'} \Delta_i^{\gamma'}\right)\Sigma_{\alpha \alpha'}}{8} \\
            = & \frac{1}{n} \frac{\bar h_{\alpha \beta} \bar h_{\alpha'\beta'} \bar h_{\omega} \bar h_{\omega'} \Sigma_{\omega \beta}  \Sigma_{\omega' \beta'} \Sigma_{\alpha \alpha'}}{16} + \frac{1}{n} \frac{\bar h_{\alpha\beta\gamma}\bar h_{\alpha'}  \bar h_{\beta'} \bar h_{\gamma'} \Sigma_{\beta \beta'} \Sigma_{\gamma' \gamma}\Sigma_{\alpha \alpha'}}{16} + \widetilde{O}\left(n^{-\frac{3}{2}}\right).
        \end{align*}
    \end{itemize}

    Combining these terms gives
    \begingroup\compactdisplaytrue
    \begin{align*}
        \Delta_{i}^\top \Sigma^{-1} \Delta_{i}
        = &\frac{1}{4} \bar h_{\alpha} \bar h_{\alpha'} \Sigma_{\alpha \alpha'} + \frac{1}{\sqrt{n}} \frac{\bar h_{\alpha'} \bar h_{\alpha \beta} \bar h_{\omega} \Sigma_{\omega \beta}\Sigma_{\alpha \alpha'}}{4}
        + \frac{1}{n} \frac{\bar h_{\alpha'} \bar h_{\alpha \beta} \bar h_{\omega\gamma} \bar h_{\gamma'} \Sigma_{\gamma \gamma'} \Sigma_{\omega \beta}\Sigma_{\alpha \alpha'}}{8}\\
        & + \frac{1}{n} \frac{\bar h_{\alpha \beta} \bar h_{\alpha'\beta'} \bar h_{\omega} \bar h_{\omega'} \Sigma_{\omega \beta}  \Sigma_{\omega' \beta'} \Sigma_{\alpha \alpha'}}{16} + \frac{1}{n} \frac{\bar h_{\alpha\beta\gamma}\bar h_{\alpha'}  \bar h_{\beta'} \bar h_{\gamma'} \Sigma_{\beta \beta'} \Sigma_{\gamma' \gamma}\Sigma_{\alpha \alpha'}}{16} + \widetilde{O}\left(n^{-\frac{3}{2}}\right).
    \end{align*}
    \endgroup

    Finally, we get
    \begin{align*}
        & \sqrt{n}\left(\bar h\left(X_i + n^{-\half} \Delta_i\right) - \bar h\left(X_i\right)\right) - \Delta_{i}^\top \Sigma^{-1} \Delta_{i}\\
        = & \frac{1}{4} \bar h_{\alpha} \bar h_{\beta} \Sigma_{\beta\alpha} + \frac{1}{\sqrt{n}} \frac{\bar h_{\alpha \alpha'} \bar h_{\beta}\bar h_{\beta'} \Sigma_{\beta \alpha}\Sigma_{\beta'\alpha'}}{8} + \frac{1}{n} \frac{\bar h_{\alpha \alpha'} \bar h_{\beta'\gamma'} \bar h_{\beta} \bar h_{\omega'} \Sigma_{\omega'\gamma'}\Sigma_{\beta \alpha}\Sigma_{\beta'\alpha'}}{16}\\
        & + \frac{1}{n}\frac{\bar h_{\alpha \alpha' \alpha''}\bar h_{\beta}\bar h_{\beta'}\bar h_{\beta''}\Sigma_{\beta \alpha}\Sigma_{\beta'\alpha'}\Sigma_{\beta'' \alpha''}}{48} + \widetilde{O}\left(n^{-\frac{3}{2}}\right).
    \end{align*}

\subsection{Expansion details in the proof of Proposition~\ref{prop:exp_Fndagger}}\label{sec:pf_prop725}
We use superscripts for the coordinates of vectors in $\RR^{d_h}$.

    A direct perturbation calculation avoids identifying the distinct first
    and last indices of the cubic coefficient. Put
    $\varepsilon=n^{-1/2}$ and write
    \begin{equation*}
        F_n^\dagger(\zeta)
        =\frac14\zeta^\top\mc V_n\zeta
         +\varepsilon P_{3,n}(\zeta)
         +\varepsilon^2P_{4,n}(\zeta),
        \qquad
        P_{3,n}(\zeta)=\frac18A_{\alpha\beta\gamma}
        \zeta^{(\alpha)}\zeta^{(\beta)}\zeta^{(\gamma)},
    \end{equation*}
    where $P_{4,n}$ is the quartic polynomial in
    Proposition~\ref{prop:expanMn2}.  Let
    $q_n(\zeta)=-\zeta^\top H_n-\frac14\zeta^\top\mc V_n\zeta$ and
    $\zeta_0=-2\xi_n$.  With
    \begin{equation*}
        g_\gamma
        \Let\{\D P_{3,n}(\zeta_0)\}_\gamma
        =S_{\gamma\omega^1\omega^2}
          \xi_n^{(\omega^1)}\xi_n^{(\omega^2)},
    \end{equation*}
    the first-order condition gives
    \begin{equation*}
        \zeta_n^{\ddagger}
        =\zeta_0-2\varepsilon\mc V_n^{-1}g
         +\widetilde{O}(\varepsilon^2).
    \end{equation*}
    Taylor expansion at $\zeta_0$ therefore yields
    \begin{align*}
        &-\zeta_n^{\ddagger\top}H_n-F_n^\dagger(\zeta_n^{\ddagger})\\
        &\quad=q_n(\zeta_0)-\varepsilon P_{3,n}(\zeta_0)
        +\varepsilon^2\left\{g^\top\mc V_n^{-1}g
        -P_{4,n}(\zeta_0)\right\}
        +\widetilde{O}(\varepsilon^3)\\
        &\quad=\brak{\mc V_n,\xi_n^{\otimes2}}
        +\frac1{\sqrt n}\brak{\mc K_n,\xi_n^{\otimes3}}
        +\frac1n\brak{\mc L_n,\xi_n^{\otimes4}}
        +\widetilde{O}(n^{-3/2}).
    \end{align*}
    Indeed, $-P_{4,n}(-2\xi_n)$ gives the first two negative
    contractions in the displayed definition of $\mc L_n$, whereas
    \begin{equation*}
        g^\top\mc V_n^{-1}g
        =S_{\gamma\omega^1\omega^2}(\mc V_n^{-1})_{\gamma\theta}
         S_{\theta\omega^3\omega^4}
         \prod_{j=1}^4\xi_n^{(\omega^j)}.
    \end{equation*}
    This is the term that cannot be collapsed by treating the first index of
    $A$ as interchangeable with its last two indices.  When $d_h=1$,
    $S=3A/2$, so the contribution reduces to $9A^2/(4\mc V_n)$, agreeing
    with the scalar Bartlett calculation.

\section{Expansion of cumulants}\label{sec:expan_cum_compute}

\subsection{Expansion of cumulants in the proof of Theorem~\MainPowerExpansionNumber{} (\ref{pf:thm4_5})}
\label{app:generic_power_correction}
\begin{lemma}[Cumulants of $\gamma^{\ddagger}$ in the proof (\ref{pf:thm4_5}) of Theorem~\MainPowerExpansionNumber]\label{lem:alter_cumulantII}
Recall that for $\eta \in \RR^{d_h}$:
\begin{align*}
        \langle \gamma^{\ddagger}, \eta \rangle =& \left\langle V_{n}^{\half}, \bar \xi_n \otimes \eta\right\rangle + \half \mc L_{W_{n},V_{n}}(\mc W_n - W_{n}, \mc V_n - V_{n}) \left\langle V_{n}^{\half}, \bar \xi_n \otimes \eta\right\rangle \\
        & - \frac{1}{2} \left\langle \mc V_n - V_{n}, \bar \xi_n \otimes \left(V_{n}^{-\half} \eta\right) \right\rangle + \frac{1}{2\sqrt{n}} \left\langle \EE_{n}[\mc K_n], \bar \xi_n \otimes \bar \xi_n \otimes \left(V_{n}^{-\half} \eta\right) \right\rangle ,
    \end{align*}
    where $\bar \xi_n = \frac{1}{\sqrt{n}}\sum_{i=1}^n  V_{n}^{-1} \mf h(X_i)$.

    Under the hypotheses of Theorem~\MainPowerExpansionNumber, including
    Assumption~\ref{a:alter_edgeworthII}, the first three cumulants of
    $\gamma^{\ddagger}$, denoted by $M_1\in\RR^{d_h}$,
    $M_2\in(\RR^{d_h})^{\otimes2}$, and
    $M_3\in(\RR^{d_h})^{\otimes3}$, admit the expansions
    \begin{align*}
        \brak{M_1, \eta} = & \brak{\tau_n, \eta} + \frac{1}{\sqrt{n}} K_{1}(\eta), \\
        \brak{M_2, \eta^{\otimes 2}} = & \brak{V_n^{-\half} W_n V_n^{-\half}, \eta^{\otimes 2}} + \frac{1}{\sqrt{n}} K_{2}(\eta) + O(n^{-1}),\\
        \brak{M_3, \eta^{\otimes 3}} = & \frac{1}{\sqrt{n}} K_{3}(\eta) + O(n^{-1}),
    \end{align*}
    where
    \begin{subequations}\label{eq:K_II}
    \begin{align}
        \tau_n = & \sqrt{n} V_{n}^{-\half} \EE_{n}\brac{\mf h(X)},\\
        K_{1}(\eta) = & \frac{1}{2}\brak{V_{n}^{\half}, \EE_n\brac{\mc L_{W_n, V_n}(\sqrt{n}\pare{\mc W_n - W_{n}}, \sqrt{n}\pare{\mc V_n - V_{n}})\bar \xi_n}\otimes \eta}\notag\\
        &- \frac{1}{2} \brak{\EE_n\brac{\sqrt{n}\left(\mc V_n - V_{n}\right)\bar \xi_n}, V_{n}^{-\half}\eta} + \frac{1}{2} \brak{\EE_n[\mc K_n], \EE_n\brac{\bar \xi_n \otimes \bar \xi_n} \otimes \left(V_{n}^{-\half} \eta\right)}, \\
        K_2(\eta)=& 2\EE_n \brac{\brak{ V_{n}^{\half}, \bar \xi_n \otimes \eta}\brak{\xi_n', \eta}} - 2\brak{\tau_n, \eta}K_{1}(\eta), \\
        K_{3}(\eta) =&\EE_n\brak{\pare{V_n^{-\half}\mf h(X)}^{\otimes 3}, \eta^{\otimes 3}} +3 \EE_n\brac{\brak{ V_{n}^{\half}, \bar \xi_n \otimes \eta}^2 \brak{\xi_n', \eta}} \notag\\
        & - 3\brak{V_n^{-\half} W_n V_n^{-\half}, \eta^{\otimes 2}}K_1(\eta) - 3\brak{\tau_n, \eta}K_2(\eta) - 3\brak{\tau_n, \eta}^2K_1(\eta),
    \end{align}
    \end{subequations}
    and $\xi_n' = \sqrt{n}\pare{\gamma^{\ddagger} - V_n^{\half}\bar \xi_n}$, $W_n = \EE_n\brac{\mc W_n}$, $V_n = \EE_{n}\brac{\mc V_n}.$
\end{lemma}

For completeness, let $v=(v_j)_{j\in[d_h]}$ and
$\kappa_{jk}=(V_n^{\half}W_n^{-1}V_n^{\half})_{jk}$, and set
$v^k=\sum_{j\in[d_h]}\kappa_{jk}v_j$.  The Hermite polynomials used in the
power correction are
\begin{equation}\label{eq:hermit_poly}
\begin{aligned}
    \mathrm h_1(v)&=(\mathrm h_1^j(v))_{j\in[d_h]},
    &\mathrm h_1^j(v)&=v^j,\\
    \mathrm h_2(v)&=(\mathrm h_2^{jk}(v))_{j,k\in[d_h]},
    &\mathrm h_2^{jk}(v)&=v^jv^k-\kappa_{jk},\\
    \mathrm h_3(v)&=(\mathrm h_3^{jk\ell}(v))_{j,k,\ell\in[d_h]},
    &\mathrm h_3^{jk\ell}(v)&=v^jv^kv^\ell-\kappa_{kj}v^\ell
      -\kappa_{\ell k}v^j-\kappa_{j\ell}v^k.
\end{aligned}
\end{equation}
For $r=1,2,3$, let $\mathsf K_r$ be the unique symmetric coefficient tensor
associated with the homogeneous polynomial $K_r$, so that
$K_r(\eta)=\brak{\mathsf K_r,\eta^{\otimes r}}$. Thus the polynomial in
Theorem~\MainPowerExpansionNumber{} is
\begin{align}\label{eq:barp_definition}
    \bar p(v)
    =\brak{\mathsf K_1,\mathrm h_1(v)}
      +\frac12\brak{\mathsf K_2,\mathrm h_2(v)}
      +\frac16\brak{\mathsf K_3,\mathrm h_3(v)},
\end{align}
with $K_1,K_2,K_3$ defined in~\eqref{eq:K_II}.

\begin{proof}[Proof of Lemma~\ref{lem:alter_cumulantII}]
    We compute the expansion of the first, second, and third cumulants of $\gamma^{\ddagger}$ up to order $n^{-\half}$. All expectations in this proof are taken under $\PP_n\opt$, so below we write $\EE$ for $\EE_{n}$.
    \begin{enumerate}
        \item For the first-order cumulants, we have
    \begingroup\compactdisplaytrue
    \begin{align*}
        \brak{M_1, \eta} =& \EE\langle\gamma^{\ddagger}, \eta\rangle \\
        =& \brak{V_{n}^{\half}, \EE\brac{\bar \xi_n} \otimes \eta} + \frac{1}{2\sqrt{n}}\brak{V_{n}^{\half}, \EE\brac{\mc L_{W_n, V_n}(\sqrt{n}\pare{\mc W_n - W_{n}}, \sqrt{n}\pare{\mc V_n - V_{n}})\bar \xi_n}\otimes \eta}\\
        &- \frac{1}{2\sqrt{n}} \brak{\EE\brac{\sqrt{n}\left(\mc V_n - V_{n}\right)\bar \xi_n}, V_{n}^{-\half}\eta} + \frac{1}{2\sqrt{n}} \brak{\EE[\mc K_n], \EE\brac{\bar \xi_n \otimes \bar \xi_n} \otimes \left(V_{n}^{-\half} \eta\right)}\\
        \Let & \brak{\tau_n, \eta} + \frac{1}{\sqrt{n}} K_{1}(\eta).
    \end{align*}
    \endgroup

    \item For the second-order cumulants, we have
    \begingroup\compactdisplaytrue
    \begin{align*}
        \brak{M_2, \eta^{\otimes 2}} =& \EE\brak{\gamma^{\ddagger} - \EE\brac{\gamma^{\ddagger}}, \eta}^2 \\
        =& \EE\langle\gamma^{\ddagger}, \eta\rangle^2 - \brak{ \EE\brac{\gamma^{\ddagger}}, \eta}^2\\
        =& \EE\brac{\left\langle V_{n}^{\half}, \bar \xi_n \otimes \eta\right\rangle^2}
        + \pare{\EE\langle\gamma^{\ddagger}, \eta\rangle^2
        - \EE\brac{\brak{ V_{n}^{\half}, \bar \xi_n \otimes \eta}^2}}
        - \pare{\brak{\tau_n, \eta} + \frac{1}{\sqrt{n}} K_{1}(\eta)}^2\\
        \Let& \brak{V_n^{-\half} W_n V_n^{-\half}, \eta^{\otimes 2}} + \frac{1}{\sqrt{n}} K_{2}(\eta) + O\pare{n^{-1}},
    \end{align*}
    \endgroup
    Defining $\xi_n' = \sqrt{n}\pare{\gamma^{\ddagger} - V_n^{\half}\bar \xi_n}$, we have
    \begin{align*}
        K_2(\eta) = 2\EE \brak{ V_{n}^{\half}, \bar \xi_n \otimes \eta}\brak{\xi_n', \eta} - 2\brak{\tau_n, \eta}K_{1}(\eta).
    \end{align*}

    \item For the third-order cumulants, we have
    \begingroup\compactdisplaytrue
    \begin{align*}
        &\brak{M_3, \eta^{\otimes 3}} \\
        =& \EE\brak{\gamma^{\ddagger} - \EE\brac{\gamma^{\ddagger}}, \eta}^3 \\
        =&  \EE\brak{\gamma^{\ddagger}, \eta}^3 - 3\pare{\EE\brak{\gamma^{\ddagger}, \eta}^2 - \brak{ \EE\brac{\gamma^{\ddagger}}, \eta}^2}\brak{\EE\brac{\gamma^{\ddagger}}, \eta} - \brak{ \EE\brac{\gamma^{\ddagger}}, \eta}^3\\
        =& \EE\pare{\brak{ V_{n}^{\half}, \bar \xi_n \otimes \eta}}^3 +\frac{3}{\sqrt{n}} \EE\brac{\brak{ V_{n}^{\half}, \bar \xi_n \otimes \eta}^2 \brak{\xi_n', \eta}}\\
        &- 3\pare{\brak{V_n^{-\half} W_n V_n^{-\half}, \eta^{\otimes 2}} + \frac{1}{\sqrt{n}} K_{2}(\eta)}\pare{\brak{\tau_n, \eta} + \frac{1}{\sqrt{n}} K_{1}(\eta)} - \pare{\brak{\tau_n, \eta} + \frac{1}{\sqrt{n}} K_{1}(\eta)}^3\\
        =& \EE\pare{\brak{ V_{n}^{\half}, \bar \xi_n \otimes \eta}}^3 +\frac{3}{\sqrt{n}} \EE\brac{\brak{ V_{n}^{\half}, \bar \xi_n \otimes \eta}^2 \brak{\xi_n', \eta}}\\
        & - 3\brak{V_n^{-\half} W_n V_n^{-\half}, \eta^{\otimes 2}}\brak{\tau_n, \eta} - \frac{3}{\sqrt{n}}\brak{V_n^{-\half} W_n V_n^{-\half}, \eta^{\otimes 2}}K_1(\eta) - \frac{3}{\sqrt{n}}\brak{\tau_n, \eta}K_2(\eta)\\
        &- \brak{\tau_n, \eta}^3 - \frac{3}{\sqrt{n}}\brak{\tau_n, \eta}^2K_1(\eta) + O\pare{n^{-1}}\\
        =& \EE\pare{\left\langle V_{n}^{\half}, \bar \xi_n \otimes \eta\right\rangle - \left\langle V_{n}^{\half}, \EE\brac{\bar \xi_n} \otimes \eta\right\rangle}^3 +\frac{3}{\sqrt{n}} \EE\brac{\brak{ V_{n}^{\half}, \bar \xi_n \otimes \eta}^2 \brak{\xi_n', \eta}}\\
        & - \frac{3}{\sqrt{n}}\brak{V_n^{-\half} W_n V_n^{-\half}, \eta^{\otimes 2}}K_1(\eta) - \frac{3}{\sqrt{n}}\brak{\tau_n, \eta}K_2(\eta) - \frac{3}{\sqrt{n}}\brak{\tau_n, \eta}^2K_1(\eta) + O\pare{n^{-1}}\\
        =& \frac{1}{\sqrt{n}}\EE\brak{\pare{V_n^{-\half}\mf h(X)}^{\otimes 3}, \eta^{\otimes 3}} +\frac{3}{\sqrt{n}} \EE\brac{\brak{ V_{n}^{\half}, \bar \xi_n \otimes \eta}^2 \brak{\xi_n', \eta}}\\
        & - \frac{3}{\sqrt{n}}\brak{V_n^{-\half} W_n V_n^{-\half}, \eta^{\otimes 2}}K_1(\eta) - \frac{3}{\sqrt{n}}\brak{\tau_n, \eta}K_2(\eta) - \frac{3}{\sqrt{n}}\brak{\tau_n, \eta}^2K_1(\eta) + O\pare{n^{-1}}\\
        \Let & \frac{1}{\sqrt{n}}K_3(\eta) + O\pare{n^{-1}}.
    \end{align*}
    \endgroup
    \end{enumerate}
\end{proof}
\subsection{Expansion of cumulants in the proof of Proposition~\ref{prop:expan_rwpi} (\ref{sec:pf_prop_5_7})}\label{app:prop_5_7}
\begin{lemma}[First three cumulants of $\bar \gamma$ under local alternatives]\label{lem:alter_seccumulant}
        Under the hypotheses of Proposition~\ref{prop:expan_rwpi}, the first
        three cumulants of $\bar \gamma$, denoted by $m_1$, $m_2$, and $m_3$,
        admit the expansions
        \begin{align*}
            m_1 = & \frac{\EE_{0}\left[\D h(X)\right] \tau_0}{\sqrt{\alpha_{2}}} + \frac{1}{\sqrt{n}} k_{1} + O(n^{-1}), \\
            m_2 = & 1 + \frac{1}{\sqrt{n}} k_{2} + O(n^{-1}),\\
            m_3 = & \frac{1}{\sqrt{n}} k_{3} + O(n^{-1}),
        \end{align*}
        where
        \begin{subequations}
        \begin{align*}
            k_{1} = & \frac{1}{2} \frac{\tau_0^\top \EE_{0}[\D^2 h(X)] \tau_0}{\sqrt{\alpha_{2}}} - \frac{\tau_0^\top \EE_{0}[\D h(X)^\top] \EE_{0}\left[h(X) \D h(X)\right] \tau_0}{\alpha_{2}^{\frac{3}{2}}} - \frac{\alpha_{3}}{2 \alpha_{2}^{\frac{3}{2}}} + \frac{\tilde \alpha_{3} \sqrt{\alpha_{2}} }{2  \tilde \alpha_{2}^2} \notag\\
            & + \frac{\tilde \alpha_{3}}{2 \sqrt{\alpha_{2}} \tilde \alpha_{2}^2} \tau_0^\top \EE_0\left[\D h(X)^\top\right] \EE_0\left[\D h(X)\right] \tau_0, \\
            k_{2} = & \left(- \frac{\alpha_{3}}{\alpha_{2}^{2}} + \frac{2 \tilde \alpha_{3}}{\tilde \alpha_{2}^2}\right) \EE_{0}[\D h(X)] \tau_0,\\
            k_{3} =& -\frac{2\alpha_{3}}{\alpha_{2}^{\frac{3}{2}}} + \frac{3\tilde \alpha_{3} \alpha_{2}^{\half}}{\tilde \alpha_{2}^2}. 
        \end{align*}
        \end{subequations}
    \end{lemma}

    \begin{proof}[Proof of Lemma~\ref{lem:alter_seccumulant}]
        Recall that $\bar \gamma = \sqrt{n} R_1 + \sqrt{n} R_2$.
        For the first-order cumulant, we have
        \begin{align*}
            m_1 & = \EE_n[\sqrt{n} R_1] + \EE_n[\sqrt{n} R_2].
        \end{align*}
        For $\EE_n[\sqrt{n} R_1]$, we have
        \begin{align*}
            \EE_n[\sqrt{n} R_1] = \EE_n\left[\frac{\sqrt{n}}{\sqrt{\alpha_{2,n}}} (A_1 + \alpha_{1,n})\right] = \frac{\sqrt{n} \alpha_{1,n}}{\sqrt{\alpha_{2,n}}},
        \end{align*}
        where $\EE_n[A_1] = 0$.

        Assumptions~2.5--2.6 and~2.9 permit Taylor expansion under $\EE_0$ and
        provide integrable bounds for the remainders.  We therefore obtain
        \begin{subequations}\label{eq:expan_alpha}
        \begin{align}
        \alpha_{1,n} &= \frac{1}{\sqrt{n}}\EE_{0}[\D h(X)] \tau_0 + \frac{1}{2 n} \tau_0^\top \EE_{0}[\D^2 h(X)] \tau_0 + O\left(n^{-\frac{3}{2}}\right),\\
        \alpha_{2,n} & = \EE_{0}\left[h(X)^2\right] + \frac{2}{\sqrt{n}}\EE_{0}\left[h(X) \D h(X)\right] \tau_0 + O\left(n^{-1}\right),\\
        \alpha_{3,n} &= \alpha_3 + O\pare{n^{-1/2}},~\tilde \alpha_{3,n} = \tilde \alpha_3 + O\pare{n^{-1/2}},~\tilde \alpha_{2,n} = \tilde \alpha_2 + O\pare{n^{-1/2}}.
        \end{align}
        \end{subequations}
        Therefore, we have
        \begingroup\compactdisplaytrue
        \begin{align}
            & \frac{\sqrt{n} \alpha_{1,n}}{\sqrt{\alpha_{2,n}}} \notag\\
            = &  \left(\EE_{0}\left[\D h(X)\right] \tau_0 + \frac{1}{2 \sqrt{n}} \tau_0^\top \EE_{0}[\D^2 h(X)] \tau_0 + O\left(\frac{1}{n}\right) \right) \times \notag\\
            &\left(\EE_{0}[h(X)^2] + \frac{2}{\sqrt{n}}\EE_{0}[h(X) \D h(X)] \tau_0 + O\left(\frac{1}{n}\right)\right)^{-\frac{1}{2}}\notag\\
            = & \frac{\EE_{0}\left[\D h(X)\right]}{\sqrt{\alpha_{2}}} \tau_0 + 
            \frac{1}{\sqrt{n}} \left(\frac{1}{2} \frac{\tau_0^\top \EE_{0}[\D^2 h(X)] \tau_0}{\sqrt{\alpha_{2}}} - \frac{\tau_0^\top \EE_{0}[\D h(X)^\top] \EE_{0}[h(X) \D h(X)] \tau_0}{\alpha_{2}^{\frac{3}{2}}} \right) + O\left(\frac{1}{n}\right). \label{eq:exp_1order_1term}
        \end{align}
        \endgroup

        For $\EE_n[\sqrt{n} R_2]$, we have
        \begin{align*}
            \EE_n[\sqrt{n} R_2] & = \EE_n\left[- \frac{1}{2 \alpha_{2,n}^{\frac{3}{2}}} \sqrt{n}A_2 (A_1 + \alpha_{1,n}) + \frac{\tilde \alpha_{3,n}}{2 \sqrt{\alpha_{2,n}} \tilde \alpha_{2,n}^2}\sqrt{n}(A_1 + \alpha_{1,n})^2\right]\\
            & = \EE_n\left[- \frac{1}{2 \alpha_{2,n}^{\frac{3}{2}}} \sqrt{n}A_2 A_1 + \frac{\tilde \alpha_{3,n}}{2 \sqrt{\alpha_{2,n}} \tilde \alpha_{2,n}^2}\left(\sqrt{n}A_1^2 + \sqrt{n}\alpha_{1,n}^2\right) \right],
        \end{align*}
        and
        \begin{align*}
            \EE_n[\sqrt{n} A_2 A_1] & = \frac{1}{\sqrt{n}} \left(\EE_n[h(X)^3] - \alpha_{2,n}\alpha_{1,n}\right) = \frac{1}{\sqrt{n}} \left(\alpha_{3,n} - \alpha_{2,n}\alpha_{1,n}\right) ,\\
            \EE_n[\sqrt{n} A_1^2] & = \frac{1}{\sqrt{n}}\left(\EE_n[h(X)^2] - \alpha_{1,n}^2\right) = \frac{1}{\sqrt{n}}\left(\alpha_{2,n} - \alpha_{1,n}^2\right).
        \end{align*}
        Thus, by \eqref{eq:expan_alpha}, we replace $\alpha_{2,n}$,
        $\alpha_{3,n}$, $\tilde\alpha_{2,n}$, and $\tilde\alpha_{3,n}$ by their
        null limits while retaining the displayed expansion of $\alpha_{1,n}$,
        and obtain
        \begingroup\compactdisplaytrue
        \begin{align}\label{eq:exp_2order_2term}
            \EE_n\brac{\sqrt{n} R_2} = \frac{1}{\sqrt{n}}\left(- \frac{\alpha_{3}}{2 \alpha_{2}^{\frac{3}{2}}} + \frac{\tilde \alpha_{3} \sqrt{\alpha_{2}} }{2  \tilde \alpha_{2}^2} + \frac{\tilde \alpha_{3}}{2 \sqrt{\alpha_{2}} \tilde \alpha_{2}^2} \tau_0^\top \EE_0\brac{\D h(X)^\top} \EE_0\brac{\D h(X)} \tau_0 \right) + O\pare{\frac{1}{n}}.
        \end{align}
        \endgroup
        Combining~\eqref{eq:exp_1order_1term} and~\eqref{eq:exp_2order_2term} gives
        \begin{align*}
            m_1 = & \frac{\EE_{0}\left[\D h(X)\right]}{\sqrt{\alpha_{2}}} \tau_0 + \frac{1}{\sqrt{n}}k_{1} + O\left(\frac{1}{n}\right),\\
            k_{1} = & \frac{1}{2} \frac{\tau_0^\top \EE_{0}[\D^2 h(X)] \tau_0}{\sqrt{\alpha_{2}}} - \frac{\tau_0^\top \EE_{0}[\D h(X)^\top] \EE_{0}\left[h(X) \D h(X)\right] \tau_0}{\alpha_{2}^{\frac{3}{2}}} - \frac{\alpha_{3}}{2 \alpha_{2}^{\frac{3}{2}}} + \frac{\tilde \alpha_{3} \sqrt{\alpha_{2}} }{2  \tilde \alpha_{2}^2}\\
            & + \frac{\tilde \alpha_{3}}{2 \sqrt{\alpha_{2}} \tilde \alpha_{2}^2} \tau_0^\top \EE_0\left[\D h(X)^\top\right] \EE_0\left[\D h(X)\right] \tau_0.
        \end{align*}

        For the second-order cumulants, we have
        \begingroup\compactdisplaytrue
        \begin{align*}
            m_2 = \text{Var}_n[\sqrt{n}R_1 + \sqrt{n}R_2] = \text{Cov}_n(\sqrt{n}R_1, \sqrt{n}R_1) + 2 \text{Cov}_n(\sqrt{n}R_1, \sqrt{n}R_2) + \text{Cov}_n(\sqrt{n} R_2, \sqrt{n} R_2),
        \end{align*}
        \endgroup
        where $\text{Var}_n$ and $\text{Cov}_n$ are the variance and the covariance under $\PP_{n}\opt$.

        For $\text{Cov}_n(\sqrt{n}R_1, \sqrt{n}R_1)$, we have
        \begin{align}\label{eq:cov11}
            \text{Cov}_n(\sqrt{n}R_1, \sqrt{n}R_1) = \frac{1}{\alpha_{2,n}} \EE_n\left[(\sqrt{n} A_1)^2\right] = \frac{\alpha_{2,n} - \alpha_{1,n}^2}{\alpha_{2,n}} = 1 + O\left(\frac{1}{n}\right).
        \end{align}

        For $\text{Cov}_n(\sqrt{n}R_1, \sqrt{n}R_2)$, we have
        \begingroup\compactdisplaytrue
        \begin{align}\label{eq:cov12}
            \text{Cov}_n(\sqrt{n}R_1, \sqrt{n}R_2) = & - \frac{\alpha_{1,n}}{2 \alpha_{2,n}^{2}} \EE_n \left[(\sqrt{n}A_1)(\sqrt{n}A_2)\right] + \frac{\alpha_{1,n}\tilde \alpha_{3,n}}{\alpha_{2,n}\tilde \alpha_{2,n}^2} \EE_n\left[(\sqrt{n}A_1)^2\right] + O\left(\frac{1}{n}\right)\notag\\
            = & - \frac{\alpha_{1,n} \alpha_{3,n}}{2 \alpha_{2,n}^{2}} + \frac{\alpha_{1,n}\tilde \alpha_{3,n}}{\tilde \alpha_{2,n}^2} + O\left(\frac{1}{n}\right)\notag\\
            = & \frac{1}{\sqrt{n}} \left(- \frac{\alpha_{3}}{2 \alpha_{2}^{2}} + \frac{\tilde \alpha_{3}}{\tilde \alpha_{2}^2}\right) \EE_{0}[\D h(X)] \tau_0 + O\left(\frac{1}{n}\right),
        \end{align}
        \endgroup
        For the first equality, we use the following identity for i.i.d.
        mean-zero random vectors
        $Z_i=(Z_i^{(1)},Z_i^{(2)},Z_i^{(3)})$ under $\PP$:
        \begin{align*}
            \EE_{\PP}\left[\left(\frac{1}{n}\sum_{i=1}^n Z_i^{(1)}\right)\left(\frac{1}{n}\sum_{i=1}^n Z_i^{(2)}\right)\left(\frac{1}{n}\sum_{i=1}^n Z_i^{(3)}\right)\right] = \frac{1}{n^2} \EE_{\PP}\left[Z_1^{(1)}Z_1^{(2)}Z_1^{(3)}\right].
        \end{align*}
        The last equality uses \eqref{eq:expan_alpha}.
        
        For $\text{Cov}_n(\sqrt{n}R_2, \sqrt{n}R_2)$, by similar reasoning as above, we get
     \begin{align}\label{eq:cov22}
            \text{Cov}_n(\sqrt{n}R_2, \sqrt{n}R_2) = O\left(\frac{1}{n}\right).
        \end{align}
        Combining~\eqref{eq:cov11}, \eqref{eq:cov12}, and~\eqref{eq:cov22} gives
        \begin{align*}
            m_2 = & 1 + \frac{1}{\sqrt{n}} k_{2} + O\left(\frac{1}{n}\right),\\
            k_{2} = & \left(- \frac{\alpha_{3}}{\alpha_{2}^{2}} + \frac{2 \tilde \alpha_{3}}{\tilde \alpha_{2}^2}\right) \EE_{0}[\D h(X)] \tau_0. 
        \end{align*}
    For the third-order cumulants, we have
    \begin{align*}
        m_3 = \EE_n\left[\left(\sqrt{n}R_1 + \sqrt{n}R_2 - \EE_n\left[\sqrt{n}R_1 + \sqrt{n}R_2\right]\right)^3\right].
    \end{align*}
    The calculation shows that
    \begin{align*}
            & \EE_n\left[\left(\sqrt{n} R_1 - \EE_n[\sqrt{n} R_1]\right)^3\right]\\
            &\quad=n^{\frac{3}{2}}\frac{\EE_n[A_1^3]}{\alpha_{2,n}^{\frac{3}{2}}}\\
            &\quad=\frac{1}{\sqrt{n}}
            \frac{\alpha_{3,n}-3\alpha_{1,n}\alpha_{2,n}+2\alpha_{1,n}^3}
            {\alpha_{2,n}^{\frac{3}{2}}}
            = \frac{1}{\sqrt{n}} \frac{\alpha_{3}}{\alpha_{2}^{\frac{3}{2}}}  + O(n^{-1}),\\
            &\EE_n\left[\left(\sqrt{n} R_1 - \EE_n[\sqrt{n} R_1]\right)^2\left(\sqrt{n} R_2 - \EE_n[\sqrt{n} R_2]\right)\right] \\
            = & n^{\frac{3}{2}}\left(-\frac{\Cov_n\left(A_1^2, A_1 A_2\right)}{2\alpha_{2,n}^{\frac{5}{2}}} + \frac{\tilde \alpha_{3,n} \Cov_n(A_1^2, A_1^2)}{2 \alpha_{2,n}^{\frac{3}{2}}\tilde \alpha_{2,n}^2}\right) + O(n^{-1})\\
            = & \frac{1}{\sqrt{n}}\left(-\frac{\alpha_{3}}{\alpha_{2}^{\frac{3}{2}}} + \frac{\tilde \alpha_{3} \alpha_{2}^{\half}}{\tilde \alpha_{2}^{2}}\right) + O(n^{-1}),\\
            & \EE_n\left[\left(\sqrt{n} R_1 - \EE_n[\sqrt{n}R_1]\right)\left(\sqrt{n}R_2 - \EE_n[\sqrt{n}R_2]\right)^2\right] =  O\left(n^{-1}\right),\\
            & \EE_n[\left(\sqrt{n} R_2 - \EE_n[\sqrt{n} R_2]\right)^3] = O\left(n^{-1}\right).
        \end{align*}
    Thus, we get
    \begin{align*}
        m_3 = \frac{1}{\sqrt{n}} k_{3} + O(n^{-1}) = \frac{1}{\sqrt{n}} \pare{-\frac{2\alpha_{3}}{\alpha_{2}^{\frac{3}{2}}} + \frac{3\tilde \alpha_{3} \alpha_{2}^{\half}}{\tilde \alpha_{2}^2}} + O(n^{-1}).
    \end{align*}
    This observation completes the proof.
    \end{proof}

\subsection{Expansion of cumulants in the proof of Theorem~\MainBartlettCoverageNumber{} (\ref{sec:pf_thm6_5})}\label{app:thm6_5}

\begin{lemma}[First to fourth cumulants of $\hat \gamma$]\label{lem:4cumuR}
    The first four cumulants of $\hat \gamma$ are
    \begin{align*}
        \tilde m_1 &= \frac{k_{11}}{\sqrt{n}} + O\left(n^{-\frac{3}{2}}\right),\\
        \tilde m_2 &= 1 + \frac{k_{22}}{n} + O\left(n^{-2}\right),\\
        \tilde m_3 &= \frac{k_{31}}{\sqrt{n}} + O\left(n^{-\frac{3}{2}}\right),\\
        \tilde m_4 &= \frac{k_{42}}{n} + O\left(n^{-2}\right),
    \end{align*}
    where
    \begin{align*}
        k_{11} =& -\frac{\alpha_3}{2 \alpha_2^{\frac{3}{2}}}  + \frac{\tilde \alpha_3 \alpha_2^{\half}}{2\tilde \alpha_2^2},\\
        k_{22} =& -\frac{3\alpha_3 \tilde \alpha_3}{2\alpha_2 \tilde \alpha_2^2} + \frac{7\alpha_3^2}{4\alpha_2^3} + \frac{-6 \tilde \alpha_3 \EE\left[h(X)\D h(X) \Sigma \D h(X)^\top\right] + 3\tilde \alpha_4 \alpha_2}{\tilde \alpha_2^3}
        \\
        & + \frac{3\EE[h(X)\D h(X) \Sigma \D^2 h(X) \Sigma\D h(X)^\top]}{\tilde \alpha_2^2}
        - \frac{\alpha_2 \tilde \alpha_3^2}{\tilde \alpha_2^4},\\
        k_{31} =& -\frac{2\alpha_3}{\alpha_2^{\frac{3}{2}}} + \frac{3\tilde \alpha_3 \alpha_2^{\half}}{\tilde \alpha_2^{2}},\\
        k_{42} =& \frac{-2 \alpha_4}{\alpha_2^2} + \frac{12 \alpha_3^2}{\alpha_2^3} + \frac{-18 \alpha_3 \tilde \alpha_3}{\alpha_2 \tilde \alpha_2^2} + \frac{12 \EE[h(X) \D h(X)\Sigma \D^2 h(X) \Sigma \D h(X)^\top]}{\tilde \alpha_2^2} \\
        &+ \frac{12 \alpha_2 \tilde \alpha_4 - 24 \tilde \alpha_3\EE[h(X)\D h(X)\Sigma \D h(X)^\top]}{\tilde \alpha_2^3}  + \frac{9\alpha_2 \tilde \alpha_3^2}{\tilde \alpha_2^4}.
    \end{align*}
\end{lemma}

\begin{proof}[Proof of Lemma~\ref{lem:4cumuR}]   
    \begin{enumerate}
        \item $\tilde m_1 = \EE[\hat \gamma] = \EE[\sqrt{n} (R_1 + R_2 + R_3)]$. Thus, we have
        \begin{align*}
            \tilde m_{1} &= \sqrt{n}\left(-\frac{\EE[A_1 A_2]}{2 \alpha_2^{\frac{3}{2}}}  + \frac{\tilde \alpha_3 \EE[A_1^2]}{2\alpha_2^{\half}\tilde \alpha_2^2}\right) + O\left(n^{-\frac{3}{2}}\right)\\
            &= \frac{1}{\sqrt{n}}\left(-\frac{\alpha_3}{2 \alpha_2^{\frac{3}{2}}}  + \frac{\tilde \alpha_3 \alpha_2^{\half}}{2\tilde \alpha_2^2}\right) + O\left(n^{-\frac{3}{2}}\right).
        \end{align*}
        \item $\tilde m_2 = \text{Var}(\sqrt{n}R) = \Cov(\sqrt{n} (R_1 + R_2 + R_3), \sqrt{n} (R_1 + R_2 + R_3))$. We have
        \begingroup\compactdisplaytrue
        \begin{align*}
            \Cov(\sqrt{n}R_1, \sqrt{n}R_2) 
            =&  n\left(-\frac{\Cov(A_1, A_1 A_2)}{2 \alpha_2^{2}} + \frac{\tilde \alpha_3 \Cov(A_1, A_1^2)}{2\alpha_2\tilde \alpha_2^2}\right) \\
            = & \frac{1}{n}\left(-\frac{\alpha_4 - \alpha_2^2}{2 \alpha_2^{2}} + \frac{\alpha_3 \tilde \alpha_3}{2\alpha_2\tilde \alpha_2^2}\right),\\
            \Cov(\sqrt{n}R_1, \sqrt{n}R_3) = & n\left(\frac{3 \Cov(A_1, A_1 A_2^2)}{8\alpha_2^{3}} - \frac{\tilde \alpha_3 \Cov(A_1, A_1^2 A_2)}{4\alpha_2^{2}\tilde \alpha_2^2}  - \frac{\tilde \alpha_3 \Cov(A_1, A_1^2 \tilde A_2)}{\alpha_2\tilde \alpha_2^{3}}\right. \\
            &\left. + \frac{\Cov(A_1, A_1^2 \tilde A_3)}{2\alpha_2\tilde \alpha_2^2} + \left( - \frac{\tilde \alpha_3^2}{8 \alpha_{2}\tilde \alpha_2^4} + \frac{\tilde \alpha_4}{2\alpha_2\tilde \alpha_2^3}\right)\Cov(A_1, A_1^3)\right)\\
            =& \frac{1}{n}\left(\frac{3 \left(\alpha_2 \alpha_4 - \alpha_2^3 + 2\alpha_3^2\right)}{8\alpha_2^{3}} - \frac{3\tilde \alpha_3  \alpha_3}{4\alpha_2\tilde \alpha_2^2} \right.\\ 
            &\left. - \frac{\tilde \alpha_3 \left(3\alpha_2 \EE[h(X)\D h(X) \Sigma \D h(X)^\top]\right)}{\alpha_2\tilde \alpha_2^{3}}\right. \\
            &\left. + \frac{3 \alpha_2 \EE[h(X)\D h(X) \Sigma \D^2 h(X) \Sigma\D h(X)^\top]}{2\alpha_2\tilde \alpha_2^2}\right.\\
            &\left. + \left( - \frac{\tilde \alpha_3^2}{8 \alpha_{2}\tilde \alpha_2^4} + \frac{\tilde \alpha_4}{2\alpha_2\tilde \alpha_2^3}\right)\left(3 \alpha_2^2\right)\right) + O(n^{-2}),\\
            \Cov(\sqrt{n}R_2, \sqrt{n}R_2) = & n\left(\frac{\Cov(A_1 A_2, A_1 A_2)}{4 \alpha_2^3} - \frac{\tilde \alpha_3 \Cov(A_1 A_2, A_1^2)}{2 \alpha_2^2 \tilde \alpha_2^2} + \frac{\tilde \alpha_3^2}{4\alpha_2 \tilde \alpha_2^4}\Cov(A_1^2, A_1^2)\right) \\
            = & \frac{1}{n}\left(\frac{\alpha_2\alpha_4 - \alpha_2^3 + \alpha_3^2}{4 \alpha_2^3} - \frac{\tilde \alpha_3 \alpha_3}{ \alpha_2 \tilde \alpha_2^2} + \frac{\tilde \alpha_3^2 \alpha_2}{2 \tilde \alpha_2^4}\right) + O\left(n^{-2}\right),\\
            \Cov(\sqrt{n}R_2, \sqrt{n}R_3) = & O\left(n^{-2}\right), \\
            \Cov(\sqrt{n}R_3, \sqrt{n}R_3) = & O\left(n^{-2}\right).
        \end{align*}     
        \endgroup
        Thus, we have
        \begingroup\compactdisplaytrue
        \begin{align*}
            \tilde m_2 
            = & \Cov(\sqrt{n}R_1, \sqrt{n}R_1) + 2\Cov(\sqrt{n}R_1, \sqrt{n}R_2) + 2\Cov(\sqrt{n}R_1, \sqrt{n}R_3) + \Cov(\sqrt{n}R_2, \sqrt{n}R_2)\\
            & + O(n^{-2})\\
            = & 1 + \frac{1}{n}\left(-\frac{\alpha_4 - \alpha_2^2}{ \alpha_2^{2}} + \frac{\alpha_3 \tilde \alpha_3}{\alpha_2\tilde \alpha_2^2}  + \frac{3 \left(\alpha_2 \alpha_4 - \alpha_2^3 + 2\alpha_3^2\right)}{4\alpha_2^{3}} - \frac{3\tilde \alpha_3  \alpha_3}{2\alpha_2\tilde \alpha_2^2}  \right. \\
            &- \frac{6\tilde \alpha_3 \left(\alpha_2 \EE[h(X)\D h(X) \Sigma \D h(X)^\top]\right)}{\alpha_2\tilde \alpha_2^{3}} \\
            & + \frac{3 \alpha_2 \EE[h(X)\D h(X) \Sigma \D^2 h(X) \Sigma\D h(X)^\top]}{\alpha_2\tilde \alpha_2^2} + \left( - \frac{\tilde \alpha_3^2}{4 \alpha_{2}\tilde \alpha_2^4} + \frac{\tilde \alpha_4}{\alpha_2\tilde \alpha_2^3}\right)\left(3 \alpha_2^2\right) \\
            &+ \frac{\alpha_2\alpha_4 - \alpha_2^3 + \alpha_3^2}{4 \alpha_2^3} - \frac{\tilde \alpha_3 \alpha_3}{ \alpha_2 \tilde \alpha_2^2} \left. + \frac{\tilde \alpha_3^2 \alpha_2}{2\tilde \alpha_2^4} \right) +  O\left(n^{-2}\right) \\
            = & 1 + \frac{1}{n}\left(-\frac{3\alpha_3 \tilde \alpha_3}{2\alpha_2 \tilde \alpha_2^2} + \frac{7\alpha_3^2}{4\alpha_2^3} + \frac{-6 \tilde \alpha_3 \EE\left[h(X)\D h(X) \Sigma \D h(X)^\top\right] + 3\tilde \alpha_4 \alpha_2}{\tilde \alpha_2^3}\right.\\ 
            &+ \frac{3\EE[h(X)\D h(X) \Sigma \D^2 h(X) \Sigma\D h(X)^\top]}{\tilde \alpha_2^2}  \left. - \frac{\alpha_2 \tilde \alpha_3^2}{\tilde \alpha_2^4}\right) +  O\left(n^{-2}\right) .
        \end{align*}
        \endgroup

        \item For the third cumulant,
        \begin{align*}
        \tilde m_3
        &=\EE\left[\left(\hat\gamma-\EE[\hat\gamma]\right)^3\right]\\
        &=\EE\left[\left(\sqrt n(R_1+R_2+R_3)
          -\EE[\sqrt n(R_1+R_2+R_3)]\right)^3\right].
        \end{align*}
        We have
        \begingroup\compactdisplaytrue
        \begin{align*}
            \EE\left[\left(\sqrt{n} R_1 - \EE[\sqrt{n} R_1]\right)^3\right] = & n^{\frac{3}{2}}\frac{\EE[A_1^3]}{\alpha_2^{\frac{3}{2}}}
            =  \frac{1}{\sqrt{n}} \frac{\alpha_3}{\alpha_2^{\frac{3}{2}}},\\
            \EE\left[\left(\sqrt{n} R_1 - \EE[\sqrt{n} R_1]\right)^2\left(\sqrt{n} R_2 - \EE[\sqrt{n} R_2]\right)\right] = & n^{\frac{3}{2}}\left(-\frac{\Cov\left(A_1^2, A_1 A_2\right)}{2\alpha_2^{\frac{5}{2}}} + \frac{\tilde \alpha_3 \Cov(A_1^2, A_1^2)}{2 \alpha_2^{\frac{3}{2}}\tilde \alpha_2^2}\right)\\
            = & \frac{1}{\sqrt{n}}\left(-\frac{\alpha_3}{\alpha_2^{\frac{3}{2}}} + \frac{\tilde \alpha_3 \alpha_2^{\half}}{\tilde \alpha_2^{2}}\right),\\
            \EE\left[\left(\sqrt{n} R_1 - \EE[\sqrt{n}R_1]\right)\left(\sqrt{n}R_2 - \EE[\sqrt{n}R_2]\right)^2\right] = & O\left(n^{-\frac{3}{2}}\right),\\
            \EE[\left(\sqrt{n} R_2 - \EE[\sqrt{n} R_2]\right)^3] = & O\left(n^{-\frac{3}{2}}\right),
        \end{align*}
        \endgroup
        and other expansion terms involved with $\sqrt{n}R_3 - \EE[\sqrt{n}R_3]$ are of order $O\left(n^{-\frac{3}{2}}\right)$. Thus,
        \begin{align*}
            \tilde m_3 = \frac{1}{\sqrt{n}}\left(-\frac{2\alpha_3}{\alpha_2^{\frac{3}{2}}} + \frac{3\tilde \alpha_3 \alpha_2^{\half}}{\tilde \alpha_2^{2}}\right) + O\left(n^{-\frac{3}{2}}\right).
        \end{align*}

        \item $\tilde m_{4} = \EE\left[\left(\hat \gamma - \EE[\hat \gamma]\right)^4\right] - 3\left(\EE\left[\left(\hat \gamma - \EE[\hat \gamma]\right)^2\right]\right)^2$.

        As for the expansion of $ \EE\left[\left(\hat \gamma - \EE[\hat \gamma]\right)^4\right]$, we have
        \[ \EE\left[\left(\sqrt{n} R_1 - \EE[\sqrt{n} R_1]\right)^4\right] = n^2 \EE\left[\frac{A_1^4}{\alpha_2^2}\right] = 3 + \frac{1}{n}\left(\frac{\alpha_4}{\alpha_2^2}-3\right),\]
        and 
        \begingroup\compactdisplaytrue
        \begin{align*}
            & \EE\left[\left(\sqrt{n} R_1 - \EE[\sqrt{n} R_1]\right)^3\left(\sqrt{n} R_2 - \EE[\sqrt{n} R_2]\right)\right] \\
            =& \frac{n^2}{\alpha_2^{\frac{3}{2}}} \left(\left(-\frac{\EE[A_1^4 A_2]}{2\alpha_{2}^{\frac{3}{2}}} + \frac{\tilde \alpha_3 \EE[A_1^5]}{2\alpha_2^{\half}\tilde \alpha_2^2}\right) - \frac{\alpha_3}{n^3}\left(-\frac{\alpha_3}{2\alpha_2^{\frac{3}{2}}} + \frac{\tilde \alpha_3 \alpha_2^{\half}}{2\tilde \alpha_2^{2}}\right)\right) + O\left(n^{-2}\right),\\
            &\EE\left[\left(\sqrt{n} R_1 - \EE[\sqrt{n} R_1]\right)^2\left(\sqrt{n} R_2 - \EE[\sqrt{n} R_2]\right)^2\right] \\
            = & \frac{n^2}{\alpha_2} \left(\frac{\EE[A_1^4 A_2^2]}{4 \alpha_2^3} - \frac{\tilde \alpha_3 \EE[A_1^5 A_2]}{2 \alpha_2^2 \tilde \alpha_2^2} + \frac{\tilde \alpha_3^2 \EE[A_1^6]}{4\alpha_2 \tilde \alpha_2^4} - 2 \EE[R_2] \left(-\frac{\EE[A_1^3 A_2]}{2 \alpha_2^{\frac{3}{2}}} + \frac{\tilde \alpha_3 \EE[A_1^4]}{2 \alpha_2^{\half}\tilde \alpha_2^2}\right) \right. \\
            &\left.+ \EE[A_1^2] (\EE[R_2])^2 \right) + O\left(n^{-2}\right),\\
            & \EE\left[\left(\sqrt{n} R_1 - \EE[\sqrt{n} R_1]\right)^3\left(\sqrt{n} R_3 - \EE[\sqrt{n} R_3]\right)\right] \\
            = & \frac{n^2}{\alpha_2^{\frac{3}{2}}} \left(\frac{3\EE[A_1^4 A_2^2]}{8 \alpha_2^{\frac{5}{2}}} - \frac{\tilde \alpha_3 \EE[A_1^5 A_2]}{4\alpha_2^{\frac{3}{2}} \tilde \alpha_2^2} - \frac{\tilde \alpha_3 \EE[A_1^5 \tilde A_2]}{\alpha_2^{\half} \tilde \alpha_2^3} + \frac{\EE[A_1^5 \tilde A_3]}{2\alpha_2^{\half}\tilde \alpha_2^2} + \left(-\frac{\tilde \alpha_3^2}{8 \alpha_2^{\half} \tilde \alpha_2^4} + \frac{\tilde \alpha_4}{2 \alpha_2^{\half} \tilde \alpha_2^3}\right) \EE[A_1^6]\right)\\
            & + O\left(n^{-2}\right),
        \end{align*}      
        \endgroup
        The empirical-average moment identities used above are
        \begin{align*}
            \EE[A_1^4 A_2] &= \frac{6(\alpha_4 - \alpha_2^2)\alpha_2 + 4\alpha_3^2}{n^3}+O(n^{-4}), \quad \EE[A_1^5] = \frac{10 \alpha_2 \alpha_3}{n^3}+O(n^{-4}),\\
            \EE[A_1^4 A_2^2] &= \frac{1}{n^3}\left(3 \alpha_2^2 (\alpha_4 - \alpha_2^2) + 12 \alpha_2 \alpha_3^2\right)+O(n^{-4}),\\
            \EE[A_1^5 A_2] &= \frac{15 \alpha_2^2 \alpha_3}{n^3}+O(n^{-4}), \quad \EE[A_1^5 \tilde A_2] = \frac{15 \alpha_2^2 \EE[h \D h\Sigma \D h^\top]}{n^3}+O(n^{-4}),\\
            \EE[A_1^5 \tilde A_3] &= \frac{15 \alpha_2^2 \EE[h \D h\Sigma \D^2 h \Sigma \D h^\top]}{n^3}+O(n^{-4}), \quad  \EE[A_1^6] = \frac{15 \alpha_2^3}{n^3}+O(n^{-4}).
        \end{align*}
        After multiplication by the $n^2$ prefactors above, the omitted terms
        contribute $O(n^{-2})$.  Moreover, every unlisted term in the expansion
        of the fourth central moment has total degree at least seven in centered
        empirical averages and is $O(n^{-2})$ by the standard partition count
        for products of centered sample averages and the available moment bounds.

        As for $\left(\EE\left[\left(\sqrt{n} R - \EE[\sqrt{n} R]\right)^2\right]\right)^2$, we have
        \begin{align*}
            &\left(\EE\left[\left(\sqrt{n} R - \EE[\sqrt{n} R]\right)^2\right]\right)^2 \\
            = & 1 + 2 \left(\EE\left[\left(\sqrt{n} R - \EE[\sqrt{n} R]\right)^2\right] - 1\right)\\
            = & 1 + \frac{2}{n}\left(-\frac{\alpha_4 - \alpha_2^2}{ \alpha_2^{2}} + \frac{\alpha_3 \tilde \alpha_3}{\alpha_2\tilde \alpha_2^2}  + \frac{3 \left(\alpha_2 \alpha_4 - \alpha_2^3 + 2\alpha_3^2\right)}{4\alpha_2^{3}} - \frac{3\tilde \alpha_3  \alpha_3}{2\alpha_2\tilde \alpha_2^2}  \right.\\
            & \left. - \frac{6\tilde \alpha_3 \left(\alpha_2 \EE\brac{h(X)\D h(X) \Sigma \D h(X)^\top}\right)}{\alpha_2\tilde \alpha_2^{3}}\right. \\
            & + \frac{3 \alpha_2 \EE\brac{h(X)\D h(X) \Sigma \D^2 h(X) \Sigma\D h(X)^\top}}{\alpha_2\tilde \alpha_2^2} + \left( - \frac{\tilde \alpha_3^2}{4 \alpha_{2}\tilde \alpha_2^4} + \frac{\tilde \alpha_4}{\alpha_2\tilde \alpha_2^3}\right)\left(3 \alpha_2^2\right) \\
            & + \frac{\alpha_2\alpha_4 - \alpha_2^3 + \alpha_3^2}{4 \alpha_2^3} - \frac{\tilde \alpha_3 \alpha_3}{ \alpha_2 \tilde \alpha_2^2} \left. + \frac{\tilde \alpha_3^2 \alpha_2}{2\tilde \alpha_2^4} \right) +  O(n^{-2}). 
        \end{align*}
        Combining these terms gives
        \begin{align*}
            \tilde m_{4} =& \frac{1}{n} \left(\frac{-2 \alpha_4}{\alpha_2^2} + \frac{12 \alpha_3^2}{\alpha_2^3} + \frac{-18 \alpha_3 \tilde \alpha_3}{\alpha_2 \tilde \alpha_2^2} + \frac{12 \EE[h(X) \D h(X)\Sigma \D^2 h(X) \Sigma \D h(X)^\top]}{\tilde \alpha_2^2}\right. \\
            &\left.+ \frac{12 \alpha_2 \tilde \alpha_4 - 24 \tilde \alpha_3\EE[h(X)\left(\D h(X)\Sigma \D h(X)^\top\right)]}{\tilde \alpha_2^3}  + \frac{9\alpha_2 \tilde \alpha_3^2}{\tilde \alpha_2^4}\right) + O\pare{n^{-2}}. 
        \end{align*}
    \end{enumerate}
    This completes the proof.
\end{proof}

\bibliographystyle{plainnat}
\bibliography{references}
\end{document}